\documentclass{article}
\usepackage{amsmath,amssymb,amsthm}
 \title{Complex structures on indecomposable 6-dimensional nilpotent real Lie Algebras
\thanks{\textit{Math. Subj. Class. [2000]}  : 17B30 (Primary), 53C15 (Secondary)}
}
\author{L. Magnin\\Mathematical Institute of Bourgogne
\thanks{UMR CNRS 5584, Universit\'{e} de Bourgogne, BP 47870, 21078 Dijon Cedex, France.}\\
\small{Email: \texttt{magnin@u-bourgogne.fr}}}
\date{\today}
\date{}
\def\Cmath{\mathbb{C}}

\def\Rmath{\mathbb{R}}

\def\T{\mathcal{T}}

\newcommand{\bv}[1]{\mathbf{#1}}

\newcommand{\myquad}{\quad}
\theoremstyle{plain}
\newtheorem{lemma}{Lemma}
\newtheorem{theorem}{Theorem}
\begin{document}
\maketitle
\begin{abstract}
We compute all complex  structures
on indecomposable 6-dimensional real Lie algebras and  their equivalence classes.
We also give for each of them a global holomorphic chart on  the
connected simply connected Lie group associated to the real Lie algebra and write down the multiplication in that chart.
\end{abstract}
%{\fontsize{8}{10}\selectfont
%\tableofcontents
%}
%%%%%%%%%%%%%%%%%%%%%%%%%%%%%%%%%%%%%%%%%%%%%%%%%%%%%%%%%%%%%%%%%%%%%%%%%%%%%%%%%%%%%%%%%%
\section{Introduction.}
%%%%%%%%%%%%%%%%%%%%%%%%%%%%%%%%%%%%%%%%%%%%%%%%%%%%%%%%%%%%%%%%%%%%%%%%%%%%%%%%%%%%%%%%%%
In the classification of nilmanifolds, a important question is to determine
the set of all integrable left invariant complex structures on a given connected
simply connected
real nilpotent finite dimensional Lie group,
or at the Lie algebra level the set
${\mathfrak{X}}_{\mathfrak{g}}$ of integrable complex structures on the nilpotent Lie algebra
${\mathfrak{g}},$
and its moduli space (\cite{Hasegawa},\cite{Cordero}, \cite{Snow}).
In the case of 6-dimensional real nilpotent Lie groups, an upper bound has been given
in  \cite{Salamon} for the dimension of
${\mathfrak{X}}_{\mathfrak{g}},$ based on a subcomplex of the Dolbeault complex.
These bounds are listed there, and Lie algebras which do not admit complex structures are specified. However, no detailed descriptions of the spaces
${\mathfrak{X}}_{\mathfrak{g}}$ are given. The aim of the present paper is to contribute in this area by supplying explicit computations of the various
${\mathfrak{X}}_{\mathfrak{g}}$  and their equivalence classes for
any indecomposable 6-dimensional real Lie algebra ${\mathfrak{g}}.$
We here are interested only in indecomposable Lie algebras, though
direct products could be processed in the same way.
%%%%%%%%%%%%%%%%%%%%%%%%%%%%%%%%%%%%%%%%%%%%%%%%%%%%%%%%%%%%%%%%%%%%%%%%%%%%%%%%%%%%%%%%%%
\section{Preliminaries.}
%%%%%%%%%%%%%%%%%%%%%%%%%%%%%%%%%%%%%%%%%%%%%%%%%%%%%%%%%%%%%%%%%%%%%%%%%%%%%%%%%%%%%%%%%%
\subsection{Labeling the algebras.}
There are 22 indecomposable nonisomorphic nilpotent real 6-dimensional Lie algebras
in the Morozov classification, labeled $M1$-$M22$ (\cite{JGP}).
Types $M14$ and $M18$ are splitted in $M14_{\pm 1}$ and $M18_{\pm 1}.$
Over $\Cmath,$ types $M14$ and $M18$ are not splitted and types $M5$ and $M10$ do not appear.
In \cite{Magnin1}, one is concerned with rank and weight systems over $\Cmath$, and a different classification is used. The correspondance with Morozov types appears there on page 130.
In the present paper, we label the algebras according to
\cite{Magnin1}, except for $M5,$  $M10,$  $M14$ and $M18.$
Note that $M5$ is the realification
$\mathfrak{n}_{\Rmath}$
of the 3-dimensional complex Heisenberg
algebra $\mathfrak{n}.$
Though $M10$ is not a realification, it appears as a subalgebra of the realification
$({\mathfrak{g}_4})_{\Rmath}$ of the complex 4-dimensional generic filiform Lie algebra
${\mathfrak{g}_4}$ in the isomorphic realisation
$[a_1,a_2]= a_3, [a_1,a_3]= a_4, [a_2,a_3]= a_4 :$ just take
$x_1= a_1, x_2 = ia_2, x_3 = ia_3, x_4=a_3, x_5= ia_4, x_6 = a_4.$
Let
$\mathfrak{g}$  be any of the labeled 6-dimensional real Lie algebras, and  let $G_0$ be the connected simply connected Lie group with Lie algebra
$\mathfrak{g}.$ From the commutation relations of the basis
$(x_j)_{1\leqslant j \leqslant6}$ of
$\mathfrak{g}$  we use,
the second kind canonical coordinates ($x \in G_0$)
\begin{equation}
\label{x6general}
x =\exp{(x^1 x_1)} \exp{(y^1 x_2)} \exp{(x^2 x_3)} \exp{(y^2 x_4)}\exp{(x^3 x_5)}\exp{(y^3 x_6)}
\end{equation}
yield a global chart for $G_0$ (see \cite{Varadarajan}, Th. 3.18.11, p.243).
We use this chart for $G_0$ in all cases  but the case of $M5$ where the natural chart is used instead.
For $1 \leqslant j \leqslant 6,$ denote by $X_j$ the left invariant vector field
on $G_0$ associated to $x_j,$
\textit{i.e.}
$$(X_j f)(x) = \left[ \frac{d}{dt} \, f(x\exp{(tx_j)}) \right]_{t=0} \; \forall \; f \in C^{\infty}(G_0).$$
Then due to the commutation relations, we have in each case
except $M5$ :
\begin{equation}
\label{X3-6}
X_3 = \frac{\partial}{\partial x^2} \; ; \;
X_4 = \frac{\partial}{\partial y^2}\; ; \;
X_5 = \frac{\partial}{\partial x^3}\; ; \;
X_6 = \frac{\partial}{\partial y^3}.
\end{equation}
\subsection{Complex structures.}
Let $\mathfrak{g}$  any finite dimensional real Lie algebra, and let
$G_0$ be the connected simply connected real Lie group with Lie algebra
$\mathfrak{g}.$
An almost complex structure on
$\mathfrak{g}$ is a linear map $J \, : \, \mathfrak{g} \rightarrow  \mathfrak{g}$
such that $J^2=-1.$ An almost complex structure on  $G_0$  is a tensor field
$x \mapsto J_x $
which at every point $x \in G_0$ is an endomorphism  of $T_x(G_0)$ such that $J_x^2=-1
.$
By definition, the almost complex structure on $G_0$ is left (resp. right) invariant if
$J_{ax}= \hat{(L_a)}_x J_x$
(resp. $J_{xa}= \hat{(R_a)}_x J_x$)
for all $a, x \in G_0,$ where
$\hat{(L_a)}_x J_x$
(resp. $\hat{(R_a)}_x J_x$)
is the endomorphism
${(L_a)_*}_x \circ J_x \circ {(L_{a^{-1}})_*}_{ax}$
(resp. ${(R_a)_*}_x \circ J_x \circ {(R_{a^{-1}})_*}_{xa})$
of $T_{ax}(G_0)$  (resp. $T_{xa}(G_0)$),
with $L_a$  (resp. $R_a$) the left (resp. right) translation
$x \mapsto ax$
(resp. $x \mapsto xa)$ and $(.)_*$ the differential.
For any almost complex structure $J$  on  $\mathfrak{g}$
there is a unique left invariant almost complex $\hat{J}$ structure on  $G_0$ such that
$\hat{J}_e = J$  ($e$ is the identity of $G_0$), and one has
$\hat{J}_{a}= \hat{(L_a)}_{*_e} J$ for all $a \in G_0.$
It is easily seen that $\hat{J}$  is right invariant if and only if
$$J \circ ad\, X = ad \, X \circ J \quad \forall X \in  \mathfrak{g},$$
that is  $(\mathfrak{g},J)$ is a \textit{complex Lie algebra}.
From the Newlander-Nirenberg theorem (\cite{NN}), $\hat{J}$ is \textit{integrable},
that is $G_0$ can be given the structure of a complex manifold with the same underlying
real structure and such that $\hat{J}$ is the canonical complex structure,
if and only if the torsion tensor of $\hat{J}$ vanishes, i.e. :
$$ [\hat{J} X, \hat{J}Y]-[X,Y]-\hat{J}[\hat{J}X,Y]-\hat{J}[X,\hat{J}Y] =0$$
for all vector fields $X,Y$ on $G_0.$ By left invariance, this is equivalent to
\begin{equation}
\label{Jinv}
 [{J} X, {J}Y]-[X,Y]-{J}[{J}X,Y]-J[X,{J}Y] =0 \quad \forall X,Y \in \mathfrak{g}.
 \end{equation}
By a complex structure  on $\mathfrak{g}$, we'll mean an
\textit{integrable} almost complex structure on  $\mathfrak{g},$ that is one satisfying
(\ref{Jinv}).
\par
Let $J$ a complex structure on $\mathfrak{g}$ and denote by $G$
the group $G_0$ endowed with the structure
of complex manifold  defined by $\hat{J}.$
Then a smooth function $f \, : \, G_0 \rightarrow G_0$ is holomorphic if and only if
its differential commutes with  $\hat{J}$ (\cite{KN}, Prop. 2.3 p. 123) :
$   \hat{J} \circ f_* =  f_* \circ \hat{J}.$
Hence left translations are holomorphic. Right translations are holomorphic, that is
$G$ is a complex Lie group,
if and only if
$\hat{J}$ is right invariant, i.e.
$(\mathfrak{g},J)$ is a complex Lie algebra.
The complexification $\mathfrak{g}_{\Cmath}$ of $\mathfrak{g}$ splits
as $\mathfrak{g}_{\Cmath}= \mathfrak{g}^{(1,0)} \oplus \mathfrak{g}^{(0,1)}$
where $\mathfrak{g}^{(1,0)}= \{X-iJX; \, X \in \mathfrak{g}\},$
$\mathfrak{g}^{(0,1)}= \{X+iJX; \, X \in \mathfrak{g}\}.$
We will denote  $\mathfrak{g}^{(1,0)}$ by $\mathfrak{m}.$ The integrability of $J$
amounts to
$\mathfrak{m}$ being a complex subalgebra of
$\mathfrak{g}_{\Cmath}.$ In that way  the set of complex structures on $\mathfrak{g}$ can be identified with the set of all complex subalgebras $\mathfrak{m}$ of
$\mathfrak{g}_{\Cmath}$ such that
$\mathfrak{g}_{\Cmath} = \mathfrak{m} \oplus \bar{\mathfrak{m}},$  bar denoting conjugation
in $\mathfrak{g}_{\Cmath}.$ This is the algebraic approach.
Our approach is more trivial since we simply fix a basis of $\mathfrak{g}$ and
compute all possible matrices in that basis for a complex structure.
From now on, we'll use the same notation $J$ for $J$ and $\hat{J}$ as well.
For any $x \in G_0,$
the complexification ${T_x(G_0)}_{\Cmath}$ of the tangent space
also splits as the direct sum of the holomorphic vectors
${T_x(G_0)}^{(1,0)}= \{X-iJX; \, X \in T_x(G_0)\}$ and the antiholomorphic vectors
${T_x(G_0)}^{(0,1)}= \{X+iJX; \, X \in T_x(G_0)\}.$
Let $H_{\Cmath}(G)$ be the space of complex valued holomorphic functions on $G.$
Then  $H_{\Cmath}(G)$ is comprised of all complex smooth functions $f$ on $G_0$
 which are annihilated by any antiholomorphic vector field. This is equivalent
to $f$ being annihilated by
all
\begin{equation}
\label{Xj-tilde}
\tilde{X}_j^{-} = X_j +i J X_j   \quad 1 \leqslant j \leqslant n
\end{equation}
with $(X_j)_{1 \leqslant j \leqslant n}$ the left invariant vector fields associated   to a basis $(x_j)_{1 \leqslant j \leqslant n}$ of
$\mathfrak{g}.$
Hence :
\begin{equation}
\label{H(G)}
H_{\Cmath}(G) =\{f \in C^{\infty}(G_0) \; ; \; \tilde{X}_j^{-} \, f = 0
\; \; \forall j\;\;  1 \leqslant j \leqslant n\}.
\end{equation}
\par
Finally, the automorphism group
$\text{Aut } \mathfrak{g}$ of $\mathfrak{g}$ acts on the set
$\mathfrak{X}_{\mathfrak{g}}$ of all complex structures on  $\mathfrak{g}$ by
$J \mapsto \Phi^{-1}\circ J \circ \Phi \quad \forall \Phi \in
\text{Aut } \mathfrak{g}.$
Two complex structures $J_1,J_2$
on   $\mathfrak{g}$  are said to be \textit{equivalent} if they are on the same
$\text{Aut } \mathfrak{g}$ orbit.
\subsection{Presentation of results.}
We consider here only indecomposable 6-dimensional nilpotent real Lie algebras which admit complex
structures. For each such $\mathfrak{g},$
we first give the commutation relations of the basis
$(x_j)_{ 1 \leqslant j \leqslant 6}$ of $\mathfrak{g}$ we use, and the
matrices $J= (J^k_j)= (\xi^k_j)$ in that
basis of the elements of ${\mathfrak{X}}_{\mathfrak{g}}$.
%The parameters are $\boxed{'boxed'}$.
The parameters are 'boxed'.
These matrices have been obtained by developping specific programs with the computer algebra system \textit{Reduce} by A. Hearn. The programs solve simultaneously
 the equation $J^2=-1$ and the torsion equations
$ij|k$ ($1 \leqslant i,j,k \leqslant 6$ ) obtained by projecting on $x_k$ the equation
 $[{J} x_i, {J}x_j]-[x_i,x_j]-{J}[{J}x_i,x_j]-J[x_i,{J}x_j] =0.$
We do not enter computational technicalities here, referring instead to
the technical report \cite{sc1}.
Let's simply say that the equations are \textit{semilinear} in the sense that
they can be solved in a succession of steps, each of which consists in solving some equation of degree 1 in some variable.
For all the Lie algebras we consider,
we prove that ${\mathfrak{X}}_{\mathfrak{g}}$ is a (smooth) submanifold of $\Rmath^{36}.$
The dimension of ${\mathfrak{X}}_{\mathfrak{g}}$ is equal to the upper bound given in  \cite{Salamon} except in the case of $M10.$
Then we give the automorphism group of $\mathfrak{g},$ representatives
of the various equivalence classes, and the commutation relations of the corresponding
algebra $\mathfrak{m} = \mathfrak{g}^{(1,0)}$ in terms of the basis
$(\tilde{x}_j)_{ 1 \leqslant j \leqslant 6}$
with $\tilde{x}_j = x_j -iJx_j.$
As $\frak{m}$ is a 3-dimensional complex Lie algebra, it is either abelian or isomorphic to the complex Heisenberg Lie algebra $\frak{n}.$
Hence, as a real Lie algebra, $\frak{m}$ is either abelian or isomorphic to $M5.$
Finally, we compute  the left invariant vector fields $X_1,X_2$ on $G_0$
in terms of the second kind canonical coordinates  (\ref{x6general})
(except in the case of $M5$, the 4 others appear in  (\ref{X3-6})), a global holomorphic chart on $G$ and the explicit look of the multiplication in $G$ in terms of that chart.
The fact that left translations are holomorphic, though the multiplication isn't except for
the canonical structure on $M5$, appears
clearly on these formulae.
There we make use of the following formula
which is easily checked from the Campbell-Hausdorff-Baker Formula :
\begin{equation}
\label{CBHcommmutateurC5}
e^X e^Y =
e^{[X,Y] +\frac{1}{2}\, \left( [X,[X,Y]] + [Y,[X,Y]] \right)
+\frac{1}{6} \,\left( [X,[X,[X,Y]]] + [Y,[Y,[X,Y]]]\right)
+\frac{1}{4}\, [X,[Y,[X,Y]]]
\;  \text{mod } \mathcal{C}^5
\mathfrak{g}}
\;
e^Y e^X
\end{equation}
where $\mathcal{C}^5  \mathfrak{g}$ denotes the $5^{\text{th}}$ central derivative.
As usual, for complex $z = x+iy \; (x,y \in \Rmath),$
$\frac{\partial}{\partial {z}}=\frac{1}{2}\left(
\frac{\partial}{\partial {x}} -i \, \frac{\partial}{\partial {y}}\right)
\; ; \; \frac{\partial}{\partial \overline{z}}=\frac{1}{2}\left(
\frac{\partial}{\partial {x}} +i \, \frac{\partial}{\partial {y}}\right) .
$
We'll abbreviate 'complex structure' to CS.
%%%%%%%%%%%%%%%%%%%%%%%%%%%%%%%%%%%%%%%%%%%%%%%%%%%%%%%%%%%%%%%%%%%%%%%%%%%%%%%%%%%%%%%%%%
\section{Lie Algebra $ {\mathcal{G}}_{6,3}$ (isomorphic to $M3$).}
%%%%%%%%%%%%%%%%%%%%%%%%%%%%%%%%%%%%%%%%%%%%%%%%%%%%%%%%%%%%%%%%%%%%%%%%%%%%%%%%%%%%%%%%%%
\smallskip  \par
Commutation relations for
$ {\mathcal{G}}_{6,3}:$
$[x_1,x_2]=x_4$;
$[x_1,x_3]=x_5$;
$[x_2,x_3]=x_6.$
\subsection{Case $ {\xi^1_6} \neq 0 .$}
\label{case31}
{\fontsize{8}{10}\selectfont
\begin{equation}
\label{63general1}
J = \begin{pmatrix}
\boxed{\xi^1_1}&
\boxed{\xi^1_2}&
\boxed{\xi^1_3}&
\xi^3_6&
-\xi^2_6&
\boxed{\xi^1_6}\\
*&
%  J^2_1=
%({\xi^5_6} {\xi^2_6}^2 + {\xi^5_5} {\xi^2_6} {\xi^1_6} - {\xi^4_6} {\xi^3_6} {\xi^2_6} - {\xi^4_5} {\xi^3_6} {\xi^1_6} + {\xi^2_6} {\xi^1_6} {\xi^1_1})/{\xi^1_6}^2&
%$  J^2_2=
( - {\xi^5_6} {\xi^2_6} - {\xi^5_5} {\xi^1_6} + {\xi^2_6} {\xi^1_2})/{\xi^1_6} &
%$  J^2_3=
({\xi^4_6} {\xi^2_6} + {\xi^4_5} {\xi^1_6} + {\xi^2_6} {\xi^1_3})/{\xi^1_6}&
%$  J^2_4=
{\xi^3_6} {\xi^2_6}/{\xi^1_6}&
%$  J^2_5=
 - {\xi^2_6}^2/{\xi^1_6}&
\boxed{\xi^2_6}\\
*&
*&
%$  J^3_3=
({\xi^5_6} {\xi^2_6} + {\xi^5_5} {\xi^1_6} + {\xi^3_6} {\xi^1_3})/{\xi^1_6}&
%$  J^3_4=
{\xi^3_6}^2/{\xi^1_6}&
%$  J^3_5=
 - {\xi^3_6} {\xi^2_6}/{\xi^1_6}&
\boxed{\xi^3_6}\\
*&
%$  J^4_2=
({\xi^5_6} {\xi^1_3} - {\xi^5_3} {\xi^1_6} + {\xi^4_6} {\xi^1_2})/{\xi^1_6}&
*&
%$  J^4_4=
( - {\xi^5_6} {\xi^2_6} - {\xi^5_5} {\xi^1_6} + {\xi^4_6} {\xi^3_6})/{\xi^1_6}
&
\boxed{\xi^4_5}&
\boxed{\xi^4_6}\\
*&
\boxed{\xi^5_2}&
\boxed{\xi^5_3}&
*&
\boxed{\xi^5_5}&
\boxed{\xi^5_6}\\
*&
*&
*&
*&
%$  J^6_5=
%({\xi^5_5} {\xi^2_6} {\xi^1_6} - {\xi^4_5} {\xi^3_6} {\xi^1_6} + {\xi^3_6} {\xi^2_6}
%{\xi^1_3} + {\xi^2_6}^2 {\xi^1_2} + {\xi^2_6} {\xi^1_6} {\xi^1_1})/{\xi^1_6}^2&
*&
*&
%$  J^6_6=
%({\xi^5_6} {\xi^2_6} - {\xi^4_6} {\xi^3_6} - {\xi^3_6} {\xi^1_3} - {\xi^2_6} {\xi^1_2} - {\xi^1_6} {\xi^1_1})/{\xi^1_6}
 \end{pmatrix}
\end{equation}
}
where
{\fontsize{6}{10}\selectfont
$ \mathbf{J^2_1}=({\xi^5_6} {\xi^2_6}^2 + {\xi^5_5} {\xi^2_6} {\xi^1_6} - {\xi^4_6} {\xi^3_6} {\xi^2_6} - {\xi^4_5} {\xi^3_6} {\xi^1_6} + {\xi^2_6} {\xi^1_6} {\xi^1_1})/{\xi^1_6}^2;$
\quad
$ \mathbf{ J^3_1}=({\xi^5_6}^2 {\xi^2_6}^3 + 2 {\xi^5_6} {\xi^5_5} {\xi^2_6}^2 {\xi^1_6} - {\xi^5_6} {\xi^4_6} {\xi^3_6} {\xi^2_6}^2 - {\xi^5_6} {\xi^4_5} {\xi^3_6} {\xi^2_6} {\xi^1_6} + {\xi^5_5}^2 {\xi^2_6} {\xi^1_6}^2 - {\xi^5_5} {\xi^4_6} {\xi^3_6} {\xi^2_6} {\xi^1_6} - {\xi^5_5} {\xi^4_5} {\xi^3_6} {\xi^1_6}^2 + {\xi^4_6} {\xi^3_6} {\xi^2_6} {\xi^1_6} {\xi^1_1} + {\xi^4_5} {\xi^3_6} {\xi^1_6}^2 {\xi^1_1} + {\xi^2_6} {\xi^1_6}^2)/({\xi^1_6}^2 ({\xi^4_6} {\xi^2_6} + {\xi^4_5} {\xi^1_6}));$
\quad
$  \mathbf{J^3_2}=( - {\xi^5_6}^2 {\xi^2_6}^2 - 2 {\xi^5_6} {\xi^5_5} {\xi^2_6} {\xi^1_6} - {\xi^5_5}^2 {\xi^1_6}^2 + {\xi^4_6} {\xi^3_6} {\xi^2_6} {\xi^1_2} + {\xi^4_5} {\xi^3_6} {\xi^1_6} {\xi^1_2} - {\xi^1_6}^2)/({\xi^1_6} ({\xi^4_6} {\xi^2_6} + {\xi^4_5} {\xi^1_6}));$
\quad
$  \mathbf{J^4_1}=( - {\xi^5_6}^3 {\xi^4_5} {\xi^2_6}^2 {\xi^1_6} - {\xi^5_6}^3 {\xi^2_6}^3 {\xi^1_3} + {\xi^5_6}^2 {\xi^5_5} {\xi^4_6} {\xi^2_6}^2 {\xi^1_6} - 2 {\xi^5_6}^2 {\xi^5_5} {\xi^4_5} {\xi^2_6} {\xi^1_6}^2 - 2 {\xi^5_6}^2 {\xi^5_5} {\xi^2_6}^2 {\xi^1_6} {\xi^1_3} + {\xi^5_6}^2 {\xi^5_3} {\xi^2_6}^3 {\xi^1_6} + 2 {\xi^5_6}^2 {\xi^4_6} {\xi^3_6} {\xi^2_6}^2 {\xi^1_3} + {\xi^5_6}^2 {\xi^4_6} {\xi^2_6}^2 {\xi^1_6} {\xi^1_1} + 2 {\xi^5_6}^2 {\xi^4_5} {\xi^3_6} {\xi^2_6} {\xi^1_6} {\xi^1_3} + 2 {\xi^5_6} {\xi^5_5}^2 {\xi^4_6} {\xi^2_6} {\xi^1_6}^2 - {\xi^5_6} {\xi^5_5}^2 {\xi^4_5} {\xi^1_6}^3 - {\xi^5_6} {\xi^5_5}^2 {\xi^2_6} {\xi^1_6}^2 {\xi^1_3} + 2 {\xi^5_6} {\xi^5_5} {\xi^5_3} {\xi^2_6}^2 {\xi^1_6}^2 + 2 {\xi^5_6} {\xi^5_5} {\xi^4_6} {\xi^3_6} {\xi^2_6} {\xi^1_6} {\xi^1_3} + 2 {\xi^5_6} {\xi^5_5} {\xi^4_6} {\xi^2_6} {\xi^1_6}^2
 {\xi^1_1} + 2 {\xi^5_6} {\xi^5_5} {\xi^4_5} {\xi^3_6} {\xi^1_6}^2 {\xi^1_3} - 2 {\xi^5_6} {\xi^5_3} {\xi^4_6} {\xi^3_6} {\xi^2_6}^2 {\xi^1_6} - 2 {\xi^5_6} {\xi^5_3} {\xi^4_5} {\xi^3_6} {\xi^2_6} {\xi^1_6}^2 + {\xi^5_6} {\xi^4_6}^2 {\xi^3_6} {\xi^2_6}^2 {\xi^1_2} + 2 {\xi^5_6} {\xi^4_6} {\xi^4_5} {\xi^3_6} {\xi^2_6} {\xi^1_6} {\xi^1_2} + {\xi^5_6} {\xi^4_5}^2 {\xi^3_6} {\xi^1_6}^
2 {\xi^1_2} - {\xi^5_6} {\xi^4_5} {\xi^1_6}^3 - {\xi^5_6} {\xi^2_6} {\xi^1_6}^2 {\xi^1_3} + {\xi^5_5}^3 {\xi^4_6} {\xi^1_6}^3 + {\xi^5_5}^2 {\xi^5_3} {\xi^2_6} {\xi^1_6}^3 + {\xi^5_5}^2
 {\xi^4_6} {\xi^1_6}^3 {\xi^1_1} - 2 {\xi^5_5} {\xi^5_3} {\xi^4_6} {\xi^3_6} {\xi^2_6} {\xi^1_6}
^2 - 2 {\xi^5_5} {\xi^5_3} {\xi^4_5} {\xi^3_6} {\xi^1_6}^3 + {\xi^5_5} {\xi^4_6} {\xi^1_6}^3
+ {\xi^5_3} {\xi^2_6} {\xi^1_6}^3 - {\xi^5_2} {\xi^4_6}^2 {\xi^3_6} {\xi^2_6}^2 {\xi^1_6} - 2
 {\xi^5_2} {\xi^4_6} {\xi^4_5} {\xi^3_6} {\xi^2_6} {\xi^1_6}^2 - {\xi^5_2} {\xi^4_5}^2 {\xi^3_6}
 {\xi^1_6}^3 + {\xi^4_6} {\xi^1_6}^3 {\xi^1_1})/
  ({\xi^1_6}^2 (({\xi^5_6} {\xi^2_6}  + {\xi^5_5} {\xi^1_6})^2 + {\xi^1_6}^2));$
\quad
$ \mathbf{J^4_3}=( - {\xi^5_6}^2 {\xi^4_6} {\xi^2_6}^2 {\xi^1_3} - 2 {\xi^5_6}^2 {\xi^4_5} {\xi^2_6} {\xi^1_6} {\xi^1_3} - 2 {\xi^5_6} {\xi^5_5} {\xi^4_5} {\xi^1_6}^2 {\xi^1_3} + 2 {\xi^5_6}
{\xi^5_3} {\xi^4_6} {\xi^2_6}^2 {\xi^1_6} + 2 {\xi^5_6} {\xi^5_3} {\xi^4_5} {\xi^2_6} {\xi^1_6}^
2 - {\xi^5_6} {\xi^4_6}^2 {\xi^2_6}^2 {\xi^1_2} - 2 {\xi^5_6} {\xi^4_6} {\xi^4_5} {\xi^2_6} {\xi^1_6} {\xi^1_2} - {\xi^5_6} {\xi^4_5}^2 {\xi^1_6}^2 {\xi^1_2} + {\xi^5_5}^2 {\xi^4_6} {\xi^1_6}^2 {\xi^1_3} + 2 {\xi^5_5} {\xi^5_3} {\xi^4_6} {\xi^2_6} {\xi^1_6}^2 + 2 {\xi^5_5} {\xi^5_3}
 {\xi^4_5} {\xi^1_6}^3 + {\xi^5_2} {\xi^4_6}^2 {\xi^2_6}^2 {\xi^1_6} + 2 {\xi^5_2} {\xi^4_6}
{\xi^4_5} {\xi^2_6} {\xi^1_6}^2 + {\xi^5_2} {\xi^4_5}^2 {\xi^1_6}^3 + {\xi^4_6} {\xi^1_6}^2
{\xi^1_3})/
 ({\xi^1_6} (({\xi^5_6} {\xi^2_6} + {\xi^5_5} {\xi^1_6})^2 + {\xi^1_6}^2));$
\quad
$ \mathbf{J^5_1}=( - {\xi^5_6}^2 {\xi^4_6} {\xi^4_5} {\xi^2_6}^2 {\xi^1_6} - {\xi^5_6}^2 {\xi^4_5}^2 {\xi^2_6} {\xi^1_6}^2 + {\xi^5_6} {\xi^5_5} {\xi^4_6}^2 {\xi^2_6}^2 {\xi^1_6} - {\xi^5_6} {\xi^5_5} {\xi^4_5}^2 {\xi^1_6}^3 + {\xi^5_6} {\xi^4_6}^2 {\xi^3_6} {\xi^2_6}^2 {\xi^1_3}
+ {\xi^5_6} {\xi^4_6}^2 {\xi^2_6}^3 {\xi^1_2} + {\xi^5_6} {\xi^4_6}^2 {\xi^2_6}^2 {\xi^1_6}
{\xi^1_1} + 2 {\xi^5_6} {\xi^4_6} {\xi^4_5} {\xi^3_6} {\xi^2_6} {\xi^1_6} {\xi^1_3} + 2 {\xi^5_6}
{\xi^4_6} {\xi^4_5} {\xi^2_6}^2 {\xi^1_6} {\xi^1_2} + 2 {\xi^5_6} {\xi^4_6} {\xi^4_5} {\xi^2_6}
{\xi^1_6}^2 {\xi^1_1} + {\xi^5_6} {\xi^4_5}^2 {\xi^3_6} {\xi^1_6}^2 {\xi^1_3} + {\xi^5_6} {\xi^4_5}^2 {\xi^2_6} {\xi^1_6}^2 {\xi^1_2} + {\xi^5_6} {\xi^4_5}^2 {\xi^1_6}^3 {\xi^1_1} + {\xi^5_5}^2 {\xi^4_6}^2 {\xi^2_6} {\xi^1_6}^2 + {\xi^5_5}^2 {\xi^4_6} {\xi^4_5} {\xi^1_6}^3 -
{\xi^5_3} {\xi^4_6}^2 {\xi^3_6} {\xi^2_6}^2 {\xi^1_6} - 2 {\xi^5_3} {\xi^4_6} {\xi^4_5} {\xi^3_6} {\xi^2_6} {\xi^1_6}^2 - {\xi^5_3} {\xi^4_5}^2 {\xi^3_6} {\xi^1_6}^3 - {\xi^5_2} {\xi^4_6}^
2 {\xi^2_6}^3 {\xi^1_6} - 2 {\xi^5_2} {\xi^4_6} {\xi^4_5} {\xi^2_6}^2 {\xi^1_6}^2 - {\xi^5_2}
 {\xi^4_5}^2 {\xi^2_6} {\xi^1_6}^3 + {\xi^4_6}^2 {\xi^2_6} {\xi^1_6}^2 + {\xi^4_6} {\xi^4_5} {\xi^1_6}^3)/
  ({\xi^1_6}^2 ({\xi^4_6} {\xi^2_6} + {\xi^4_5} {\xi^1_6})^2);$
\quad
$  \mathbf{J^5_4}=( - {\xi^5_6}^2 {\xi^2_6}^2 - 2 {\xi^5_6} {\xi^5_5} {\xi^2_6} {\xi^1_6} + {\xi^5_6} {\xi^4_6} {\xi^3_6} {\xi^2_6} + {\xi^5_6} {\xi^4_5} {\xi^3_6} {\xi^1_6} - {\xi^5_5}^2 {\xi^1_6}^2 - {\xi^1_6}^2)/({\xi^1_6} ({\xi^4_6} {\xi^2_6} + {\xi^4_5} {\xi^1_6}));$
\quad
$  \mathbf{J^6_1}=( - {\xi^5_6}^4 {\xi^4_5} {\xi^2_6}^4 {\xi^1_6} - {\xi^5_6}^4 {\xi^2_6}^5 {\xi^1_3} + {\xi^5_6}^3 {\xi^5_5} {\xi^4_6} {\xi^2_6}^4 {\xi^1_6} - 3 {\xi^5_6}^3 {\xi^5_5} {\xi^4_5} {\xi^2_6}^3 {\xi^1_6}^2 - 4 {\xi^5_6}^3 {\xi^5_5} {\xi^2_6}^4 {\xi^1_6} {\xi^1_3} + {\xi^5_6}^3 {\xi^4_6} {\xi^4_5} {\xi^3_6} {\xi^2_6}^3 {\xi^1_6} + 3 {\xi^5_6}^3 {\xi^4_6} {\xi^3_6} {\xi^2_6}^4 {\xi^1_3} + {\xi^5_6}^3 {\xi^4_6} {\xi^2_6}^4 {\xi^1_6} {\xi^1_1} + {\xi^5_6}
^3 {\xi^4_5}^2 {\xi^3_6} {\xi^2_6}^2 {\xi^1_6}^2 + 3 {\xi^5_6}^3 {\xi^4_5} {\xi^3_6} {\xi^2_6}^3 {\xi^1_6} {\xi^1_3} + {\xi^5_6}^3 {\xi^4_5} {\xi^2_6}^3 {\xi^1_6}^2 {\xi^1_1} + 3 {\xi^5_6}^2 {\xi^5_5}^2 {\xi^4_6} {\xi^2_6}^3 {\xi^1_6}^2 - 3 {\xi^5_6}^2 {\xi^5_5}^2 {\xi^4_5} {\xi^2_6}^2 {\xi^1_6}^3 - 6 {\xi^5_6}^2 {\xi^5_5}^2 {\xi^2_6}^3 {\xi^1_6}^2 {\xi^1_3} - {\xi^5_6}^2 {\xi^5_5} {\xi^4_6}^2 {\xi^3_6} {\xi^2_6}^3 {\xi^1_6} + {\xi^5_6}^2 {\xi^5_5} {\xi^4_6} {\xi^4_5} {\xi^3_6} {\xi^2_6}^2 {\xi^1_6}^2 + 7 {\xi^5_6}^2 {\xi^5_5} {\xi^4_6}
{\xi^3_6} {\xi^2_6}^3 {\xi^1_6} {\xi^1_3} - {\xi^5_6}^2 {\xi^5_5} {\xi^4_6} {\xi^2_6}^4 {\xi^1_6} {\xi^1_2} + 2 {\xi^5_6}^2 {\xi^5_5} {\xi^4_6} {\xi^2_6}^3 {\xi^1_6}^2 {\xi^1_1} + 2 {\xi^5_6}^2 {\xi^5_5} {\xi^4_5}^2 {\xi^3_6} {\xi^2_6} {\xi^1_6}^3 + 7 {\xi^5_6}^2 {\xi^5_5} {\xi^4_5} {\xi^3_6} {\xi^2_6}^2 {\xi^1_6}^2 {\xi^1_3} - {\xi^5_6}^2 {\xi^5_5} {\xi^4_5} {\xi^2_6}^
3 {\xi^1_6}^2 {\xi^1_2} + 2 {\xi^5_6}^2 {\xi^5_5} {\xi^4_5} {\xi^2_6}^2 {\xi^1_6}^3 {\xi^1_1} - 2 {\xi^5_6}^2 {\xi^5_3} {\xi^4_6} {\xi^3_6} {\xi^2_6}^4 {\xi^1_6} - 2 {\xi^5_6}^2 {\xi^5_3} {\xi^4_5} {\xi^3_6} {\xi^2_6}^3 {\xi^1_6}^2 - {\xi^5_6}^2 {\xi^5_2} {\xi^4_6} {\xi^2_6}^5
 {\xi^1_6} - {\xi^5_6}^2 {\xi^5_2} {\xi^4_5} {\xi^2_6}^4 {\xi^1_6}^2 - 2 {\xi^5_6}^2 {\xi^4_6}^2 {\xi^3_6}^2 {\xi^2_6}^3 {\xi^1_3} + {\xi^5_6}^2 {\xi^4_6}^2 {\xi^3_6} {\xi^2_6}^4
{\xi^1_2} - {\xi^5_6}^2 {\xi^4_6}^2 {\xi^3_6} {\xi^2_6}^3 {\xi^1_6} {\xi^1_1} - 4 {\xi^5_6}^
2 {\xi^4_6} {\xi^4_5} {\xi^3_6}^2 {\xi^2_6}^2 {\xi^1_6} {\xi^1_3} + 2 {\xi^5_6}^2 {\xi^4_6}
{\xi^4_5} {\xi^3_6} {\xi^2_6}^3 {\xi^1_6} {\xi^1_2} - {\xi^5_6}^2 {\xi^4_6} {\xi^4_5} {\xi^3_6}
{\xi^2_6}^2 {\xi^1_6}^2 {\xi^1_1} - {\xi^5_6}^2 {\xi^4_6} {\xi^3_6} {\xi^2_6}^3 {\xi^1_6} {\xi^1_3} {\xi^1_1} - {\xi^5_6}^2 {\xi^4_6} {\xi^2_6}^4 {\xi^1_6} {\xi^1_2} {\xi^1_1} - {\xi^5_6}
^2 {\xi^4_6} {\xi^2_6}^3 {\xi^1_6}^2 {\xi^1_1}^2 - 2 {\xi^5_6}^2 {\xi^4_5}^2 {\xi^3_6}^
2 {\xi^2_6} {\xi^1_6}^2 {\xi^1_3} + {\xi^5_6}^2 {\xi^4_5}^2 {\xi^3_6} {\xi^2_6}^2 {\xi^1_6}
^2 {\xi^1_2} - {\xi^5_6}^2 {\xi^4_5} {\xi^3_6} {\xi^2_6}^2 {\xi^1_6}^2 {\xi^1_3} {\xi^1_1} -
 {\xi^5_6}^2 {\xi^4_5} {\xi^2_6}^3 {\xi^1_6}^2 {\xi^1_2} {\xi^1_1} - {\xi^5_6}^2 {\xi^4_5}
{\xi^2_6}^2 {\xi^1_6}^3 {\xi^1_1}^2 - 2 {\xi^5_6}^2 {\xi^4_5} {\xi^2_6}^2 {\xi^1_6}^3 -
2 {\xi^5_6}^2 {\xi^2_6}^3 {\xi^1_6}^2 {\xi^1_3} + 3 {\xi^5_6} {\xi^5_5}^3 {\xi^4_6} {\xi^2_6}^2 {\xi^1_6}^3 - {\xi^5_6} {\xi^5_5}^3 {\xi^4_5} {\xi^2_6} {\xi^1_6}^4 - 4 {\xi^5_6} {\xi^5_5}^3 {\xi^2_6}^2 {\xi^1_6}^3 {\xi^1_3} - 2 {\xi^5_6} {\xi^5_5}^2 {\xi^4_6}^2 {\xi^3_6}
{\xi^2_6}^2 {\xi^1_6}^2 - {\xi^5_6} {\xi^5_5}^2 {\xi^4_6} {\xi^4_5} {\xi^3_6} {\xi^2_6} {\xi^1_6}^3 + 5 {\xi^5_6} {\xi^5_5}^2 {\xi^4_6} {\xi^3_6} {\xi^2_6}^2 {\xi^1_6}^2 {\xi^1_3} - 2
{\xi^5_6} {\xi^5_5}^2 {\xi^4_6} {\xi^2_6}^3 {\xi^1_6}^2 {\xi^1_2} + {\xi^5_6} {\xi^5_5}^2 {\xi^4_6} {\xi^2_6}^2 {\xi^1_6}^3 {\xi^1_1} + {\xi^5_6} {\xi^5_5}^2 {\xi^4_5}^2 {\xi^3_6} {\xi^1_6}^4 + 5 {\xi^5_6} {\xi^5_5}^2 {\xi^4_5} {\xi^3_6} {\xi^2_6} {\xi^1_6}^3 {\xi^1_3} - 2 {\xi^5_6} {\xi^5_5}^2 {\xi^4_5} {\xi^2_6}^2 {\xi^1_6}^3 {\xi^1_2} + {\xi^5_6} {\xi^5_5}^2 {\xi^4_5} {\xi^2_6} {\xi^1_6}^4 {\xi^1_1} - 4 {\xi^5_6} {\xi^5_5} {\xi^5_3} {\xi^4_6} {\xi^3_6} {\xi^2_6}^3 {\xi^1_6}^2 - 4 {\xi^5_6} {\xi^5_5} {\xi^5_3} {\xi^4_5} {\xi^3_6} {\xi^2_6}^2 {\xi^1_6}^
3 - 2 {\xi^5_6} {\xi^5_5} {\xi^5_2} {\xi^4_6} {\xi^2_6}^4 {\xi^1_6}^2 - 2 {\xi^5_6} {\xi^5_5}
{\xi^5_2} {\xi^4_5} {\xi^2_6}^3 {\xi^1_6}^3 - 2 {\xi^5_6} {\xi^5_5} {\xi^4_6}^2 {\xi^3_6}^2
{\xi^2_6}^2 {\xi^1_6} {\xi^1_3} + 2 {\xi^5_6} {\xi^5_5} {\xi^4_6}^2 {\xi^3_6} {\xi^2_6}^3 {\xi^1_6} {\xi^1_2} - 2 {\xi^5_6} {\xi^5_5} {\xi^4_6}^2 {\xi^3_6} {\xi^2_6}^2 {\xi^1_6}^2 {\xi^1_1} - 4 {\xi^5_6} {\xi^5_5} {\xi^4_6} {\xi^4_5} {\xi^3_6}^2 {\xi^2_6} {\xi^1_6}^2 {\xi^1_3} + 4
{\xi^5_6} {\xi^5_5} {\xi^4_6} {\xi^4_5} {\xi^3_6} {\xi^2_6}^2 {\xi^1_6}^2 {\xi^1_2} - 2 {\xi^5_6} {\xi^5_5} {\xi^4_6} {\xi^4_5} {\xi^3_6} {\xi^2_6} {\xi^1_6}^3 {\xi^1_1} - 2 {\xi^5_6} {\xi^5_5}
 {\xi^4_6} {\xi^3_6} {\xi^2_6}^2 {\xi^1_6}^2 {\xi^1_3} {\xi^1_1} - 2 {\xi^5_6} {\xi^5_5} {\xi^4_6} {\xi^2_6}^3 {\xi^1_6}^2 {\xi^1_2} {\xi^1_1} - 2 {\xi^5_6} {\xi^5_5} {\xi^4_6} {\xi^2_6}^2
{\xi^1_6}^3 {\xi^1_1}^2 + {\xi^5_6} {\xi^5_5} {\xi^4_6} {\xi^2_6}^2 {\xi^1_6}^3 - 2 {\xi^5_6} {\xi^5_5} {\xi^4_5}^2 {\xi^3_6}^2 {\xi^1_6}^3 {\xi^1_3} + 2 {\xi^5_6} {\xi^5_5} {\xi^4_5}^
2 {\xi^3_6} {\xi^2_6} {\xi^1_6}^3 {\xi^1_2} - 2 {\xi^5_6} {\xi^5_5} {\xi^4_5} {\xi^3_6} {\xi^2_6}
 {\xi^1_6}^3 {\xi^1_3} {\xi^1_1} - 2 {\xi^5_6} {\xi^5_5} {\xi^4_5} {\xi^2_6}^2 {\xi^1_6}^3 {\xi^1_2} {\xi^1_1} - 2 {\xi^5_6} {\xi^5_5} {\xi^4_5} {\xi^2_6} {\xi^1_6}^4 {\xi^1_1}^2 - 3 {\xi^5_6} {\xi^5_5} {\xi^4_5} {\xi^2_6} {\xi^1_6}^4 - 4 {\xi^5_6} {\xi^5_5} {\xi^2_6}^2 {\xi^1_6}^3
 {\xi^1_3} + 2 {\xi^5_6} {\xi^5_3} {\xi^4_6}^2 {\xi^3_6}^2 {\xi^2_6}^3 {\xi^1_6} + 4 {\xi^5_6} {\xi^5_3} {\xi^4_6} {\xi^4_5} {\xi^3_6}^2 {\xi^2_6}^2 {\xi^1_6}^2 + 2 {\xi^5_6} {\xi^5_3}
{\xi^4_5}^2 {\xi^3_6}^2 {\xi^2_6} {\xi^1_6}^3 - {\xi^5_6} {\xi^4_6}^3 {\xi^3_6}^2 {\xi^2_6}
^3 {\xi^1_2} - 3 {\xi^5_6} {\xi^4_6}^2 {\xi^4_5} {\xi^3_6}^2 {\xi^2_6}^2 {\xi^1_6} {\xi^1_2}
 - 3 {\xi^5_6} {\xi^4_6} {\xi^4_5}^2 {\xi^3_6}^2 {\xi^2_6} {\xi^1_6}^2 {\xi^1_2} + {\xi^5_6}
{\xi^4_6} {\xi^4_5} {\xi^3_6} {\xi^2_6} {\xi^1_6}^3 + 3 {\xi^5_6} {\xi^4_6} {\xi^3_6} {\xi^2_6}^
2 {\xi^1_6}^2 {\xi^1_3} + {\xi^5_6} {\xi^4_6} {\xi^2_6}^2 {\xi^1_6}^3 {\xi^1_1} - {\xi^5_6}
{\xi^4_5}^3 {\xi^3_6}^2 {\xi^1_6}^3 {\xi^1_2} + {\xi^5_6} {\xi^4_5}^2 {\xi^3_6} {\xi^1_6}^4
 + 3 {\xi^5_6} {\xi^4_5} {\xi^3_6} {\xi^2_6} {\xi^1_6}^3 {\xi^1_3} + {\xi^5_6} {\xi^4_5} {\xi^2_6} {\xi^1_6}^4 {\xi^1_1} + {\xi^5_5}^4 {\xi^4_6} {\xi^2_6} {\xi^1_6}^4 - {\xi^5_5}^4 {\xi^2_6} {\xi^1_6}^4 {\xi^1_3} - {\xi^5_5}^3 {\xi^4_6}^2 {\xi^3_6} {\xi^2_6} {\xi^1_6}^3 - {\xi^5_5}^3 {\xi^4_6} {\xi^4_5} {\xi^3_6} {\xi^1_6}^4 + {\xi^5_5}^3 {\xi^4_6} {\xi^3_6} {\xi^2_6} {\xi^1_6}^3 {\xi^1_3} - {\xi^5_5}^3 {\xi^4_6} {\xi^2_6}^2 {\xi^1_6}^3 {\xi^1_2} + {\xi^5_5}^3
{\xi^4_5} {\xi^3_6} {\xi^1_6}^4 {\xi^1_3} - {\xi^5_5}^3 {\xi^4_5} {\xi^2_6} {\xi^1_6}^4 {\xi^1_2} - 2 {\xi^5_5}^2 {\xi^5_3} {\xi^4_6} {\xi^3_6} {\xi^2_6}^2 {\xi^1_6}^3 - 2 {\xi^5_5}^2
{\xi^5_3} {\xi^4_5} {\xi^3_6} {\xi^2_6} {\xi^1_6}^4 - {\xi^5_5}^2 {\xi^5_2} {\xi^4_6} {\xi^2_6}
^3 {\xi^1_6}^3 - {\xi^5_5}^2 {\xi^5_2} {\xi^4_5} {\xi^2_6}^2 {\xi^1_6}^4 + {\xi^5_5}^2
{\xi^4_6}^2 {\xi^3_6} {\xi^2_6}^2 {\xi^1_6}^2 {\xi^1_2} - {\xi^5_5}^2 {\xi^4_6}^2 {\xi^3_6}
 {\xi^2_6} {\xi^1_6}^3 {\xi^1_1} + 2 {\xi^5_5}^2 {\xi^4_6} {\xi^4_5} {\xi^3_6} {\xi^2_6} {\xi^1_6}^3 {\xi^1_2} - {\xi^5_5}^2 {\xi^4_6} {\xi^4_5} {\xi^3_6} {\xi^1_6}^4 {\xi^1_1} - {\xi^5_5}
^2 {\xi^4_6} {\xi^3_6} {\xi^2_6} {\xi^1_6}^3 {\xi^1_3} {\xi^1_1} - {\xi^5_5}^2 {\xi^4_6} {\xi^2_6}^2 {\xi^1_6}^3 {\xi^1_2} {\xi^1_1} - {\xi^5_5}^2 {\xi^4_6} {\xi^2_6} {\xi^1_6}^4 {\xi^1_1}^2 + {\xi^5_5}^2 {\xi^4_6} {\xi^2_6} {\xi^1_6}^4 + {\xi^5_5}^2 {\xi^4_5}^2 {\xi^3_6} {\xi^1_6}^4 {\xi^1_2} - {\xi^5_5}^2 {\xi^4_5} {\xi^3_6} {\xi^1_6}^4 {\xi^1_3} {\xi^1_1} - {\xi^5_5}^2 {\xi^4_5} {\xi^2_6} {\xi^1_6}^4 {\xi^1_2} {\xi^1_1} - {\xi^5_5}^2 {\xi^4_5} {\xi^1_6}^5
{\xi^1_1}^2 - {\xi^5_5}^2 {\xi^4_5} {\xi^1_6}^5 - 2 {\xi^5_5}^2 {\xi^2_6} {\xi^1_6}^4 {\xi^1_3} + 2 {\xi^5_5} {\xi^5_3} {\xi^4_6}^2 {\xi^3_6}^2 {\xi^2_6}^2 {\xi^1_6}^2 + 4 {\xi^5_5}
 {\xi^5_3} {\xi^4_6} {\xi^4_5} {\xi^3_6}^2 {\xi^2_6} {\xi^1_6}^3 + 2 {\xi^5_5} {\xi^5_3} {\xi^4_5}^2 {\xi^3_6}^2 {\xi^1_6}^4 - {\xi^5_5} {\xi^4_6}^2 {\xi^3_6} {\xi^2_6} {\xi^1_6}^3 - {\xi^5_5} {\xi^4_6} {\xi^4_5} {\xi^3_6} {\xi^1_6}^4 + {\xi^5_5} {\xi^4_6} {\xi^3_6} {\xi^2_6} {\xi^1_6}^3 {\xi^1_3} - {\xi^5_5} {\xi^4_6} {\xi^2_6}^2 {\xi^1_6}^3 {\xi^1_2} + {\xi^5_5} {\xi^4_5}
{\xi^3_6} {\xi^1_6}^4 {\xi^1_3} - {\xi^5_5} {\xi^4_5} {\xi^2_6} {\xi^1_6}^4 {\xi^1_2} - 2 {\xi^5_3} {\xi^4_6} {\xi^3_6} {\xi^2_6}^2 {\xi^1_6}^3 - 2 {\xi^5_3} {\xi^4_5} {\xi^3_6} {\xi^2_6} {\xi^1_6}^4 + {\xi^5_2} {\xi^4_6}^3 {\xi^3_6}^2 {\xi^2_6}^3 {\xi^1_6} + 3 {\xi^5_2} {\xi^4_6}
^2 {\xi^4_5} {\xi^3_6}^2 {\xi^2_6}^2 {\xi^1_6}^2 + 3 {\xi^5_2} {\xi^4_6} {\xi^4_5}^2 {\xi^3_6}^2 {\xi^2_6} {\xi^1_6}^3 - {\xi^5_2} {\xi^4_6} {\xi^2_6}^3 {\xi^1_6}^3 + {\xi^5_2} {\xi^4_5}^3 {\xi^3_6}^2 {\xi^1_6}^4 - {\xi^5_2} {\xi^4_5} {\xi^2_6}^2 {\xi^1_6}^4 + {\xi^4_6}^
2 {\xi^3_6} {\xi^2_6}^2 {\xi^1_6}^2 {\xi^1_2} - {\xi^4_6}^2 {\xi^3_6} {\xi^2_6} {\xi^1_6}^3
{\xi^1_1} + 2 {\xi^4_6} {\xi^4_5} {\xi^3_6} {\xi^2_6} {\xi^1_6}^3 {\xi^1_2} - {\xi^4_6} {\xi^4_5}
 {\xi^3_6} {\xi^1_6}^4 {\xi^1_1} - {\xi^4_6} {\xi^3_6} {\xi^2_6} {\xi^1_6}^3 {\xi^1_3} {\xi^1_1}
 - {\xi^4_6} {\xi^2_6}^2 {\xi^1_6}^3 {\xi^1_2} {\xi^1_1} - {\xi^4_6} {\xi^2_6} {\xi^1_6}^4 {\xi^1_1}^2 + {\xi^4_5}^2 {\xi^3_6} {\xi^1_6}^4 {\xi^1_2} - {\xi^4_5} {\xi^3_6} {\xi^1_6}^4 {\xi^1_3} {\xi^1_1} - {\xi^4_5} {\xi^2_6} {\xi^1_6}^4 {\xi^1_2} {\xi^1_1} - {\xi^4_5} {\xi^1_6}^5
{\xi^1_1}^2 - {\xi^4_5} {\xi^1_6}^5 - {\xi^2_6} {\xi^1_6}^4 {\xi^1_3})
 /({\xi^1_6}^3 ({\xi^4_6} {\xi^2_6} + {\xi^4_5} {\xi^1_6})(( {\xi^5_6} {\xi^2_6}
  + {\xi^5_5} {\xi^1_6})^2 +      {\xi^1_6}^2 ))
  ;$
\quad
$  \mathbf{J^6_2}=({\xi^5_6}^2 {\xi^2_6}^2 {\xi^1_3} + 2 {\xi^5_6} {\xi^5_5} {\xi^2_6} {\xi^1_6} {\xi^1_3} - {\xi^5_6} {\xi^4_6} {\xi^3_6} {\xi^2_6} {\xi^1_3} + {\xi^5_6} {\xi^4_6} {\xi^2_6}^2 {\xi^1_2} - {\xi^5_6} {\xi^4_5} {\xi^3_6} {\xi^1_6} {\xi^1_3} + {\xi^5_6} {\xi^4_5} {\xi^2_6} {\xi^1_6}
 {\xi^1_2} + {\xi^5_5}^2 {\xi^1_6}^2 {\xi^1_3} + {\xi^5_5} {\xi^4_6} {\xi^2_6} {\xi^1_6} {\xi^1_2} + {\xi^5_5} {\xi^4_5} {\xi^1_6}^2 {\xi^1_2} + {\xi^5_3} {\xi^4_6} {\xi^3_6} {\xi^2_6} {\xi^1_6} + {\xi^5_3} {\xi^4_5} {\xi^3_6} {\xi^1_6}^2 + {\xi^5_2} {\xi^4_6} {\xi^2_6}^2 {\xi^1_6} + {\xi^5_2} {\xi^4_5} {\xi^2_6} {\xi^1_6}^2 - {\xi^4_6}^2 {\xi^3_6} {\xi^2_6} {\xi^1_2} - {\xi^4_6}
{\xi^4_5} {\xi^3_6} {\xi^1_6} {\xi^1_2} - {\xi^4_6} {\xi^3_6} {\xi^2_6} {\xi^1_3} {\xi^1_2} - {\xi^4_6} {\xi^2_6}^2 {\xi^1_2}^2 - {\xi^4_6} {\xi^2_6} {\xi^1_6} {\xi^1_2} {\xi^1_1} - {\xi^4_5} {\xi^3_6} {\xi^1_6} {\xi^1_3} {\xi^1_2} - {\xi^4_5} {\xi^2_6} {\xi^1_6} {\xi^1_2}^2 - {\xi^4_5} {\xi^1_6}^2 {\xi^1_2} {\xi^1_1} + {\xi^1_6}^2 {\xi^1_3})/({\xi^1_6}^2 ({\xi^4_6} {\xi^2_6} + {\xi^4_5} {\xi^1_6}));$
\quad
$  \mathbf{J^6_3}=( - {\xi^5_6}^3 {\xi^2_6}^3 {\xi^1_3} - 3 {\xi^5_6}^2 {\xi^5_5} {\xi^2_6}^2
{\xi^1_6} {\xi^1_3} + {\xi^5_6}^2 {\xi^5_3} {\xi^2_6}^3 {\xi^1_6} + {\xi^5_6}^2 {\xi^4_6} {\xi^3_6} {\xi^2_6}^2 {\xi^1_3} - {\xi^5_6}^2 {\xi^4_6} {\xi^2_6}^3 {\xi^1_2} + 2 {\xi^5_6}^2
{\xi^4_5} {\xi^3_6} {\xi^2_6} {\xi^1_6} {\xi^1_3} - {\xi^5_6}^2 {\xi^4_5} {\xi^2_6}^2 {\xi^1_6}
{\xi^1_2} - {\xi^5_6}^2 {\xi^3_6} {\xi^2_6}^2 {\xi^1_3}^2 - {\xi^5_6}^2 {\xi^2_6}^3 {\xi^1_3} {\xi^1_2} - {\xi^5_6}^2 {\xi^2_6}^2 {\xi^1_6} {\xi^1_3} {\xi^1_1} - 3 {\xi^5_6} {\xi^5_5}^
2 {\xi^2_6} {\xi^1_6}^2 {\xi^1_3} + 2 {\xi^5_6} {\xi^5_5} {\xi^5_3} {\xi^2_6}^2 {\xi^1_6}^2 -
 2 {\xi^5_6} {\xi^5_5} {\xi^4_6} {\xi^2_6}^2 {\xi^1_6} {\xi^1_2} + 2 {\xi^5_6} {\xi^5_5} {\xi^4_5} {\xi^3_6} {\xi^1_6}^2 {\xi^1_3} - 2 {\xi^5_6} {\xi^5_5} {\xi^4_5} {\xi^2_6} {\xi^1_6}^2 {\xi^1_2} - 2 {\xi^5_6} {\xi^5_5} {\xi^3_6} {\xi^2_6} {\xi^1_6} {\xi^1_3}^2 - 2 {\xi^5_6} {\xi^5_5}
{\xi^2_6}^2 {\xi^1_6} {\xi^1_3} {\xi^1_2} - 2 {\xi^5_6} {\xi^5_5} {\xi^2_6} {\xi^1_6}^2 {\xi^1_3} {\xi^1_1} - 2 {\xi^5_6} {\xi^5_3} {\xi^4_6} {\xi^3_6} {\xi^2_6}^2 {\xi^1_6} - 2 {\xi^5_6} {\xi^5_3} {\xi^4_5} {\xi^3_6} {\xi^2_6} {\xi^1_6}^2 + {\xi^5_6} {\xi^4_6}^2 {\xi^3_6} {\xi^2_6}^2
{\xi^1_2} + 2 {\xi^5_6} {\xi^4_6} {\xi^4_5} {\xi^3_6} {\xi^2_6} {\xi^1_6} {\xi^1_2} + {\xi^5_6} {\xi^4_5}^2 {\xi^3_6} {\xi^1_6}^2 {\xi^1_2} - {\xi^5_6} {\xi^2_6} {\xi^1_6}^2 {\xi^1_3} - {\xi^5_5}^3 {\xi^1_6}^3 {\xi^1_3} + {\xi^5_5}^2 {\xi^5_3} {\xi^2_6} {\xi^1_6}^3 - {\xi^5_5}^2 {\xi^4_6} {\xi^3_6} {\xi^1_6}^2 {\xi^1_3} - {\xi^5_5}^2 {\xi^4_6} {\xi^2_6} {\xi^1_6}^2 {\xi^1_2}
 - {\xi^5_5}^2 {\xi^4_5} {\xi^1_6}^3 {\xi^1_2} - {\xi^5_5}^2 {\xi^3_6} {\xi^1_6}^2 {\xi^1_3}
^2 - {\xi^5_5}^2 {\xi^2_6} {\xi^1_6}^2 {\xi^1_3} {\xi^1_2} - {\xi^5_5}^2 {\xi^1_6}^3 {\xi^1_3} {\xi^1_1} - 2 {\xi^5_5} {\xi^5_3} {\xi^4_6} {\xi^3_6} {\xi^2_6} {\xi^1_6}^2 - 2 {\xi^5_5}
{\xi^5_3} {\xi^4_5} {\xi^3_6} {\xi^1_6}^3 - {\xi^5_5} {\xi^1_6}^3 {\xi^1_3} + {\xi^5_3} {\xi^2_6} {\xi^1_6}^3 - {\xi^5_2} {\xi^4_6}^2 {\xi^3_6} {\xi^2_6}^2 {\xi^1_6} - 2 {\xi^5_2} {\xi^4_6}
 {\xi^4_5} {\xi^3_6} {\xi^2_6} {\xi^1_6}^2 - {\xi^5_2} {\xi^4_5}^2 {\xi^3_6} {\xi^1_6}^3 - {\xi^4_6} {\xi^3_6} {\xi^1_6}^2 {\xi^1_3} - {\xi^4_6} {\xi^2_6} {\xi^1_6}^2 {\xi^1_2} - {\xi^4_5}
{\xi^1_6}^3 {\xi^1_2} - {\xi^3_6} {\xi^1_6}^2 {\xi^1_3}^2 - {\xi^2_6} {\xi^1_6}^2 {\xi^1_3}
{\xi^1_2} - {\xi^1_6}^3 {\xi^1_3} {\xi^1_1})
 /({\xi^1_6}^2 (({\xi^5_6} {\xi^2_6} + {\xi^5_5} {\xi^1_6})^2 + {\xi^1_6}^2));$
\quad
$  \mathbf{J^6_4}=( - {\xi^5_6}^2 {\xi^2_6}^3 - 2 {\xi^5_6} {\xi^5_5} {\xi^2_6}^2 {\xi^1_6} + 2
{\xi^5_6} {\xi^4_6} {\xi^3_6} {\xi^2_6}^2 + 2 {\xi^5_6} {\xi^4_5} {\xi^3_6} {\xi^2_6} {\xi^1_6} -
 {\xi^5_5}^2 {\xi^2_6} {\xi^1_6}^2 + {\xi^5_5} {\xi^4_6} {\xi^3_6} {\xi^2_6} {\xi^1_6} + {\xi^5_5} {\xi^4_5} {\xi^3_6} {\xi^1_6}^2 - {\xi^4_6}^2 {\xi^3_6}^2 {\xi^2_6} - {\xi^4_6} {\xi^4_5}
{\xi^3_6}^2 {\xi^1_6} - {\xi^4_6} {\xi^3_6}^2 {\xi^2_6} {\xi^1_3} - {\xi^4_6} {\xi^3_6} {\xi^2_6}^2 {\xi^1_2} - {\xi^4_6} {\xi^3_6} {\xi^2_6} {\xi^1_6} {\xi^1_1} - {\xi^4_5} {\xi^3_6}^2 {\xi^1_6} {\xi^1_3} - {\xi^4_5} {\xi^3_6} {\xi^2_6} {\xi^1_6} {\xi^1_2} - {\xi^4_5} {\xi^3_6} {\xi^1_6}
^2 {\xi^1_1} - {\xi^2_6} {\xi^1_6}^2)/({\xi^1_6}^2 ({\xi^4_6} {\xi^2_6} + {\xi^4_5} {\xi^1_6}
));$
\quad
$  \mathbf{J^6_5}=({\xi^5_5} {\xi^2_6} {\xi^1_6} - {\xi^4_5} {\xi^3_6} {\xi^1_6} + {\xi^3_6} {\xi^2_6}
{\xi^1_3} + {\xi^2_6}^2 {\xi^1_2} + {\xi^2_6} {\xi^1_6} {\xi^1_1})/{\xi^1_6}^2;$
\quad
$  \mathbf{J^6_6}=({\xi^5_6} {\xi^2_6} - {\xi^4_6} {\xi^3_6} - {\xi^3_6} {\xi^1_3} - {\xi^2_6} {\xi^1_2} - {\xi^1_6} {\xi^1_1})/{\xi^1_6};$
}
and the parameters are subject to the condition
\begin{equation}
\label{cond63_1general}
{\xi^1_6} ({\xi^4_6} {\xi^2_6} + {\xi^4_5} {\xi^1_6}) \neq 0 .
 \end{equation}
Now, equivalence by a suitable automorphism of the form
\begin{equation}
{\fontsize{8}{10}\selectfont
\label{63specialautom1}
\Phi = \begin{pmatrix}
1&
0&
0&
0&
0&
0\\
0&
1&
0&
0&
0&
0\\
0&
0&
1&
0&
0&
0\\
{b^4_1}&
{b^4_2}&
{b^4_3}&
1&
0&
0\\
{b^5_1}&
{b^5_2}&
{b^5_3}&
0&
1&
0\\
{b^6_1}&
{b^6_2}&
{b^6_3}&
0&
0&
1\end{pmatrix}
}
\end{equation}
reduces to the case
$ {\xi^1_1}={\xi^1_2}={\xi^1_3}={\xi^4_6}={\xi^5_2}={\xi^5_3}={\xi^5_6}= 0 .$
Then from (\ref{cond63_1general}), $\xi^4_5\xi^1_6 \neq 0$ and
applying equivalence by
${\fontsize{8}{10}\selectfont
 \Psi = \text{diag}\left( \begin{pmatrix}
1&
0&
0\\
{\xi^2_6}/{\xi^1_6}&
0&
 - {\xi^4_5} c\\
{\xi^3_6}/{\xi^1_6}&
c&
 - {\xi^5_5} c \end{pmatrix},\begin{pmatrix}
0&
 - {\xi^4_5} c&
0\\
c&
 - {\xi^5_5} c&

0\\
{\xi^2_6} c/{\xi^1_6}&
c ( - {\xi^5_5} {\xi^2_6} + {\xi^4_5} {\xi^3_6})/{\xi^1_6}&
{\xi^4_5} c^2\end{pmatrix} \right) }$
where
$c=|\xi^4_5\xi^1_6|^{-\frac{1}{2}}$, we get reduced according to the sign of $\xi^4_5\xi^1_6$ to either
\begin{equation}
{\fontsize{8}{10}\selectfont
\label{53-J_1}
 J_1 = \begin{pmatrix}
0&
0&
0&
0&
0&
1\\
0&
0&
1&
0&
0&
0\\
0&
-1&
0&
0&
0&
0\\
0&
0&
0&
0&
1&
0\\
0&
0&
0&
-1&
0&
0\\
-1&
0&
0&
0&
0&
0\end{pmatrix}
}
\text{ or }
{\fontsize{8}{10}\selectfont
J_1^{-} = \begin{pmatrix}
0&0&0&0&0&-1\\0&0&1&0&0&0\\0&-1&0&0&0&0\\0&0&0&0&1&0\\0&0&0&-1&0&0\\1&0&0&0&0&0
\end{pmatrix}.
}
\end{equation}
Now $ J_1^{-} $ is  equivalent to $ J_1 ,$
hence any  CS with $ {\xi^1_6} \neq 0 $
is equivalent to $ J_1.$
%%%%%%%%%%%%%%%%%%%%%%%%%%%%%%%%%%%%%%%%%%%%%%%%%%%%%%%%%%%%%%%%%%%%%%%%%%%%%%%%%%%%%%%
\subsection{Case $ {\xi^1_6} = 0, {\xi^2_5} \neq 0 .$}
{\fontsize{8}{10}\selectfont
\begin{equation}
\label{63general2}
J = \begin{pmatrix}
-\xi^6_6&
*&
*&
0&
0&
0\\
\boxed{\xi^2_1}&
*&
%$  J^2_2=( - {\xi^5_5} {\xi^2_5} - {\xi^4_5} {\xi^2_4} + {\xi^2_4} {\xi^2_3})/{\xi^2_5}
%( - {\xi^5_5} {\xi^2_5} - {\xi^4_5} {\xi^2_4} + {\xi^2_4} {\xi^2_3})/{\xi^2_5}&
\boxed{\xi^2_3}&
\boxed{\xi^2_4}&
\boxed{\xi^2_5}&
0\\
\boxed{\xi^3_1}&
*&
%$  J^3_3=({\xi^6_6} {\xi^2_5} - {\xi^2_4} {\xi^2_3})/{\xi^2_5};$\\
({\xi^6_6} {\xi^2_5} - {\xi^2_4} {\xi^2_3})/{\xi^2_5}&
%$  J^3_4=( - {\xi^2_4}^2)/{\xi^2_5};$\\
 - {\xi^2_4}^2/{\xi^2_5}&
%$  J^3_5= - {\xi^2_4};$\\
- {\xi^2_4}&
0\\
%$  J^4_1=( - {\xi^6_5} {\xi^2_3} + {\xi^6_3} {\xi^2_5} + {\xi^4_5} {\xi^2_1})/{\xi^2_5}
( - {\xi^6_5} {\xi^2_3} + {\xi^6_3} {\xi^2_5} + {\xi^4_5} {\xi^2_1})/{\xi^2_5} &
*&
\boxed{\xi^4_3}&
%$  J^4_4=( - {\xi^6_6} {\xi^2_5} + {\xi^4_5} {\xi^2_4})/{\xi^2_5};$\\
( - {\xi^6_6} {\xi^2_5} + {\xi^4_5} {\xi^2_4})/{\xi^2_5}&
\boxed{\xi^4_5}&
%$  J^4_6=({\xi^2_5} ({\xi^6_6}^2 + 1))/({\xi^3_1} {\xi^2_5} + {\xi^2_4} {\xi^2_1});$\\
{\xi^2_5} ({\xi^6_6}^2 + 1)/({\xi^3_1} {\xi^2_5} + {\xi^2_4} {\xi^2_1})\\
*&
*&
*&
%$  J^5_4=({\xi^2_4} ({\xi^6_6} + {\xi^5_5}))/{\xi^2_5};$\\
{\xi^2_4} ({\xi^6_6} + {\xi^5_5})/{\xi^2_5}&
\boxed{\xi^5_5}&
%$  J^5_6= - ({\xi^2_4} ({\xi^6_6}^2 + 1))/({\xi^3_1} {\xi^2_5} + {\xi^2_4} {\xi^2_1})
- {\xi^2_4} ({\xi^6_6}^2 + 1)/({\xi^3_1} {\xi^2_5} + {\xi^2_4} {\xi^2_1})\\
*&
*&
\boxed{\xi^6_3}&
%$  J^6_4=({\xi^6_5} {\xi^2_4} - {\xi^3_1} {\xi^2_5} - {\xi^2_4} {\xi^2_1})/{\xi^2_5};$\\
({\xi^6_5} {\xi^2_4} - {\xi^3_1} {\xi^2_5} - {\xi^2_4} {\xi^2_1})/{\xi^2_5}&
\boxed{\xi^6_5}&
\boxed{\xi^6_6}
 \end{pmatrix}
\end{equation}
}
where
{\fontsize{6}{10}\selectfont
$ \mathbf{J^1_3}= - ({\xi^2_5} ({\xi^6_6}^2 + 1))/({\xi^3_1} {\xi^2_5} + {\xi^2_4} {\xi^2_1})
;$\quad
$  \mathbf{J^2_2}=( - {\xi^5_5} {\xi^2_5} - {\xi^4_5} {\xi^2_4} + {\xi^2_4} {\xi^2_3})/{\xi^2_5}
;$\quad
$  \mathbf{J^3_2}=({\xi^2_4} ({\xi^6_6} {\xi^2_5} + {\xi^5_5} {\xi^2_5} + {\xi^4_5} {\xi^2_4} - {\xi^2_4}
 {\xi^2_3}))/{\xi^2_5}^2;$
  \quad
$  \mathbf{J^4_2}=( - {\xi^6_6}^2 {\xi^6_5} {\xi^2_5}^2 + {\xi^6_6} {\xi^4_5} {\xi^3_1} {\xi^2_5}^
2 + {\xi^6_6} {\xi^4_5} {\xi^2_5} {\xi^2_4} {\xi^2_1} - {\xi^6_5} {\xi^2_5}^2 - {\xi^5_5} {\xi^4_5}
 {\xi^3_1} {\xi^2_5}^2 - {\xi^5_5} {\xi^4_5} {\xi^2_5} {\xi^2_4} {\xi^2_1} - {\xi^4_5}^2 {\xi^3_1}
 {\xi^2_5} {\xi^2_4} - {\xi^4_5}^2 {\xi^2_4}^2 {\xi^2_1} + {\xi^4_3} {\xi^3_1} {\xi^2_5}^
2 {\xi^2_4} + {\xi^4_3} {\xi^2_5} {\xi^2_4}^2 {\xi^2_1})/({\xi^2_5}^2 ({\xi^3_1} {\xi^2_5} +
{\xi^2_4} {\xi^2_1}));$
  \quad
$  \mathbf{J^5_1}=({\xi^6_6} {\xi^2_5} {\xi^2_1} + {\xi^6_5} {\xi^2_4} {\xi^2_3} - {\xi^6_3} {\xi^2_5}
{\xi^2_4} + {\xi^5_5} {\xi^2_5} {\xi^2_1} - {\xi^3_1} {\xi^2_5} {\xi^2_3} - {\xi^2_4} {\xi^2_3}
{\xi^2_1})/{\xi^2_5}^2;$
  \quad
$  \mathbf{J^5_2}=({\xi^6_6}^2 {\xi^6_5} {\xi^2_5}^2 {\xi^2_4} + {\xi^6_6}^2 {\xi^2_5}^2 {\xi^2_4}
 {\xi^2_1} - {\xi^6_6} {\xi^4_5} {\xi^3_1} {\xi^2_5}^2 {\xi^2_4} - {\xi^6_6} {\xi^4_5} {\xi^2_5}
 {\xi^2_4}^2 {\xi^2_1} - {\xi^6_6} {\xi^3_1} {\xi^2_5}^2 {\xi^2_4} {\xi^2_3} - {\xi^6_6} {\xi^2_5}
 {\xi^2_4}^2 {\xi^2_3} {\xi^2_1} + {\xi^6_5} {\xi^2_5}^2 {\xi^2_4} - {\xi^5_5}^2 {\xi^3_1}
{\xi^2_5}^3 - {\xi^5_5}^2 {\xi^2_5}^2 {\xi^2_4} {\xi^2_1} - {\xi^5_5} {\xi^4_5} {\xi^3_1} {\xi^2_5}^2
 {\xi^2_4} - {\xi^5_5} {\xi^4_5} {\xi^2_5} {\xi^2_4}^2 {\xi^2_1} + {\xi^5_5} {\xi^3_1}
{\xi^2_5}^2 {\xi^2_4} {\xi^2_3} + {\xi^5_5} {\xi^2_5} {\xi^2_4}^2 {\xi^2_3} {\xi^2_1} + {\xi^4_5}
 {\xi^3_1} {\xi^2_5} {\xi^2_4}^2 {\xi^2_3} + {\xi^4_5} {\xi^2_4}^3 {\xi^2_3} {\xi^2_1} - {\xi^4_3}
 {\xi^3_1} {\xi^2_5}^2 {\xi^2_4}^2 - {\xi^4_3} {\xi^2_5} {\xi^2_4}^3 {\xi^2_1} - {\xi^3_1}
 {\xi^2_5}^3)/({\xi^2_5}^3 ({\xi^3_1} {\xi^2_5} + {\xi^2_4} {\xi^2_1}));$
  \quad
$  \mathbf{J^5_3}=({\xi^6_6}^2 {\xi^2_5}^2 {\xi^2_1} - {\xi^6_6} {\xi^3_1} {\xi^2_5}^2 {\xi^2_3} -
 {\xi^6_6} {\xi^2_5} {\xi^2_4} {\xi^2_3} {\xi^2_1} + {\xi^5_5} {\xi^3_1} {\xi^2_5}^2 {\xi^2_3} +
{\xi^5_5} {\xi^2_5} {\xi^2_4} {\xi^2_3} {\xi^2_1} + {\xi^4_5} {\xi^3_1} {\xi^2_5} {\xi^2_4} {\xi^2_3}
+ {\xi^4_5} {\xi^2_4}^2 {\xi^2_3} {\xi^2_1} - {\xi^4_3} {\xi^3_1} {\xi^2_5}^2 {\xi^2_4} - {\xi^4_3}
 {\xi^2_5} {\xi^2_4}^2 {\xi^2_1} + {\xi^2_5}^2 {\xi^2_1})/({\xi^2_5}^2 ({\xi^3_1} {\xi^2_5}
+ {\xi^2_4} {\xi^2_1}));$
  \quad
$  \mathbf{J^6_1}=({\xi^6_6}^2 {\xi^6_5} {\xi^2_5}^2 {\xi^2_1} - 2 {\xi^6_6} {\xi^6_5} {\xi^3_1} {\xi^2_5}^2
 {\xi^2_3} - 2 {\xi^6_6} {\xi^6_5} {\xi^2_5} {\xi^2_4} {\xi^2_3} {\xi^2_1} + 2 {\xi^6_6}
 {\xi^6_3} {\xi^3_1} {\xi^2_5}^3 + 2 {\xi^6_6} {\xi^6_3} {\xi^2_5}^2 {\xi^2_4} {\xi^2_1} + {\xi^6_5}
 {\xi^2_5}^2 {\xi^2_1} + {\xi^4_5} {\xi^3_1}^2 {\xi^2_5}^2 {\xi^2_3} + 2 {\xi^4_5} {\xi^3_1}
 {\xi^2_5} {\xi^2_4} {\xi^2_3} {\xi^2_1} + {\xi^4_5} {\xi^2_4}^2 {\xi^2_3} {\xi^2_1}^2 - {\xi^4_3}
 {\xi^3_1}^2 {\xi^2_5}^3 - 2 {\xi^4_3} {\xi^3_1} {\xi^2_5}^2 {\xi^2_4} {\xi^2_1} - {\xi^4_3}
 {\xi^2_5} {\xi^2_4}^2 {\xi^2_1}^2)/({\xi^2_5}^3 ({\xi^6_6}^2 + 1));$
  \quad
$  \mathbf{J^6_2}=( - {\xi^6_6} {\xi^6_5} {\xi^2_5} - {\xi^6_5} {\xi^5_5} {\xi^2_5} - {\xi^6_5} {\xi^4_5}
 {\xi^2_4} + {\xi^6_3} {\xi^2_5} {\xi^2_4} + {\xi^4_5} {\xi^3_1} {\xi^2_5} + {\xi^4_5} {\xi^2_4}
 {\xi^2_1})/{\xi^2_5}^2;$
}
and the parameters are subject to the condition
\begin{equation}
\label{cond63_2general}
\xi^2_5( {\xi^2_5} {\xi^3_1} + {\xi^2_4} {\xi^2_1}) \neq 0 .
 \end{equation}
Now, equivalence by a suitable automorphism of the form
(\ref{63specialautom1})
reduces to the case
$ {\xi^2_1}={\xi^2_3}={\xi^4_5}={\xi^5_5}={\xi^6_5}={\xi^4_3}={\xi^6_3}= 0 $.
Applying then equivalence by
$
{\fontsize{8}{10}\selectfont
\Psi = \text{diag} \left( \begin{pmatrix}
 - {\xi^6_6}/{\xi^2_5}&
 - {\xi^6_6} {\xi^2_4}/({\xi^3_1} {\xi^2_5})&
 - 1/{\xi^3_1}&
\\
0&
1&
0\\
{\xi^3_1}/{\xi^2_5}&
0&
0
\end{pmatrix}, \begin{pmatrix}
 - {\xi^6_6}/{\xi^2_5}&
0&
1/{\xi^3_1}\\
{\xi^6_6} {\xi^2_4}/{\xi^2_5}^2&
1/{\xi^2_5}&
0\\
 - {\xi^3_1}/{\xi^2_5}&
0&
0\end{pmatrix}\right)
}
,$
$ J2 = \Psi^{-1}J\Psi $ is
$$
{\fontsize{8}{10}\selectfont
J2 = \begin{pmatrix}
0&
0&
 - {\xi^2_5}/{\xi^3_1}&
0&
 - {\xi^2_4}/{\xi^3_1}&
 - {\xi^2_4}^2/{\xi^3_1}^2\\
0&
0&
0&
0&
1&
{\xi^2_4}/{\xi^3_1}\\
{\xi^3_1}/{\xi^2_5}&
{\xi^2_4}/{\xi^2_5}&
0&
0&
0&
0\\
0&
0&
0&
0&
0&
{\xi^2_5}/{\xi^3_1}\\
0&
-1&
0&
{\xi^2_4}/{\xi^2_5}&
0&
0\\
0&
0&
0&
 - {\xi^3_1}/{\xi^2_5}&
0&
0\end{pmatrix}\, .
}
$$
Suppose first that $ {\xi^2_4} \neq 0 $.
Then this $ J2 $ is a CS
belonging in the case \ref{case31}.
Suppose now that $ {\xi^2_4} = 0 $.
Applying equivalence by the automorphism
$ \Lambda = \text{diag}(1,{\xi^3_1}/{\xi^2_5},{\xi^3_1}/{\xi^2_5},{\xi^3_1}/{\xi^2_5}
,{\xi^3_1}/{\xi^2_5},{\xi^3_1}^2/{\xi^2_5}^2),$
%$$ \Lambda = \begin{pmatrix}1&0&0&0&0&0\\0&{\xi^3_1}/{\xi^2_5}&0&0&0&0\\0&0&{\xi^3_1}/{\xi^2_5}&0&0&0\\0&0&0&{\xi^3_1}/{\xi^2_5}&0&0\\0&0&0&0&{\xi^3_1}/{\xi^2_5}&0\\0&0&0&0&0&{\xi^3_1}^2/{\xi^2_5}^2\end{pmatrix}$$
$ \Lambda^{-1} J2 \Lambda$ is the matrix
\begin{equation}
{\fontsize{8}{10}\selectfont
\label{63J_2}
J_2= \begin{pmatrix}
0&
0&
-1&
0&
0&
0\\
0&
0&
0&
0&
1&
0\\
1&
0&
0&
0&
0&
0\\
0&
0&
0&
0&
0&
1\\
0&
-1&
0&
0&
0&
0\\
0&
0&
0&
-1&
0&
0\end{pmatrix}
}
\end{equation}
Hence, from the result of the case \ref{case31}, any  CS with $ {\xi^1_6} = 0 , \xi^2_5 \neq 0$ is equivalent
to  either $J_1$ in
(\ref{53-J_1})
or  $J_2$ in
(\ref{63J_2}).
Since $ J_2 $ is equivalent to $ J_1 $ by the automorphism
$ M = \text{diag} \left(
\left( \begin{smallmatrix} 0&1\\-1&0\end{smallmatrix}\right),
-1,0,0,1,\left(
\begin{smallmatrix} 0&-1\\1&0\end{smallmatrix}\right) \right)$
\textit{i.e. } $ M^{-1}J_2 M = J_1 $,
we get that  any  CS with
$ {\xi^1_6} =0, {\xi^2_5} \neq 0 $
is equivalent to  $ J_1. $
%%%%%%%%%%%%%%%%%%%%%%%%%%%%%%%%%%%%%%%%%%%%%%%%%%%%%%%%%%%%%%%%%%%%%%%%%%%%%%%%%%%%%%%%%
\subsection{Case $ {\xi^1_6} = 0, {\xi^2_5} = 0.$}
%{\fontsize{8}{10}\selectfont
\begin{equation}
\label{63general3}
J = \begin{pmatrix}
\boxed{\xi^1_1}&
\boxed{\xi^1_2}&
0&
0&
0&
0\\
- \frac{{\xi^1_1}^2 + 1}{\xi^1_2}&
- {\xi^1_1}&
0&
0&
0&
0\\
\boxed{\xi^3_1}&
%$  J^3_2=({\xi^1_2} ({\xi^4_1} {\xi^3_4} + {\xi^3_3} {\xi^3_1} + {\xi^3_1} {\xi^1_1}))/({\xi^1_1}^2 + 1);$\\
\frac{{\xi^1_2} ({\xi^4_1} {\xi^3_4} + {\xi^3_3} {\xi^3_1} + {\xi^3_1} {\xi^1_1})}
{{\xi^1_1}^2 + 1}&
\boxed{\xi^3_3}&
\boxed{\xi^3_4}&
0&
0\\
\boxed{\xi^4_1}&
%$  J^4_2=({\xi^1_2} ( - {\xi^4_1} {\xi^3_4} {\xi^3_3} + {\xi^4_1} {\xi^3_4} {\xi^1_1} - {\xi^3_3}^2 {\xi^3_1} - {\xi^3_1}))/({\xi^3_4} ({\xi^1_1}^2 + 1));$\\
\frac{{\xi^1_2} ( - {\xi^4_1} {\xi^3_4} {\xi^3_3} + {\xi^4_1} {\xi^3_4} {\xi^1_1} - {\xi^3_3}^2 {\xi^3_1} - {\xi^3_1})}{{\xi^3_4} ({\xi^1_1}^2 + 1)}&
%*&
%$  J^4_3= - \frac{{\xi^3_3}^2 + 1}{\xi^3_4};$\\
- \frac{{\xi^3_3}^2 + 1}{\xi^3_4}&
%$  J^4_4= - {\xi^3_3};$\\
- {\xi^3_3}&
0&
0\\
%$  J^5_1=({\xi^6_4} {\xi^4_1} {\xi^1_2} + {\xi^6_3} {\xi^3_1} {\xi^1_2} - {\xi^6_2} {\xi^1_1}^2 - {\xi^6_2})/({\xi^1_1}^2 + 1);$\\
\frac{{\xi^6_4} {\xi^4_1} {\xi^1_2} + {\xi^6_3} {\xi^3_1} {\xi^1_2} - {\xi^6_2} {\xi^1_1}^2 - {\xi^6_2}}{{\xi^1_1}^2 + 1}&
*&
*&
%$  J^5_4=({\xi^1_2} ( - {\xi^6_4} {\xi^3_3} - {\xi^6_4} {\xi^1_1} + {\xi^6_3} {\xi^3_4}))/%({\xi^1_1}^2 + 1);$\\
\frac{{\xi^1_2} ( - {\xi^6_4} {\xi^3_3} - {\xi^6_4} {\xi^1_1} + {\xi^6_3} {\xi^3_4})}{
{\xi^1_1}^2 + 1}&
%$  J^5_5={\xi^1_1};$\\
{\xi^1_1}&
%$  J^5_6={\xi^1_2};$\\
{\xi^1_2}\\
\boxed{\xi^6_1}&
\boxed{\xi^6_2}&
\boxed{\xi^6_3}&
\boxed{\xi^6_4}&
%$  J^6_5= - ({\xi^1_1}^2 + 1)/{\xi^1_2};$\\
- \frac{{\xi^1_1}^2 + 1}{\xi^1_2}&
%$  J^6_6= - {\xi^1_1};$\\
- {\xi^1_1}
 \end{pmatrix}
\end{equation}
%}
where
{\fontsize{6}{10}\selectfont
$  \mathbf{J^5_2}=({\xi^1_2} ( - {\xi^6_4} {\xi^4_1} {\xi^3_4} {\xi^3_3} {\xi^1_2} + {\xi^6_4} {\xi^4_1} {\xi^3_4} {\xi^1_2} {\xi^1_1} - {\xi^6_4} {\xi^3_3}^2 {\xi^3_1} {\xi^1_2} - {\xi^6_4} {\xi^3_1}
 {\xi^1_2} + {\xi^6_3} {\xi^4_1} {\xi^3_4}^2 {\xi^1_2} + {\xi^6_3} {\xi^3_4} {\xi^3_3} {\xi^3_1}
{\xi^1_2} + {\xi^6_3} {\xi^3_4} {\xi^3_1} {\xi^1_2} {\xi^1_1} - 2 {\xi^6_2} {\xi^3_4} {\xi^1_1}^3
 - 2 {\xi^6_2} {\xi^3_4} {\xi^1_1} + {\xi^6_1} {\xi^3_4} {\xi^1_2} {\xi^1_1}^2 + {\xi^6_1} {\xi^3_4} {\xi^1_2}))/({\xi^3_4} ({\xi^1_1}^2 + 1)^2);$\quad
$  \mathbf{J^5_3}=({\xi^1_2} ( - {\xi^6_4} {\xi^3_3}^2 - {\xi^6_4} + {\xi^6_3} {\xi^3_4} {\xi^3_3} -
 {\xi^6_3} {\xi^3_4} {\xi^1_1}))/({\xi^3_4} ({\xi^1_1}^2 + 1));$
}
and the parameters are subject to the condition
\begin{equation}
\label{cond63_3general}
 {\xi^1_2} {\xi^3_4}  \neq 0 .
 \end{equation}
Now, equivalence by a suitable automorphism of the form
(\ref{63specialautom1})
reduces to the case
$ {\xi^3_1}={\xi^3_3}={\xi^4_1}={\xi^6_1}={\xi^6_2}={\xi^6_3}={\xi^6_4}=  0. $
Applying then equivalence by the automorphism
${\fontsize{8}{10}\selectfont
 \Psi = \text{diag} \left( \begin{pmatrix}
{\xi^1_2}&0&0\\- {\xi^1_1}&1&0\\0&0&{\xi^3_4} {\xi^1_2}\end{pmatrix},
\begin{pmatrix} {\xi^1_2}&0&0\\0&{\xi^3_4} {\xi^1_2}^2&0\\
0& - {\xi^3_4} {\xi^1_2} {\xi^1_1}&{\xi^3_4} {\xi^1_2}\end{pmatrix}\right)
}
,$
we get into the case of a CS $ J $ where
all parameters but $  {\xi^1_2} = {\xi^3_4} = 1  $ vanish
and $  {\xi^1_2} = {\xi^3_4} = 1 ,$ that is
$ J = \text{diag} \left(
 \left(\begin{smallmatrix} 0&1\\-1&0\end{smallmatrix} \right),
 \left(\begin{smallmatrix} 0&1\\-1&0\end{smallmatrix} \right),
 \left(\begin{smallmatrix} 0&1\\-1&0\end{smallmatrix} \right)
 \right)
 ,$
which is equivalent to
the matrix $J_1 $ in (\ref{53-J_1}).
Now, with the automorphism
$\text{diag}(1,-1,1,-1,1,-1),$
$J$ is equivalent to its opposite
\begin{equation}
\label{63equivfinal}
 J_0
 = \text{diag} \left(
 \left(\begin{smallmatrix} 0&-1\\1&0\end{smallmatrix} \right),
 \left(\begin{smallmatrix} 0&-1\\1&0\end{smallmatrix} \right),
 \left(\begin{smallmatrix} 0&-1\\1&0\end{smallmatrix} \right)
 \right)\, .
\end{equation}
\par Hence, any  CS with
$ {\xi^1_6} =0, {\xi^2_5} = 0 $
is equivalent to $ J_0 $ in
(\ref{63equivfinal}).
%%%%%%%%%%%%%%%%%%%%%%%%%%%%%%%%%%%%%%%%%%%%%%%%%%%%%%%%%%%%%%%%%%%%%%%%%%%%%%%%%%%%%%%%%%
Commutation relations of  $ \mathfrak{m}$ for $J_0:$
$[\tilde{x}_1,\tilde{x}_3]=\tilde{x}_5$;
$[\tilde{x}_1,\tilde{x}_4]=\tilde{x}_6$;
$[\tilde{x}_2,\tilde{x}_3]=\tilde{x}_6$;
$[\tilde{x}_2,\tilde{x}_4]= - \tilde{x}_5.$
%%%%%%%%%%%%%%%%%%%%%%%%%%%%%%%%%%%%%%%%%%%%%%%%%%%%%%%%%%%%%%%%%%%%%%%%%%%%%%%%%%%%%%%%%%
\subsection{Conclusions.}
Any $J \in \mathfrak{X}_{6,3}$ is equivalent to  $J_0$ defined by (\ref{63equivfinal}).
Hence  $\mathfrak{X}_{6,3}$ is comprised of the single $\text{Aut }  {\mathcal{G}}_{6,3}$ orbit of $J_0.$
Now
$\text{Aut } {\mathcal{G}}_{6,3}$ consists of the matrices
$${\fontsize{8}{10}\selectfont
 \Phi = \begin{pmatrix}
b^1_1&
b^1_2&
b^1_3&
0&
0&
0\\
b^2_1&
b^2_2&
b^2_3&
0&
0&
0\\
b^3_1&
b^3_2&
b^3_3&
0&
0&
0\\
b^4_1&
b^4_2&
b^4_3&
b^2_2 b^1_1 -b^2_1 b^1_2&
b^2_3 b^1_1 -b^2_1 b^1_3&
b^2_3 b^1_2 -b^2_2 b^1_3\\
b^5_1&
b^5_2&
b^5_3&
b^3_2 b^1_1 -b^3_1 b^1_2&
b^3_3 b^1_1 -b^3_1 b^1_3&
b^3_3 b^1_2 -b^3_2 b^1_3\\
b^6_1&
b^6_2&
b^6_3&
b^3_2 b^2_1 -b^3_1 b^2_2&
b^3_3 b^2_1 -b^3_1 b^2_3&
b^3_3 b^2_2 -b^3_2 b^2_3
\end{pmatrix}
}$$
where the $b^i_j$'s are arbitrary reals with the condition $\det \Phi \neq 0,$
and  the stabilizer of $J_0$ is  6-dimensional.
Hence
$\mathfrak{X}_{6,3}$
is a submanifold of dimension 12 of $\Rmath^{36}$ (\cite{Bourbaki}, Chap. 3, par. 1,
Prop. 14).
 We also remark that
$\mathfrak{X}_{6,3}$
 is the zero set of a polynomial map
 $F \; : \; \Rmath^{36} \rightarrow \Rmath^{81} \times \Rmath^{36};$
 however this map is not a subimmersion, that is its rank is not locally constant.
\subsection{}
%\subsubsection{}
\begin{equation*}
X_1 = \frac{\partial}{\partial x^1}- x^2\frac{\partial}{\partial x^3}- y^1\frac{\partial}{\partial y^2} \quad , \quad
X_2 = \frac{\partial}{\partial y^1} - x^2\frac{\partial}{\partial y^3} \, .
\end{equation*}
%\subsubsection{Holomorphic functions for $J_0.$}
Let $G$ denote the group $G_0$ endowed with the left invariant  structure
of complex manifold  defined by
$J_0$ (\ref{63equivfinal}).
Then
$ H_{\Cmath}(G) =\{f \in C^{\infty}(G_0) \; ; \; \tilde{X}_j^{-} \, f = 0
\; \forall j=1,3,5\}.$
One has
\begin{equation*}
\tilde{X}_1^{-}  =2\, \frac{\partial}{\partial \overline{z^1}} -2\, x^2 \frac{\partial}{\partial \overline{z^3}}
- y^1\frac{\partial}{\partial y^2} \; ; \;
\tilde{X}_3^{-}  =2\, \frac{\partial}{\partial \overline{z^2}} \; ; \;
\tilde{X}_5^{-}  =2\, \frac{\partial}{\partial \overline{z^3}}
\end{equation*}
where
$z^j = x^j +i y^j \; (1 \leqslant j \leqslant 3). $
Then $f \in C^{\infty}(G_0)$ is in $ H_{\Cmath}(G)$ if and only if it is holomorphic
with respect to $z^2$ and $z^3$ and satisfies
$$ 2\, \frac{\partial f}{\partial \overline{z^1}} =
\frac{z^1 -\overline{z^1}}{2} \, \frac{\partial f}{\partial {z^2}}
\, .
$$
Hence the 3 functions $w^1 =z^1, w^2= z^2 +\frac{|z^1|^2}{4} - \frac{(\overline{z^1})^2}{8} , w^3=z^3$ are holomorphic.
Let $F : G \rightarrow \Cmath^3$ defined
by $F=(w^1,w^2,w^3).$
$F$ is a biholomorphic bijection, hence a global chart on $G.$
We determine now how the multiplication of $G$ looks like in that chart.
Let $a,x \in G $ with respective second kind canonical coordinates
$(x^1,y^1,x^2,y^2,x^3,y^3), (\alpha^1, \beta^1, \alpha^2, \beta^2, \alpha^3, \beta^3)$
as in  (\ref{x6general}).
With obvious notations,
$a =[w^1_a, w^2_a,w^3_a],$
$x =[w^1_x, w^2_x,w^3_x],$
$a \, x =[w^1_{a x}, w^2_{a x},w^3_{a x}].$
Computations yield :
\begin{eqnarray*}
w^1_{a x} &=& w^1_a + w^1_x\\
w^2_{a x} &=& w^2_a + w^2_x - \frac{w^1_{a} - \overline{w^1_{a}}}{4}\; w^1_x\\
w^3_{a x} &=& w^3_a + w^3_x +
\frac{1}{2} (w^2_{a} - \overline{w^2_{a}}
+\frac{(\overline{w^1_{a}})^2 -(w^1_{a})^2}{4} - |w^1_{a}|^2)\;  w^1_x .
\end{eqnarray*}
%%%%%%%%%%%%%%%%%%%%%%%%%%%%%%%%%%%%%%%%%%%%%%%%%%%%%%%%%%%%%%%%%%%%%%%%%%%%%%%%%%%%%%%%%
%%%%%%%%%%%%%%%%%%%%%%%%%%%%%%%%%%%%%%%%%%%%%%%%%%%%%%%%%%%%%%%%%%%%%%%%%%%%%%%%%%%%%%%%%
\section{Lie Algebra $ {\mathcal{G}}_{6,7}$ (isomorphic to $M6$).}
%%%%%%%%%%%%%%%%%%%%%%%%%%%%%%%%%%%%%%%%%%%%%%%%%%%%%%%%%%%%%%%%%%%%%%%%%%%%%%%%%%%%%%%%%%
Commutation relations for $ {\mathcal{G}}_{6,7}:$
$[x_1,x_2]=x_4$;
$[x_1,x_3]=x_5$;
$[x_1,x_4]=x_6$;
$[x_2,x_3]= - x_6.$
%{\fontsize{8}{10}\selectfont
\begin{equation}
\label{67general}
J = \begin{pmatrix}
\boxed{\xi^1_1}&
\boxed{\xi^1_2}&
0&
0&
0&
0\\
%$  J^2_1= - ({\xi^1_1}^2 + 1)/{\xi^1_2};$\\
- \frac{{\xi^1_1}^2 + 1}{\xi^1_2}&
%$  J^2_2= - {\xi^1_1};$\\
- {\xi^1_1}&
0&
0&
0&
0\\
*&
\boxed{\xi^3_2}&
%$  J^3_3=( - {\xi^5_5} {\xi^3_4} + {\xi^5_5} {\xi^1_2} + {\xi^3_4} {\xi^1_1})/{\xi^1_2}
c&
\boxed{\xi^3_4}&
0&
0\\
*&
\boxed{\xi^4_2}&
-\frac{c^2+1}{\xi^3_4}&
-c&
0&
0\\
*&
*&
*&
\boxed{\xi^5_4}&
\boxed{\xi^5_5}&
%$  J^5_6= - ({\xi^3_4} {\xi^1_2})/({\xi^3_4} - {\xi^1_2});$\\
- \frac{{\xi^3_4} {\xi^1_2}}{{\xi^3_4} - {\xi^1_2}}\\
\boxed{\xi^6_1}&
\boxed{\xi^6_2}&
*&
\boxed{\xi^6_4}&
%$  J^6_5=({\xi^5_5}^2 {\xi^3_4} - {\xi^5_5}^2 {\xi^1_2} + {\xi^3_4} - {\xi^1_2})/({\xi^3_4} {\xi^1_2});$\\
\frac{({\xi^5_5}^2+1)( {\xi^3_4} - {\xi^1_2})
}{{\xi^3_4} {\xi^1_2}}&
%$  J^6_6= - {\xi^5_5};$\\
- {\xi^5_5}
 \end{pmatrix}
\end{equation}
%}
where
{\fontsize{6}{10}\selectfont
$ \mathbf{J^3_1}=({\xi^5_5} {\xi^3_4} {\xi^3_2} - {\xi^5_5} {\xi^3_2} {\xi^1_2} - {\xi^4_2} {\xi^3_4}
{\xi^1_2} - {\xi^3_4} {\xi^3_2} {\xi^1_1} + {\xi^3_2} {\xi^1_2} {\xi^1_1})/{\xi^1_2}^2;$
 \quad
$ \mathbf{c}= \mathbf{J^3_3}=( - {\xi^5_5} {\xi^3_4} + {\xi^5_5} {\xi^1_2} + {\xi^3_4} {\xi^1_1})/{\xi^1_2}
;$
 \quad
$  \mathbf{J^4_1}=({\xi^5_5}^2 {\xi^3_4}^2 {\xi^3_2} - 2 {\xi^5_5}^2 {\xi^3_4} {\xi^3_2} {\xi^1_2}
 + {\xi^5_5}^2 {\xi^3_2} {\xi^1_2}^2 - {\xi^5_5} {\xi^4_2} {\xi^3_4}^2 {\xi^1_2} + {\xi^5_5}
{\xi^4_2} {\xi^3_4} {\xi^1_2}^2 - 2 {\xi^5_5} {\xi^3_4}^2 {\xi^3_2} {\xi^1_1} + 2 {\xi^5_5} {\xi^3_4} {\xi^3_2} {\xi^1_2} {\xi^1_1} + {\xi^4_2} {\xi^3_4}^2 {\xi^1_2} {\xi^1_1} + {\xi^4_2} {\xi^3_4} {\xi^1_2}^2 {\xi^1_1} + {\xi^3_4}^2 {\xi^3_2} {\xi^1_1}^2 + {\xi^3_2} {\xi^1_2}^2)/(
{\xi^3_4} {\xi^1_2}^3);$
 \quad
$  \mathbf{J^5_1}=( - {\xi^6_4} {\xi^5_5}^2 {\xi^3_4}^2 {\xi^3_2} {\xi^1_2} + {\xi^6_4} {\xi^5_5}^
2 {\xi^3_4} {\xi^3_2} {\xi^1_2}^2 + {\xi^6_4} {\xi^5_5} {\xi^4_2} {\xi^3_4}^2 {\xi^1_2}^2 + 2
 {\xi^6_4} {\xi^5_5} {\xi^3_4}^2 {\xi^3_2} {\xi^1_2} {\xi^1_1} - 2 {\xi^6_4} {\xi^5_5} {\xi^3_4}
{\xi^3_2} {\xi^1_2}^2 {\xi^1_1} - {\xi^6_4} {\xi^4_2} {\xi^3_4}^2 {\xi^1_2}^2 {\xi^1_1} - {\xi^6_4} {\xi^3_4}^2 {\xi^3_2} {\xi^1_2} {\xi^1_1}^2 - {\xi^6_4} {\xi^3_4} {\xi^3_2} {\xi^1_2}^2
+ {\xi^6_2} {\xi^3_4}^2 {\xi^1_2}^2 {\xi^1_1}^2 + {\xi^6_2} {\xi^3_4}^2 {\xi^1_2}^2 + {\xi^6_1} {\xi^5_5} {\xi^3_4}^2 {\xi^1_2}^3 - {\xi^6_1} {\xi^3_4}^2 {\xi^1_2}^3 {\xi^1_1} + {\xi^5_5}^3 {\xi^5_4} {\xi^3_4}^2 {\xi^3_2} - 2 {\xi^5_5}^3 {\xi^5_4} {\xi^3_4} {\xi^3_2} {\xi^1_2} + {\xi^5_5}^3 {\xi^5_4} {\xi^3_2} {\xi^1_2}^2 - {\xi^5_5}^2 {\xi^5_4} {\xi^4_2} {\xi^3_4}^
2 {\xi^1_2} + {\xi^5_5}^2 {\xi^5_4} {\xi^4_2} {\xi^3_4} {\xi^1_2}^2 - {\xi^5_5}^2 {\xi^5_4}
{\xi^3_4}^2 {\xi^3_2} {\xi^1_1} + 2 {\xi^5_5}^2 {\xi^5_4} {\xi^3_4} {\xi^3_2} {\xi^1_2} {\xi^1_1} - {\xi^5_5}^2 {\xi^5_4} {\xi^3_2} {\xi^1_2}^2 {\xi^1_1} + {\xi^5_5} {\xi^5_4} {\xi^3_4}^2
{\xi^3_2} - 2 {\xi^5_5} {\xi^5_4} {\xi^3_4} {\xi^3_2} {\xi^1_2} + {\xi^5_5} {\xi^5_4} {\xi^3_2} {\xi^1_2}^2 - {\xi^5_4} {\xi^4_2} {\xi^3_4}^2 {\xi^1_2} + {\xi^5_4} {\xi^4_2} {\xi^3_4} {\xi^1_2}
^2 - {\xi^5_4} {\xi^3_4}^2 {\xi^3_2} {\xi^1_1} + 2 {\xi^5_4} {\xi^3_4} {\xi^3_2} {\xi^1_2} {\xi^1_1} - {\xi^5_4} {\xi^3_2} {\xi^1_2}^2 {\xi^1_1})/({\xi^3_4} {\xi^1_2}^2 ({\xi^5_5}^2 +1) ( {\xi^3_4} - {\xi^1_2}));$
 \quad
$  \mathbf{J^5_2}=({\xi^6_4} {\xi^5_5} {\xi^3_4}^2 {\xi^3_2} - 2 {\xi^6_4} {\xi^5_5} {\xi^3_4} {\xi^3_2} {\xi^1_2} - {\xi^6_4} {\xi^4_2} {\xi^3_4}^2 {\xi^1_2} - {\xi^6_4} {\xi^3_4}^2 {\xi^3_2} {\xi^1_1} + {\xi^6_2} {\xi^5_5} {\xi^3_4}^2 {\xi^1_2} + {\xi^6_2} {\xi^3_4}^2 {\xi^1_2} {\xi^1_1} -
 {\xi^6_1} {\xi^3_4}^2 {\xi^1_2}^2 + {\xi^5_5}^2 {\xi^5_4} {\xi^3_4} {\xi^3_2} - {\xi^5_5}^2
{\xi^5_4} {\xi^3_2} {\xi^1_2} + {\xi^5_4} {\xi^3_4} {\xi^3_2} - {\xi^5_4} {\xi^3_2} {\xi^1_2})/(
{\xi^3_4} ({\xi^5_5}^2  +1)( {\xi^3_4} - {\xi^1_2}));$
 \quad
$ \mathbf{J^5_3}=({\xi^6_4} {\xi^1_2}^2 - {\xi^5_5} {\xi^5_4} {\xi^3_4} + {\xi^5_5} {\xi^5_4} {\xi^1_2} + {\xi^5_4} {\xi^3_4} {\xi^1_1} - {\xi^5_4} {\xi^1_2} {\xi^1_1})/({\xi^1_2} ({\xi^3_4} - {\xi^1_2}));$
 \quad
$  \mathbf{J^6_3}=( - {\xi^6_4} {\xi^5_5} {\xi^3_4}^2 + 2 {\xi^6_4} {\xi^5_5} {\xi^3_4} {\xi^1_2} +
{\xi^6_4} {\xi^3_4}^2 {\xi^1_1} - {\xi^5_5}^2 {\xi^5_4} {\xi^3_4} + {\xi^5_5}^2 {\xi^5_4} {\xi^1_2} - {\xi^5_4} {\xi^3_4} + {\xi^5_4} {\xi^1_2})/({\xi^3_4}^2 {\xi^1_2});$\\
}
and the parameters are subject to the condition
\begin{equation}
\label{cond67general}
{\xi^1_2} \xi^3_4 ({\xi^1_2}- {\xi^3_4}) \neq 0 .
 \end{equation}
\par Now the automorphism group of
$ {\mathcal{G}}_{6,7}$ is comprised of the matrices
$${\fontsize{8}{10}\selectfont
 \Phi = \begin{pmatrix}
{b^1_1}&
0&
0&
0&
0&
0\\
{b^2_1}&
{b^2_2}&
0&
0&
0&
0\\
{b^3_1}&
{b^3_2}&
{b^1_1}^2&
0&
0&
0\\
{b^4_1}&
{b^4_2}&
{b^4_3}&
{b^2_2} {b^1_1}&
0&
0\\
{b^5_1}&
{b^5_2}&
{b^5_3}&
{b^3_2} {b^1_1}&
{b^1_1}^3&
0\\
{b^6_1}&
{b^6_2}&
{b^6_3}&
{b^4_2} {b^1_1} - {b^3_2} {b^2_1} + {b^3_1} {b^2_2}&
{b^1_1} ({b^4_3} - {b^2_1} {b^1_1})&
{b^2_2} {b^1_1}^2\end{pmatrix}
 }$$
where $ {b^2_2}{b^1_1} \neq 0.$
Taking suitable values for the $b^i_j$'s,
equivalence by $\Phi$
leads to the case where
$ {\xi^1_1}={\xi^3_2}={\xi^4_2}={\xi^5_4}={\xi^5_5}={\xi^6_1}={\xi^6_2}={\xi^6_4}= 0 $
and $ {\xi^1_2}=1, {\xi^3_4}=\alpha,$ where $ \alpha  \neq 0,1 $:
\begin{equation}
\label{67equivfinal}
J_{\alpha} =
\text{diag} \left(
\left( \begin{smallmatrix} 0&1\\-1&0\end{smallmatrix} \right),
\left( \begin{smallmatrix}0&\alpha\\  - 1/\alpha&0
\end{smallmatrix} \right),
\left( \begin{smallmatrix}
0&  - \alpha/(\alpha - 1)\\
(\alpha - 1)/\alpha
\end{smallmatrix} \right) \right)
\quad
(\alpha  \neq 0,1)
\end{equation}
\par Hence, any CS
is equivalent to $ J_{\alpha} $ in
(\ref{67equivfinal}).
It is easily seen that  the $J_{\alpha}$'s corresponding to different values
of $\alpha$ are not equivalent.
%%%%%%%%%%%%%%%%%%%%%%%%%%%%%%%%%%%%%%%%%%%%%%%%%%%%%%%%%%%%%%%%%%%%%%%%%%%%%%%%%%%%%%%%%%
Commutation relations of $ \mathfrak{m}$ for $J_{\alpha}:$
$
[\tilde{x}_1,\tilde{x}_3]=\tilde{x}_5 \, ; \,
[\tilde{x}_1,\tilde{x}_4]=(1-\alpha)\tilde{x}_6 \, ; \,
[\tilde{x}_2,\tilde{x}_3]=\frac{1-\alpha}{\alpha} \, \tilde{x}_6 \, ; \,
[\tilde{x}_2,\tilde{x}_4]= - \alpha \tilde{x}_5.$
%%%%%%%%%%%%%%%%%%%%%%%%%%%%%%%%%%%%%%%%%%%%%%%%%%%%%%%%%%%%%%%%%%%%%%%%%%%%%%%%%%%%%%%%%%
%\subsection{Conclusions.}
\par
From (\ref{67general}), $\mathfrak{X}_{6,7}$
is a submanifold of dimension 10 in $\Rmath^{36}.$
It is  the disjoint union of the continuously many orbits of the
$J_{\alpha}$'s in (\ref{67equivfinal}).
%%%%%%%%%%%%%%%%%%%%%%%%%%%%%%%%%%%%%%%%%%%%%%%%%%%%%%%%%%%%%%%%%%%%%%%%%%%%%%%%%%%%%%%%%
\subsection{}
%\subsubsection{}
\begin{equation*}
X_1 = \frac{\partial}{\partial x^1}- y^1\frac{\partial}{\partial y^2}
- x^2\frac{\partial}{\partial x^3} - y^2\frac{\partial}{\partial y^3}
\quad , \quad
X_2 = \frac{\partial}{\partial y^1} + x^2\frac{\partial}{\partial y^3}.
\end{equation*}
%\subsubsection{Holomorphic functions for $J_{\alpha}.$}
Let $G$ denote the group $G_0$ endowed with the left invariant  structure
of complex manifold  defined by
$J_{\alpha}$ in (\ref{67equivfinal}).
Then $ H_{\Cmath}(G) =\{f \in C^{\infty}(G_0) \; ; \; \tilde{X}_j^{-} \, f = 0
\; \forall j=1,3,5\}.$
%%%%%%%%%%%%%%%%%%%%%%%%%%%%%%%%%%%%%%%%%%%%%%%%%%%%%%%%%%%%%%%%%%%%%%%%%%%%%%%%%%%%%
One has
\begin{equation*}
\tilde{X}_1^{-}  =2\, \frac{\partial}{\partial \overline{w^1}}
-y^1\,  \frac{\partial}{\partial {y^2}}
-x^2\,  \frac{\partial}{\partial {x^3}}
-(y^2+ix^2)\,  \frac{\partial}{\partial {y^3}}
\quad , \quad
\tilde{X}_3^{-}  =2\, \frac{\partial}{\partial \overline{w^2}}
\quad , \quad
\tilde{X}_5^{-}  =2\, \frac{\partial}{\partial \overline{w^3}}
\quad ,
\end{equation*}
where
%begin{equation*}
$
w^1 = x^1 -i y^1
\quad , \quad
w^2 = x^2 -i\alpha y^2
\quad , \quad
w^3 = x^2 +i\frac{\alpha}{\alpha-1}\,  y^3.
$
%\end{equation*}
%%%%%%%%%%%%%%%%%%%%%%%%%%%%%%%%%%%%%%%%%%%%%%%%%%%%%%%%%%%%%%%%%%%%%%%%%%%%%%%%%%%%%
Then $f \in C^{\infty}(G_0)$ is in $ H_{\Cmath}(G)$ if and only if it is holomorphic
with respect to $w^2$ and $w^3$ and satisfies
$$2\, \frac{\partial f}{\partial \overline{w^1}} -\alpha \, \frac{w^1-\overline{w^1}}{2} \,  \frac{\partial f}{\partial {w^2}}
+\frac{w^2}{\alpha-1} \,\frac{\partial f}{\partial w^3}=0.$$
The 3 functions
\begin{eqnarray*}
\varphi^1 &=&w^1 \quad , \quad
\varphi^2 =w^2 +\alpha \, \left(-\frac{(\overline{w^1})^2}{8} + \frac{|w^1|^2}{4}\right)\quad ,\\
\varphi^3 &=&w^3 +\frac{\alpha}{8(1-\alpha)}  \, (\overline{w^1})^2
\left(\frac{w^1}{2} -\frac{\overline{w^1}}{3} \right)
+ \frac{\overline{w^1} w^2}{2(1-\alpha)}
\end{eqnarray*}
are holomorphic.
Let $F : G \rightarrow \Cmath^3$ defined
by
$F=(\varphi^1,\varphi^2,\varphi^3).$
$F$ is a biholomorphic bijection, hence a global chart on $G.$
We determine now how the multiplication of $G$ looks like in that chart.
Let $a,x \in G $ with respective second kind canonical coordinates
$(x^1,y^1,x^2,y^2,x^3,y^3), (\alpha^1, \beta^1, \alpha^2, \beta^2, \alpha^3, \beta^3)$
as in  (\ref{x6general}).
With obvious notations,
$a =[w^1_a, w^2_a,w^3_a],$
$x =[w^1_x, w^2_x,w^3_x],$
$a \, x =[w^1_{a x}, w^2_{a x},w^3_{a x}],$
$a =[\varphi^1_a, \varphi^2_a,\varphi^3_a],$
$x =[\varphi^1_x, \varphi^2_x,\varphi^3_x],$
$a \, x =[\varphi^1_{a x}, \varphi^2_{a x},\varphi^3_{a x}].$
Computations yield :
\begin{eqnarray*}
w^1_{a x} &=& w^1_a + w^1_x \quad , \quad
w^2_{a x} = w^2_a + w^2_x + i \alpha \beta^1 x^1 \quad ,\\
w^3_{a x} &=& w^3_a + w^3_x - \alpha^2x^1+ i \frac{\alpha}{\alpha-1}\, \left(-\beta^1 x^1
+\alpha^2 y^1 +\frac{1}{2}\, \beta^1 (x^1)^2\right).
\end{eqnarray*}
We then get
\begin{equation*}
\varphi^1_{a x} = \varphi^1_a + \varphi^1_x
\quad , \quad
\varphi^2_{a x} = \varphi^2_a + \varphi^2_x
+ \frac{\alpha}{4} \, \left(2 \overline{\varphi^1_a} - \varphi^1_a\right)\, \varphi^1_x
\quad , \quad
\varphi^3_{a x} = \varphi^3_a +\varphi^3_x + \chi(a,x)
\end{equation*}
where
\begin{multline*}
\chi(a,x) = \varphi^1_x \left( -\frac{1}{2}\, \varphi^2_a
+  \frac{\alpha}{2(1-\alpha)}\, \overline{\varphi^2_a}
+\frac{\alpha(2+\alpha)}{16(1-\alpha)}\,(\overline{\varphi^1_a})^2
+\frac{\alpha^2}{16(1-\alpha)}\,({\varphi^1_a})^2
-\frac{\alpha^2}{4(1-\alpha)}\,|\varphi^1_a|^2\right)
\\
+\frac{\alpha}{16(1-\alpha)}\, \left( {\varphi^1_a} -\overline{\varphi^1_a} \right)
\,(\varphi^1_x)^2
+\frac{1}{2(1-\alpha)}\, \overline{\varphi^1_a} \,\varphi^2_x .
\end{multline*}
%%%%%%%%%%%%%%%%%%%%%%%%%%%%%%%%%%%%%%%%%%%%%%%%%%%%%%%%%%%%%%%%%%%%%%%%%%%%%%%%%%%%%%%%%
\section{Lie Algebra $ {\mathcal{G}}_{6,4}$ (isomorphic to $M7$).}
%%%%%%%%%%%%%%%%%%%%%%%%%%%%%%%%%%%%%%%%%%%%%%%%%%%%%%%%%%%%%%%%%%%%%%%%%%%%%%%%%%%%%%%%%%
Commutation relations for
$ {\mathcal{G}}_{6,4}:$
$[x_1,x_2]=x_4$;
$[x_1,x_3]=x_6$;
$[x_2,x_4]=x_5.$
%{\fontsize{8}{10}\selectfont
\begin{equation}
\label{64general}
J = \begin{pmatrix}
\boxed{\xi^1_1}&
\boxed{\xi^1_2}&
0&
0&
0&
0\\
%$  J^2_1= - ({\xi^1_1}^2 + 1)/{\xi^1_2};$\\
- \frac{({\xi^1_1}^2 + 1)}{\xi^1_2}&
-{\xi^1_1}&
0&
0&
0&
0\\
*&
\boxed{\xi^3_2}&
b&
\boxed{\xi^3_4}&
0&
0\\
*&
\boxed{\xi^4_2}&
-\frac{b^2+1}{\xi^3_4}&
-b&
0&
0\\
*&
*&
*&
*&
\boxed{\xi^5_5}&
-\frac{{\xi^5_5}^2+1}{c}\\
\boxed{\xi^6_1}&
\boxed{\xi^6_2}&
\boxed{\xi^6_3}&
\boxed{\xi^6_4}&
c&
%$  J^6_6= - {\xi^5_5};$\\
- {\xi^5_5}
 \end{pmatrix}
\end{equation}
%}
where
{\fontsize{6}{10}\selectfont
$  \mathbf{J^3_1}=( - {\xi^5_5} {\xi^4_2} {\xi^3_4} + 2 {\xi^5_5} {\xi^3_2} {\xi^1_1} - {\xi^4_2} {\xi^3_4} {\xi^1_1} + {\xi^3_2} {\xi^1_1}^2 - {\xi^3_2})/({\xi^1_2} ({\xi^5_5} + {\xi^1_1}));$
\quad
$\mathbf{b}=  \mathbf{J^3_3}=( - {\xi^5_5} {\xi^1_1} + 1)/({\xi^5_5} + {\xi^1_1});$
\quad
$  \mathbf{J^4_1}=({\xi^5_5}^2 {\xi^3_2} {\xi^1_1}^2 + {\xi^5_5}^2 {\xi^3_2} + {\xi^5_5} {\xi^4_2}
 {\xi^3_4} {\xi^1_1}^2 + {\xi^5_5} {\xi^4_2} {\xi^3_4} + {\xi^4_2} {\xi^3_4} {\xi^1_1}^3 + {\xi^4_2} {\xi^3_4} {\xi^1_1} + {\xi^3_2} {\xi^1_1}^2 + {\xi^3_2})
  /({\xi^3_4} {\xi^1_2} ({\xi^5_5}+ {\xi^1_1})^2);$
\quad
$  \mathbf{J^5_1}=({\xi^6_4} {\xi^5_5}^2 {\xi^3_2} {\xi^1_1}^4 + 2 {\xi^6_4} {\xi^5_5}^2 {\xi^3_2}
 {\xi^1_1}^2 + {\xi^6_4} {\xi^5_5}^2 {\xi^3_2} + {\xi^6_4} {\xi^5_5} {\xi^4_2} {\xi^3_4} {\xi^1_1}^4 + 2 {\xi^6_4} {\xi^5_5} {\xi^4_2} {\xi^3_4} {\xi^1_1}^2 + {\xi^6_4} {\xi^5_5} {\xi^4_2}
{\xi^3_4} + {\xi^6_4} {\xi^4_2} {\xi^3_4} {\xi^1_1}^5 + 2 {\xi^6_4} {\xi^4_2} {\xi^3_4} {\xi^1_1}
^3 + {\xi^6_4} {\xi^4_2} {\xi^3_4} {\xi^1_1} + {\xi^6_4} {\xi^3_2} {\xi^1_1}^4 + 2 {\xi^6_4}
{\xi^3_2} {\xi^1_1}^2 + {\xi^6_4} {\xi^3_2} - {\xi^6_3} {\xi^5_5}^2 {\xi^4_2} {\xi^3_4}^2 {\xi^1_1}^2 - {\xi^6_3} {\xi^5_5}^2 {\xi^4_2} {\xi^3_4}^2 + 2 {\xi^6_3} {\xi^5_5}^2 {\xi^3_4}
{\xi^3_2} {\xi^1_1}^3 + 2 {\xi^6_3} {\xi^5_5}^2 {\xi^3_4} {\xi^3_2} {\xi^1_1} - 2 {\xi^6_3} {\xi^5_5} {\xi^4_2} {\xi^3_4}^2 {\xi^1_1}^3 - 2 {\xi^6_3} {\xi^5_5} {\xi^4_2} {\xi^3_4}^2 {\xi^1_1} + 3 {\xi^6_3} {\xi^5_5} {\xi^3_4} {\xi^3_2} {\xi^1_1}^4 + 2 {\xi^6_3} {\xi^5_5} {\xi^3_4} {\xi^3_2} {\xi^1_1}^2 - {\xi^6_3} {\xi^5_5} {\xi^3_4} {\xi^3_2} - {\xi^6_3} {\xi^4_2} {\xi^3_4}^2
{\xi^1_1}^4 - {\xi^6_3} {\xi^4_2} {\xi^3_4}^2 {\xi^1_1}^2 + {\xi^6_3} {\xi^3_4} {\xi^3_2} {\xi^1_1}^5 - {\xi^6_3} {\xi^3_4} {\xi^3_2} {\xi^1_1} - {\xi^6_2} {\xi^5_5}^2 {\xi^3_4} {\xi^1_1}^
4 - 2 {\xi^6_2} {\xi^5_5}^2 {\xi^3_4} {\xi^1_1}^2 - {\xi^6_2} {\xi^5_5}^2 {\xi^3_4} - 2 {\xi^6_2} {\xi^5_5} {\xi^3_4} {\xi^1_1}^5 - 4 {\xi^6_2} {\xi^5_5} {\xi^3_4} {\xi^1_1}^3 - 2 {\xi^6_2} {\xi^5_5} {\xi^3_4} {\xi^1_1} - {\xi^6_2} {\xi^3_4} {\xi^1_1}^6 - 2 {\xi^6_2} {\xi^3_4} {\xi^1_1}^4 - {\xi^6_2} {\xi^3_4} {\xi^1_1}^2 - {\xi^6_1} {\xi^5_5}^3 {\xi^3_4} {\xi^1_2} {\xi^1_1}
^2 - {\xi^6_1} {\xi^5_5}^3 {\xi^3_4} {\xi^1_2} - {\xi^6_1} {\xi^5_5}^2 {\xi^3_4} {\xi^1_2} {\xi^1_1}^3 - {\xi^6_1} {\xi^5_5}^2 {\xi^3_4} {\xi^1_2} {\xi^1_1} + {\xi^6_1} {\xi^5_5} {\xi^3_4}
{\xi^1_2} {\xi^1_1}^4 + {\xi^6_1} {\xi^5_5} {\xi^3_4} {\xi^1_2} {\xi^1_1}^2 + {\xi^6_1} {\xi^3_4} {\xi^1_2} {\xi^1_1}^5 + {\xi^6_1} {\xi^3_4} {\xi^1_2} {\xi^1_1}^3)
 /({\xi^3_4}^2 {\xi^1_2}^2 ({\xi^5_5}+ {\xi^1_1})^3);$
\quad
$  \mathbf{J^5_2}=({\xi^6_4} {\xi^4_2} {\xi^1_1}^2 + {\xi^6_4} {\xi^4_2} + {\xi^6_3} {\xi^3_2} {\xi^1_1}^2 + {\xi^6_3} {\xi^3_2} - {\xi^6_2} {\xi^5_5} {\xi^1_1}^2 - {\xi^6_2} {\xi^5_5} - {\xi^6_2}
 {\xi^1_1}^3 - {\xi^6_2} {\xi^1_1} + {\xi^6_1} {\xi^1_2} {\xi^1_1}^2 + {\xi^6_1} {\xi^1_2})/(
{\xi^3_4} {\xi^1_2} ({\xi^5_5} + {\xi^1_1}));$
\quad
$  \mathbf{J^5_3}=( - {\xi^6_4} {\xi^5_5}^2 {\xi^1_1}^4 - 2 {\xi^6_4} {\xi^5_5}^2 {\xi^1_1}^2 -
 {\xi^6_4} {\xi^5_5}^2 - {\xi^6_4} {\xi^1_1}^4 - 2 {\xi^6_4} {\xi^1_1}^2 - {\xi^6_4} - {\xi^6_3} {\xi^5_5}^3 {\xi^3_4} {\xi^1_1}^2 - {\xi^6_3} {\xi^5_5}^3 {\xi^3_4} - 3 {\xi^6_3} {\xi^5_5}^2 {\xi^3_4} {\xi^1_1}^3 - 3 {\xi^6_3} {\xi^5_5}^2 {\xi^3_4} {\xi^1_1} - 2 {\xi^6_3} {\xi^5_5} {\xi^3_4} {\xi^1_1}^4 - {\xi^6_3} {\xi^5_5} {\xi^3_4} {\xi^1_1}^2 + {\xi^6_3} {\xi^5_5} {\xi^3_4} + {\xi^6_3} {\xi^3_4} {\xi^1_1}^3 + {\xi^6_3} {\xi^3_4} {\xi^1_1})
  /({\xi^3_4}^2 {\xi^1_2} ({\xi^5_5}  + {\xi^1_1})^3);$
\quad
$  \mathbf{J^5_4}=( - {\xi^6_4} {\xi^5_5}^2 {\xi^1_1}^2 - {\xi^6_4} {\xi^5_5}^2 - {\xi^6_4} {\xi^1_1}^2 - {\xi^6_4} + {\xi^6_3} {\xi^5_5} {\xi^3_4} {\xi^1_1}^2 + {\xi^6_3} {\xi^5_5} {\xi^3_4} + {\xi^6_3} {\xi^3_4} {\xi^1_1}^3 + {\xi^6_3} {\xi^3_4} {\xi^1_1})/({\xi^3_4} {\xi^1_2} ({\xi^5_5}+ {\xi^1_1})^2);$
\quad
}
and the parameters are subject to the condition
\begin{equation}
\label{cond64general}
{\xi^1_2} \xi^3_4 ({\xi^1_1}+ {\xi^5_5}) \neq 0 .
 \end{equation}
\par Now the automorphism group of
$ {\mathcal{G}}_{6,4}$ is comprised of the matrices
$${\fontsize{8}{10}\selectfont
 \Phi = \begin{pmatrix}
{b^1_1}&
0&
0&
0&
0&
0\\
0&
{b^2_2}&
0&
0&
0&
0\\
{b^3_1}&
{b^3_2}&
{b^3_3}&
0&
0&
0\\
{b^4_1}&
{b^4_2}&
0&
{b^2_2} {b^1_1}&
0&
0\\
{b^5_1}&
{b^5_2}&
{b^5_3}&
 - {b^4_1} {b^2_2}&
{b^2_2}^2 {b^1_1}&
0\\
{b^6_1}&
{b^6_2}&
{b^6_3}&
{b^3_2} {b^1_1}&
0&
{b^3_3} {b^1_1}\end{pmatrix}
 }$$
where
$ {b^3_3} {b^2_2} {b^1_1} \neq 0.$
Taking suitable values for the $b^i_j$'s,
equivalence by $\Phi$
leads to the case where
${\xi^3_2}={\xi^4_2}={\xi^6_2}={\xi^6_1}={\xi^6_3}={\xi^6_4}= 0, $
$ {\xi^1_2}={\xi^3_4}=1$
and ${\xi^1_1}=\alpha,{\xi^5_5} = \beta , \alpha \neq -\beta :$
\begin{equation}
\label{64equivfinal}
 J_{\alpha,\beta} =
 \text{diag} \left(
 \begin{pmatrix}
\alpha&1\\
 - (\alpha^2 + 1)& - \alpha
 \end{pmatrix} ,
 \begin{pmatrix}
\frac{ - \alpha \beta + 1}{\alpha + \beta}&
1\\
 - \frac{(\alpha^2 +1)(\beta^2 + 1)}{(\alpha +\beta)^2}&
\frac{\alpha \beta - 1}{\alpha + \beta}&
 \end{pmatrix} ,
 \begin{pmatrix}
\beta&
\frac{(\alpha^2 + 1)(  \beta^2 + 1)}{\alpha + \beta}\\
 - \frac{\alpha + \beta}{\alpha^2 + 1}&
 - \beta\end{pmatrix}
 \right)
 \quad (\alpha \neq -\beta)
 \end{equation}
\par Hence, any  CS
is equivalent to  $ J_{\alpha,\beta} $ in
(\ref{64equivfinal}).
The $J_{\alpha,\beta}$'s corresponding to different couples $(\alpha, \beta)$ are not equivalent.
%%%%%%%%%%%%%%%%%%%%%%%%%%%%%%%%%%%%%%%%%%%%%%%%%%%%%%%%%%%%%%%%%%%%%%%%%%%%%%%%%%%%%%%%%%
Commutation relations of  $ \mathfrak{m}$  for  $J_{\alpha,\beta}:$
$[\tilde{x}_1,\tilde{x}_3]= -
\frac{
(\beta^2+1)(\alpha^2+1)^2}
{(\alpha+\beta)^2}\,
\tilde{x}_5
-\frac{\beta (1+\alpha^2)}{\alpha +\beta} \,
\tilde{x}_6
$;\quad
$[\tilde{x}_1,\tilde{x}_4]=
\frac{\alpha \beta -1}{\alpha+\beta}\, (1+\alpha^2)
\tilde{x}_5
- \alpha \tilde{x}_6$;
$[\tilde{x}_2,\tilde{x}_3]= -
\frac{\alpha}{(\alpha+\beta)^2} \,(1+ \alpha^2 )(1+ \beta^2) \tilde{x}_5
+
\frac{\alpha\beta -1}{\alpha +\beta} \, \tilde{x}_6$;
\quad
$[\tilde{x}_2,\tilde{x}_4]= \beta \, \frac{\alpha^2+1}{\alpha +\beta} \, \tilde{x}_5
-  \tilde{x}_6$.
%\subsection{Conclusions.}
\par
From (\ref{64general}), $\mathfrak{X}_{6,4}$
is a submanifold of dimension 10 in $\Rmath^{36}.$
It is  the disjoint union of the continuously many orbits of the
$J_{\alpha,\beta}$'s in (\ref{64equivfinal}).
%%%%%%%%%%%%%%%%%%%%%%%%%%%%%%%%%%%%%%%%%%%%%%%%%%%%%%%%%%%%%%%%%%%%%%%%%%%%%%%%%%%%%%%%%
%%%%%%%%%%%%%%%%%%%%%%%%%%%%%%%%%%%%%%%%%%%%%%%%%%%%%%%%%%%%%%%%%%%%%%%%%%%%%%%%%%%%%%%%%%
\subsection{}
%\subsubsection{}
\begin{equation*}
X_1 = \frac{\partial}{\partial x^1}+ \frac{1}{2} (y^1)^2 \, \frac{\partial}{\partial x^3}
- x^2\frac{\partial}{\partial y^3} - y^1\frac{\partial}{\partial y^2}
\quad , \quad
X_2 = \frac{\partial}{\partial y^1} + y^2\frac{\partial}{\partial x^3}.
\end{equation*}
%\subsubsection{Holomorphic functions for $J_{\alpha,\beta}.$}
Let $G$ denote the group $G_0$ endowed with the left invariant  structure
of complex manifold  defined by
$J_{\alpha,\beta}$ (\ref{64equivfinal}).
Then $ H_{\Cmath}(G) =\{f \in C^{\infty}(G_0) \; ; \; \tilde{X}_j^{-} \, f = 0
\; \forall j=1,3,5\}.$
One has
\begin{eqnarray*}
\tilde{X}_1^{-} & =&
(1+i\alpha)\left[
2\, \frac{\partial}{\partial \overline{w^1}}
+ \left( \frac{1}{2} (y^1)^2 +iy^2 (1-i\alpha) \right) \frac{\partial}{\partial x^3}
-x^2 \frac{\partial}{\partial x^3}
-y^1 \frac{\partial}{\partial y^2} \right]
\quad , \\
\tilde{X}_3^{-} & =&2\, \frac{\partial}{\partial \overline{w^2}}
\quad , \quad
\tilde{X}_5^{-}  =2\, \frac{\partial}{\partial \overline{w^3}}
\quad , \quad
\end{eqnarray*}
where
\begin{eqnarray*}
w^1 &=& x^1 -i (\alpha x^1 + y^1)\\
w^2 &=& x^2 + \frac{(1-\alpha \beta) (\alpha+\beta)}{(1+\alpha^2)(1+\beta^2)}y^2
-i  \frac{(\alpha+\beta)^2}{(1+\alpha^2)(1+\beta^2)} y^2
\\
w^3 &=& x^3 +\frac{\beta(\alpha^2+1)}{\alpha+\beta}\,  y^3
-i \frac{\alpha^2+1}{\alpha+\beta}\,  y^3 .
\end{eqnarray*}
Then $f \in C^{\infty}(G_0)$ is in $ H_{\Cmath}(G)$ if and only if it is holomorphic
with respect to $w^2$ and $w^3$ and satisfies the equation
$$2\, \frac{\partial f}{\partial \overline{w^1}}
+ \left( \frac{1}{4} ((\alpha -i) w^1 + (\alpha +i) \overline{w^1})^2
- \frac{(\alpha^2+1)(\beta-i)}{\alpha+\beta} w^2 \right)
\,  \frac{\partial f}{\partial w^3}
+ \frac{A}{2}          ((\alpha -i) w^1 + (\alpha +i) \overline{w^1})
\,  \frac{\partial f}{\partial w^2}
=0$$
where
$$ A=\frac{(\alpha+\beta)(1-\alpha\beta -i(\alpha+\beta))}{(1+\alpha^2)(1+\beta^2)}.$$
The 3 functions
\begin{equation*}
\varphi^1 =w^1
\quad , \quad
\varphi^2 = w^2 - \frac{A}{4}  \, \left( \frac{\alpha + i}{2} \,
(\overline{w^1})^2
+ (\alpha-i) \, |w^1|^2\right) \quad ,
\end{equation*}
\begin{multline*}
\varphi^3 =w^3 - \frac{1}{16} (\alpha-i)^2 \, (w^1)^2 \overline{w^1}
-\frac{\alpha^2+1}{16}\,
\left( 1+ A \frac{(\beta-i)(\alpha -i)}{\alpha + \beta} \right) w^1 (\overline{w^1})^2
\\
-\frac{\alpha+i}{48}\,
\left( \alpha +i+ 2 A \frac{(\beta-i)(\alpha^2 +1)}{\alpha + \beta} \right)  (\overline{w^1})^3
+\frac{(\alpha^2+1)(\beta-i)}{2(\alpha + \beta)}\,   \overline{w^1}  w^2
\end{multline*}
are holomorphic.
Let $F : G \rightarrow \Cmath^3$ defined
by
$F=(\varphi^1,\varphi^2,\varphi^3).$
$F$ is a biholomorphic bijection, hence a global chart on $G.$
We determine now how the multiplication of $G$ looks like in that chart.
Let $a,x \in G $ with respective second kind canonical coordinates
$(x^1,y^1,x^2,y^2,x^3,y^3), (\alpha^1, \beta^1, \alpha^2, \beta^2, \alpha^3, \beta^3)$
as in  (\ref{x6general}).
With obvious notations, computations yield:
\begin{eqnarray*}
w^1_{a x} &=& w^1_a + w^1_x\\
w^2_{a x} &=& w^2_a + w^2_x - \frac{\alpha + \beta}{(1+\alpha^2)(1+\beta^2)}
(1- \alpha \beta -i(\alpha +\beta)\,  \beta^1 x^1
\\
w^3_{a x} &=& w^3_a + w^3_x
+\frac{1}{2}\, (\beta^1)^2 x^1  -\beta^2 y^1   +\beta^1 x^1 y^1
-\frac{(\alpha^2+1)(\beta-i)}{\alpha+\beta}\,  \alpha^2x^1
.
\end{eqnarray*}
We then get
\begin{eqnarray*}
\varphi^1_{a x} &=& \varphi^1_a + \varphi^1_x\\
\varphi^2_{a x} &=& \varphi^2_a + \varphi^2_x + \frac{\alpha+\beta}{4(\alpha^2+1)(\beta-i)} \,
\left(2 (1-i\alpha)\overline{\varphi^1_a} - (\alpha^2+1)\varphi^1_a\right) \, \varphi^1_x
\\
\varphi^3_{a x} &=& \varphi^3_a +\varphi^3_x + \chi(a,x)
\end{eqnarray*}
where
\begin{multline*}
\chi(a,x) =
\frac{1}{16(\alpha+\beta)^2}
\left(
-(\alpha+\beta)(\alpha+i)^2(\beta+i) \,(\overline{\varphi^1_a})^2
+(\alpha+\beta)(\alpha+2\beta+i)(\alpha-i)^2  \,({\varphi^1_a})^2
\right.\\ \left.
+4(\alpha+\beta)(\alpha+\beta +i(\alpha \beta -1)) \,|\varphi^1_a|^2
-8 (1+\alpha^2)(1+\beta^2)   \, \overline{\varphi^2_a}
-8 (1+\alpha^2)(\alpha -i)(\beta-i)   \, \varphi^2_a
\right)
\, \varphi^1_x
\\
+
\frac{1}{16}
\left( 2   (\alpha^2 -1-2i\alpha)\,\varphi^1_a
+(\alpha^2 +3+2i)\,\overline{\varphi^1_a}\right)
\, (\varphi^1_x)^2
+\frac{(\alpha^2+1)(\beta-i)}{2(\alpha+\beta)}\, \overline{\varphi^1_a} \,\varphi^2_x .
\end{multline*}
%%%%%%%%%%%%%%%%%%%%%%%%%%%%%%%%%%%%%%%%%%%%%%%%%%%%%%%%%%%%%%%%%%%%%%%%%%%%%%%%%%%%%%%%%%
%%%%%%%%%%%%%%%%%%%%%%%%%%%%%%%%%%%%%%%%%%%%%%%%%%%%%%%%%%%%%%%%%%%%%%%%%%%%%%%%%%%%%%%%%%
\section{Lie Algebra ${\mathcal{G}}_{6,1}$ (isomorphic to $M4$).}
%%%%%%%%%%%%%%%%%%%%%%%%%%%%%%%%%%%%%%%%%%%%%%%%%%%%%%%%%%%%%%%%%%%%%%%%%%%%%%%%%%%%%%%%%%
Commutation relations for
$ {\mathcal{G}}_{6,1}:$
$[x_1,x_2]=x_5$;
$[x_1,x_4]=x_6$;
$[x_2,x_3]=x_6.$
\subsection{Case $ {\xi^2_1} \neq 0 , \xi^2_3 \neq 0, \xi^1_3 \neq \xi^2_4.$}
\label{61_1}
%{\fontsize{8}{10}\selectfont
\begin{equation}
\label{61general1}
J = \begin{pmatrix}
\boxed{\xi^1_1}&
*&
\boxed{\xi^1_3}&
*&
0&
0\\
\boxed{\xi^2_1}&
\boxed{\xi^2_2}&
\boxed{\xi^2_3}&
\boxed{\xi^2_4}&
0&
0\\
*&
*&
*&
*&
0&
0\\
*&
*&
*&
*&
0&
0\\
*&
*&
*&
*&
\boxed{\xi^5_5}&
%$  J^5_6=( - ({\xi^5_5}^2 + 1))/{\xi^6_5};$\\
 - ({\xi^5_5}^2 + 1)/{\xi^6_5}\\
\boxed{\xi^6_1}&
\boxed{\xi^6_2}&
\boxed{\xi^6_3}&
\boxed{\xi^6_4}&
\boxed{\xi^6_5}&
-{\xi^5_5}
 \end{pmatrix}
\end{equation}
%}
where
{\fontsize{6}{10}\selectfont
$  \mathbf{J^1_2}=( - (({\xi^2_4} {\xi^2_3} {\xi^1_1} - {\xi^2_4} {\xi^2_1} {\xi^1_3} - {\xi^2_3} {\xi^2_2} {\xi^1_3} - ({\xi^2_4} - {\xi^1_3}) {\xi^5_5} {\xi^2_3}) {\xi^6_5} + ({\xi^5_5}^2 + 1) (
{\xi^2_4} - {\xi^1_3}) {\xi^2_1}))/({\xi^6_5} {\xi^2_3}^2);$
 \quad
$  \mathbf{J^1_4}=( - (({\xi^5_5}^2 + 1) ({\xi^2_4} - {\xi^1_3}) - {\xi^6_5} {\xi^2_4} {\xi^1_3}))
/({\xi^6_5} {\xi^2_3});$
 \quad
$  \mathbf{J^3_1}=( - ({\xi^6_5} {\xi^5_5} {\xi^2_4}^2 {\xi^2_3} {\xi^2_1} - {\xi^6_5} {\xi^5_5} {\xi^2_4} {\xi^2_3} {\xi^2_1} {\xi^1_3} - {\xi^6_5} {\xi^2_4}^2 {\xi^2_3} {\xi^2_1} {\xi^1_1} + {\xi^6_5} {\xi^2_4}^2 {\xi^2_1}^2 {\xi^1_3} + {\xi^6_5} {\xi^2_4} {\xi^2_3}^2 {\xi^1_1}^2 + {\xi^6_5} {\xi^2_4} {\xi^2_3}^2 - {\xi^6_5} {\xi^2_4} {\xi^2_3} {\xi^2_1} {\xi^1_3} {\xi^1_1} - {\xi^5_5}^2 {\xi^2_4}^2 {\xi^2_1}^2 + {\xi^5_5}^2 {\xi^2_4} {\xi^2_3} {\xi^2_2} {\xi^2_1} + {\xi^5_5}^2 {\xi^2_4} {\xi^2_3} {\xi^2_1} {\xi^1_1} + {\xi^5_5}^2 {\xi^2_4} {\xi^2_1}^2 {\xi^1_3} -
{\xi^5_5}^2 {\xi^2_3} {\xi^2_2} {\xi^2_1} {\xi^1_3} - {\xi^5_5}^2 {\xi^2_3} {\xi^2_1} {\xi^1_3}
{\xi^1_1} - {\xi^2_4}^2 {\xi^2_1}^2 + {\xi^2_4} {\xi^2_3} {\xi^2_2} {\xi^2_1} + {\xi^2_4} {\xi^2_3} {\xi^2_1} {\xi^1_1} + {\xi^2_4} {\xi^2_1}^2 {\xi^1_3} - {\xi^2_3} {\xi^2_2} {\xi^2_1} {\xi^1_3} - {\xi^2_3} {\xi^2_1} {\xi^1_3} {\xi^1_1}))/(({\xi^5_5}^2 + 1) ({\xi^2_4} - {\xi^1_3}) {\xi^2_3}^2);$
 \quad
$  \mathbf{J^3_2}=( - ({\xi^6_5}^2 {\xi^5_5} {\xi^2_4}^2 {\xi^2_3}^2 {\xi^2_2} + {\xi^6_5}^2 {\xi^5_5} {\xi^2_4}^2 {\xi^2_3}^2 {\xi^1_1} - {\xi^6_5}^2 {\xi^5_5} {\xi^2_4}^2 {\xi^2_3} {\xi^2_1} {\xi^1_3} - {\xi^6_5}^2 {\xi^5_5} {\xi^2_4} {\xi^2_3}^2 {\xi^2_2} {\xi^1_3} - {\xi^6_5}^2
 {\xi^5_5} {\xi^2_4} {\xi^2_3}^2 {\xi^1_3} {\xi^1_1} + {\xi^6_5}^2 {\xi^5_5} {\xi^2_4} {\xi^2_3}
 {\xi^2_1} {\xi^1_3}^2 - {\xi^6_5}^2 {\xi^2_4}^2 {\xi^2_3}^2 {\xi^2_2} {\xi^1_1} - {\xi^6_5}
^2 {\xi^2_4}^2 {\xi^2_3}^2 {\xi^1_1}^2 + {\xi^6_5}^2 {\xi^2_4}^2 {\xi^2_3} {\xi^2_2} {\xi^2_1} {\xi^1_3} + 2 {\xi^6_5}^2 {\xi^2_4}^2 {\xi^2_3} {\xi^2_1} {\xi^1_3} {\xi^1_1} - {\xi^6_5}
^2 {\xi^2_4}^2 {\xi^2_1}^2 {\xi^1_3}^2 + {\xi^6_5}^2 {\xi^2_4} {\xi^2_3}^2 {\xi^2_2} {\xi^1_3} {\xi^1_1} - {\xi^6_5}^2 {\xi^2_4} {\xi^2_3}^2 {\xi^1_3} - {\xi^6_5}^2 {\xi^2_4} {\xi^2_3} {\xi^2_2} {\xi^2_1} {\xi^1_3}^2 + {\xi^6_5} {\xi^5_5}^3 {\xi^2_4}^2 {\xi^2_3} {\xi^2_1} - 2
 {\xi^6_5} {\xi^5_5}^3 {\xi^2_4} {\xi^2_3} {\xi^2_1} {\xi^1_3} + {\xi^6_5} {\xi^5_5}^3 {\xi^2_3}
 {\xi^2_1} {\xi^1_3}^2 - {\xi^6_5} {\xi^5_5}^2 {\xi^2_4}^2 {\xi^2_3} {\xi^2_2} {\xi^2_1} - 2
{\xi^6_5} {\xi^5_5}^2 {\xi^2_4}^2 {\xi^2_3} {\xi^2_1} {\xi^1_1} + 2 {\xi^6_5} {\xi^5_5}^2 {\xi^2_4}^2 {\xi^2_1}^2 {\xi^1_3} + {\xi^6_5} {\xi^5_5}^2 {\xi^2_4} {\xi^2_3}^2 {\xi^2_2}^2 +
{\xi^6_5} {\xi^5_5}^2 {\xi^2_4} {\xi^2_3}^2 + 2 {\xi^6_5} {\xi^5_5}^2 {\xi^2_4} {\xi^2_3} {\xi^2_2} {\xi^2_1} {\xi^1_3} + 2 {\xi^6_5} {\xi^5_5}^2 {\xi^2_4} {\xi^2_3} {\xi^2_1} {\xi^1_3} {\xi^1_1} - 2 {\xi^6_5} {\xi^5_5}^2 {\xi^2_4} {\xi^2_1}^2 {\xi^1_3}^2 - {\xi^6_5} {\xi^5_5}^2 {\xi^2_3}^2 {\xi^2_2}^2 {\xi^1_3} - {\xi^6_5} {\xi^5_5}^2 {\xi^2_3}^2 {\xi^1_3} - {\xi^6_5} {\xi^5_5}^2 {\xi^2_3} {\xi^2_2} {\xi^2_1} {\xi^1_3}^2 + {\xi^6_5} {\xi^5_5} {\xi^2_4}^2 {\xi^2_3}
 {\xi^2_1} - 2 {\xi^6_5} {\xi^5_5} {\xi^2_4} {\xi^2_3} {\xi^2_1} {\xi^1_3} + {\xi^6_5} {\xi^5_5}
{\xi^2_3} {\xi^2_1} {\xi^1_3}^2 - {\xi^6_5} {\xi^2_4}^2 {\xi^2_3} {\xi^2_2} {\xi^2_1} - 2 {\xi^6_5} {\xi^2_4}^2 {\xi^2_3} {\xi^2_1} {\xi^1_1} + 2 {\xi^6_5} {\xi^2_4}^2 {\xi^2_1}^2 {\xi^1_3}
 + {\xi^6_5} {\xi^2_4} {\xi^2_3}^2 {\xi^2_2}^2 + {\xi^6_5} {\xi^2_4} {\xi^2_3}^2 + 2 {\xi^6_5} {\xi^2_4} {\xi^2_3} {\xi^2_2} {\xi^2_1} {\xi^1_3} + 2 {\xi^6_5} {\xi^2_4} {\xi^2_3} {\xi^2_1} {\xi^1_3} {\xi^1_1} - 2 {\xi^6_5} {\xi^2_4} {\xi^2_1}^2 {\xi^1_3}^2 - {\xi^6_5} {\xi^2_3}^2 {\xi^2_2}^2 {\xi^1_3} - {\xi^6_5} {\xi^2_3}^2 {\xi^1_3} - {\xi^6_5} {\xi^2_3} {\xi^2_2} {\xi^2_1}
{\xi^1_3}^2 - {\xi^5_5}^4 {\xi^2_4}^2 {\xi^2_1}^2 + 2 {\xi^5_5}^4 {\xi^2_4} {\xi^2_1}^2
{\xi^1_3} - {\xi^5_5}^4 {\xi^2_1}^2 {\xi^1_3}^2 - 2 {\xi^5_5}^2 {\xi^2_4}^2 {\xi^2_1}^2
+ 4 {\xi^5_5}^2 {\xi^2_4} {\xi^2_1}^2 {\xi^1_3} - 2 {\xi^5_5}^2 {\xi^2_1}^2 {\xi^1_3}^2 -
 {\xi^2_4}^2 {\xi^2_1}^2 + 2 {\xi^2_4} {\xi^2_1}^2 {\xi^1_3} - {\xi^2_1}^2 {\xi^1_3}^2))/
(({\xi^5_5}^2 + 1) ({\xi^2_4} - {\xi^1_3}) {\xi^6_5} {\xi^2_3}^3);$
 \quad
$  \mathbf{J^3_3}=( - ({\xi^6_5} {\xi^5_5} {\xi^2_4} {\xi^2_3} - {\xi^6_5} {\xi^2_4} {\xi^2_3} {\xi^1_1} + {\xi^6_5} {\xi^2_4} {\xi^2_1} {\xi^1_3} - {\xi^5_5}^2 {\xi^2_4} {\xi^2_1} + {\xi^5_5}^2 {\xi^2_3} {\xi^2_2} + {\xi^5_5}^2 {\xi^2_1} {\xi^1_3} - {\xi^2_4} {\xi^2_1} + {\xi^2_3} {\xi^2_2} +
 {\xi^2_1} {\xi^1_3}))/(({\xi^5_5}^2 + 1) {\xi^2_3});$
 \quad
$  \mathbf{J^3_4}=( - ({\xi^6_5}^2 {\xi^5_5} {\xi^2_4}^2 {\xi^2_3} - {\xi^6_5}^2 {\xi^2_4}^2 {\xi^2_3} {\xi^1_1} + {\xi^6_5}^2 {\xi^2_4}^2 {\xi^2_1} {\xi^1_3} - {\xi^6_5} {\xi^5_5}^2 {\xi^2_4}^2 {\xi^2_1} + {\xi^6_5} {\xi^5_5}^2 {\xi^2_4} {\xi^2_3} {\xi^2_2} - {\xi^6_5} {\xi^5_5}^2
{\xi^2_4} {\xi^2_3} {\xi^1_1} + 2 {\xi^6_5} {\xi^5_5}^2 {\xi^2_4} {\xi^2_1} {\xi^1_3} - {\xi^6_5}
 {\xi^2_4}^2 {\xi^2_1} + {\xi^6_5} {\xi^2_4} {\xi^2_3} {\xi^2_2} - {\xi^6_5} {\xi^2_4} {\xi^2_3}
{\xi^1_1} + 2 {\xi^6_5} {\xi^2_4} {\xi^2_1} {\xi^1_3} - {\xi^5_5}^4 {\xi^2_4} {\xi^2_1} + {\xi^5_5}^4 {\xi^2_1} {\xi^1_3} - 2 {\xi^5_5}^2 {\xi^2_4} {\xi^2_1} + 2 {\xi^5_5}^2 {\xi^2_1} {\xi^1_3} - {\xi^2_4} {\xi^2_1} + {\xi^2_1} {\xi^1_3}))/(({\xi^5_5}^2 + 1) {\xi^6_5} {\xi^2_3}^2)
;$
 \quad
$  \mathbf{J^4_1}=( - ((({\xi^5_5}^2 + 1) {\xi^2_1} - {\xi^6_5} {\xi^5_5} {\xi^2_3}) ({\xi^2_4} -
{\xi^1_3}) {\xi^2_1} - (({\xi^2_3} {\xi^1_1}^2 + {\xi^2_3} - {\xi^2_1} {\xi^1_3} {\xi^1_1}) {\xi^2_3} - ({\xi^2_3} {\xi^1_1} - {\xi^2_1} {\xi^1_3}) {\xi^2_4} {\xi^2_1}) {\xi^6_5}))/(({\xi^5_5}
^2 + 1) ({\xi^2_4} - {\xi^1_3}) {\xi^2_3});$
 \quad
$  \mathbf{J^4_2}=( - ((({\xi^5_5}^2 + 1) {\xi^2_1} - {\xi^6_5} {\xi^5_5} {\xi^2_3}) ({\xi^2_4} -
{\xi^1_3}) ({\xi^2_3} {\xi^2_2} + {\xi^2_3} {\xi^1_1} - {\xi^2_1} {\xi^1_3}) - (({\xi^2_3} {\xi^2_2} {\xi^1_1} - {\xi^2_3} - {\xi^2_2} {\xi^2_1} {\xi^1_3}) {\xi^2_3} {\xi^1_3} - ({\xi^2_3} {\xi^2_2} + {\xi^2_3} {\xi^1_1} - {\xi^2_1} {\xi^1_3}) ({\xi^2_3} {\xi^1_1} - {\xi^2_1} {\xi^1_3}) {\xi^2_4}) {\xi^6_5}))/(({\xi^5_5}^2 + 1) ({\xi^2_4} - {\xi^1_3}) {\xi^2_3}^2);$
 \quad
$  \mathbf{J^4_3}= - (({\xi^2_3} {\xi^1_1} - {\xi^2_1} {\xi^1_3} - {\xi^5_5} {\xi^2_3}) {\xi^6_5} + (
{\xi^5_5}^2 + 1) {\xi^2_1})/({\xi^5_5}^2 + 1);$
 \quad
$  \mathbf{J^4_4}=( - (({\xi^2_3} {\xi^1_1} - {\xi^2_1} {\xi^1_3} - {\xi^5_5} {\xi^2_3}) {\xi^6_5} {\xi^2_4} + ({\xi^2_3} {\xi^1_1} - {\xi^2_1} {\xi^1_3} + {\xi^2_4} {\xi^2_1}) ({\xi^5_5}^2 + 1)))
/(({\xi^5_5}^2 + 1) {\xi^2_3});$
 \quad
$  \mathbf{J^5_1}=((({\xi^5_5} - {\xi^1_1}) {\xi^6_1} - {\xi^6_2} {\xi^2_1}) ({\xi^5_5}^2 + 1) ({\xi^2_4} - {\xi^1_3}) {\xi^2_3}^2 + ({\xi^6_5} {\xi^5_5} {\xi^2_4}^2 {\xi^2_3} {\xi^2_1} - 2 {\xi^6_5} {\xi^5_5} {\xi^2_4} {\xi^2_3} {\xi^2_1} {\xi^1_3} + {\xi^6_5} {\xi^5_5} {\xi^2_3} {\xi^2_1}
{\xi^1_3}^2 - {\xi^6_5} {\xi^2_4}^2 {\xi^2_3} {\xi^2_1} {\xi^1_1} + {\xi^6_5} {\xi^2_4}^2 {\xi^2_1}^2 {\xi^1_3} + {\xi^6_5} {\xi^2_4} {\xi^2_3}^2 {\xi^1_1}^2 + {\xi^6_5} {\xi^2_4} {\xi^2_3}^2 - {\xi^6_5} {\xi^2_4} {\xi^2_3} {\xi^2_1} {\xi^1_3} {\xi^1_1} - {\xi^5_5}^2 {\xi^2_4}^2
{\xi^2_1}^2 + {\xi^5_5}^2 {\xi^2_4} {\xi^2_3} {\xi^2_2} {\xi^2_1} + {\xi^5_5}^2 {\xi^2_4} {\xi^2_3} {\xi^2_1} {\xi^1_1} + 2 {\xi^5_5}^2 {\xi^2_4} {\xi^2_1}^2 {\xi^1_3} - {\xi^5_5}^2 {\xi^2_3} {\xi^2_2} {\xi^2_1} {\xi^1_3} - {\xi^5_5}^2 {\xi^2_3} {\xi^2_1} {\xi^1_3} {\xi^1_1} - {\xi^5_5}^2 {\xi^2_1}^2 {\xi^1_3}^2 - {\xi^2_4}^2 {\xi^2_1}^2 + {\xi^2_4} {\xi^2_3} {\xi^2_2} {\xi^2_1} + {\xi^2_4} {\xi^2_3} {\xi^2_1} {\xi^1_1} + 2 {\xi^2_4} {\xi^2_1}^2 {\xi^1_3} - {\xi^2_3}
 {\xi^2_2} {\xi^2_1} {\xi^1_3} - {\xi^2_3} {\xi^2_1} {\xi^1_3} {\xi^1_1} - {\xi^2_1}^2 {\xi^1_3}
^2) {\xi^6_3} - (({\xi^2_3} {\xi^1_1}^2 + {\xi^2_3} - {\xi^2_1} {\xi^1_3} {\xi^1_1}) {\xi^2_3}
 - ({\xi^2_3} {\xi^1_1} - {\xi^2_1} {\xi^1_3}) {\xi^2_4} {\xi^2_1}) {\xi^6_5} {\xi^6_4} {\xi^2_3}
+ (({\xi^5_5}^2 + 1) {\xi^2_1} - {\xi^6_5} {\xi^5_5} {\xi^2_3}) ({\xi^6_4} {\xi^2_3} - {\xi^6_3} {\xi^1_3}) ({\xi^2_4} - {\xi^1_3}) {\xi^2_1})/(({\xi^5_5}^2 + 1) ({\xi^2_4} - {\xi^1_3}) {\xi^6_5} {\xi^2_3}^2);$
 \quad
$  \mathbf{J^5_2}=(((({\xi^5_5}^2 + 1) ({\xi^2_4} - {\xi^1_3}) {\xi^6_1} {\xi^2_1} + {\xi^6_5} {\xi^6_2} {\xi^5_5} {\xi^2_3}^2 - {\xi^6_5} {\xi^6_2} {\xi^2_3}^2 {\xi^2_2}) ({\xi^5_5}^2 + 1) (
{\xi^2_4} - {\xi^1_3}) - ({\xi^6_5} {\xi^5_5} {\xi^2_4} {\xi^2_3}^2 {\xi^2_2} + {\xi^6_5} {\xi^5_5} {\xi^2_4} {\xi^2_3}^2 {\xi^1_1} - {\xi^6_5} {\xi^5_5} {\xi^2_3}^2 {\xi^2_2} {\xi^1_3} - {\xi^6_5} {\xi^5_5} {\xi^2_3}^2 {\xi^1_3} {\xi^1_1} - {\xi^6_5} {\xi^2_4} {\xi^2_3}^2 {\xi^2_2} {\xi^1_1} - {\xi^6_5} {\xi^2_4} {\xi^2_3}^2 {\xi^1_1}^2 + {\xi^6_5} {\xi^2_4} {\xi^2_3} {\xi^2_2}
{\xi^2_1} {\xi^1_3} + 2 {\xi^6_5} {\xi^2_4} {\xi^2_3} {\xi^2_1} {\xi^1_3} {\xi^1_1} - {\xi^6_5} {\xi^2_4} {\xi^2_1}^2 {\xi^1_3}^2 + {\xi^6_5} {\xi^2_3}^2 {\xi^2_2} {\xi^1_3} {\xi^1_1} - {\xi^6_5} {\xi^2_3}^2 {\xi^1_3} - {\xi^6_5} {\xi^2_3} {\xi^2_2} {\xi^2_1} {\xi^1_3}^2 - {\xi^5_5}^2
{\xi^2_4} {\xi^2_3} {\xi^2_2} {\xi^2_1} - {\xi^5_5}^2 {\xi^2_4} {\xi^2_3} {\xi^2_1} {\xi^1_1} +
{\xi^5_5}^2 {\xi^2_3} {\xi^2_2} {\xi^2_1} {\xi^1_3} + {\xi^5_5}^2 {\xi^2_3} {\xi^2_1} {\xi^1_3}
{\xi^1_1} - {\xi^2_4} {\xi^2_3} {\xi^2_2} {\xi^2_1} - {\xi^2_4} {\xi^2_3} {\xi^2_1} {\xi^1_1} + {\xi^2_3} {\xi^2_2} {\xi^2_1} {\xi^1_3} + {\xi^2_3} {\xi^2_1} {\xi^1_3} {\xi^1_1}) {\xi^6_5} {\xi^6_4}
) {\xi^2_3} + ({\xi^6_5} {\xi^5_5}^2 {\xi^2_4}^2 {\xi^2_3}^2 - 2 {\xi^6_5} {\xi^5_5}^2 {\xi^2_4} {\xi^2_3}^2 {\xi^1_3} + {\xi^6_5} {\xi^5_5}^2 {\xi^2_3}^2 {\xi^1_3}^2 + {\xi^6_5} {\xi^5_5} {\xi^2_4}^2 {\xi^2_3}^2 {\xi^2_2} - {\xi^6_5} {\xi^5_5} {\xi^2_4} {\xi^2_3}^2 {\xi^2_2}
{\xi^1_3} + {\xi^6_5} {\xi^5_5} {\xi^2_4} {\xi^2_3}^2 {\xi^1_3} {\xi^1_1} - {\xi^6_5} {\xi^5_5}
{\xi^2_3}^2 {\xi^1_3}^2 {\xi^1_1} - {\xi^6_5} {\xi^2_4}^2 {\xi^2_3}^2 {\xi^2_2} {\xi^1_1} -
{\xi^6_5} {\xi^2_4}^2 {\xi^2_3}^2 {\xi^1_1}^2 + {\xi^6_5} {\xi^2_4}^2 {\xi^2_3} {\xi^2_2} {\xi^2_1} {\xi^1_3} + 2 {\xi^6_5} {\xi^2_4}^2 {\xi^2_3} {\xi^2_1} {\xi^1_3} {\xi^1_1} - {\xi^6_5}
{\xi^2_4}^2 {\xi^2_1}^2 {\xi^1_3}^2 + {\xi^6_5} {\xi^2_4} {\xi^2_3}^2 {\xi^2_2} {\xi^1_3} {\xi^1_1} - {\xi^6_5} {\xi^2_4} {\xi^2_3}^2 {\xi^1_3} - {\xi^6_5} {\xi^2_4} {\xi^2_3} {\xi^2_2} {\xi^2_1} {\xi^1_3}^2 - {\xi^5_5}^3 {\xi^2_4}^2 {\xi^2_3} {\xi^2_1} + 2 {\xi^5_5}^3 {\xi^2_4}
{\xi^2_3} {\xi^2_1} {\xi^1_3} - {\xi^5_5}^3 {\xi^2_3} {\xi^2_1} {\xi^1_3}^2 - {\xi^5_5}^2 {\xi^2_4}^2 {\xi^2_3} {\xi^2_2} {\xi^2_1} - {\xi^5_5}^2 {\xi^2_4}^2 {\xi^2_3} {\xi^2_1} {\xi^1_1}
+ {\xi^5_5}^2 {\xi^2_4}^2 {\xi^2_1}^2 {\xi^1_3} + {\xi^5_5}^2 {\xi^2_4} {\xi^2_3}^2 {\xi^2_2}^2 + {\xi^5_5}^2 {\xi^2_4} {\xi^2_3}^2 + 2 {\xi^5_5}^2 {\xi^2_4} {\xi^2_3} {\xi^2_2} {\xi^2_1} {\xi^1_3} - {\xi^5_5}^2 {\xi^2_4} {\xi^2_1}^2 {\xi^1_3}^2 - {\xi^5_5}^2 {\xi^2_3}^2
{\xi^2_2}^2 {\xi^1_3} - {\xi^5_5}^2 {\xi^2_3}^2 {\xi^1_3} - {\xi^5_5}^2 {\xi^2_3} {\xi^2_2}
{\xi^2_1} {\xi^1_3}^2 + {\xi^5_5}^2 {\xi^2_3} {\xi^2_1} {\xi^1_3}^2 {\xi^1_1} - {\xi^5_5} {\xi^2_4}^2 {\xi^2_3} {\xi^2_1} + 2 {\xi^5_5} {\xi^2_4} {\xi^2_3} {\xi^2_1} {\xi^1_3} - {\xi^5_5} {\xi^2_3} {\xi^2_1} {\xi^1_3}^2 - {\xi^2_4}^2 {\xi^2_3} {\xi^2_2} {\xi^2_1} - {\xi^2
 _4}^2 {\xi^2_3} {\xi^2_1} {\xi^1_1} + {\xi^2_4}^2 {\xi^2_1}^2 {\xi^1_3} + {\xi^2_4} {\xi^2_3}^2 {\xi^2_2}
^2 + {\xi^2_4} {\xi^2_3}^2 + 2 {\xi^2_4} {\xi^2_3} {\xi^2_2} {\xi^2_1} {\xi^1_3} - {\xi^2_4}
{\xi^2_1}^2 {\xi^1_3}^2 - {\xi^2_3}^2 {\xi^2_2}^2 {\xi^1_3} - {\xi^2_3}^2 {\xi^1_3} - {\xi^2_3} {\xi^2_2} {\xi^2_1} {\xi^1_3}^2 + {\xi^2_3} {\xi^2_1} {\xi^1_3}^2 {\xi^1_1}) {\xi^6_5} {\xi^6_3} + ({\xi^2_4} {\xi^2_3} {\xi^1_1} - {\xi^2_4} {\xi^2_1} {\xi^1_3} - {\xi^2_3} {\xi^2_2} {\xi^1_3} - ({\xi^2_4} - {\xi^1_3}) {\xi^5_5} {\xi^2_3}) ({\xi^5_5}^2 + 1) ({\xi^2_4} - {\xi^1_3})
 {\xi^6_5} {\xi^6_1} {\xi^2_3} - (({\xi^5_5}^2 + 1) {\xi^2_1} - {\xi^6_5} {\xi^5_5} {\xi^2_3})
({\xi^6_5} {\xi^6_4} {\xi^2_3} {\xi^2_1} {\xi^1_3} - {\xi^6_5} {\xi^6_3} {\xi^5_5} {\xi^2_4} {\xi^2_3} + {\xi^6_5} {\xi^6_3} {\xi^5_5} {\xi^2_3} {\xi^1_3} + {\xi^6_5} {\xi^6_3} {\xi^2_4} {\xi^2_3}
{\xi^1_1} - {\xi^6_5} {\xi^6_3} {\xi^2_4} {\xi^2_1} {\xi^1_3} - {\xi^6_5} {\xi^6_3} {\xi^2_3} {\xi^1_3} {\xi^1_1} + {\xi^6_3} {\xi^5_5}^2 {\xi^2_4} {\xi^2_1} - {\xi^6_3} {\xi^5_5}^2 {\xi^2_1} {\xi^1_3} + {\xi^6_3} {\xi^2_4} {\xi^2_1} - {\xi^6_3} {\xi^2_1} {\xi^1_3}) ({\xi^2_4} - {\xi^1_3}))/
(({\xi^5_5}^2 + 1) ({\xi^2_4} - {\xi^1_3}) {\xi^6_5}^2 {\xi^2_3}^3);$
 \quad
$  \mathbf{J^5_3}=( - ((({\xi^6_2} {\xi^2_3} + {\xi^6_1} {\xi^1_3} - {\xi^6_3} {\xi^5_5} - {\xi^6_4}
{\xi^2_1}) ({\xi^5_5}^2 + 1) - ({\xi^2_3} {\xi^1_1} - {\xi^2_1} {\xi^1_3} - {\xi^5_5} {\xi^2_3}
) {\xi^6_5} {\xi^6_4}) {\xi^2_3} + ({\xi^2_4} {\xi^2_3} {\xi^1_1} - {\xi^2_4} {\xi^2_1} {\xi^1_3}
- {\xi^2_3} {\xi^2_2} {\xi^1_3} - ({\xi^2_4} - {\xi^1_3}) {\xi^5_5} {\xi^2_3}) {\xi^6_5} {\xi^6_3}
 - ((({\xi^5_5}^2 + 1) {\xi^2_2} + ({\xi^5_5} - {\xi^2_2}) {\xi^6_5} {\xi^1_3}) {\xi^2_3} - ({\xi^5_5}^2 + 1) ({\xi^2_4} - {\xi^1_3}) {\xi^2_1}) {\xi^6_3}))/(({\xi^5_5}^2 + 1) {\xi^6_5}
 {\xi^2_3});$
 \quad
$  \mathbf{J^5_4}=( - ({\xi^6_5}^2 {\xi^6_4} {\xi^5_5} {\xi^2_4} {\xi^2_3}^2 - {\xi^6_5}^2 {\xi^6_4} {\xi^2_4} {\xi^2_3}^2 {\xi^1_1} + {\xi^6_5}^2 {\xi^6_4} {\xi^2_4} {\xi^2_3} {\xi^2_1} {\xi^1_3} - {\xi^6_5}^2 {\xi^6_3} {\xi^5_5} {\xi^2_4}^2 {\xi^2_3} + {\xi^6_5}^2 {\xi^6_3} {\xi^2_4}
^2 {\xi^2_3} {\xi^1_1} - {\xi^6_5}^2 {\xi^6_3} {\xi^2_4}^2 {\xi^2_1} {\xi^1_3} - {\xi^6_5} {\xi^6_4} {\xi^5_5}^3 {\xi^2_3}^2 - {\xi^6_5} {\xi^6_4} {\xi^5_5}^2 {\xi^2_4} {\xi^2_3} {\xi^2_1}
 - {\xi^6_5} {\xi^6_4} {\xi^5_5}^2 {\xi^2_3}^2 {\xi^1_1} + {\xi^6_5} {\xi^6_4} {\xi^5_5}^2 {\xi^2_3} {\xi^2_1} {\xi^1_3} - {\xi^6_5} {\xi^6_4} {\xi^5_5} {\xi^2_3}^2 - {\xi^6_5} {\xi^6_4} {\xi^2_4} {\xi^2_3} {\xi^2_1} - {\xi^6_5} {\xi^6_4} {\xi^2_3}^2 {\xi^1_1} + {\xi^6_5} {\xi^6_4} {\xi^2_3} {\xi^2_1} {\xi^1_3} + {\xi^6_5} {\xi^6_3} {\xi^5_5}^2 {\xi^2_4}^2 {\xi^2_1} - {\xi^6_5} {\xi^6_3} {\xi^5_5}^2 {\xi^2_4} {\xi^2_3} {\xi^2_2} + {\xi^6_5} {\xi^6_3} {\xi^5_5}^2 {\xi^2_4} {\xi^2_3} {\xi^1_1} - 2 {\xi^6_5} {\xi^6_3} {\xi^5_5}^2 {\xi^2_4} {\xi^2_1} {\xi^1_3} + {\xi^6_5}
{\xi^6_3} {\xi^2_4}^2 {\xi^2_1} - {\xi^6_5} {\xi^6_3} {\xi^2_4} {\xi^2_3} {\xi^2_2} + {\xi^6_5}
{\xi^6_3} {\xi^2_4} {\xi^2_3} {\xi^1_1} - 2 {\xi^6_5} {\xi^6_3} {\xi^2_4} {\xi^2_1} {\xi^1_3} + {\xi^6_5} {\xi^6_2} {\xi^5_5}^2 {\xi^2_4} {\xi^2_3}^2 + {\xi^6_5} {\xi^6_2} {\xi^2_4} {\xi^2_3}^2
 + {\xi^6_5} {\xi^6_1} {\xi^5_5}^2 {\xi^2_4} {\xi^2_3} {\xi^1_3} + {\xi^6_5} {\xi^6_1} {\xi^2_4}
{\xi^2_3} {\xi^1_3} + {\xi^6_3} {\xi^5_5}^4 {\xi^2_4} {\xi^2_1} - {\xi^6_3} {\xi^5_5}^4 {\xi^2_1} {\xi^1_3} + 2 {\xi^6_3} {\xi^5_5}^2 {\xi^2_4} {\xi^2_1} - 2 {\xi^6_3} {\xi^5_5}^2 {\xi^2_1}
{\xi^1_3} + {\xi^6_3} {\xi^2_4} {\xi^2_1} - {\xi^6_3} {\xi^2_1} {\xi^1_3} - {\xi^6_1} {\xi^5_5}^4
 {\xi^2_4} {\xi^2_3} + {\xi^6_1} {\xi^5_5}^4 {\xi^2_3} {\xi^1_3} - 2 {\xi^6_1} {\xi^5_5}^2 {\xi^2_4} {\xi^2_3} + 2 {\xi^6_1} {\xi^5_5}^2 {\xi^2_3} {\xi^1_3} - {\xi^6_1} {\xi^2_4} {\xi^2_3} +
{\xi^6_1} {\xi^2_3} {\xi^1_3}))/(({\xi^5_5}^2 + 1) {\xi^6_5}^2 {\xi^2_3}^2);$
}
and the parameters are subject to the condition
\begin{equation}
\label{cond61general1}
{\xi^2_1} {\xi^2_3} \xi^6_5 ( {\xi^1_3} -{\xi^2_4})  \neq 0 .
 \end{equation}
\par Now the automorphism group of
$ {\mathcal{G}}_{6,1}$ is comprised of the matrices
\begin{equation}
{\fontsize{8}{10}\selectfont
\label{aut61}
\Phi = \begin{pmatrix}
{b^1_1}&
{b^1_2}&
0&
0&
0&
0\\
{b^2_1}&
{b^2_2}&
0&
0&
0&
0\\
{b^3_1}&
{b^3_2}&
{b^1_1} u&
 - {b^1_2} u&
0&
0\\
{b^4_1}&
{b^4_2}&
 - {b^2_1} u&
{b^2_2} u&
0&
0\\
{b^5_1}&
{b^5_2}&
{b^5_3}&
{b^5_4}&
{b^2_2} {b^1_1} - {b^2_1} {b^1_2}&
0\\
{b^6_1}&
{b^6_2}&
{b^6_3}&
{b^6_4}&
{b^3_2} {b^2_1} - {b^3_1} {b^2_2} - {b^4_1} {b^1_2} + {b^4_2} {b^1_1}&
({b^2_2} {b^1_1} - {b^2_1} {b^1_2}) u\end{pmatrix}
}
\end{equation}
where  $u \in \Rmath, \, u \neq 0$ and
${b^2_2} {b^1_1} - {b^2_1} {b^1_2} \neq 0.$
Taking suitable values for $u$ and the $b^i_j$'s,
equivalence by $\Phi$ leads to the case
$ {\xi^1_1}={\xi^1_3}={\xi^2_2}={\xi^5_5}={\xi^6_1}={\xi^6_2}={\xi^6_3}={\xi^6_4} =0 $
and $ {\xi^2_1}={\xi^2_3}={\xi^6_5}=1, \xi^2_4 = \alpha \neq 0$:
\begin{equation}
{\fontsize{8}{10}\selectfont
\label{61equivfinal1}
J_{\alpha} = \begin{pmatrix}
0&
 - \alpha&
0&
 - \alpha&
0&
0\\
1&
0&
1&
\alpha&
0&
0\\
\alpha - 1&
\alpha - 1&
\alpha&
(\alpha + 1) \alpha&
0&
0\\
 - (\alpha - 1)/\alpha&
0&
-1&
 - \alpha&
0&
0\\
0&
0&
0&
0&
0&
-1\\
0&
0&
0&
0&
1&
0\end{pmatrix}
\quad
\quad
(\alpha \neq 0).
}
\end{equation}
The $ J_{\alpha}$'s corresponding to distinct $\alpha$'s are not equivalent.
\par Commutation relations of $ \mathfrak{m} : $
$  [ \tilde{x}_1,\tilde{x}_2 ]= (1-\alpha) \tilde{x}_5 $;
$ [\tilde{x}_1,\tilde{x}_3]= - \tilde{x}_6$;
$  [\tilde{x}_1,\tilde{x}_4]= - \alpha (\tilde{x}_5 + \tilde{x}_6)$;
$  [\tilde{x}_2,\tilde{x}_3]=\alpha \tilde{x}_5 $;
$  [\tilde{x}_2,\tilde{x}_4]=\alpha (\alpha \tilde{x}_5 - \tilde{x}_6)$;
$  [\tilde{x}_3,\tilde{x}_4]= - \alpha \tilde{x}_5 .$
\subsection{Case $ {\xi^2_1} \neq 0 , \xi^2_3 = 0, \xi^1_3 \neq \xi^2_4.$}
%{\fontsize{8}{10}\selectfont
\begin{equation}
\label{61general2}
J = \begin{pmatrix}
%$  J^1_1=( - ({\xi^4_1} {\xi^2_4} + {\xi^2_2} {\xi^2_1}))/{\xi^2_1};$\\
  - \frac{{\xi^4_1} {\xi^2_4} + {\xi^2_2} {\xi^2_1}}{{\xi^2_1}}&
%$  J^1_2=( - ({\xi^2_2}^2 + 1 + {\xi^4_2} {\xi^2_4}))/{\xi^2_1};$\\
- \frac{{\xi^2_2}^2 + 1 + {\xi^4_2} {\xi^2_4}}{{\xi^2_1}}&
\boxed{\xi^1_3}&
*&
0&
0\\
\boxed{\xi^2_1}&
\boxed{\xi^2_2}&
0&
\boxed{\xi^2_4}&
0&
0\\
*&*&
%$  J^3_3=( - ({\xi^5_5} {\xi^2_4} - {\xi^5_5} {\xi^1_3} + {\xi^2_2} {\xi^1_3}))/{\xi^2_4}
 - \frac{{\xi^5_5} {\xi^2_4} - {\xi^5_5} {\xi^1_3} + {\xi^2_2} {\xi^1_3}}{{\xi^2_4}} &
*&
0&
0\\
\boxed{\xi^4_1}&
\boxed{\xi^4_2}&
%$  J^4_3=( - {\xi^2_1} {\xi^1_3})/{\xi^2_4};$\\
 - \frac{{\xi^2_1} {\xi^1_3}}{{\xi^2_4}}&
*&
%$  J^4_4=({\xi^4_1} {\xi^2_4}^2 + {\xi^2_2} {\xi^2_1} {\xi^1_3} + ({\xi^2_4} - {\xi^1_3})
%{\xi^5_5} {\xi^2_1})/({\xi^2_4} {\xi^2_1});$\\
0&
0\\
*&
*&
*&
*&
\boxed{\xi^5_5}&
%$  J^5_6= - ({\xi^2_4} {\xi^1_3})/({\xi^2_4} - {\xi^1_3});$\\
\frac{{\xi^2_4} {\xi^1_3}}{{\xi^2_4} - {\xi^1_3}}\\
\boxed{\xi^6_1}&
\boxed{\xi^6_2}&
\boxed{\xi^6_3}&
\boxed{\xi^6_4}&
%$  J^6_5=(({\xi^5_5}^2 + 1) ({\xi^2_4} - {\xi^1_3}))/({\xi^2_4} {\xi^1_3});$\\
\frac{({\xi^5_5}^2 + 1) ({\xi^2_4} - {\xi^1_3})}{{\xi^2_4} {\xi^1_3}}&
-{\xi^5_5}
 \end{pmatrix}
\end{equation}
%}
where
{\fontsize{6}{10}\selectfont
$  \mathbf{J^1_4}=( - (({\xi^2_4} + {\xi^1_3}) {\xi^2_2} {\xi^2_1} + {\xi^4_1} {\xi^2_4}^2 + ({\xi^2_4} - {\xi^1_3}) {\xi^5_5} {\xi^2_1}))/{\xi^2_1}^2;$
\quad
$  \mathbf{J^3_1}=( - (({\xi^2_4} - {\xi^1_3}) {\xi^4_1} {\xi^2_2} - {\xi^4_2} {\xi^2_4} {\xi^2_1} -
({\xi^2_4} - {\xi^1_3}) {\xi^5_5} {\xi^4_1}))/({\xi^2_1} {\xi^1_3});$
\quad
$ \mathbf{J^3_2}=(({\xi^2_4} + {\xi^1_3}) {\xi^4_2} {\xi^2_2} {\xi^2_1} - ({\xi^2_2}^2 + 1) {\xi^4_1} {\xi^2_4} + ({\xi^2_4} - {\xi^1_3}) {\xi^5_5} {\xi^4_2} {\xi^2_1})/({\xi^2_1}^2 {\xi^1_3})
;$
\quad
$  \mathbf{J^3_4}=(((({\xi^2_4} - {\xi^1_3}) {\xi^5_5}^2 - ({\xi^2_4} - {\xi^1_3}) {\xi^5_5} {\xi^2_2} + ({\xi^2_2}^2 + 1) {\xi^2_4} + {\xi^4_2} {\xi^2_4}^2) ({\xi^2_4} - {\xi^1_3}) + ({\xi^2_2}^2 + 1 + {\xi^4_2} {\xi^2_4}) {\xi^2_4} {\xi^1_3}) {\xi^2_1} + (({\xi^2_4} + {\xi^1_3}) {\xi^2_2} {\xi^2_1} + {\xi^4_1} {\xi^2_4}^2) ({\xi^5_5} - {\xi^2_2}) ({\xi^2_4} - {\xi^1_3}))/({\xi^2_4} {\xi^2_1}^2 {\xi^1_3});$
\quad
$  \mathbf{J^4_4}=({\xi^4_1} {\xi^2_4}^2 + {\xi^2_2} {\xi^2_1} {\xi^1_3} + ({\xi^2_4} - {\xi^1_3})
{\xi^5_5} {\xi^2_1})/({\xi^2_4} {\xi^2_1});$
\quad
$  \mathbf{J^5_1}=( - ({\xi^6_4} {\xi^4_1} {\xi^2_1} {\xi^1_3} + {\xi^6_3} {\xi^5_5} {\xi^4_1} {\xi^2_4} - {\xi^6_3} {\xi^5_5} {\xi^4_1} {\xi^1_3} + {\xi^6_3} {\xi^4_2} {\xi^2_4} {\xi^2_1} - {\xi^6_3}
{\xi^4_1} {\xi^2_4} {\xi^2_2} + {\xi^6_3} {\xi^4_1} {\xi^2_2} {\xi^1_3} + {\xi^6_2} {\xi^2_1}^2
{\xi^1_3} - {\xi^6_1} {\xi^5_5} {\xi^2_1} {\xi^1_3} - {\xi^6_1} {\xi^4_1} {\xi^2_4} {\xi^1_3} - {\xi^6_1} {\xi^2_2} {\xi^2_1} {\xi^1_3}) {\xi^2_4})/(({\xi^5_5}^2 + 1) ({\xi^2_4} - {\xi^1_3}) {\xi^2_1});$
\quad
$  \mathbf{J^5_2}=( - (({\xi^6_4} {\xi^4_2} {\xi^2_1} - {\xi^6_2} {\xi^5_5} {\xi^2_1} + {\xi^6_2} {\xi^2_2} {\xi^2_1} - {\xi^6_1} {\xi^4_2} {\xi^2_4} - {\xi^6_1} {\xi^2_2}^2 - {\xi^6_1}) {\xi^2_1}
{\xi^1_3} + ({\xi^5_5} {\xi^4_2} {\xi^2_1} + {\xi^4_2} {\xi^2_2} {\xi^2_1} - {\xi^4_1} {\xi^2_2}^
2 - {\xi^4_1}) ({\xi^2_4} - {\xi^1_3}) {\xi^6_3} - (({\xi^2_2}^2 + 1) {\xi^4_1} - 2 {\xi^4_2}
 {\xi^2_2} {\xi^2_1}) {\xi^6_3} {\xi^1_3}) {\xi^2_4})/(({\xi^5_5}^2 + 1) ({\xi^2_4} - {\xi^1_3}
) {\xi^2_1}^2);$
\quad
$  \mathbf{J^5_3}=(({\xi^6_4} {\xi^2_1} {\xi^1_3} + 2 {\xi^6_3} {\xi^5_5} {\xi^2_4} - {\xi^6_3} {\xi^5_5} {\xi^1_3} + {\xi^6_3} {\xi^2_2} {\xi^1_3} - {\xi^6_1} {\xi^2_4} {\xi^1_3}) {\xi^1_3})/(({\xi^5_5}^2 + 1) ({\xi^2_4} - {\xi^1_3}));$
\quad
$  \mathbf{J^5_4}=((({\xi^2_4} - {\xi^1_3}) {\xi^6_4} {\xi^2_1} + {\xi^6_1} {\xi^2_4} {\xi^1_3} - ({\xi^5_5} - {\xi^2_2}) ({\xi^2_4} - {\xi^1_3}) {\xi^6_3}) (({\xi^2_4} + {\xi^1_3}) {\xi^2_2} {\xi^2_1} + {\xi^4_1} {\xi^2_4}^2) + (({\xi^6_4} {\xi^5_5} - {\xi^6_2} {\xi^2_4}) {\xi^2_1} + ({\xi^2_4} - {\xi^1_3}) {\xi^6_1} {\xi^5_5}) {\xi^2_4} {\xi^2_1} {\xi^1_3} - (({\xi^4_1} {\xi^2_4} + {\xi^2_2} {\xi^2_1}) ({\xi^2_4} - {\xi^1_3}) {\xi^2_4} + ({\xi^4_1} {\xi^2_4} + {\xi^2_2} {\xi^2_1})
{\xi^2_4} {\xi^1_3} + ({\xi^2_4} - {\xi^1_3}) {\xi^5_5} {\xi^2_1} {\xi^1_3}) {\xi^6_4} {\xi^2_1} -
 ((({\xi^2_4} - {\xi^1_3}) {\xi^5_5}^2 - ({\xi^2_4} - {\xi^1_3}) {\xi^5_5} {\xi^2_2} + ({\xi^2_2}^2 + 1) {\xi^2_4} + {\xi^4_2} {\xi^2_4}^2) ({\xi^2_4} - {\xi^1_3}) + ({\xi^2_2}^2 + 1 +
 {\xi^4_2} {\xi^2_4}) {\xi^2_4} {\xi^1_3}) {\xi^6_3} {\xi^2_1})/(({\xi^5_5}^2 + 1) ({\xi^2_4} -
 {\xi^1_3}) {\xi^2_1}^2);$
}
and the parameters are subject to the condition
\begin{equation}
\label{cond61general2}
{\xi^2_1} {\xi^2_4} \xi^1_3 ( {\xi^1_3} -{\xi^2_4})  \neq 0 .
 \end{equation}
Taking $u=1$ and suitable values for the $b^i_j$'s in (\ref{aut61}),
equivalence by $\Phi$
 switches to the case \ref{61_5}
${\xi^2_1}={\xi^2_3}=0$
below.
\subsection{Case $ {\xi^2_1} \neq 0 , \xi^1_3 = \xi^2_4.$}
\label{six3}
%{\fontsize{8}{10}\selectfont
\begin{equation}
\label{61general3}
J = \begin{pmatrix}
\boxed{\xi^1_1}&
\boxed{\xi^1_2}&
0&
0&
0&
0\\
%$  J^2_1=( - ({\xi^1_1}^2 + 1))/{\xi^1_2};$\\
 - \frac{{\xi^1_1}^2 + 1}{{\xi^1_2}}&
-\xi^1_1&
0&
0&
0&
0\\
%$  J^3_1={\xi^4_2};$\\
{\xi^4_2}&
%$  J^3_2=((2 {\xi^4_2} {\xi^1_1} - {\xi^4_1} {\xi^1_2}) {\xi^1_2})/({\xi^1_1}^2 + 1)
%;$\\
\frac{(2 {\xi^4_2} {\xi^1_1} - {\xi^4_1} {\xi^1_2}) {\xi^1_2}}{{\xi^1_1}^2 + 1} &
%$  J^3_3={\xi^1_1};$\\
{\xi^1_1}&
%$  J^3_4= - {\xi^1_2};$\\
- {\xi^1_2}&
0&
0\\
\boxed{\xi^4_1}&
\boxed{\xi^4_2}&
%$  J^4_3=({\xi^1_1}^2 + 1)/{\xi^1_2};$\\
\frac{{\xi^1_1}^2 + 1}{{\xi^1_2}}&
-\xi^1_1&
0&
0\\
*&
*&
%$  J^5_3=(({\xi^5_5} - {\xi^1_1}) {\xi^6_3} {\xi^1_2} - ({\xi^1_1}^2 + 1) {\xi^6_4})/({\xi^6_5} {\xi^1_2});$\\
\frac{({\xi^5_5} - {\xi^1_1}) {\xi^6_3} {\xi^1_2} - ({\xi^1_1}^2 + 1) {\xi^6_4}}{{\xi^6_5} {\xi^1_2}}&
%$  J^5_4=({\xi^6_4} {\xi^5_5} + {\xi^6_4} {\xi^1_1} + {\xi^6_3} {\xi^1_2})/{\xi^6_5};$\\
\frac{{\xi^6_4} {\xi^5_5} + {\xi^6_4} {\xi^1_1} + {\xi^6_3} {\xi^1_2}}{{\xi^6_5}}&
\boxed{\xi^5_5}&
%$  J^5_6=( - ({\xi^5_5}^2 + 1))/{\xi^6_5};$\\
 - \frac{{\xi^5_5}^2 + 1}{{\xi^6_5}}
\\
\boxed{\xi^6_1}&
\boxed{\xi^6_2}&
\boxed{\xi^6_3}&
\boxed{\xi^6_4}&
\boxed{\xi^6_5}&
-{\xi^5_5}
 \end{pmatrix}
\end{equation}
%}
where
{\fontsize{6}{10}\selectfont
$  \mathbf{J^5_1}=(({\xi^5_5} - {\xi^1_1}) {\xi^6_1} {\xi^1_2} + ({\xi^1_1}^2 + 1) {\xi^6_2} - {\xi^6_3} {\xi^4_2} {\xi^1_2} - {\xi^6_4} {\xi^4_1} {\xi^1_2})/({\xi^6_5} {\xi^1_2});$
\quad
$  \mathbf{J^5_2}=( - ((2 {\xi^4_2} {\xi^1_1} - {\xi^4_1} {\xi^1_2}) {\xi^6_3} {\xi^1_2} + ({\xi^1_1}
^2 + 1) {\xi^6_4} {\xi^4_2} - ({\xi^6_2} {\xi^5_5} + {\xi^6_2} {\xi^1_1} - {\xi^6_1} {\xi^1_2})
 ({\xi^1_1}^2 + 1)))/(({\xi^1_1}^2 + 1) {\xi^6_5});$
}
and the parameters are subject to the condition
\begin{equation}
\label{cond61general3}
{\xi^1_2} {\xi^6_5}   \neq 0 .
 \end{equation}
Taking $u= \frac{\xi^6_5}{1+{\xi^5_5}^2}$ and suitable values for the $b^i_j$'s in (\ref{aut61}),
equivalence by $\Phi$
leads to the case where
$ {\xi^1_1}={\xi^4_1}={\xi^4_2}={\xi^5_5}={\xi^6_1}={\xi^6_2}={\xi^6_3}={\xi^6_4} =0 $
and $ {\xi^1_2}={\xi^6_5}=1$:
\begin{equation}
\label{6_1equivfinal2}
J = \text{diag} \left(
\left( \begin{smallmatrix} 0&1\\-1&0 \end{smallmatrix}\right),
\left( \begin{smallmatrix} 0&-1\\1&0 \end{smallmatrix}\right),
\left( \begin{smallmatrix} 0&-1\\1&0 \end{smallmatrix}\right)\right)
\,
.\end{equation}
$ J$  is not equivalent to any $J_{\alpha}.$
$ \mathfrak{m}$  is here an abelian algebra.
\subsection{Case $ {\xi^2_1} = 0 , \xi^2_3 \neq 0.$}
%{\fontsize{8}{10}\selectfont
\begin{equation}
\label{61general4}
J = \begin{pmatrix}
\boxed{\xi^1_1}&
*&
*&
*&
0&
0\\
0&
\boxed{\xi^2_2}&
\boxed{\xi^2_3}&
\boxed{\xi^2_4}&
0&
0\\
%$  J^3_1=( - {\xi^4_1} {\xi^2_4})/{\xi^2_3};$\\
 - \frac{{\xi^4_1} {\xi^2_4}}{{\xi^2_3}}&
*&
%$  J^3_3= - ({\xi^6_5} {\xi^5_5} {\xi^2_4} - {\xi^6_5} {\xi^2_4} {\xi^1_1} + {\xi^5_5}^2 {\xi^2_2} + {\xi^2_2})/({\xi^5_5}^2 + 1);$\\
- \frac{{\xi^6_5} {\xi^5_5} {\xi^2_4} - {\xi^6_5} {\xi^2_4} {\xi^1_1} + {\xi^5_5}^2 {\xi^2_2} + {\xi^2_2}}{{\xi^5_5}^2 + 1}&
*&
0&
0\\
\boxed{\xi^4_1}&
*&
%$  J^4_3=(({\xi^5_5} - {\xi^1_1}) {\xi^6_5} {\xi^2_3})/({\xi^5_5}^2 + 1);$\\
\frac{({\xi^5_5} - {\xi^1_1}) {\xi^6_5} {\xi^2_3}}{{\xi^5_5}^2 + 1}&
%$  J^4_4= - (({\xi^5_5}^2 + 1) {\xi^1_1} - ({\xi^5_5} - {\xi^1_1}) {\xi^6_5} {\xi^2_4})/(
%{\xi^5_5}^2 + 1);$\\
- \frac{({\xi^5_5}^2 + 1) {\xi^1_1} - ({\xi^5_5} - {\xi^1_1}) {\xi^6_5} {\xi^2_4}}{
{\xi^5_5}^2 + 1}&
0&
0\\
*&
*&
*&
*&
\boxed{\xi^5_5}&
%$  J^5_6=( - ({\xi^5_5}^2 + 1))/{\xi^6_5};$\\
 - \frac{{\xi^5_5}^2 + 1}{{\xi^6_5}}
\\
\boxed{\xi^6_1}&
\boxed{\xi^6_2}&
\boxed{\xi^6_3}&
\boxed{\xi^6_4}&
\boxed{\xi^6_5}&
-{\xi^5_5}
 \end{pmatrix}
\end{equation}
%}
where
{\fontsize{6}{10}\selectfont
$  \mathbf{J^1_2}=(({\xi^5_5}^2 + 1) ({\xi^2_2} - {\xi^1_1}) {\xi^4_1} {\xi^2_4} + ({\xi^5_5} - {\xi^2_2}) ({\xi^1_1}^2 + 1) {\xi^6_5} {\xi^2_3})/(({\xi^5_5}^2 + 1) {\xi^4_1} {\xi^2_3})
;$
\myquad
$  \mathbf{J^1_3}=(({\xi^5_5}^2 + 1) {\xi^4_1} {\xi^2_4} - ({\xi^1_1}^2 + 1) {\xi^6_5} {\xi^2_3})
/(({\xi^5_5}^2 + 1) {\xi^4_1});$\myquad
$  \mathbf{J^1_4}=( - (({\xi^5_5}^2 + 1) ({\xi^1_1}^2 + 1) {\xi^2_3} - ({\xi^5_5}^2 + 1) {\xi^4_1} {\xi^2_4}^2 + ({\xi^1_1}^2 + 1) {\xi^6_5} {\xi^2_4} {\xi^2_3}))/(({\xi^5_5}^2 + 1)
{\xi^4_1} {\xi^2_3});$\myquad
$  \mathbf{J^3_2}=( - ((({\xi^2_2}^2 + 1) {\xi^2_3} - {\xi^4_1} {\xi^2_4}^2) ({\xi^5_5}^2 + 1)
 - ({\xi^2_2} {\xi^1_1} - 1 - ({\xi^2_2} + {\xi^1_1}) {\xi^5_5}) {\xi^6_5} {\xi^2_4} {\xi^2_3}))
/(({\xi^5_5}^2 + 1) {\xi^2_3}^2);$\myquad
$  \mathbf{J^3_4}=( - ({\xi^6_5} {\xi^5_5} {\xi^2_4} - {\xi^6_5} {\xi^2_4} {\xi^1_1} + {\xi^5_5}^2
{\xi^2_2} - {\xi^5_5}^2 {\xi^1_1} + {\xi^2_2} - {\xi^1_1}) {\xi^2_4})/(({\xi^5_5}^2 + 1) {\xi^2_3});$\myquad
$  \mathbf{J^4_2}=( - (({\xi^2_2} {\xi^1_1} - 1 - ({\xi^2_2} + {\xi^1_1}) {\xi^5_5}) {\xi^6_5} {\xi^2_3} + ({\xi^5_5}^2 + 1) {\xi^4_1} {\xi^2_4}))/(({\xi^5_5}^2 + 1) {\xi^2_3});$
\myquad
$  \mathbf{J^5_1}=( - ({\xi^6_4} {\xi^4_1} {\xi^2_3} - {\xi^6_3} {\xi^4_1} {\xi^2_4} - {\xi^6_1} {\xi^5_5} {\xi^2_3} + {\xi^6_1} {\xi^2_3} {\xi^1_1}))/({\xi^6_5} {\xi^2_3});$\myquad
$  \mathbf{J^5_2}=((({\xi^5_5}^2 + 1) ({\xi^5_5} - {\xi^2_2}) {\xi^6_2} {\xi^4_1} {\xi^2_3} - ({\xi^5_5}^2 + 1) ({\xi^2_2} - {\xi^1_1}) {\xi^6_1} {\xi^4_1} {\xi^2_4} - ({\xi^5_5} - {\xi^2_2}) ({\xi^1_1}^2 + 1) {\xi^6_5} {\xi^6_1} {\xi^2_3} - (({\xi^2_2} + {\xi^1_1}) {\xi^5_5} + {\xi^1_1}
^2 + 1) {\xi^6_5} {\xi^6_4} {\xi^4_1} {\xi^2_3} - ({\xi^2_2} {\xi^1_1} - 1 - ({\xi^2_2} + {\xi^1_1}) {\xi^5_5}) {\xi^6_5} {\xi^6_3} {\xi^4_1} {\xi^2_4}) {\xi^2_3} + (({\xi^2_2}^2 + 1) {\xi^2_3} - {\xi^4_1} {\xi^2_4}^2) ({\xi^5_5}^2 + 1) {\xi^6_3} {\xi^4_1} + (({\xi^5_5}^2 + 1)
{\xi^4_1} {\xi^2_4} + ({\xi^2_2} + {\xi^1_1}) {\xi^6_5} {\xi^2_3} {\xi^1_1}) {\xi^6_4} {\xi^4_1}
{\xi^2_3})/(({\xi^5_5}^2 + 1) {\xi^6_5} {\xi^4_1} {\xi^2_3}^2);$\myquad
$  \mathbf{J^5_3}=( - (((({\xi^5_5} - {\xi^1_1}) {\xi^6_4} {\xi^4_1} - ({\xi^1_1}^2 + 1) {\xi^6_1}
) {\xi^2_3} - ({\xi^2_2} - {\xi^1_1}) {\xi^6_3} {\xi^4_1} {\xi^2_4}) {\xi^6_5} + ({\xi^6_2} {\xi^2_3} + {\xi^6_1} {\xi^2_4} - {\xi^6_3} {\xi^5_5}) ({\xi^5_5}^2 + 1) {\xi^4_1} - (({\xi^5_5}^2
 + 1) {\xi^2_2} + ({\xi^5_5} - {\xi^2_2}) {\xi^6_5} {\xi^2_4}) {\xi^6_3} {\xi^4_1}))/(({\xi^5_5}
^2 + 1) {\xi^6_5} {\xi^4_1});$\myquad
$  \mathbf{J^5_4}=( - (((({\xi^5_5} - {\xi^1_1}) {\xi^6_4} {\xi^4_1} - ({\xi^1_1}^2 + 1) {\xi^6_1}
) {\xi^6_5} {\xi^2_4} - ({\xi^5_5}^2 + 1) ({\xi^1_1}^2 + 1) {\xi^6_1}) {\xi^2_3} + (({\xi^6_2} {\xi^2_3} + {\xi^6_1} {\xi^2_4}) {\xi^2_4} - ({\xi^5_5} + {\xi^1_1}) {\xi^6_4} {\xi^2_3}) ({\xi^5_5}^2 + 1) {\xi^4_1} + (({\xi^1_1}^2 + 1) {\xi^2_3} - {\xi^4_1} {\xi^2_4}^2) ({\xi^5_5}
 - {\xi^2_2}) {\xi^6_5} {\xi^6_3} - (({\xi^5_5}^2 + 1) ({\xi^2_2} - {\xi^1_1}) {\xi^4_1} {\xi^2_4} + ({\xi^5_5} - {\xi^2_2}) ({\xi^1_1}^2 + 1) {\xi^6_5} {\xi^2_3} + ({\xi^2_2} - {\xi^1_1})
 {\xi^6_5} {\xi^4_1} {\xi^2_4}^2) {\xi^6_3}))/(({\xi^5_5}^2 + 1) {\xi^6_5} {\xi^4_1} {\xi^2_3}
);$
}
and the parameters are subject to the condition
\begin{equation}
\label{cond61general4}
{\xi^2_3}{\xi^4_1} {\xi^6_5}   \neq 0 .
 \end{equation}
Taking suitable values for $u$ and the $b^i_j$'s in (\ref{aut61}),
equivalence by $\Phi$ switches to the case \ref{61_1},
more precisely $ {\xi^2_1}={\xi^2_3}={\xi^2_4}=1, {\xi^1_3}=0.$
Hence $J$ is equivalent to $J_{\alpha}$ in (\ref{61equivfinal1}) with $\alpha=1.$
\subsection{Case $ {\xi^2_1} = 0 , \xi^2_3 = 0.$}
\label{61_5}
%{\fontsize{8}{10}\selectfont
\begin{equation}
\label{61general5}
J = \begin{pmatrix}
%$  J^1_1=({\xi^5_5} {\xi^2_4} - {\xi^5_5} {\xi^1_3} + {\xi^2_2} {\xi^1_3})/{\xi^2_4};$\\
\frac{{\xi^5_5} {\xi^2_4} - {\xi^5_5} {\xi^1_3} + {\xi^2_2} {\xi^1_3}}{{\xi^2_4}}&
\boxed{\xi^1_2}&
\boxed{\xi^1_3}&
%$  J^1_4=(({\xi^2_4} + {\xi^1_3}) {\xi^2_2} {\xi^1_2} + {\xi^3_2} {\xi^2_4} {\xi^1_3} + ({\xi^2_4} - {\xi^1_3}) {\xi^5_5} {\xi^1_2})/({\xi^2_2}^2 + 1);$\\
\frac{({\xi^2_4} + {\xi^1_3}) {\xi^2_2} {\xi^1_2} + {\xi^3_2} {\xi^2_4} {\xi^1_3} + ({\xi^2_4} - {\xi^1_3}) {\xi^5_5} {\xi^1_2}}{{\xi^2_2}^2 + 1}&
0&
0\\
0&
\boxed{\xi^2_2}&
0&
\boxed{\xi^2_4}&
0&
0\\
*&
\boxed{\xi^3_2}&
%$  J^3_3=( - ({\xi^5_5} {\xi^2_4} - {\xi^5_5} {\xi^1_3} + {\xi^2_2} {\xi^1_3}))/{\xi^2_4}
-\frac{{\xi^5_5} {\xi^2_4} - {\xi^5_5} {\xi^1_3} + {\xi^2_2} {\xi^1_3}}{{\xi^2_4}}&
*&
0&
0\\
0&
%$  J^4_2=( - ({\xi^2_2}^2 + 1))/{\xi^2_4};$\\
- \frac{{\xi^2_2}^2 + 1}{{\xi^2_4}}&
0&
%$  J^4_4= - {\xi^2_2};$\\
- {\xi^2_2}&
0&
0\\
*&
*&
*&
*&
\boxed{\xi^5_5}&
%$  J^5_6= - ({\xi^2_4} {\xi^1_3})/({\xi^2_4} - {\xi^1_3});$\\
- \frac{{\xi^2_4} {\xi^1_3}}{{\xi^2_4} - {\xi^1_3}}
\\
\boxed{\xi^6_1}&
\boxed{\xi^6_2}&
\boxed{\xi^6_3}&
\boxed{\xi^6_4}&
%$  J^6_5=(({\xi^5_5}^2 + 1) ({\xi^2_4} - {\xi^1_3}))/({\xi^2_4} {\xi^1_3});$\\
\frac{({\xi^5_5}^2 + 1) ({\xi^2_4} - {\xi^1_3})}{{\xi^2_4} {\xi^1_3}}&
-{\xi^5_5}
 \end{pmatrix}
\end{equation}
%}
where
{\fontsize{6}{10}\selectfont
$  \mathbf{J^3_1}=( - (({\xi^2_4} - {\xi^1_3}) {\xi^5_5}^2 {\xi^2_4} - ({\xi^2_4} - {\xi^1_3}) {\xi^5_5}^2 {\xi^1_3} + 2 ({\xi^2_4} - {\xi^1_3}) {\xi^5_5} {\xi^2_2} {\xi^1_3} + ({\xi^2_2}^2 +
1) {\xi^2_4} {\xi^1_3} + ({\xi^2_4} - {\xi^2_2}^2 {\xi^1_3}) ({\xi^2_4} - {\xi^1_3})))/({\xi^2_4}^2 {\xi^1_3});$\myquad
$  \mathbf{J^3_4}=( - ((({\xi^2_4} + {\xi^1_3}) {\xi^2_2} {\xi^1_2} + {\xi^3_2} {\xi^2_4} {\xi^1_3})
({\xi^2_4} - {\xi^1_3}) {\xi^5_5} - (({\xi^2_4} + {\xi^1_3}) {\xi^2_2} {\xi^1_2} + {\xi^3_2} {\xi^2_4} {\xi^1_3}) ({\xi^2_4} - {\xi^1_3}) {\xi^2_2} + (({\xi^2_4} - {\xi^1_3})^2 {\xi^5_5}^2 -
 ({\xi^2_4} - {\xi^1_3})^2 {\xi^5_5} {\xi^2_2} + ({\xi^2_4} - {\xi^1_3}) ({\xi^2_2}^2 + 1)
{\xi^2_4} + ({\xi^2_2}^2 + 1) {\xi^2_4} {\xi^1_3}) {\xi^1_2}))/(({\xi^2_2}^2 + 1) {\xi^2_4}
{\xi^1_3});$\myquad
$  \mathbf{J^5_1}=(({\xi^5_5} - {\xi^2_2}) {\xi^6_1} {\xi^2_4} {\xi^1_3}^2 + ({\xi^2_4} - {\xi^2_2}
^2 {\xi^1_3}) ({\xi^2_4} - {\xi^1_3}) {\xi^6_3} + (({\xi^2_4} - {\xi^1_3}) {\xi^5_5}^2 {\xi^2_4} - ({\xi^2_4} - {\xi^1_3}) {\xi^5_5}^2 {\xi^1_3} + 2 ({\xi^2_4} - {\xi^1_3}) {\xi^5_5} {\xi^2_2} {\xi^1_3} + ({\xi^2_2}^2 + 1) {\xi^2_4} {\xi^1_3}) {\xi^6_3})/(({\xi^5_5}^2 + 1) ({\xi^2_4} - {\xi^1_3}) {\xi^2_4});$\myquad
$  \mathbf{J^5_2}=((({\xi^6_2} {\xi^5_5} - {\xi^6_2} {\xi^2_2} - {\xi^6_1} {\xi^1_2} - {\xi^6_3} {\xi^3_2}) {\xi^2_4} + ({\xi^2_2}^2 + 1) {\xi^6_4}) {\xi^1_3})/(({\xi^5_5}^2 + 1) ({\xi^2_4} -
{\xi^1_3}));$\myquad
$  \mathbf{J^5_3}=((2 {\xi^6_3} {\xi^5_5} {\xi^2_4} - {\xi^6_3} {\xi^5_5} {\xi^1_3} + {\xi^6_3} {\xi^2_2} {\xi^1_3} - {\xi^6_1} {\xi^2_4} {\xi^1_3}) {\xi^1_3})/(({\xi^5_5}^2 + 1) ({\xi^2_4} - {\xi^1_3}));$\myquad
$  \mathbf{J^5_4}=((({\xi^2_4} - {\xi^1_3})^2 {\xi^5_5}^2 - ({\xi^2_4} - {\xi^1_3})^2 {\xi^5_5}
{\xi^2_2} + ({\xi^2_4} - {\xi^1_3}) ({\xi^2_2}^2 + 1) {\xi^2_4} + ({\xi^2_2}^2 + 1) {\xi^2_4} {\xi^1_3}) {\xi^6_3} {\xi^1_2} + ({\xi^6_4} {\xi^5_5} {\xi^2_2}^2 + {\xi^6_4} {\xi^5_5} + {\xi^6_4} {\xi^2_2}^3 + {\xi^6_4} {\xi^2_2} - {\xi^6_2} {\xi^2_4} {\xi^2_2}^2 - {\xi^6_2} {\xi^2_4}
 - {\xi^6_1} {\xi^5_5} {\xi^2_4} {\xi^1_2} + {\xi^6_1} {\xi^5_5} {\xi^1_3} {\xi^1_2}) {\xi^2_4} {\xi^1_3} + (({\xi^2_4} + {\xi^1_3}) {\xi^2_2} {\xi^1_2} + {\xi^3_2} {\xi^2_4} {\xi^1_3}) (({\xi^2_4} - {\xi^1_3}) {\xi^6_3} {\xi^5_5} - ({\xi^2_4} - {\xi^1_3}) {\xi^6_3} {\xi^2_2} - {\xi^6_1} {\xi^2_4} {\xi^1_3}))/(({\xi^5_5}^2 + 1) ({\xi^2_4} - {\xi^1_3}) ({\xi^2_2}^2 + 1));$
}
and the parameters are subject to the condition
\begin{equation}
\label{cond61general5}
{\xi^1_3}{\xi^2_4} (\xi^2_4 -\xi^1_3)   \neq 0 .
 \end{equation}
$\bullet$ Suppose first $\xi^2_4 \neq -\xi^1_3,$
taking suitable values for $u$ and the $b^i_j$'s in (\ref{aut61}),
equivalence by $\Phi$ leads to the case where
$ {\xi^1_2}={\xi^2_2}={\xi^3_2}={\xi^5_5}={\xi^6_1}={\xi^6_2}={\xi^6_3}={\xi^6_4} =0 $
and $ {\xi^1_3}=1, {\xi^2_4}= \beta \neq 0, \pm1$:
\begin{equation}
{\fontsize{8}{10}\selectfont
\label{61equiv51}
J'_{\beta} = \begin{pmatrix}
0&
0&
1&
0&
0&
0\\
0&
0&
0&
\beta&
0&
0\\
-1&
0&
0&
0&
0&
0\\
0&
 - 1/\beta&
0&
0&
0&
0\\
0&
0&
0&
0&
0&
 - \beta/(\beta - 1)\\
0&
0&
0&
0&
(\beta - 1)/\beta&
0\end{pmatrix}
\quad \quad (\beta \neq 0,\pm 1).
}
\end{equation}
The  $ J$'s corresponding to different $\beta,\beta'$
are not equivalent unless $ \beta'=1/\beta.$
$J'_{\beta}$ is equivalent to $J_{\alpha}$ in
(\ref{61equivfinal1})
if and only if $ \alpha=\frac{(\beta-1)^2}{\beta}.$
\\$\bullet$ Suppose now  ${\xi^2_4} = -{\xi^1_3}$ in (\ref{61general5}).
Taking suitable values for $u$ and the $b^i_j$'s in (\ref{aut61}),
equivalence by $\Phi$ leads to the case where
$ {\xi^1_2}={\xi^2_2}={\xi^5_5}={\xi^6_1}={\xi^6_2}={\xi^6_3}={\xi^6_4} =0 $
and $ {\xi^1_3}=1, {\xi^2_4}=-1, \xi^3_2 = \gamma :$
\begin{equation}
{\fontsize{8}{10}\selectfont
\label{61equiv52}
J''_{\gamma} = \begin{pmatrix}
0&
0&
1&
 - \gamma&
0&
0\\
0&
0&
0&
-1&
0&
0\\
-1&
\gamma&
0&
0&
0&
0\\
0&
1&
0&
0&
0&
0\\
0&
0&
0&
0&
0&
 - 1/2\\
0&
0&
0&
0&
2&
0\end{pmatrix}
\quad \quad (\gamma \in \Rmath).
}
\end{equation}
\par Commutation relations of $ \mathfrak{m} : $
$  [\tilde{x}_1,\tilde{x}_2]=\tilde{x}_5;$
$  [\tilde{x}_1,\tilde{x}_4]=2 \tilde{x}_6;$
$  [\tilde{x}_2,\tilde{x}_3]=2 \tilde{x}_6;$
$ [\tilde{x}_2,\tilde{x}_4]= - 2 \gamma \tilde{x}_6 ;$
$  [\tilde{x}_3,\tilde{x}_4]=\tilde{x}_5.$\par
$J''_{\gamma}$
is not equivalent to any $J_{\alpha}$
in (\ref{61equivfinal1}) nor to (\ref{6_1equivfinal2}), and
each $J''_{\gamma} \, (\gamma \neq 0)$
is equivalent to $J''_{1}.$
$J''_{1}$ is not equivalent to $J''_{0}.$
%\end{document}
\subsection{Conclusions.}
One has with obvious notations
\begin{equation}
\label{decompx6_1}
\mathfrak{X}_{6,1} =
\mathfrak{X}_{\xi^2_1 \neq 0} \cup \mathfrak{X}_{\xi^2_3 \neq 0}
\cup \mathfrak{X}_{\xi^2_4 \neq 0}
\end{equation}
and $$  \mathfrak{X}_{\xi^2_3 \neq 0} \subset \mathfrak{X}_{\xi^1_3 \neq \xi^2_4}.$$
It can be seen that the formula (\ref{61general1}), which still makes sense for
$\xi^6_5 \neq 0$ under the  only assumption that
$\xi^2_3 \neq 0,$ yields all  of  $\mathfrak{X}_{\xi^2_3 \neq 0}.$
Hence $\mathfrak{X}_{\xi^2_3 \neq 0}$ is a 12-dimensional submanifold of $\Rmath^{36}$ with a global chart.
The automorphism
$ \Phi = \text{diag} \left(1,1,
\left( \begin{smallmatrix} 0&1\\1&0 \end{smallmatrix}\right),
1,1\right)  $
 switches $\xi^2_3$ and $\xi^2_4,$ hence
$\mathfrak{X}_{\xi^2_4 \neq 0} = \Phi \, \mathfrak{X}_{\xi^2_3 \neq 0}\,  {\Phi}^{-1}$ is
also  a 12-dimensional submanifold of $\Rmath^{36}.$
Consider now $\mathfrak{X}_{\xi^2_1 \neq 0}.$ There are 3 subcases:
\begin{eqnarray}
\label{caselocalsmnfld1}
\xi^2_3 & \neq 0& \quad (\text{which implies } \xi^1_3 \neq \xi^2_4 \text{ from section \ref{six3}})\\
\label{caselocalsmnfld2}
\xi^2_3 &=0&  \quad \text{and } \quad \xi^1_3 \neq \xi^2_4\\
\label{caselocalsmnfld3}
\xi^2_3 &=0&  \quad \text{and } \quad \xi^1_3 = \xi^2_4.
\end{eqnarray}
To prove that
$\mathfrak{X}_{\xi^2_1 \neq 0}$
is a 12-dimensional submanifold of $\Rmath^{36},$ it is sufficient to prove that it is a local
submanifold in the neighborhood of any of its points.
Take any $K \in \mathfrak{X}_{\xi^2_1 \neq 0} :$  $K =(\xi^i_j(K)).$
In case (\ref{caselocalsmnfld1}),
$K \in \mathfrak{X}_{\xi^2_3 \neq 0}$
and then, from the section \ref{61_1},
$\mathfrak{X}_{\xi^2_1 \neq 0}$ is
a local
12-dimensional submanifold of $\Rmath^{36}.$
Suppose now $K$ belongs in case (\ref{caselocalsmnfld2}) or (\ref{caselocalsmnfld3}).
To solve the initial system comprised of all the torsion equations and the equation $J^2=-1$ in $\Rmath^{36}$ in the neighborhood of $K,$
one has to complete first a set of common steps,  and then we are left with solving the
system $S$ of the remaining equations in the 15 variables
$\xi^1_1,\xi^1_2,\xi^1_3,\xi^2_1,\xi^2_2,\xi^2_3,\xi^2_4,\xi^4_1,\xi^4_2,\xi^5_5,
\xi^6_1,\xi^6_2,\xi^6_3,\xi^6_4,\xi^6_5$ in the open subset
$\xi^6_5  \neq 0$ of $\Rmath^{15}.$
Among these equations, we single out the 3 following equations :
 \begin{equation}
 \label{f,g,h}
 \begin{cases}
 f =0\\ g=0 \\h=0
 \end{cases}
 \end{equation}
 where :
$ f= {J^2}^2_1=( - {\xi^6_5} {\xi^5_5} {\xi^2_3} {\xi^2_2} - {\xi^6_5} {\xi^5_5} {\xi^2_3} {\xi^1_1} + {\xi^6_5} {\xi^2_3} {\xi^2_2} {\xi^1_1} - {\xi^6_5} {\xi^2_3} {\xi^2_1} {\xi^1_2} - {\xi^6_5} {\xi^2_3} + {\xi^5_5}^2 {\xi^4_2} {\xi^2_3} + {\xi^5_5}^2 {\xi^4_1} {\xi^2_4} + {\xi^5_5}
^2 {\xi^2_2} {\xi^2_1} + {\xi^5_5}^2 {\xi^2_1} {\xi^1_1} + {\xi^4_2} {\xi^2_3} + {\xi^4_1} {\xi^2_4} + {\xi^2_2} {\xi^2_1} + {\xi^2_1} {\xi^1_1})/({\xi^5_5}^2 + 1);$
\quad
$g=  {J^2}^2_2=( - {\xi^6_5} {\xi^5_5} {\xi^2_3} {\xi^2_2}^2 - {\xi^6_5} {\xi^5_5} {\xi^2_3}
 {\xi^2_1} {\xi^1_2} - {\xi^6_5} {\xi^5_5} {\xi^2_3} + {\xi^6_5} {\xi^2_3} {\xi^2_2}^2 {\xi^1_1}
- {\xi^6_5} {\xi^2_3} {\xi^2_2} {\xi^2_1} {\xi^1_2} + {\xi^6_5} {\xi^2_3} {\xi^1_1} + {\xi^5_5}^2
 {\xi^4_2} {\xi^2_4} {\xi^2_1} + {\xi^5_5}^2 {\xi^4_2} {\xi^2_3} {\xi^2_2} - {\xi^5_5}^2 {\xi^4_2} {\xi^2_3} {\xi^1_1} + {\xi^5_5}^2 {\xi^4_1} {\xi^2_3} {\xi^1_2} + {\xi^5_5}^2 {\xi^2_2}^2
{\xi^2_1} + {\xi^5_5}^2 {\xi^2_1}^2 {\xi^1_2} + {\xi^4_2} {\xi^2_4} {\xi^2_1} + {\xi^4_2} {\xi^2_3} {\xi^2_2} - {\xi^4_2} {\xi^2_3} {\xi^1_1} + {\xi^4_1} {\xi^2_3} {\xi^1_2} + {\xi^2_2}^2 {\xi^2_1} + {\xi^2_1}^2 {\xi^1_2})/({\xi^2_1} ({\xi^5_5}^2 + 1));$\quad
$ h= {J^2}^2_3=({\xi^6_5} {\xi^5_5} {\xi^2_4} {\xi^2_3} - {\xi^6_5} {\xi^5_5} {\xi^2_3} {\xi^1_3} - {\xi^6_5} {\xi^2_4} {\xi^2_3} {\xi^1_1} + {\xi^6_5} {\xi^2_4} {\xi^2_1} {\xi^1_3} - {\xi^6_5}
 {\xi^2_3}^2 {\xi^1_2} + {\xi^6_5} {\xi^2_3} {\xi^2_2} {\xi^1_3} - {\xi^5_5}^2 {\xi^2_4} {\xi^2_1} + {\xi^5_5}^2 {\xi^2_1} {\xi^1_3} - {\xi^2_4} {\xi^2_1} + {\xi^2_1} {\xi^1_3})/({\xi^5_5}^2
 + 1).$
 \\
The solution $J$ that's looked for is of type $\xi^2_1 \neq 0$ hence belongs in one of the 3 cases
(\ref{caselocalsmnfld1}),
(\ref{caselocalsmnfld2}),
(\ref{caselocalsmnfld3}).
If $J$ belongs in
 case  (\ref{caselocalsmnfld2}) or (\ref{caselocalsmnfld3}), the system $S$ is equivalent to the 3 equations
 (\ref{f,g,h}).
If $J$ belongs in the case  (\ref{caselocalsmnfld1}), the system $S$ is equivalent to the 3 equations
 (\ref{f,g,h}) if and only if $c(J) \neq 0$ where
$ c(J)=
 ({\xi^5_5}^2 +1)\xi^2_1 +\xi^2_1\xi^2_4 \xi^6_5 -(\xi^5_5 -\xi^2_2)\xi^2_3\xi^6_5.
 $
Now, if $K$ belongs in case (\ref{caselocalsmnfld2}),
$c(K)=
 \xi^2_1(K)({\xi^5_5(K)}^2 +1)+\xi^2_4(K) \xi^6_5(K)) =
 \frac{1}{\xi^1_3(K)}\, \xi^2_1(K)\xi^2_4(K)({\xi^5_5(K)}^2 +1)\neq 0 $
 since in that case $\xi^6_5(K) = \frac{ {\xi^5_5(K)}^2 +1}{\xi^2_4(K) \xi^1_3(K)} \,
 (\xi^2_4(K)-\xi^1_3(K))$
 (see (\ref{61general2})).
If $K$ belongs in case (\ref{caselocalsmnfld3}),
$c(K)= \xi^2_1(K)({\xi^5_5(K)}^2 +1)\neq 0$
 (see (\ref{61general3})).
 Hence in both cases, one has $c(J)\neq 0$ in some neighborhood of $K$ and
the remaining system is equivalent in that neighborhood to the 3 equations
 (\ref{f,g,h}).
 We will now show that the system
 (\ref{f,g,h}) is of maximal rank 3 at $K,$  that is some 3-jacobian doesn't vanish.
\\ $\bullet$
Suppose $K$ belongs in case (\ref{caselocalsmnfld2}). Then
$\frac{D(f,g,h)}{D(\xi^1_1,\xi^1_2,\xi^6_5)}(K) = -\frac{\xi^2_4(K)\xi^2_1(K)^3\xi^1_3(K)}{{\xi^5_5(K)}^2+1} \neq 0.$
\\ $\bullet$
Suppose $K$ belongs in case (\ref{caselocalsmnfld3}). Then
$\frac{D(f,g,h)}{D(\xi^1_1,\xi^1_2,\xi^2_4)}(K) = -\xi^2_1(K)^3 \neq 0.$
\\ Hence  the system
 (\ref{f,g,h}) is of maximal rank 3 at $K,$ and it follows that
$\mathfrak{X}_{\xi^2_1 \neq 0}$ is
is a local
submanifold in the neighborhood of $K.$
\par
Hence
$\mathfrak{X}_{\xi^2_1 \neq 0}$ is
is a 12-dimensional submanifold of $\Rmath^{36},$
and so is
$\mathfrak{X}_{6,1}$ from (\ref{decompx6_1}).
Any element of
$\mathfrak{X}_{6,1}$ is
is equivalent to either $J$ in
(\ref{6_1equivfinal2}),
or   $J_{\alpha} (\alpha \neq 0)$ in (\ref{61equivfinal1}),
or   $J''_{1},$ or   $J''_{0} $
in (\ref{61equiv52}).
%%%%%%%%%%%%%%%%%%%%%%%%%%%%%%%%%%%%%%%%%%%%%%%%%%%%%%%%%%%%%%%%%%%%%%%%%%%%%%%%%%%%%%%%%%
%%%%%%%%%%%%%%%%%%%%%%%%%%%%%%%%%%%%%%%%%%%%%%%%%%%%%%%%%%%%%%%%%%%%%%%%%%%%%%%%%%%%%%%%%%
\subsection{}
%\subsubsection{}
\begin{equation*}
X_1 = \frac{\partial}{\partial x^1} - y^1 \, \frac{\partial}{\partial x^3}
- y^2\frac{\partial}{\partial y^3} \quad , \quad
X_2 = \frac{\partial}{\partial y^1} - x^2\frac{\partial}{\partial y^3}.
\end{equation*}
\subsubsection{Holomorphic functions for $J_{\alpha}.$}
Let $G$ denote the group $G_0$ endowed with the left invariant  structure
of complex manifold  defined by
$J_{\alpha}$
(\ref{61equivfinal1}).
Then $ H_{\Cmath}(G) =\{f \in C^{\infty}(G_0) \; ; \; \tilde{X}_j^{-} \, f = 0
\; \forall j\; 1 \leqslant j \leqslant 6\}.$
As
$\tilde{X}_2^{-} = -i( \alpha \tilde{X}_1^{-} +(1-\alpha)\tilde{X}_3^{-}),$
$\tilde{X}_4^{-} = -i\alpha \tilde{X}_1^{-} +(1-i)\tilde{X}_3^{-},$
$\tilde{X}_6^{-} = -i\tilde{X}_5^{-},$ one has
$ H_{\Cmath}(G) =\{f \in C^{\infty}(G_0) \; ; \; \tilde{X}_j^{-} \, f = 0
\; \forall j=1,3,5\}.$ Now
\begin{eqnarray*}
\tilde{X}_1^{-}  &=&
2\, \frac{\partial}{\partial \overline{z^1}}
+i(\alpha-1) \left(\frac{\partial}{\partial x^2}  -\frac{1}{\alpha}
\frac{\partial}{\partial y^2} \right)
-y^1 \frac{\partial}{\partial x^3}
-(y^2+ix^2) \frac{\partial}{\partial y^3} \quad , \\
\tilde{X}_3^{-}  &=& i  \frac{\partial}{\partial y^1}
+(1+i\alpha)   \frac{\partial}{\partial x^2}
-i   \frac{\partial}{\partial y^2}
-i x^2  \frac{\partial}{\partial y^3}\quad , \quad
\tilde{X}_5^{-}  =2\, \frac{\partial}{\partial \overline{z^3}} \quad,
\end{eqnarray*}
where
$z^1 = x^1 +i y^1,\;$
$z^3 = x^3 +i  y^3.$
%\end{eqnarray*}
Then $f \in C^{\infty}(G_0)$ is in $ H_{\Cmath}(G)$ if and only if it is holomorphic
with respect to  $z^3$ and satisfies the  2 equations
\begin{eqnarray}
2\, \frac{\partial f}{\partial \overline{w^2}}
+ \frac{\partial f}{\partial \overline{z^1}}
- \frac{\partial f}{\partial {z^1}}
+  \frac{1}{2} \left(
(1-i\alpha) w^2 + (1+i\alpha) \overline{w^2} \right)
\frac{\partial f}{\partial {z^3}}
&=&0\\
2\, \frac{\partial f}{\partial \overline{z^1}}
+\frac{\alpha-1}{\alpha} \left(
\frac{\partial f}{\partial \overline{w^2}} - \frac{\partial f}{\partial w^2} \right)
+ \left( \frac{\overline{z^1} -z^1}{2i} +
\left(1-\frac{i\alpha}{2} \right)
w^2
+\frac{i\alpha}{2}\,  \overline{w^2} \right)
\,  \frac{\partial f}{\partial z^3}
&=&0
\end{eqnarray}
where
$w^2 = x^2 +\alpha y^2 -iy^2.$
We set $w^1=z^1, w^3=z^3.$
The 3 functions
$$
\varphi^1 =2w^1 +\overline{w^2}+ w^2 \quad , \quad
\varphi^2 = 2w^2 +  \frac{\alpha -1}{\alpha} (\overline{w^1} + w^1)\quad ,
$$
$$
\varphi^3 =w^3 + \frac{1}{32} \left(
4i(\overline{w^1})^2
-8i\overline{w^1} \overline{w^2}
-8i\overline{w^1} {w^1}
-8(2-i)\overline{w^1} {w^2}
-(4+i)(\overline{w^2})^2
+4i\overline{w^2} {w^1}
+4i\overline{w^2} {w^2}
\right)$$
if $\alpha =1$ and if not
\begin{multline*}
\varphi^3 =w^3
-\frac{\alpha i }{2(\alpha-1)}\, \overline{w^1} {w^2}
-\frac{1+i\alpha}{8}\, (\overline{w^2})^2
-\frac{1-\alpha+i\alpha(1+\alpha)}{4(\alpha-1)}\, \overline{w^2} {w^2}
\\
+\frac{i\alpha}{2(\alpha-1)}{w^1} {w^2}
-\frac{3\alpha^2 -2\alpha -1 -i\alpha(\alpha+1)^2}{8(\alpha-1)^2}\, ({w^2})^2
\end{multline*}
are holomorphic.
Let $F : G \rightarrow \Cmath^3$ defined
by
$F=(\varphi^1,\varphi^2,\varphi^3).$
$F$ is a a global chart on $G.$
We determine now how the multiplication of $G$ looks like in that chart.
Let $a,x \in G $ with respective second kind canonical coordinates
$(x^1,y^1,x^2,y^2,x^3,y^3), (\alpha^1, \beta^1, \alpha^2, \beta^2, \alpha^3, \beta^3)$
as in  (\ref{x6general}).
With obvious notations,
%$a =[w^1_a, w^2_a,w^3_a],$
%$x =[w^1_x, w^2_x,w^3_x],$
%$a \, x =[w^1_{a x}, w^2_{a x},w^3_{a x}],$
%$a =[\varphi^1_a, \varphi^2_a,\varphi^3_a],$
%$x =[\varphi^1_x, \varphi^2_x,\varphi^3_x],$
%$a \, x =[\varphi^1_{a x}, \varphi^2_{a x},\varphi^3_{a x}].$
computations yield:
\begin{eqnarray}
w^1_{a x} &=& w^1_a + w^1_x\\
w^2_{a x} &=& w^2_a + w^2_x\\
\label{w3ax}
w^3_{a x} &=& w^3_a + w^3_x
-b^1 x^1  -i(b^2 x^1  +a^2y^1).
\end{eqnarray}
We then get
\begin{equation*}
\varphi^1_{a x} = \varphi^1_a + \varphi^1_x
\quad , \quad
\varphi^2_{a x} = \varphi^2_a + \varphi^2_x
\quad , \quad
\varphi^3_{a x} = \varphi^3_a +\varphi^3_x + \chi(a,x)
\quad ,
\end{equation*}
where for $\alpha \neq 1$
\begin{multline*}
\chi(a,x) =
\frac{1}{8}\, \varphi^1_x
\left( (\overline{\varphi^1_a} +\varphi^1_a) ((1 -i) \alpha - 1)
  - \alpha \, \overline{\varphi^2_a}
+\frac{\alpha (1-\alpha+2i\alpha)}{\alpha-1}
\, \varphi^2_a
\right)
\\
+\frac{1}{8}\,  \varphi^2_x \left(2 \frac{i\alpha}{1-i} \, \overline{\varphi^1_a}
+ \alpha \overline{\varphi^2_a}  +
2\frac{\alpha(1+i\alpha)}{(\alpha-1)(1-i)} \, \varphi^1_a
+
\alpha \frac{(1 - 2 i) \alpha^2 - 1}{ (\alpha- 1)^2}
\,\varphi^2_a
\right)
\end{multline*}
and for $\alpha =1$
\begin{multline*}
\chi(a,x) =
\frac{1}{32}\, \varphi^1_x
\left( -4i \overline{\varphi^1_a} +2i\varphi^1_a
  - 4 \overline{\varphi^2_a}
+3i \varphi^2_a
\right)
+\frac{1}{64}\,  \varphi^2_x \left(8(i-1) \, \overline{\varphi^1_a}
+ 8 \overline{\varphi^2_a}  -2i\, \varphi^1_a + (4-5i)\,\varphi^2_a \right)
.
\end{multline*}
\subsubsection{Holomorphic functions for $J.$}
Now $J$ is defined in (\ref{6_1equivfinal2}).
Then $f \in C^{\infty}(G_0)$ is in $ H_{\Cmath}(G)$ if and only if it is holomorphic
with respect to  $z^2$ and $z^3$ and satisfies the  equation
\begin{equation*}
2\, \frac{\partial f}{\partial \overline{w^1}}
=(z^2+y^1)
\, \frac{\partial f}{\partial \overline{z^3}}
\end{equation*}
where
$w^1 = x^1 -iy^1,$
$z^2 = x^2 +i  y^2, $
$z^3 = x^3 +i  y^3. $
We set $w^2=z^2, w^3=z^3.$
The 3 functions
\begin{equation*}
\varphi^1 =w^1
\quad , \quad
\varphi^2 = w^2
\quad , \quad
\varphi^3 = w^3 + \frac{1}{2} w^2 \overline{w^1}
+ \frac{i}{4} w^1 \overline{w^1}
- \frac{i}{8}  {\overline{w^1}}^2
\end{equation*}
are holomorphic.
Let $F : G \rightarrow \Cmath^3$ defined
by
$F=(\varphi^1,\varphi^2,\varphi^3).$
$F$ is a global chart on $G,$ and
\begin{equation*}
\varphi^1_{a x} = \varphi^1_a + \varphi^1_x
\quad , \quad
\varphi^2_{a x} = \varphi^2_a + \varphi^2_x
\quad , \quad
\varphi^3_{a x} = \varphi^3_a +\varphi^3_x + \chi(a,x)
\quad ,
\end{equation*}
where from (\ref{w3ax})
\begin{equation*}
\chi(a,x) =
\frac{1}{4} \, \varphi^1_x
\left( 2i \overline{\varphi^1_a} -i\varphi^1_a   +2 \overline{\varphi^2_a} \right)
+\frac{1}{2}\,  \varphi^2_x \overline{\varphi^1_a}
.
\end{equation*}
\subsubsection{Holomorphic functions for $J''_{\gamma}.$}
$J''_{\gamma}$ is defined in
(\ref{61equiv52})
for any real $\gamma.$
Here
$\tilde{X}_2^{-} = i\tilde{X}_1^{-},$
$\tilde{X}_4^{-} = -i \tilde{X}_2^{-} -i\gamma\tilde{X}_1^{-},$
$\tilde{X}_6^{-} = -\frac{i}{2}\, \tilde{X}_5^{-},$
hence
$ H_{\Cmath}(G) =\{f \in C^{\infty}(G_0) \; ; \; \tilde{X}_j^{-} \, f = 0
\; \forall j=1,2,5\}.$
One has
\begin{eqnarray*}
\tilde{X}_1^{-}  &=&
 2\,\frac{\partial}{\partial \overline{w^1}}
-i\frac{\partial}{\partial x^2}
-y^1 \frac{\partial}{\partial x^3}
-y^2 \frac{\partial}{\partial y^3} \quad , \\
\tilde{X}_2^{-}
 &=&2\,\frac{\partial}{\partial \overline{w^2}}
+i\gamma \frac{\partial}{\partial x^2}
-x^2  \frac{\partial}{\partial y^3}
\quad , \quad
\tilde{X}_5^{-}  =2\, \frac{\partial}{\partial \overline{w^3}}
\quad ,
\end{eqnarray*}
where  $
w^1 = x^1 -i x^2
\quad , \quad
w^1 = y^1 +i y^2
\quad , \quad
w^3 = x^3 +\frac{i}{2}  y^3.$
Then $f \in C^{\infty}(G_0)$ is in $ H_{\Cmath}(G)$ if and only if it is holomorphic
with respect to  $w^3$ and satisfies the  2 equations
\begin{eqnarray}
2\, \frac{\partial f}{\partial \overline{w^2}}
-\frac{\gamma}{2}\left(  \frac{\partial f}{\partial \overline{w^1}}
- \frac{\partial f}{\partial {w^1}} \right)
+  \frac{1}{4} \left(
 w^1 - \overline{w^1} \right)
\frac{\partial f}{\partial {w^3}}
&=&0\\
2\, \frac{\partial f}{\partial \overline{w^1}}
-\frac{1}{4} \left( 3w^2 +\overline{w^2} \right)
\,  \frac{\partial f}{\partial w^3}
&=&0\, .
\end{eqnarray}
The 3 functions
\begin{eqnarray*}
\varphi^1 &=&\gamma \overline{w^2}-4w^1 \quad , \quad
\varphi^2 = w^2 \quad , \\
\varphi^3 &=&w^3 + \frac{1}{32} \left(4
\overline{w^1} \overline{w^2}
+12\overline{w^1} {w^2}
+\gamma (\overline{w^2})^2
-4\overline{w^2} {w^1}
+12{w^1}{w^2}\right)
\end{eqnarray*}
are holomorphic.
Let $F : G \rightarrow \Cmath^3$ defined
by
$F=(\varphi^1,\varphi^2,\varphi^3).$
$F$ is a global chart on $G.$
Instead of (\ref{w3ax}), we here have
$$w^3_{a x} = w^3_a + w^3_x
-b^1 x^1  -\frac{i}{2}\, (b^2 x^1  +a^2y^1), $$
whence
\begin{equation*}
\varphi^1_{a x} = \varphi^1_a + \varphi^1_x
\quad , \quad
\varphi^2_{a x} = \varphi^2_a + \varphi^2_x
\quad , \quad
\varphi^3_{a x} = \varphi^3_a +\varphi^3_x + \chi(a,x)
\end{equation*}
with
\begin{equation*}
\chi(a,x) =
\frac{1}{16} \, \varphi^1_x \overline{\varphi^2_a}
+\frac{1}{16}\,  \varphi^2_x \left( -\overline{\varphi^1_a}
-2{\varphi^1_a}
+2\gamma \overline{\varphi^2_a}
+\gamma {\varphi^2_a}\right).
\end{equation*}
%%%%%%%%%%%%%%%%%%%%%%%%%%%%%%%%%%%%%%%%%%%%%%%%%%%%%%%%%%%%%%%%%%%%%%%%%%%%%%%%%%%%%%%%%%
\section{Lie Algebra $ {\mathcal{G}}_{6,6}$ (isomorphic to $M1$).}
%%%%%%%%%%%%%%%%%%%%%%%%%%%%%%%%%%%%%%%%%%%%%%%%%%%%%%%%%%%%%%%%%%%%%%%%%%%%%%%%%%%%%%%%%%
Commutation relations for
$ {\mathcal{G}}_{6,6}:$
$[x_1,x_2]=x_4$;
$[x_2,x_3]=x_6$;
$[x_2,x_4]=x_5$.
%\subsection{}
\begin{equation}
\label{66general}
J = \begin{pmatrix}
\boxed{\xi^1_1}&
 - \frac{{\xi^1_1}^2 + 1}{\xi^2_1}&
0&
0&
0&
0\\
\boxed{\xi^2_1}&
 - {\xi^1_1}&
0&
0&
0&
0\\
*&
%(({\xi^3_3} - {\xi^1_1}) {\xi^4_1} - {\xi^4_2} {\xi^2_1})/{\xi^4_3}&
*&
%(({\xi^3_3} + {\xi^1_1}) {\xi^4_2} {\xi^2_1} + ({\xi^1_1}^2 + 1) {\xi^4_1})/({\xi^4_3} {\xi^2_1})&
\boxed{\xi^3_3}&
 - \frac{{\xi^3_3}^2 + 1}{\xi^4_3}&
0&
0\\
\boxed{\xi^4_1}&
\boxed{\xi^4_2}&
\boxed{\xi^4_3}&
 - {\xi^3_3}&
0&
0\\
\boxed{\xi^5_1}&
*&
%( - (({\xi^3_3} - {\xi^1_1}) {\xi^5_3} {\xi^4_1} - ({\xi^3_3} - {\xi^1_1}) {\xi^5_1} {\xi^4_3} -
 %{\xi^5_3} {\xi^4_2} {\xi^2_1} + {\xi^5_4} {\xi^4_3} {\xi^4_1} + {\xi^6_1} {\xi^4_3}^2))/({\xi^4_3} {\xi^2_1})&
\boxed{\xi^5_3}&
\boxed{\xi^5_4}&
 - {\xi^3_3}&
{\xi^4_3}\\
\boxed{\xi^6_1}&
*&
%( - ((({\xi^3_3} + {\xi^1_1}) {\xi^4_1} + {\xi^4_2} {\xi^2_1}) {\xi^5_4} {\xi^4_3} + ({\xi^5_3}
%{\xi^4_1} - {\xi^5_1} {\xi^4_3}) ({\xi^3_3}^2 + 1) + ({\xi^3_3} + {\xi^1_1}) {\xi^6_1} {\xi^4_3}^2))/({\xi^4_3}^2 {\xi^2_1})&
 - {\xi^5_4}&
*&
% (({\xi^3_3}^2 + 1) {\xi^5_3} + 2 {\xi^5_4} {\xi^4_3} {\xi^3_3})/{\xi^4_3}^2&
- \frac{{\xi^3_3}^2 + 1}{\xi^4_3}&
{\xi^3_3}\end{pmatrix}
\end{equation}
where
{\fontsize{6}{10}\selectfont
$ \mathbf{J^3_1}=(({\xi^3_3} - {\xi^1_1}) {\xi^4_1} - {\xi^4_2} {\xi^2_1})/{\xi^4_3};\myquad$
$ \mathbf{J^3_2}=(({\xi^3_3} + {\xi^1_1}) {\xi^4_2} {\xi^2_1} + ({\xi^1_1}^2 + 1) {\xi^4_1})/({\xi^4_3} {\xi^2_1});\myquad$
$ \mathbf{J^5_2}=( - (({\xi^3_3} - {\xi^1_1}) {\xi^5_3} {\xi^4_1} - ({\xi^3_3} - {\xi^1_1}) {\xi^5_1} {\xi^4_3} - {\xi^5_3} {\xi^4_2} {\xi^2_1} + {\xi^5_4} {\xi^4_3} {\xi^4_1} + {\xi^6_1} {\xi^4_3}^2))/({\xi^4_3} {\xi^2_1});\myquad$
$ \mathbf{J^6_2}=( - ((({\xi^3_3} + {\xi^1_1}) {\xi^4_1} + {\xi^4_2} {\xi^2_1}) {\xi^5_4} {\xi^4_3} + ({\xi^5_3} {\xi^4_1} - {\xi^5_1} {\xi^4_3}) ({\xi^3_3}^2 + 1) + ({\xi^3_3} + {\xi^1_1})
{\xi^6_1} {\xi^4_3}^2))/({\xi^4_3}^2 {\xi^2_1});$
\myquad$\mathbf{J^6_4}= (({\xi^3_3}^2 + 1) {\xi^5_3} + 2 {\xi^5_4} {\xi^4_3} {\xi^3_3})/{\xi^4_3}^2;$
}
and the parameters are subject to the condition
\begin{equation}
\label{cond616general1}
\xi^2_1\xi^4_3 \neq 0.
 \end{equation}
%\par The matrix $ J $ is :\\
%$  J^1_1={\xi^1_1};$\\
%$  J^1_2=( - ({\xi^1_1}^2 + 1))/{\xi^2_1};$\\
%$  J^1_3=0;$\\
%$  J^1_4=0;$\\
%$  J^1_5=0;$\\
%$  J^1_6=0;$\\
%$  J^2_1={\xi^2_1};$\\
%$  J^2_2= - {\xi^1_1};$\\
%$  J^2_3=0;$\\
%$  J^2_4=0;$\\
%$  J^2_5=0;$\\
%$  J^2_6=0;$\\
%$  J^3_1=(({\xi^3_3} - {\xi^1_1}) {\xi^4_1} - {\xi^4_2} {\xi^2_1})/{\xi^4_3};$\\
%$  J^3_2=(({\xi^3_3} + {\xi^1_1}) {\xi^4_2} {\xi^2_1} + ({\xi^1_1}^2 + 1) {\xi^4_1})/({\xi%^4_3} {\xi^2_1});$\\
%$  J^3_3={\xi^3_3};$\\
%$  J^3_4=( - ({\xi^3_3}^2 + 1))/{\xi^4_3};$\\
%$  J^3_5=0;$\\
%$  J^3_6=0;$\\
%$  J^4_1={\xi^4_1};$\\
%$  J^4_2={\xi^4_2};$\\
%$  J^4_3={\xi^4_3};$\\
%$  J^4_4= - {\xi^3_3};$\\
%$  J^4_5=0;$\\
%$  J^4_6=0;$\\
%$  J^5_1={\xi^5_1};$\\
%$  J^5_2=( - (({\xi^3_3} - {\xi^1_1}) {\xi^5_3} {\xi^4_1} - ({\xi^3_3} - {\xi^1_1}) {\xi^5%_1} {\xi^4_3} - {\xi^5_3} {\xi^4_2} {\xi^2_1} + {\xi^5_4} {\xi^4_3} {\xi^4_1} + {\xi^6_1} %{\xi^4_3}
%^2))/({\xi^4_3} {\xi^2_1});$\\
%$  J^5_3={\xi^5_3};$\\
%$  J^5_4={\xi^5_4};$\\
%$  J^5_5= - {\xi^3_3};$\\
%$  J^5_6={\xi^4_3};$\\
%$  J^6_1={\xi^6_1};$\\
%$  J^6_2=( - ((({\xi^3_3} + {\xi^1_1}) {\xi^4_1} + {\xi^4_2} {\xi^2_1}) {\xi^5_4} {\xi^4_3%}
%+ ({\xi^5_3} {\xi^4_1} - {\xi^5_1} {\xi^4_3}) ({\xi^3_3}^2 + 1) + ({\xi^3_3} + {\xi^1_1}) {\xi^6_1} {\xi^4_3}^2))/({\xi^4_3}^2 {\xi^2_1});$\\
%$  J^6_3= - {\xi^5_4};$\\
%$  J^6_4=(({\xi^3_3}^2 + 1) {\xi^5_3} + 2 {\xi^5_4} {\xi^4_3} {\xi^3_3})/{\xi^4_3}^2
%;$\\
%$  J^6_5=( - ({\xi^3_3}^2 + 1))/{\xi^4_3};$\\
%$  J^6_6={\xi^3_3};$\\
\par Now the automorphism group of
$ {\mathcal{G}}_{6,6}$ is comprised of the matrices
$$
{\fontsize{8}{10}\selectfont
\Phi = \begin{pmatrix}
{b^1_1}&
{b^1_2}&
0&
0&
0&
0\\
0&
{b^2_2}&
{b^2_3}&
0&
0&
0\\
{b^3_1}&
{b^3_2}&
{b^3_3}&
0&
0&
0\\
{b^4_1}&
{b^4_2}&
{b^4_3}&
{b^2_2} {b^1_1}&
0&
0\\
{b^5_1}&
{b^5_2}&
{b^5_3}&
 - {b^4_1} {b^2_2}&
{b^2_2}^2 {b^1_1}&
{b^4_3} {b^2_2}\\
{b^6_1}&
{b^6_2}&
{b^6_3}&
 - {b^3_1} {b^2_2}&
0&
{b^3_3} {b^2_2}\end{pmatrix}
 } $$
where $ b^1_1 {b^2_2} {b^3_3} \neq 0.$
Taking
$$
{\fontsize{8}{10}\selectfont
\Phi = \begin{pmatrix}
1&
0&
0&
0&
0&
0\\
0&
1&
0&
0&
0&
0\\
0&
0&
1&
0&
0&
0\\
0&
0&
0&
1&
0&
0\\
0&
0&
 - ({\xi^5_4} {\xi^4_3})/({\xi^3_3}^2 + 1)&
0&
1&
0\\
{b^6_1}&
%(({b^5_2} {\xi^3_3}^2 + {b^5_2} + {\xi^5_4} {\xi^4_2}) {\xi^2_1} - ({\xi^3_3} - {\xi^1_1}) xi
%(5_4} {\xi^4_1} + ({b^5_1} {\xi^3_3} + {b^5_1} {\xi^1_1} - {\xi^5_1}) ({\xi^3_3}^2 + 1))/((
%{\xi^3_3}^2 + 1) {\xi^4_3})&
{b^6_2}&
%( - ((({\xi^3_3} + {\xi^1_1}) {\xi^4_2} {\xi^2_1} - 2 ({\xi^3_3} - {\xi^1_1}) {\xi^4_1} {\xi^3,3
%)) {\xi^5_4} {\xi^4_3} - (({\xi^3_3} - {\xi^1_1}) {\xi^5_3} {\xi^4_1} - ({\xi^3_3} - {\xi^1_1})
%{\xi^5_1} {\xi^4_3} - {\xi^5_3} {\xi^4_2} {\xi^2_1}) ({\xi^3_3}^2 + 1) - ({\xi^3_3}^2 + 1)
%{\xi^6_1} {\xi^4_3}^2 + ({\xi^3_3}^2 + 1) ({\xi^1_1}^2 + 1) {b^5_1} {\xi^4_3} - ({\xi^3_3}
%^2 + 1) ({\xi^3_3} - {\xi^1_1}) {b^5_2} {\xi^4_3} {\xi^2_1} + ({\xi^3_3} - {\xi^1_1})^2 {\xi^
%5_4} {\xi^4_3} {\xi^4_1}))/(({\xi^3_3}^2 + 1) {\xi^4_3}^2 {\xi^2_1})&
%{b^6_3}&
- (({\xi^3_3}^2 + 1) {\xi^5_3} + 2 {\xi^5_4} {\xi^4_3} {\xi^3_3})/
(({\xi^3_3}^2 + 1)\xi^4_3)&
0&
0&
1\end{pmatrix}
}
$$
with suitable values for  $b^6_1,b^6_2,$
equivalence by $\Phi$
leads to the case where
$ \xi^5_1=\xi^5_3=\xi^5_4={\xi^6_1} = 0. $
Then equivalence by
$${\fontsize{8}{10}\selectfont
 \Phi = \begin{pmatrix}
1&
0&
0&
0&
0&
0\\
0&
1&
0&
0&
0&
0\\
0&
(
 - {\xi^4_2} {\xi^2_1} + {\xi^4_1} {\xi^3_3} - {\xi^4_1} {\xi^1_1})/({\xi^4_3} {\xi^2_1})&
1&
0&
0&
0\\
0&
{\xi^4_1}/{\xi^2_1}&
0&
1&
0&
0\\
0&
0&
0&
0&
1&
0\\
0&
0&
0&
0&
0&
1\end{pmatrix}
}$$
leads to the case where
moreover $ {\xi^4_1}={\xi^4_2}=0 .$
Finally equivalence by
${\fontsize{8}{10}\selectfont
 \Phi = \text{diag}
\left( \begin{pmatrix}
{\xi^2_1}&  - {\xi^1_1}\\0&1
 \end{pmatrix},
\begin{pmatrix}
({\xi^4_3} {\xi^2_1})/({\xi^3_3}^2 + 1)&
0\\
 - ({\xi^4_3} {\xi^3_3} {\xi^2_1})/({\xi^3_3}^2 + 1)&
{\xi^2_1}
 \end{pmatrix},
\begin{pmatrix}
{\xi^2_1}&
 - ({\xi^4_3} {\xi^3_3} {\xi^2_1})/({\xi^3_3}^2 + 1)\\
0&
({\xi^4_3} {\xi^2_1})/({\xi^3_3}^2 + 1)\end{pmatrix} \right)
}
$
leads to
\begin{equation}
\label{66equivfinal}
J = \text{diag} \left(
\begin{pmatrix}0&-1\\1&0 \end{pmatrix},
\begin{pmatrix}0&-1\\1&0 \end{pmatrix},
\begin{pmatrix}0&1\\-1&0 \end{pmatrix}
\right)
\, .
\end{equation}
Commutation relations of $ \mathfrak{m} : $
$  [\tilde{x}_1,\tilde{x}_3]= - \tilde{x}_5;$
 $  [\tilde{x}_1,\tilde{x}_4]=\tilde{x}_6;$
$  [\tilde{x}_2,\tilde{x}_3]=\tilde{x}_6;$
 $  [\tilde{x}_2,\tilde{x}_4]=\tilde{x}_5.$
%\subsection{Conclusions.}
\par
From (\ref{66general}),  $\frak{X}_{6,6}$  is a 10-dimensional submanifold of $\Rmath^{36}.$
There is only one $\text{Aut }  {\mathcal{G}}_{6,6}$ orbit, and any element of
$\frak{X}_{6,6}$     is equivalent to $J$ in (\ref{66equivfinal}).
%%%%%%%%%%%%%%%%%%%%%%%%%%%%%%%%%%%%%%%%%%%%%%%%%%%%%%%%%%%%%%%%%%%%%%%%%%%%%%%%%%%%%%%%%%
\subsection{}
%\subsubsection{}
\begin{equation*}
X_1 = \frac{\partial}{\partial x^1} - y^1 \, \frac{\partial}{\partial y^2}
+ \frac{(y^1)^2}{2} \, \frac{\partial}{\partial x^3}
\quad , \quad
X_2 = \frac{\partial}{\partial y^1} - y^2\frac{\partial}{\partial x^3}
- x^2\frac{\partial}{\partial y^3}.
\end{equation*}
%\subsubsection{Holomorphic functions for $J.$}
Let $G$ denote the group $G_0$ endowed with the left invariant  structure
of complex manifold  defined by
$J$ in (\ref{66equivfinal}).
Then $ H_{\Cmath}(G) =\{f \in C^{\infty}(G_0) \; ; \; \tilde{X}_j^{-} \, f = 0
\; \forall j=1,3,5\}.$
One has
\begin{eqnarray*}
\tilde{X}_1^{-}  &=&
2\, \frac{\partial}{\partial \overline{w^1}}
-y^1 \frac{\partial}{\partial y^2}
+\left( \frac{(y^1)^2}{2} -iy^2 \right)
\frac{\partial}{\partial x^3}
-ix^2\frac{\partial}{\partial y^3} \, ,
\\
\tilde{X}_3^{-}  &=& 2  \frac{\partial}{\partial \overline{w^2}}
\quad , \quad
\tilde{X}_5^{-}  = 2\, \frac{\partial}{\partial \overline{w^3}}
\quad ,
\end{eqnarray*}
where
$w^1 = x^1 +i y^1
\quad , \quad
w^2 = x^2 +i y^2
\quad , \quad
w^3 = x^3 -i  y^3.$
Then $f \in C^{\infty}(G_0)$ is in $ H_{\Cmath}(G)$ if and only if it is holomorphic
with respect to  $w^2$ and $w^3$ and satisfies the equation
\begin{eqnarray}
2\, \frac{\partial f}{\partial \overline{w^1}}
-\frac{w^1-\overline{w^1}}{2}\, \frac{\partial f}{\partial {w^2}}
-\left(
\frac{(w^1-\overline{w^1})^2}{8}
+w^2 \right)
\frac{\partial f}{\partial {w^3}}
&=&0.
\end{eqnarray}
The 3 functions
\begin{eqnarray*}
\varphi^1 &=&w^1
\quad , \quad
\varphi^2 = w^2 +  \frac{1}{4} \left(w^1\overline{w^1} -\frac{(\overline{w^1})^2}{2}\right)\quad , \\
\varphi^3 &=&w^3 + \frac{1}{48} \left(
-(\overline{w^1})^3
+3\overline{w^1} (w^1)^2
+24\overline{w^1} w^2 \right)
\end{eqnarray*}
are holomorphic.
Let $F : G \rightarrow \Cmath^3$ defined
by
$F=(\varphi^1,\varphi^2,\varphi^3).$
%$F$ is a biholomorphic bijection, hence a global chart on $G.$
$F$ is a global chart on $G.$
We determine now how the multiplication of $G$ looks like in that chart.
Let $a,x \in G $ with respective second kind canonical coordinates
$(x^1,y^1,x^2,y^2,x^3,y^3), (\alpha^1, \beta^1, \alpha^2, \beta^2, \alpha^3, \beta^3)$
as in  (\ref{x6general}).
With obvious notations,  computations yield:
$$
w^1_{a x} = w^1_a + w^1_x\quad , \quad
w^2_{a x} = w^2_a + w^2_x-ib^1x^1 \quad, \quad
w^3_{a x} = w^3_a + w^3_x
+\frac{(b^1)^2}{2} x^1  -(b^2 -b^1x^1) y^1  +ia^2y^1.
$$
We then get
$$
\varphi^1_{a x} = \varphi^1_a + \varphi^1_x\quad , \quad
\varphi^2_{a x} = \varphi^2_a + \varphi^2_x +\frac{1}{4}\, \varphi^1_x (2 \overline{\varphi^1_a} -\varphi^1_a)
\quad , \quad
\varphi^3_{a x} = \varphi^3_a +\varphi^3_x + \chi(a,x)
\quad ,
$$
where
\begin{equation*}
\chi(a,x) =
\frac{1}{16}\, (\varphi^1_x)^2
\left( 3\overline{\varphi^1_a} -2\varphi^1_a
\right)
+\frac{1}{16}\, \varphi^1_x
\left( 2(\overline{\varphi^1_a})^2 -(\varphi^1_a)^2
+8\overline{\varphi^2_a} \right)
+\frac{1}{2}\,\varphi^2_x \overline{\varphi^1_a} .
\end{equation*}
%%%%%%%%%%%%%%%%%%%%%%%%%%%%%%%%%%%%%%%%%%%%%%%%%%%%%%%%%%%%%%%%%%%%%%%%%%%%%%%%%%%%%%%%%%
%%%%%%%%%%%%%%%%%%%%%%%%%%%%%%%%%%%%%%%%%%%%%%%%%%%%%%%%%%%%%%%%%%%%%%%%%%%%%%%%%%%%%%%%%%
\section{Lie Algebra $ {\mathcal{G}}_{6,5}$ (isomorphic to $M8$).}
%%%%%%%%%%%%%%%%%%%%%%%%%%%%%%%%%%%%%%%%%%%%%%%%%%%%%%%%%%%%%%%%%%%%%%%%%%%%%%%%%%%%%%%%%%
Commutation relations for
$ {\mathcal{G}}_{6,5}:$
$[x_1,x_2]=x_4$;
$[x_1,x_4]=x_5$;
$[x_2,x_3]=x_6$;
$[x_2,x_4]=x_6$.
%\subsection{}
%{\fontsize{8}{10}\selectfont
\begin{equation}
\label{6-5general}
J= \begin{pmatrix}
a&-\frac{a^2+1}{\xi^2_1}&0&0&0&0\\
\boxed{\xi^2_1}&-a&0&0&0&0\\
\boxed{\xi^3_1}&
*&
b&
-\frac{b^2+1}{\xi^4_3}&
0&0\\
\boxed{\xi^4_1}&*&\boxed{\xi^4_3}&
%$  J^4_4=(({\xi^5_5} + {\xi^4_3}) {\xi^2_1} - {\xi^6_5} {\xi^4_3})/{\xi^2_1};$\\
-b&
0&0\\
\boxed{\xi^5_1}&*&
%$  J^5_3=(((2 {\xi^5_5} + {\xi^4_3}) {\xi^6_3} - {\xi^6_4} {\xi^4_3}) {\xi^2_1} - {\xi^6_5}
%{\xi^6_3} {\xi^4_3})/({\xi^6_5} {\xi^2_1});$\\
\frac{((2 {\xi^5_5} + {\xi^4_3}) {\xi^6_3} - {\xi^6_4} {\xi^4_3}) {\xi^2_1} - {\xi^6_5}
{\xi^6_3} {\xi^4_3}}{{\xi^6_5} {\xi^2_1}}&
*&\boxed{\xi^5_5}&-\frac{{\xi^5_5}^2+1}{\xi^6_5}\\
\boxed{\xi^6_1}&*&\boxed{\xi^6_3}&\boxed{\xi^6_4}&\boxed{\xi^6_5}&-\xi^5_5
\end{pmatrix}
\end{equation}
%}
where
{\fontsize{6}{10}\selectfont
$ \mathbf{a}= \mathbf{J^1_1}=( - (({\xi^5_5}^2 + 1) {\xi^2_1} - {\xi^6_5} {\xi^5_5} {\xi^4_3}))/({\xi^6_5} {\xi^4_3});$
\myquad
$  \mathbf{J^3_2}=(((({\xi^5_5} + 2 {\xi^4_3}) {\xi^5_5} {\xi^4_1} + {\xi^4_3}^2 {\xi^4_1} + {\xi^4_3}^2 {\xi^3_1} + {\xi^4_1}) {\xi^6_5} + ({\xi^5_5}^2 + 1) {\xi^3_1} {\xi^2_1}) {\xi^2_1}^2
+ (({\xi^6_5} {\xi^4_3} - 2 {\xi^5_5} {\xi^2_1}) {\xi^4_1} - (2 {\xi^4_1} + {\xi^3_1}) {\xi^4_3}
 {\xi^2_1}) {\xi^6_5}^2 {\xi^4_3})/({\xi^6_5} {\xi^4_3} {\xi^2_1}^3);$
\myquad
$\mathbf{b}=  \mathbf{J^3_3}=( - (({\xi^5_5} + {\xi^4_3}) {\xi^2_1} - {\xi^6_5} {\xi^4_3}))/{\xi^2_1};$
\myquad
$  \mathbf{J^4_2}=((({\xi^5_5}^2 + 1) {\xi^2_1}^2 + {\xi^6_5}^2 {\xi^4_3}^2) {\xi^4_1} - (({\xi^4_1} + {\xi^3_1}) {\xi^4_3} + 2 {\xi^5_5} {\xi^4_1}) {\xi^6_5} {\xi^4_3} {\xi^2_1})/({\xi^6_5}
{\xi^4_3} {\xi^2_1}^2);$
\myquad
$  \mathbf{J^5_2}=( - (((({\xi^5_5} {\xi^4_1} + 2 {\xi^4_3} {\xi^4_1} + 2 {\xi^4_3} {\xi^3_1}) {\xi^5_5} + {\xi^4_3}^2 {\xi^4_1} + {\xi^4_3}^2 {\xi^3_1} + {\xi^4_1}) {\xi^6_3} - ({\xi^6_1} {\xi^4_3} + {\xi^5_1} {\xi^2_1}) ({\xi^5_5}^2 + 1) - ({\xi^4_1} + {\xi^3_1}) {\xi^6_4} {\xi^4_3}^2
) {\xi^2_1}^2 - ((2 ({\xi^6_3} {\xi^4_1} - {\xi^5_1} {\xi^2_1}) {\xi^5_5} + (2 {\xi^4_1} + {\xi^3_1}) {\xi^6_3} {\xi^4_3}) {\xi^2_1} - ({\xi^6_5} {\xi^6_3} + {\xi^6_4} {\xi^2_1}) {\xi^4_3} {\xi^4_1}) {\xi^6_5} {\xi^4_3}))/({\xi^6_5} {\xi^4_3} {\xi^2_1}^3);$
\myquad
$  \mathbf{J^5_4}=( - ((({\xi^6_4} {\xi^4_3}^2 - {\xi^6_3}) {\xi^2_1} - {\xi^6_5} {\xi^6_4} {\xi^4_3}^2) {\xi^2_1} - (({\xi^5_5} + {\xi^4_3}) {\xi^2_1} - {\xi^6_5} {\xi^4_3})^2 {\xi^6_3}))/({\xi^6_5} {\xi^4_3} {\xi^2_1}^2);$
\myquad
$  \mathbf{J^6_2}=( - (({\xi^6_5} {\xi^5_1} + {\xi^6_4} {\xi^4_1} + {\xi^6_3} {\xi^3_1}) {\xi^6_5} {\xi^4_3} - ({\xi^5_5}^2 + 1) {\xi^6_1} {\xi^2_1}))/({\xi^6_5} {\xi^4_3} {\xi^2_1});$
}
and the parameters are subject to the condition
\begin{equation}
\label{cond65general}
\xi^2_1\xi^4_3\xi^6_5 \neq 0.
 \end{equation}
\par
The automorphisms of $ {\mathcal{G}}_{6,5}$  fall into 2 kinds
The first kind  is comprised of the matrices
\begin{equation}
{\fontsize{8}{10}\selectfont
\label{autom6_5firstkind}
\Phi = \begin{pmatrix}
{b^1_1}&
0&
0&
0&
0&
0\\
0&
{b^2_2}&
0&
0&
0&
0\\
{b^3_1}&
{b^3_2}&
{b^2_2} {b^1_1}&
0&
0&
0\\
{b^4_1}&
{b^4_2}&
0&
{b^2_2} {b^1_1}&
0&
0\\
{b^5_1}&
{b^5_2}&
{b^5_3}&
{b^4_2} {b^1_1}&
{b^2_2} {b^1_1}^2&
0\\
{b^6_1}&
{b^6_2}&
{b^6_3}&
 - ({b^4_1} + {b^3_1}) {b^2_2}&
0&
{b^2_2}^2 {b^1_1}\end{pmatrix}
 }
 \end{equation}
 where $b^1_1b^2_2 \neq 0.$
The second kind  is comprised of the matrices
\begin{equation}
{\fontsize{8}{10}\selectfont
\label{autom6_5secondkind}
\Phi = \begin{pmatrix}
0&
{b^1_2}&
0&
0&
0&
0\\
{b^2_1}&
0&
0&
0&
0&
0\\
{b^3_1}&
{b^3_2}&
{b^2_1} {b^1_2}&
0&
0&
0\\
{b^4_1}&
{b^4_2}&
 - {b^2_1} {b^1_2}&
 - {b^2_1} {b^1_2}&
0&
0\\
{b^5_1}&
{b^5_2}&
{b^5_3}&
 - {b^4_1} {b^1_2}&
0&
 - {b^2_1} {b^1_2}^2\\
{b^6_1}&
{b^6_2}&
{b^6_3}&
{b^2_1} ({b^4_2} + {b^3_2})&
 - {b^2_1}^2 {b^1_2}&
0\end{pmatrix}
}
 \end{equation}
 where $b^1_2b^2_1 \neq 0.$
Taking suitable values for the $b^i_j$'s,
equivalence by $\Phi$ in
(\ref{autom6_5firstkind})
leads to the case where
$ {\xi^5_1}=\xi^5_5={\xi^6_1}=\xi^6_3=\xi^6_4=\xi^6_5=
{\xi^3_1}={\xi^4_1}=0 $ and moreover $ {\xi^2_1}=1: $
{\fontsize{8}{10} \selectfont
\begin{equation}
\label{6-5}
 J({\xi^4_3},{\xi^5_5},{\xi^6_5}) = \begin{pmatrix}
%( - ({\xi^5_5}^2 + 1 - {\xi^6_5} {\xi^5_5} {\xi^4_3}))/({\xi^6_5} {\xi^4_3})&
-\frac{{\xi^5_5}^2 +1 -\xi^6_5\xi^5_5\xi^4_3}{\xi^6_5\xi^4_3} &
%( - ({\xi^6_5}^2 {\xi^4_3}^2 - 2 {\xi^6_5} {\xi^5_5} {\xi^4_3} + {\xi^5_5}^2 + 1) ({\xi^5_5}^2 + 1))/({\xi^6_5}^2 {\xi^4_3}^2)&
%b&
-\frac{({\xi^5_5}^2 +1 -\xi^6_5\xi^5_5\xi^4_3)^2}{(\xi^6_5\xi^4_3)^2}-1 &
0&
0&
0&
0\\
1&
\frac{{\xi^5_5}^2 +1 -\xi^6_5\xi^5_5\xi^4_3}{\xi^6_5\xi^4_3} &
%({\xi^5_5}^2 + 1 - {\xi^6_5} {\xi^5_5} {\xi^4_3})/({\xi^6_5} {\xi^4_3})&
%c&
0&
0&
0&
0\\
0&
0&
 - ({\xi^5_5} + {\xi^4_3} - {\xi^6_5} {\xi^4_3})&
 - \frac{({\xi^5_5} + {\xi^4_3} - {\xi^6_5} {\xi^4_3})^2 + 1}{\xi^4_3}&
0&
0\\
0&
0&
{\xi^4_3}&
{\xi^5_5} + {\xi^4_3} - {\xi^6_5} {\xi^4_3}&
0&
0\\
0&
0&
0&
0&
{\xi^5_5}&
 - \frac{{\xi^5_5}^2 + 1}{\xi^6_5}\\
0&
0&
0&
0&
{\xi^6_5}&
 - {\xi^5_5}\end{pmatrix}
 \end{equation}
}
where $\xi^4_3\xi^6_5 \neq 0.$
\\
\par Commutation relations of $ \mathfrak{m} : $
 $  [\tilde{x}_1,\tilde{x}_3]=
\frac{1}{\xi^6_5} \,\left(
(-{\xi^5_5} {\xi^4_3} +{\xi^5_5}^2 + 1) \tilde{x}_5
+(-{\xi^6_5} {\xi^4_3} + {\xi^5_5}) \tilde{x}_6  \right);
$
%( - (({\xi^6_5} {\xi^4_3} \tilde{x}_6 + {\xi^5_5} {\xi^4_3} \tilde{x}_5 - {\xi^5_5} \tilde{x}_6) {\xi^6_5} - ({\xi^5_5}^2 + 1) \tilde{x}_5))/{\xi^6_5}
%;\\$
\\ $  [\tilde{x}_1,\tilde{x}_4]=
\frac{1}{{\xi^6_5} {\xi^4_3}}\,
(({\xi^5_5}^2 + 1) ({\xi^5_5} + {\xi^4_3})
+ {\xi^6_5} {\xi^5_5} {\xi^4_3}({\xi^6_5} {\xi^4_3} -2{\xi^5_5}- {\xi^4_3}))\,\tilde{x}_5
+\frac{1}{{\xi^4_3}}\,
(({\xi^6_5} {\xi^4_3} - {\xi^5_5})^2 +
{\xi^4_3}(\xi^5_5 - \xi^6_5 {\xi^4_3})+1)\, \tilde{x}_6;
 $
%\\ $  [\tilde{x}_1,\tilde{x}_4]=(\tilde{x}_5 (({\xi^5_5}^2 + 1) ({\xi^5_5} + {\xi^4_3}
%) + {\xi^6_5} {\xi^4_3}))/({\xi^6_5} {\xi^4_3}) + ({\xi^6_5}^2 {\xi^4_3}^2 \tilde{x}_6 + {\xi^6_5} {\xi^5_5} {\xi^4_3}^2 \tilde{x}_5 - 2 {\xi^6_5} {\xi^5_5} {\xi^4_3} \tilde{x}_6 - {\xi^6_5}
%{\xi^4_3}^2 \tilde{x}_6 - 2 {\xi^5_5}^2 {\xi^4_3} \tilde{x}_5 + {\xi^5_5}^2 \tilde{x}_6 - {\xi^5_5} {\xi^4_3}^2 \tilde{x}_5 + {\xi^5_5} {\xi^4_3} \tilde{x}_6 - ({\xi^4_3} \tilde{x}_5 -
%\tilde{x}_6))/{\xi^4_3};\\$
\\ $  [\tilde{x}_2,\tilde{x}_3]=
\frac{({\xi^6_5} {\xi^4_3} -  {\xi^5_5} )^2 + 1}{{\xi^6_5}^2 {\xi^4_3}}
\, \left(
(1+{\xi^5_5}^2 )\tilde{x}_5 + {\xi^6_5} {\xi^5_5} \tilde{x}_6
\right);
$
%\\ $  [\tilde{x}_2,\tilde{x}_3]=(({\xi^6_5}^2 {\xi^4_3}^2 - 2 {\xi^6_5} {\xi^5_5} {\xi^4_3} + {\xi^5_5}^2 + 1) ({\xi^6_5} {\xi^5_5} \tilde{x}_6 + {\xi^5_5}^2 \tilde{x}_5 +
%\tilde{x}_5))/({\xi^6_5}^2 {\xi^4_3});\\$
\\ $  [\tilde{x}_2,\tilde{x}_4]=
\frac{({\xi^6_5} {\xi^4_3} -  {\xi^5_5} )^2 + 1}{{\xi^6_5}^2 {\xi^4_3}^2}
\, \left(
-( {\xi^5_5}^2+1)(\xi^4_3\xi^6_5-\xi^5_5 -\xi^4_3)  \tilde{x}_5
+\xi^6_5 ( - {\xi^6_5}{\xi^5_5} {\xi^4_3} + {\xi^5_5}^2+1   +  {\xi^5_5} \xi^4_3)
\tilde{x}_6 \right).$
%\\ $  [\tilde{x}_2,\tilde{x}_4]=( - ({\xi^6_5}^2 {\xi^5_5} {\xi^4_3} \tilde{x}_6 + {\xi^6_5} {\xi^5_5}^2 {\xi^4_3} \tilde{x}_5 - {\xi^6_5} {\xi^5_5}^2 \tilde{x}_6 - {\xi^6_5} {\xi^5_5}
% {\xi^4_3} \tilde{x}_6 + {\xi^6_5} {\xi^4_3} \tilde{x}_5 - {\xi^6_5} \tilde{x}_6 - {\xi^5_5}^3
%\tilde{x}_5 - {\xi^5_5}^2 {\xi^4_3} \tilde{x}_5 - {\xi^5_5} \tilde{x}_5 - {\xi^4_3} \tilde{x}_5) (
%{\xi^6_5}^2 {\xi^4_3}^2 - 2 {\xi^6_5} {\xi^5_5} {\xi^4_3} + {\xi^5_5}^2 + 1))/({\xi^6_5}^2
% {\xi^4_3}^2)\\$

$ J({\xi^4_3},{\xi^5_5},{\xi^6_5}) ,  J({\eta^4_3},{\eta^5_5},{\eta^6_5}) $
as in (\ref{6-5}) are equivalent under some first kind automorphism if and only if
$\eta^4_3=\xi^4_3,
\eta^5_5=\xi^5_5,
\eta^6_5=\xi^6_5.$
They are equivalent
under some second kind automorphism
if and only if
$ \eta^4_3= - (({\xi^6_5} {\xi^4_3} -  {\xi^5_5})^2 +1)/{\xi^4_3},$
$\eta^5_5= - {\xi^5_5},$
$\eta^6_5= {\xi^6_5} {\xi^4_3}^2/
(({\xi^6_5} {\xi^4_3} -  {\xi^5_5})^2 +1).$
%\subsection{Conclusions.}
\par
From (\ref{6-5general}), $\mathfrak{X}_{6,5}$
is a submanifold of dimension 10 in $\Rmath^{36}.$
Each CS is equivalent to  some
$ J({\xi^4_3},{\xi^5_5},{\xi^6_5})$ in (\ref{6-5}).
%%%%%%%%%%%%%%%%%%%%%%%%%%%%%%%%%%%%%%%%%%%%%%%%%%%%%%%%%%%%%%%%%%%%%%%%%%%%%%%%%%%%%%%%%%
\subsection{}
%\subsubsection{}
\begin{equation*}
X_1 = \frac{\partial}{\partial x^1} - y^1 \, \frac{\partial}{\partial y^2}
 -y^2\, \frac{\partial}{\partial x^3}
+ \frac{(y^1)^2}{2} \, \frac{\partial}{\partial y^3}
\quad , \quad
X_2 = \frac{\partial}{\partial y^1} - (x^2+y^2) \, \frac{\partial}{\partial y^3}.
\end{equation*}
%\subsubsection{Holomorphic functions for $ J({\xi^4_3},{\xi^5_5},{\xi^6_5}).$}
Let $G$ denote the group $G_0$ endowed with the left invariant  structure
of complex manifold  defined by
$J(\xi^4_3,\xi^5_5,\xi^6_5)$  in (\ref{6-5}) where $\xi^4_3\xi^6_5 \neq 0.$
Then $ H_{\Cmath}(G) =\{f \in C^{\infty}(G_0) \; ; \; \tilde{X}_j^{-} \, f = 0
\; \forall j=1,3,5\}.$
One has
\begin{eqnarray*}
\tilde{X}_1^{-}  &=&
2\, \frac{\partial}{\partial \overline{w^1}}
-y^1 (1+iA)\, \frac{\partial}{\partial y^2}
-y^2 (1+iA)\, \frac{\partial}{\partial x^3}
+\left( \frac{(y^1)^2}{2} (1+iA)  -i(x^2+y^2) \right) \,\frac{\partial}{\partial y^3}
\quad , \\
\tilde{X}_3^{-}  &=& 2  \frac{\partial}{\partial \overline{w^2}}
\quad , \quad
\tilde{X}_5^{-}  =2\, \frac{\partial}{\partial \overline{w^3}}
\quad ,
\end{eqnarray*}
where
\begin{eqnarray*}
w^1 &=& x^1-Ay^1 +i y^1\\
w^2 &=& x^2 +\frac{\xi^5_5 +\xi^4_3 -\xi^6_5\xi^4_3}{\xi^4_3} \, y^2  + \frac{i}{\xi^4_3}\, y^2\\
w^3 &=& x^3 - \frac{\xi^5_5}{\xi^6_5} \, y^3 +\frac{i}{\xi^6_5} \,   y^3 \\
 A&=& - \frac{{\xi^5_5}^2 + 1 - {\xi^6_5} {\xi^5_5} {\xi^4_3}}{{\xi^6_5} {\xi^4_3}}.
\end{eqnarray*}
Then $f \in C^{\infty}(G_0)$ is in $ H_{\Cmath}(G)$ if and only if it is holomorphic
with respect to  $w^2$ and $w^3$ and satisfies the equation
\begin{multline*}
2\, \frac{\partial f}{\partial \overline{w^1}}
-\frac{w^1-\overline{w^1}}{2i}\,
\frac{(1+iA)(\xi^5_5 +\xi^4_3 -\xi^6_5\xi^4_3 +i)}{\xi^4_3}\, \frac{\partial f}{\partial {w^2}}          \\
-\left[ (1+iA)\xi^4_3 \, \frac{w^2-\overline{w^2}}{2i}
+\left(
\frac{(w^1-\overline{w^1})^2}{8}\,(1+iA)
+i \,\left( \frac{w^2+\overline{w^2}}{2} - (\xi^5_5  -\xi^6_5 \xi^4_3)
\, \frac{w^2-\overline{w^2}}{2i}
\right) \right) \frac{i-\xi^5_5}{\xi^6_5}
\right] \,
\frac{\partial f}{\partial {w^3}}
=0.
\end{multline*}
The 3 functions
\begin{equation*}
\varphi^1 =w^1 \quad  , \quad
\varphi^2 = w^2 +  \frac{1+iA}{4i\xi^4_3} (\xi^5_5 +\xi^4_3 - \xi^6_5\xi^4_3 +i)
\left(w^1\overline{w^1} -\frac{(\overline{w^1})^2}{2}\right)
\quad ,
\end{equation*}
\begin{multline*}
\varphi^3 =w^3 +
\frac{1}{48{\xi^6_5}^2 {\xi^4_3}^2} \,
\overline{w^1}^3 ( - 2 i {\xi^6_5}^2 {\xi^5_5}^2 {\xi^4_3}^2 - 4 {\xi^6_5}^2 {\xi^5_5} {\xi^4_3}^2 + 2 i {\xi^6_5}^2 {\xi^4_3}^2 + 4 i {\xi^6_5} {\xi^5_5}^3
{\xi^4_3} + i {\xi^6_5} {\xi^5_5}^2 {\xi^4_3}^2 + 4 {\xi^6_5} {\xi^5_5}^2 {\xi^4_3}
 + 2 \xi^6_5 {\xi^5_5} {\xi^4_3}^2 \\
 + 4 i {\xi^6_5} {\xi^5_5} {\xi^4_3} - i {\xi^6_5} {\xi^4_3}^2 + 4
 {\xi^6_5} {\xi^4_3} - 2 i {\xi^5_5}^4 - i {\xi^5_5}^3 {\xi^4_3} - {\xi^5_5}^2 {\xi^4_3} -
4 i {\xi^5_5}^2 - i {\xi^5_5} {\xi^4_3} - {\xi^4_3} - 2 i)
\\
+
\frac{1}{16{\xi^6_5}^2 {\xi^4_3}^2} \,
\overline{w^1}^2 w^1 (i {\xi^6_5}^2 {\xi^5_5}^2 {\xi^4_3}^2 + 2 {\xi^6_5}^2 {\xi^5_5} {\xi^4_3}
^2 - i {\xi^6_5}^2 {\xi^4_3}^2 - 2 i {\xi^6_5} {\xi^5_5}^3 {\xi^4_3} - 2 {\xi^6_5}
{\xi^5_5}^2 {\xi^4_3}
 - 2 i {\xi^6_5} {\xi^5_5} {\xi^4_3} - 2 {\xi^6_5} {\xi^4_3} + i {\xi^5_5}^4 +
 2 i {\xi^5_5}^2 + i)
 \\+
 \frac{1}{16 {\xi^6_5}^2 {\xi^4_3}} \,
 \overline{w^1} {w^1}^2 ( - i {\xi^6_5} {\xi^5_5}^2 {\xi^4_3} - 2 {\xi^6_5} {\xi^5_5} {\xi^4_3} + i {\xi^6_5} {\xi^4_3} + i {\xi^5_5}^3 +
{\xi^5_5}^2 + i {\xi^5_5} + 1)
 -
 \frac{ i {\xi^5_5} + 1}{2 {\xi^6_5}} \,  \overline{w^1} w^2
 \end{multline*}
are holomorphic.
Let $F : G \rightarrow \Cmath^3$ defined
by
$F=(\varphi^1,\varphi^2,\varphi^3).$
$F$ is a global chart on $G.$
We determine now how the multiplication of $G$ looks like in that chart.
Let $a,x \in G $ with respective second kind canonical coordinates
$(x^1,y^1,x^2,y^2,x^3,y^3), (\alpha^1, \beta^1, \alpha^2, \beta^2, \alpha^3, \beta^3)$
as in  (\ref{x6general}).
With obvious notations, computations yield:
\begin{eqnarray*}
w^1_{a x} &=& w^1_a + w^1_x\\
w^2_{a x} &=& w^2_a + w^2_x-\frac{b^1x^1}{\xi^4_3} \, (\xi^5_5 +\xi^4_3 -\xi^6_5 \xi^4_3 +i)\\
w^3_{a x} &=& w^3_a + w^3_x
-b^2 x^1 +\frac{1}{2} b^1(x^1)^2
+\frac{i-\xi^5_5}{\xi^6_5} \,
\left(
\frac{1}{2} \, (b^1)^2 x^1  -b^2y^1  + b^1x^1 y^1 -a^2y^1 \right).
\end{eqnarray*}
We then get
\begin{equation*}
\varphi^1_{a x} = \varphi^1_a + \varphi^1_x \quad , \quad
\varphi^2_{a x} = \varphi^2_a + \varphi^2_x + C \varphi^1_x \quad , \quad
\varphi^3_{a x} = \varphi^3_a +\varphi^3_x
+ D_1 (\varphi^1_x)^2
+ D_2 \varphi^1_x
+ D_3 \varphi^2_x
\end{equation*}
where
$C =
(\overline{\varphi^1_a} (i {{\xi^6_5}} {{\xi^4_3}} - i {{\xi^5_5}} - i {{\xi^4_3}} + 1))/(2 {{\xi^4_3}}) + (\varphi^1_a ( - {{\xi^6_5}}^2 {{\xi^5_5}} {{\xi^4_3}}^2 - i {{\xi^6_5}}^2 {{\xi^4_3}}^2 +
2 {{\xi^6_5}} {{\xi^5_5}}^2 {{\xi^4_3}}
+ {{\xi^6_5}} {{\xi^5_5}} {{\xi^4_3}}^2 + 2 i {{\xi^6_5}} {{\xi^5_5}} {{\xi^4_3}} + i {{\xi^6_5}} {{\xi^4_3}}^2 - {{\xi^5_5}}^3 - {{\xi^5_5}}^2 {{\xi^4_3}} - i {{\xi^5_5}}^2 - {{\xi^5_5}} - {{\xi^4_3}} - i)
)/(4 {{\xi^6_5}} {{\xi^4_3}}^2) ;     \\
D_1=(\overline{\varphi^1_a} ( - i {{\xi^6_5}}^2 {{\xi^5_5}}^2 {{\xi^4_3}}^2 + 2 {{\xi^6_5}}^2
 {{\xi^5_5}} {{\xi^4_3}}^2 + i {{\xi^6_5}}^2 {{\xi^4_3}}^2 + 2 i {{\xi^6_5}} {{\xi^5_5}}^3 {{\xi^4_3}} +
 i {{\xi^6_5}} {{\xi^5_5}}^2 {{\xi^4_3}}^2 - 2 {{\xi^6_5}} {{\xi^5_5}}^2 {{\xi^4_3}} - 2 {{\xi^6_5}} {{\xi^5_5}} {{\xi^4_3}}^2 + 2 i {{\xi^6_5}} {{\xi^5_5}} {{\xi^4_3}} + 3 i {{\xi^6_5}} {{\xi^4_3}}^2 - 2 {{\xi^6_5}} {{\xi^4_3}} - i {{\xi^5_5}}^4 - i {{\xi^5_5}}^3 {{\xi^4_3}} - {{\xi^5_5}}^2 {{\xi^4_3}} - 2 i {{\xi^5_5}}^2 - i {{\xi^5_5}} {{\xi^4_3}} - {{\xi^4_3}} - i))/(16 {{\xi^6_5}}^2 {{\xi^4_3}}^2) + (\varphi^1_a (i
 {{\xi^6_5}}^2 {{\xi^5_5}}^2 {{\xi^4_3}}^2 - 2 {{\xi^6_5}}^2 {{\xi^5_5}} {{\xi^4_3}}^2 - i {{\xi^6_5}}
^2 {{\xi^4_3}}^2 - 2 i {{\xi^6_5}} {{\xi^5_5}}^3 {{\xi^4_3}} - 2 i {{\xi^6_5}} {{\xi^5_5}}^2 {{\xi^4_3}}
^2 + 2 {{\xi^6_5}} {{\xi^5_5}}^2 {{\xi^4_3}} - 2 i {{\xi^6_5}} {{\xi^5_5}} {{\xi^4_3}} - 2 i {{\xi^6_5}}
{{\xi^4_3}}^2 + 2 {{\xi^6_5}} {{\xi^4_3}} + i {{\xi^5_5}}^4 + 2 i {{\xi^5_5}}^3 {{\xi^4_3}} + 2 {{\xi^5_5}}^2 {{\xi^4_3}} + 2 i {{\xi^5_5}}^2 + 2 i {{\xi^5_5}} {{\xi^4_3}} + 2 {{\xi^4_3}} + i))/(16 {{\xi^6_5}}^2 {{\xi^4_3}}^2)
;\\
D_2=(\overline{\varphi^1_a}^2 ( - i {{\xi^6_5}}^2 {{\xi^5_5}} {{\xi^4_3}}^2 - {{\xi^6_5}}^2 {\xi^4_3}^2 + 2 i {{\xi^6_5}} {{\xi^5_5}}^2 {{\xi^4_3}} + i {{\xi^6_5}} {{\xi^5_5}} {{\xi^4_3}}^2 + 4 {{\xi^6_5}} {{\xi^5_5}} {{\xi^4_3}} + {{\xi^6_5}} {{\xi^4_3}}^2 - 2 i {{\xi^6_5}} {{\xi^4_3}} - i {{\xi^5_5}}^3 - i
 {{\xi^5_5}}^2 {{\xi^4_3}} - 3 {{\xi^5_5}}^2 - 2 {{\xi^5_5}} {{\xi^4_3}} - i {{\xi^5_5}} + i {{\xi^4_3}} -
 3))/(16 {{\xi^6_5}} {{\xi^4_3}}) + (\overline{\varphi^1_a} \varphi^1_a ({{\xi^6_5}} {{\xi^4_3}} - {{\xi^5_5}} - {{\xi^4_3}} +
 i))/4 + ( - i \overline{\varphi^2_a} {{\xi^4_3}})/2 + ({\varphi^1_a}^2 (i {{\xi^6_5}}^3 {{\xi^5_5}} {{\xi^4_3}}^2 -
{{\xi^6_5}}^3 {{\xi^4_3}}^2 - 2 i {{\xi^6_5}}^2 {{\xi^5_5}}^2 {{\xi^4_3}} - i {{\xi^6_5}}^2 {{\xi^5_5}}
 {{\xi^4_3}}^2 + {{\xi^6_5}}^2 {{\xi^4_3}}^2 - 2 i {{\xi^6_5}}^2 {{\xi^4_3}} + i {{\xi^6_5}} {{\xi^5_5}}
^3 + {{\xi^6_5}} {{\xi^5_5}}^2 + i {{\xi^6_5}} {{\xi^5_5}} + {{\xi^6_5}} + i {{\xi^5_5}}^3 + {{\xi^5_5}}
^2 + i {{\xi^5_5}} + 1))/(16 {{\xi^6_5}}^2 {{\xi^4_3}})
+ (\varphi^2_a (i {{\xi^6_5}} {{\xi^4_3}} - i {{\xi^5_5}} - 1))/(2 {{\xi^6_5}});\\
D_3 =(\overline{\varphi^1_a} ( - i {{\xi^5_5}} - 1))/(2 {{\xi^6_5}}).$
%%%%%%%%%%%%%%%%%%%%%%%%%%%%%%%%%%%%%%%%%%%%%%%%%%%%%%%%%%%%%%%%%%%%%%%%%%%%%%%%%%%%%%%%%%
%%%%%%%%%%%%%%%%%%%%%%%%%%%%%%%%%%%%%%%%%%%%%%%%%%%%%%%%%%%%%%%%%%%%%%%%%%%%%%%%%%%%%%%%%%
\section{Lie Algebra $ {\mathcal{G}}_{6,8}$ (isomorphic to $M9$).}
%%%%%%%%%%%%%%%%%%%%%%%%%%%%%%%%%%%%%%%%%%%%%%%%%%%%%%%%%%%%%%%%%%%%%%%%%%%%%%%%%%%%%%%%%%
Commutation relations for
$ {\mathcal{G}}_{6,8}:$
$[x_1,x_2]=x_4$;
$[x_1,x_4]=x_5$;
$[x_2,x_3]=x_5$;
$[x_2,x_4]=x_6.$
%\subsection{}
{\fontsize{8}{10} \selectfont
\begin{equation} \label{68general}
J = \begin{pmatrix}
\frac{{\xi^5_5} {\xi^4_3} + {\xi^5_5} {\xi^2_1} - {\xi^3_3} {\xi^2_1}}{{\xi^4_3}}&
 - \frac{({\xi^5_5} {\xi^4_3} + {\xi^5_5} {\xi^2_1} - {\xi^3_3} {\xi^2_1})^2 + {\xi^4_3}^2
 }{{\xi^4_3}^2 {\xi^2_1}}&
0&
0&
0&
0\\
\boxed{\xi^2_1}&
 - \frac{{\xi^5_5} {\xi^4_3} + {\xi^5_5} {\xi^2_1} - {\xi^3_3} {\xi^2_1}}{{\xi^4_3}}&
0&
0&
0&
0\\
\boxed{\xi^3_1}&
\boxed{\xi^3_2}&
\boxed{\xi^3_3}&
 - \frac{{\xi^3_3}^2 + 1}{{\xi^4_3}}&
0&
0\\
 *&
 *&
\boxed{\xi^4_3}&
 - {\xi^3_3}&
0&
0\\
*&
*&
 - \frac{({\xi^6_4} {\xi^4_3} - {\xi^6_3} {\xi^5_5} + {\xi^6_3} {\xi^3_3}) ({\xi^4_3} + {\xi^2_1})}{
{\xi^4_3} {\xi^2_1}}&
*&
\boxed{\xi^5_5}&
 - \frac{({\xi^5_5}^2 + 1) ({\xi^4_3} + {\xi^2_1})}{{\xi^4_3} {\xi^2_1}}\\
\boxed{\xi^6_1}&
\boxed{\xi^6_2}&
\boxed{\xi^6_3}&
\boxed{\xi^6_4}&
\frac{{\xi^4_3} {\xi^2_1}}{{\xi^4_3} + {\xi^2_1}}&
 - {\xi^5_5}\end{pmatrix}
\end{equation}
where
$ \mathbf{J^4_1}=({\xi^4_3} {\xi^3_3} {\xi^3_1} + {\xi^4_3} {\xi^3_2} {\xi^2_1} - {\xi^3_3} {\xi^3_1}
{\xi^2_1} + ({\xi^4_3} + {\xi^2_1}) {\xi^5_5} {\xi^3_1})/({\xi^3_3}^2 + 1);$
 \myquad
$  \mathbf{J^4_2}=( - ((({\xi^5_5} - {\xi^3_3}) ({\xi^4_3} {\xi^3_3} {\xi^3_1} + {\xi^4_3} {\xi^3_2}
{\xi^2_1} - {\xi^3_3} {\xi^3_1} {\xi^2_1}) - ({\xi^4_3} + {\xi^2_1}) {\xi^5_5} {\xi^3_3} {\xi^3_1}
) ({\xi^4_3} + {\xi^2_1}) + (({\xi^4_3} + {\xi^2_1})^2 {\xi^5_5}^2 + ({\xi^3_3}^2 + 1) {\xi^4_3}^2) {\xi^3_1}))/(({\xi^3_3}^2 + 1) {\xi^4_3} {\xi^2_1});$
 \myquad
$  \mathbf{J^5_1}=( - ({\xi^6_4} {\xi^5_5} {\xi^4_3}^2 {\xi^3_1} + {\xi^6_4} {\xi^5_5} {\xi^4_3} {\xi^3_1} {\xi^2_1} + {\xi^6_4} {\xi^4_3}^2 {\xi^3_3} {\xi^3_1} + {\xi^6_4} {\xi^4_3}^2 {\xi^3_2}
{\xi^2_1} - {\xi^6_4} {\xi^4_3} {\xi^3_3} {\xi^3_1} {\xi^2_1} + {\xi^6_3} {\xi^4_3} {\xi^3_3}^2
{\xi^3_1} + {\xi^6_3} {\xi^4_3} {\xi^3_1} + {\xi^6_2} {\xi^4_3} {\xi^3_3}^2 {\xi^2_1} + {\xi^6_2}
 {\xi^4_3} {\xi^2_1} + {\xi^6_1} {\xi^5_5} {\xi^3_3}^2 {\xi^2_1} + {\xi^6_1} {\xi^5_5} {\xi^2_1}
- {\xi^6_1} {\xi^3_3}^3 {\xi^2_1} - {\xi^6_1} {\xi^3_3} {\xi^2_1}) ({\xi^4_3} + {\xi^2_1}))/((
{\xi^3_3}^2 + 1) {\xi^4_3}^2 {\xi^2_1});$
 \myquad
$  \mathbf{J^5_2}=(({\xi^6_4} {\xi^5_5}^2 {\xi^4_3}^3 {\xi^3_1} + 2 {\xi^6_4} {\xi^5_5}^2 {\xi^4_3}^2 {\xi^3_1} {\xi^2_1} + {\xi^6_4} {\xi^5_5}^2 {\xi^4_3} {\xi^3_1} {\xi^2_1}^2 + {\xi^6_4}
{\xi^5_5} {\xi^4_3}^3 {\xi^3_2} {\xi^2_1} - 2 {\xi^6_4} {\xi^5_5} {\xi^4_3}^2 {\xi^3_3} {\xi^3_1} {\xi^2_1} + {\xi^6_4} {\xi^5_5} {\xi^4_3}^2 {\xi^3_2} {\xi^2_1}^2 - 2 {\xi^6_4} {\xi^5_5} {\xi^4_3} {\xi^3_3} {\xi^3_1} {\xi^2_1}^2 - {\xi^6_4} {\xi^4_3}^3 {\xi^3_3} {\xi^3_2} {\xi^2_1} +
{\xi^6_4} {\xi^4_3}^3 {\xi^3_1} - {\xi^6_4} {\xi^4_3}^2 {\xi^3_3} {\xi^3_2} {\xi^2_1}^2 + {\xi^6_4} {\xi^4_3} {\xi^3_3}^2 {\xi^3_1} {\xi^2_1}^2 - {\xi^6_3} {\xi^4_3}^2 {\xi^3_3}^2 {\xi^3_2} {\xi^2_1} - {\xi^6_3} {\xi^4_3}^2 {\xi^3_2} {\xi^2_1} + 2 {\xi^6_2} {\xi^5_5} {\xi^4_3}^2
{\xi^3_3}^2 {\xi^2_1} + 2 {\xi^6_2} {\xi^5_5} {\xi^4_3}^2 {\xi^2_1} + {\xi^6_2} {\xi^5_5} {\xi^4_3} {\xi^3_3}^2 {\xi^2_1}^2 + {\xi^6_2} {\xi^5_5} {\xi^4_3} {\xi^2_1}^2 - {\xi^6_2} {\xi^4_3}
 {\xi^3_3}^3 {\xi^2_1}^2 - {\xi^6_2} {\xi^4_3} {\xi^3_3} {\xi^2_1}^2 + {\xi^6_1} {\xi^5_5}^2
 {\xi^4_3}^2 {\xi^3_3}^2 + {\xi^6_1} {\xi^5_5}^2 {\xi^4_3}^2 + 2 {\xi^6_1} {\xi^5_5}^2 {\xi^4_3} {\xi^3_3}^2 {\xi^2_1} + 2 {\xi^6_1} {\xi^5_5}^2 {\xi^4_3} {\xi^2_1} + {\xi^6_1} {\xi^5_5}^2 {\xi^3_3}^2 {\xi^2_1}^2 + {\xi^6_1} {\xi^5_5}^2 {\xi^2_1}^2 - 2 {\xi^6_1} {\xi^5_5}
{\xi^4_3} {\xi^3_3}^3 {\xi^2_1} - 2 {\xi^6_1} {\xi^5_5} {\xi^4_3} {\xi^3_3} {\xi^2_1} - 2 {\xi^6_1} {\xi^5_5} {\xi^3_3}^3 {\xi^2_1}^2 - 2 {\xi^6_1} {\xi^5_5} {\xi^3_3} {\xi^2_1}^2 + {\xi^6_1} {\xi^4_3}^2 {\xi^3_3}^2 + {\xi^6_1} {\xi^4_3}^2 + {\xi^6_1} {\xi^3_3}^4 {\xi^2_1}^2 +
{\xi^6_1} {\xi^3_3}^2 {\xi^2_1}^2) ({\xi^4_3} + {\xi^2_1}))/(({\xi^3_3}^2 + 1) {\xi^4_3}^3
 {\xi^2_1}^2);$
 \myquad
$  \mathbf{J^5_4}=(({\xi^6_4} {\xi^5_5} {\xi^4_3} + {\xi^6_4} {\xi^4_3} {\xi^3_3} + {\xi^6_3} {\xi^3_3}
^2 + {\xi^6_3}) ({\xi^4_3} + {\xi^2_1}))/({\xi^4_3}^2 {\xi^2_1});$
 }
and the parameters are subject to the condition
\begin{equation}
\label{condM9general}
\xi^2_1 \xi^4_3(\xi^4_3+\xi^2_1)  \neq 0.
\end{equation}
\par Now the automorphism group of
$ {\mathcal{G}}_{6,8}$ is comprised of the matrices
$$
{\fontsize{8}{10} \selectfont
\Phi = \begin{pmatrix}
{b^1_1}&
{b^1_2}&
0&
0&
0&
0\\
0&
{b^2_2}&
0&
0&
0&
0\\
{b^3_1}&
{b^3_2}&
{b^1_1}^2&
0&
0&
0\\
{b^4_1}&
{b^4_2}&
0&
{b^2_2} {b^1_1}&
0&
0\\
{b^5_1}&
{b^5_2}&
{b^5_3}&
 - ({b^4_1} {b^1_2} + {b^3_1} {b^2_2} - {b^4_2} {b^1_1})&
{b^2_2} {b^1_1}^2&
{b^2_2} {b^1_2} {b^1_1}\\
{b^6_1}&
{b^6_2}&
{b^6_3}&
 - {b^4_1} {b^2_2}&
0&
{b^2_2}^2 {b^1_1}\end{pmatrix}
 }
 $$
where $ {b^2_2} {b^1_1}\neq 0. $
Taking suitable values for the $b^i_j$'s,
we are led to the case where
 $ {\xi^2_1}=1, {\xi^3_1}= {\xi^3_2}= {\xi^5_5}= {\xi^6_1}= {\xi^6_3}=0$
and moreover
$ {\xi^6_2}= {\xi^6_4}=0.$
Hence any   $ J $ in (\ref{68general})
is equivalent to :
\begin{equation}
{\fontsize{8}{10} \selectfont
\label{68final}
 J({\xi^3_3},{\xi^4_3}) = \begin{pmatrix}
 - {\xi^3_3}/{\xi^4_3}&
 - ({\xi^4_3}^2 + {\xi^3_3}^2)/{\xi^4_3}^2&
0&
0&
0&
0\\
1&
{\xi^3_3}/{\xi^4_3}&
0&
0&
0&
0\\
0&
0&
{\xi^3_3}&
 - ({\xi^3_3}^2 + 1)/{\xi^4_3}&
0&
0\\
0&
0&
{\xi^4_3}&
 - {\xi^3_3}&
0&
0\\
0&
0&
0&
0&
0&
 - ({\xi^4_3} + 1)/{\xi^4_3}\\
0&
0&
0&
0&
{\xi^4_3}/({\xi^4_3} + 1)&
0\end{pmatrix}
}
\end{equation}
 where $ {\xi^4_3} \neq 0,- 1.$
 $ J({\xi^3_3},{\xi^4_3})\cong J(\eta^3_3,\eta^4_3)$
if and only if $ {\xi^3_3} = \eta^3_3$ and ${\xi^4_3}=\eta^4_3.$
Commutation relations of $ \mathfrak{m} : $
$  [\tilde{x}_1,\tilde{x}_3]= -{\xi^4_3} \tilde{x}_6  ;$
$  [\tilde{x}_2,\tilde{x}_3]= ({\xi^4_3} + 1) \tilde{x}_5 - {\xi^3_3} \tilde{x}_6;$
 $  [\tilde{x}_2,\tilde{x}_4]=
 \frac{{\xi^3_3} ( - {\xi^4_3}^2 + 1)}{{\xi^4_3}^2} \,  \tilde{x}_5
 +  \frac{{\xi^4_3} + {\xi^3_3}^2}{\xi^4_3} \,  \tilde{x}_6  .$
%\subsection{Conclusions.}
\par
From (\ref{68general}), $\mathfrak{X}_{6,8}$
is a submanifold of dimension 10 in $\Rmath^{36}.$
It is  the disjoint union of the continuously many orbits of the
$J(\xi^3_3,\xi^4_3)$ defined in (\ref{68final}) where $\xi^4_3 \neq 0,-1.$
%%%%%%%%%%%%%%%%%%%%%%%%%%%%%%%%%%%%%%%%%%%%%%%%%%%%%%%%%%%%%%%%%%%%%%%%%%%%%%%%%%%%%%%%%%
\subsection{}
%\subsubsection{}
\begin{equation*}
X_1 = \frac{\partial}{\partial x^1} - y^1 \, \frac{\partial}{\partial y^2}
+ y^2\, \frac{\partial}{\partial x^3}
+ \frac{(y^1)^2}{2} \, \frac{\partial}{\partial y^3}
\quad , \quad
X_2 = \frac{\partial}{\partial y^1} - x^2\frac{\partial}{\partial x^3}
- y^2\frac{\partial}{\partial y^3}.
\end{equation*}
%\subsubsection{Holomorphic functions for $J(\xi^3_3,\xi^4_3).$}
Let $G$ denote the group $G_0$ endowed with the left invariant  structure
of complex manifold  defined by
$J(\xi^3_3,\xi^4_3)$ defined in (\ref{68final}) where $\xi^4_3 \neq 0,-1.$
Let $H_{\Cmath}(G)$ the space of complex valued holomorphic functions on $G.$
Then $ H_{\Cmath}(G) =\{f \in C^{\infty}(G_0) \; ; \; \tilde{X}_j^{-} \, f = 0
\; \forall j=1,3,5\}.$
One has
\begin{eqnarray*}
\tilde{X}_1^{-}  &=&
2\, \frac{\partial}{\partial \overline{w^1}}
-y^1 \left(1-\frac{i\xi^3_3}{\xi^4_3}\right)\, \frac{\partial}{\partial y^2}
-\left(y^2 \left(1-\frac{i\xi^3_3}{\xi^4_3}\right)+ix^2\right)\, \frac{\partial}{\partial x^3}
+\left(\frac{(y^1)^2}{2} \left(1-\frac{i\xi^3_3}{\xi^4_3}\right)-iy^2\right)\, \frac{\partial}{\partial y^3}
\quad , \\
\tilde{X}_3^{-}  &=& 2  \frac{\partial}{\partial \overline{w^2}}
\quad , \quad
\tilde{X}_5^{-}  =2\, \frac{\partial}{\partial \overline{w^3}}
\quad ,
\end{eqnarray*}
where
$
w^1 = x^1 +\frac{\xi^3_3}{\xi^4_3}\, y^1 +iy^1
\quad , \quad
w^2 = x^2 -\frac{\xi^3_3}{\xi^4_3}\, y^2 +\frac{i}{\xi^4_3}\, y^2
\quad , \quad
w^3 = x^3 +i \, \frac{\xi^4_3+1}{\xi^4_3}\,   y^3.
$
Then $f \in C^{\infty}(G_0)$ is in $ H_{\Cmath}(G)$ if and only if it is holomorphic
with respect to  $w^2$ and $w^3$ and satisfies the equation
\begin{eqnarray*}
2\, \frac{\partial f}{\partial \overline{w^1}}
-i A\frac{w^1-\overline{w^1}}{2}\, \frac{\partial f}{\partial {w^2}}
-\left( iw^2 + B \, \frac{(w^1-\overline{w^1})^2}{8} \right)
\frac{\partial f}{\partial {w^3}}
&=&0
\end{eqnarray*}
where
\begin{equation*}
A=\frac{1}{\xi^4_3} \left( \xi^3_3\left(1-\frac{1}{\xi^4_3}\right)
-i \left(1+\frac{{\xi^3_3}^2}{\xi^4_3}\right) \right)
\quad , \quad
B=\frac{\xi^4_3+1}{{\xi^4_3}^2}\, \left( \xi^3_3+i\xi^4_3\right)
\, .
\end{equation*}
The 3 functions
\begin{eqnarray*}
\varphi^1 &=&w^1\quad , \quad
\varphi^2 = w^2 +  \frac{iA}{4} \left(w^1\overline{w^1} -\frac{(\overline{w^1})^2}{2}\right)\quad , \\
\varphi^3 &=&w^3 +
\frac{i}{2}\, \overline{w^1} w^2 -\frac{B}{48}\, (w^1-\overline{w^1})^3
-\frac{A}{16}\, w^1(\overline{w^1})^2
+\frac{A}{24}\, (\overline{w^1})^3
\end{eqnarray*}
are holomorphic.
Let $F : G \rightarrow \Cmath^3$ defined
by
$F=(\varphi^1,\varphi^2,\varphi^3).$
$F$ is a global chart on $G.$
We determine now how the multiplication of $G$ looks like in that chart.
Let $a,x \in G $ with respective second kind canonical coordinates
$(x^1,y^1,x^2,y^2,x^3,y^3), (\alpha^1, \beta^1, \alpha^2, \beta^2, \alpha^3, \beta^3)$
as in  (\ref{x6general}).
With obvious notations,
%$a =[w^1_a, w^2_a,w^3_a],$
%$x =[w^1_x, w^2_x,w^3_x],$
%$a \, x =[w^1_{a x}, w^2_{a x},w^3_{a x}],$
%$a =[\varphi^1_a, \varphi^2_a,\varphi^3_a],$
%$x =[\varphi^1_x, \varphi^2_x,\varphi^3_x],$
%$a \, x =[\varphi^1_{a x}, \varphi^2_{a x},\varphi^3_{a x}].$
computations yield:
\begin{eqnarray*}
w^1_{a x} &=& w^1_a + w^1_x \quad , \quad
w^2_{a x} = w^2_a + w^2_x+b^1x^1\, \frac{\xi^3_3-i}{\xi^4_3} \quad ,\\
w^3_{a x} &=& w^3_a + w^3_x
-b^2 x^1
+b^1\frac{(x^1)^2}{2}   -a^2y^1
+i \, \frac{\xi^4_3+1}{\xi^4_3}\left( \frac{(b^1)^2}{2} x^1  -(b^2 -b^1x^1) y^1 \right).
\end{eqnarray*}
We then get
\begin{equation*}
\varphi^1_{a x} = \varphi^1_a + \varphi^1_x
\quad , \quad
\varphi^2_{a x} = \varphi^2_a + \varphi^2_x + \chi^2(a,x)
\quad , \quad
\varphi^3_{a x} = \varphi^3_a +\varphi^3_x + \chi^3(a,x)
\quad ,
\end{equation*}
where
\begin{equation*}
\chi^2(a,x) =
\frac{\xi^3_3-i}{4{\xi^4_3}^2}\, \varphi^1_x
\left( 2i\xi^4_3 \overline{\varphi^1_a} +(\xi^3_3-i\xi^4_3) \varphi^1_a  \right)
\quad ,
\end{equation*}
\begin{multline*}
\chi^3(a,x) =
\frac{i}{2}\,\varphi^2_x \overline{\varphi^1_a}
+\frac{1}{16{\xi^4_3}^2}\, (\varphi^1_x)^2  (4i{\xi^4_3}^2-(\xi^3_3-i)(3\xi^4_3+i\xi^3_3))
\left(\overline{\varphi^1_a} -\varphi^1_a\right)
\\
+\varphi^1_x
\left( \frac{1}{16{\xi^4_3}} \, (\overline{\varphi^1_a})^2 (-\xi^3_3(\xi^4_3-i\xi^3_3)-i\xi^4_3-3\xi^3_3+2i)
\right. \\       \left.
+\frac{1}{16{\xi^4_3}} \, ({\varphi^1_a})^2 (-\xi^3_3(\xi^4_3-i\xi^3_3)-i\xi^4_3+\xi^3_3-2i)
\right. \\       \left.
+\frac{\xi^3_3(\xi^4_3-i\xi^3_3)}{4{\xi^4_3}} \, \overline{\varphi^1_a} \varphi^1_a
-\frac{i\xi^4_3}{2} \, \overline{\varphi^2_a}
+\frac{i(\xi^4_3+1)}{2} \, {\varphi^2_a}
\right).
\end{multline*}
%%%%%%%%%%%%%%%%%%%%%%%%%%%%%%%%%%%%%%%%%%%%%%%%%%%%%%%%%%%%%%%%%%%%%%%%%%%%%%%%%%%%%%%%%%
%%%%%%%%%%%%%%%%%%%%%%%%%%%%%%%%%%%%%%%%%%%%%%%%%%%%%%%%%%%%%%%%%%%%%%%%%%%%%%%%%%%%%%%%%%
\section{Lie Algebra $ M10.$}
%%%%%%%%%%%%%%%%%%%%%%%%%%%%%%%%%%%%%%%%%%%%%%%%%%%%%%%%%%%%%%%%%%%%%%%%%%%%%%%%%%%%%%%%%%
Commutation relations for
$ {M10}:$
$[x_1,x_2]=x_3$;
$[x_1,x_3]=x_5$;
$[x_1,x_4]=x_6$;
$[x_2,x_3]= - x_6$;
$[x_2,x_4]=x_5$.
\subsection{Case $ {\xi^4_3} \neq \xi^2_1 .$}
\label{M10case1}
%{\fontsize{8}{10}\selectfont
\begin{equation}
\label{M10general1}
J = \begin{pmatrix}
\boxed{\xi^1_1}&
-\frac{{\xi^1_1}^2+1}{\xi^2_1}&
0&
0&
0&
0\\
\boxed{\xi^2_1}&
-\xi^1_1&
0&
0&
0&
0\\
\boxed{\xi^3_1}&
%$  J^3_2=(({\xi^3_3}^2 + 1) {\xi^4_1} - ({\xi^3_3} + {\xi^1_1}) {\xi^4_3} {\xi^3_1})/({\xi^4_3} {\xi^2_1});$\\
\frac{({\xi^3_3}^2 + 1) {\xi^4_1} - ({\xi^3_3} + {\xi^1_1}) {\xi^4_3} {\xi^3_1}}{{\xi^4_3} {\xi^2_1}}&
\boxed{\xi^3_3}&
-\frac{{\xi^3_3}^2+1}{\xi^4_3}&
0&
0\\
\boxed{\xi^4_1}&
%$  J^4_2=(({\xi^3_3} - {\xi^1_1}) {\xi^4_1} - {\xi^4_3} {\xi^3_1})/{\xi^2_1};$\\
\frac{({\xi^3_3} - {\xi^1_1}) {\xi^4_1} - {\xi^4_3} {\xi^3_1}}{{\xi^2_1}}&
\boxed{\xi^4_3}&
-{\xi^3_3}&
0&
0\\
*&
*&
\boxed{\xi^5_3}&
*&
r&
-\frac{r^2+1}{b}\\
\boxed{\xi^6_1}&
\boxed{\xi^6_2}&
\boxed{\xi^6_3}&
*&
b&
-r
 \end{pmatrix}
\end{equation}
%}
where
{\fontsize{6}{10}\selectfont
$  \mathbf{J^5_1}=(((2 {\xi^4_3} {\xi^1_1} - {\xi^3_3} {\xi^2_1}) {\xi^2_1} + ({\xi^1_1}^2 + 1) {\xi^3_3} - ({\xi^1_1}^2 + 1 + {\xi^2_1}^2) {\xi^1_1}) {\xi^6_1} {\xi^4_3}^2 + (({\xi^4_3} {\xi^1_1}^2 + {\xi^4_3} - 2 {\xi^3_3} {\xi^2_1} {\xi^1_1} - 2 {\xi^2_1}) {\xi^4_3} + ({\xi^3_3}^
2 + 1) {\xi^2_1}^2) {\xi^5_3} {\xi^4_1} + ({\xi^4_3}^2 {\xi^2_1} - {\xi^4_3} {\xi^2_1}^2 -
{\xi^4_3} {\xi^1_1}^2 - {\xi^4_3} + {\xi^3_3}^2 {\xi^2_1} + {\xi^2_1}) {\xi^6_2} {\xi^4_3} {\xi^2_1} - (((({\xi^3_3} - {\xi^1_1}) {\xi^4_1} - {\xi^4_3} {\xi^3_1}) ({\xi^3_3}^2 + 1) - ({\xi^4_3}^2 {\xi^3_1} - {\xi^4_3} {\xi^4_1} {\xi^3_3} - {\xi^4_3} {\xi^4_1} {\xi^1_1} + 2 {\xi^4_1}
{\xi^3_3} {\xi^2_1}) {\xi^4_3}) {\xi^2_1} + ({\xi^1_1}^2 + 1 + {\xi^2_1}^2) {\xi^4_3}^2 {\xi^3_1}) {\xi^6_3})/
 (({\xi^4_3} {\xi^1_1} -{\xi^3_3} {\xi^2_1})^2 +({\xi^4_3} -{\xi^2_1})^2 ) {\xi^4_3});$
 \myquad
$  \mathbf{J^5_2}=({\xi^6_3} {\xi^4_3}^2 {\xi^4_1} {\xi^2_1} {\xi^1_1}^2 + {\xi^6_3} {\xi^4_3}^2
{\xi^4_1} {\xi^2_1} - {\xi^6_3} {\xi^4_3}^2 {\xi^3_3} {\xi^3_1} {\xi^2_1}^2 + {\xi^6_3} {\xi^4_3}^2 {\xi^3_3} {\xi^3_1} {\xi^1_1}^2 + {\xi^6_3} {\xi^4_3}^2 {\xi^3_3} {\xi^3_1} + {\xi^6_3}
{\xi^4_3}^2 {\xi^3_1} {\xi^2_1}^2 {\xi^1_1} + {\xi^6_3} {\xi^4_3}^2 {\xi^3_1} {\xi^1_1}^3 +
{\xi^6_3} {\xi^4_3}^2 {\xi^3_1} {\xi^1_1} + {\xi^6_3} {\xi^4_3} {\xi^4_1} {\xi^3_3}^2 {\xi^2_1}
^2 - {\xi^6_3} {\xi^4_3} {\xi^4_1} {\xi^3_3}^2 {\xi^1_1}^2 - {\xi^6_3} {\xi^4_3} {\xi^4_1} {\xi^3_3}^2 - 2 {\xi^6_3} {\xi^4_3} {\xi^4_1} {\xi^3_3} {\xi^2_1}^2 {\xi^1_1} - {\xi^6_3} {\xi^4_3} {\xi^4_1} {\xi^2_1}^2 - {\xi^6_3} {\xi^4_3} {\xi^4_1} {\xi^1_1}^2 - {\xi^6_3} {\xi^4_3} {\xi^4_1} - 2 {\xi^6_3} {\xi^4_3} {\xi^3_3}^2 {\xi^3_1} {\xi^2_1} {\xi^1_1} - 2 {\xi^6_3} {\xi^4_3}
{\xi^3_1} {\xi^2_1} {\xi^1_1} + 2 {\xi^6_3} {\xi^4_1} {\xi^3_3}^3 {\xi^2_1} {\xi^1_1} - {\xi^6_3}
 {\xi^4_1} {\xi^3_3}^2 {\xi^2_1} {\xi^1_1}^2 + {\xi^6_3} {\xi^4_1} {\xi^3_3}^2 {\xi^2_1} + 2
{\xi^6_3} {\xi^4_1} {\xi^3_3} {\xi^2_1} {\xi^1_1} - {\xi^6_3} {\xi^4_1} {\xi^2_1} {\xi^1_1}^2 +
{\xi^6_3} {\xi^4_1} {\xi^2_1} - {\xi^6_2} {\xi^4_3}^2 {\xi^3_3} {\xi^2_1}^3 + {\xi^6_2} {\xi^4_3}^2 {\xi^3_3} {\xi^2_1} {\xi^1_1}^2 + {\xi^6_2} {\xi^4_3}^2 {\xi^3_3} {\xi^2_1} + {\xi^6_2}
{\xi^4_3}^2 {\xi^2_1}^3 {\xi^1_1} + {\xi^6_2} {\xi^4_3}^2 {\xi^2_1} {\xi^1_1}^3 + {\xi^6_2}
{\xi^4_3}^2 {\xi^2_1} {\xi^1_1} - 2 {\xi^6_2} {\xi^4_3} {\xi^3_3}^2 {\xi^2_1}^2 {\xi^1_1} - 2
 {\xi^6_2} {\xi^4_3} {\xi^2_1}^2 {\xi^1_1} - {\xi^6_1} {\xi^4_3}^3 {\xi^2_1} {\xi^1_1}^2 - {\xi^6_1} {\xi^4_3}^3 {\xi^2_1} + {\xi^6_1} {\xi^4_3}^2 {\xi^2_1}^2 {\xi^1_1}^2 + {\xi^6_1} {\xi^4_3}^2 {\xi^2_1}^2 + {\xi^6_1} {\xi^4_3}^2 {\xi^1_1}^4 + 2 {\xi^6_1} {\xi^4_3}^2 {\xi^1_1}^2 + {\xi^6_1} {\xi^4_3}^2 - {\xi^6_1} {\xi^4_3} {\xi^3_3}^2 {\xi^2_1} {\xi^1_1}^2 - {\xi^6_1} {\xi^4_3} {\xi^3_3}^2 {\xi^2_1} - {\xi^6_1} {\xi^4_3} {\xi^2_1} {\xi^1_1}^2 - {\xi^6_1}
{\xi^4_3} {\xi^2_1} - {\xi^5_3} {\xi^4_3}^3 {\xi^3_1} {\xi^1_1}^2 - {\xi^5_3} {\xi^4_3}^3 {\xi^3_1} + {\xi^5_3} {\xi^4_3}^2 {\xi^4_1} {\xi^3_3} {\xi^1_1}^2 + {\xi^5_3} {\xi^4_3}^2 {\xi^4_1} {\xi^3_3} - {\xi^5_3} {\xi^4_3}^2 {\xi^4_1} {\xi^1_1}^3 - {\xi^5_3} {\xi^4_3}^2 {\xi^4_1}
{\xi^1_1} + 2 {\xi^5_3} {\xi^4_3}^2 {\xi^3_3} {\xi^3_1} {\xi^2_1} {\xi^1_1} + 2 {\xi^5_3} {\xi^4_3}^2 {\xi^3_1} {\xi^2_1} - 2 {\xi^5_3} {\xi^4_3} {\xi^4_1} {\xi^3_3}^2 {\xi^2_1} {\xi^1_1} + 2
 {\xi^5_3} {\xi^4_3} {\xi^4_1} {\xi^3_3} {\xi^2_1} {\xi^1_1}^2 - 2 {\xi^5_3} {\xi^4_3} {\xi^4_1}
{\xi^3_3} {\xi^2_1} + 2 {\xi^5_3} {\xi^4_3} {\xi^4_1} {\xi^2_1} {\xi^1_1} - {\xi^5_3} {\xi^4_3} {\xi^3_3}^2 {\xi^3_1} {\xi^2_1}^2 - {\xi^5_3} {\xi^4_3} {\xi^3_1} {\xi^2_1}^2 + {\xi^5_3} {\xi^4_1} {\xi^3_3}^3 {\xi^2_1}^2 - {\xi^5_3} {\xi^4_1} {\xi^3_3}^2 {\xi^2_1}^2
 {\xi^1_1} + {\xi^5_3} {\xi^4_1} {\xi^3_3} {\xi^2_1}^2 - {\xi^5_3} {\xi^4_1} {\xi^2_1}^2 {\xi^1_1})/
 (({\xi^4_3} {\xi^1_1} -{\xi^3_3} {\xi^2_1})^2 +({\xi^4_3} -{\xi^2_1})^2 ) {\xi^4_3}\xi^2_1);$
 \myquad
$ \mathbf{J^5_4}=( - ({\xi^6_3} {\xi^4_3}^2 {\xi^2_1}^2 + 2 {\xi^6_3} {\xi^4_3} {\xi^3_3} {\xi^2_1} {\xi^1_1} - 2 {\xi^6_3} {\xi^4_3} {\xi^2_1} + {\xi^6_3} {\xi^3_3}^2 {\xi^1_1}^2 + {\xi^6_3}
{\xi^3_3}^2 + {\xi^6_3} {\xi^1_1}^2 + {\xi^6_3} + {\xi^5_3} {\xi^4_3}^2 {\xi^3_3} {\xi^2_1} -
 {\xi^5_3} {\xi^4_3}^2 {\xi^2_1} {\xi^1_1} - 2 {\xi^5_3} {\xi^4_3} {\xi^3_3} {\xi^1_1}^2 - 2
{\xi^5_3} {\xi^4_3} {\xi^3_3} + {\xi^5_3} {\xi^3_3}^3 {\xi^2_1} + {\xi^5_3} {\xi^3_3}^2 {\xi^2_1} {\xi^1_1} + {\xi^5_3} {\xi^3_3} {\xi^2_1} + {\xi^5_3} {\xi^2_1} {\xi^1_1}))/
 (((({\xi^4_3} -{\xi^2_1})\xi^2_1 -  ({\xi^1_1}^2+1))\xi^4_3
 + {\xi^2_1} ( {\xi^3_3}^2 +1))  {\xi^4_3});$
 \myquad
$\mathbf{r}=  \mathbf{J^5_5}= - ((({\xi^4_3} {\xi^1_1} - {\xi^3_3} {\xi^2_1}) {\xi^2_1} + ({\xi^1_1}^2 + 1)
{\xi^3_3}) {\xi^4_3} - ({\xi^3_3}^2 + 1) {\xi^2_1} {\xi^1_1})/((({\xi^4_3} - {\xi^2_1}) {\xi^2_1} - ({\xi^1_1}^2 + 1)) {\xi^4_3} + ({\xi^3_3}^2 + 1) {\xi^2_1});$
 \myquad
$  \mathbf{J^6_4}=( - ((({\xi^4_3} {\xi^3_3} + {\xi^4_3} {\xi^1_1} - 2 {\xi^3_3} {\xi^2_1}) {\xi^4_3}
 + ({\xi^3_3}^2 + 1) ({\xi^3_3} - {\xi^1_1})) {\xi^6_3} {\xi^2_1} - (({\xi^4_3} {\xi^1_1}^2
+ {\xi^4_3} - 2 {\xi^3_3} {\xi^2_1} {\xi^1_1} - 2 {\xi^2_1}) {\xi^4_3} + ({\xi^3_3}^2 + 1) {\xi^2_1}^2) {\xi^5_3}))/((({\xi^4_3}^2 - {\xi^4_3} {\xi^2_1} + {\xi^3_3}^2 + 1) {\xi^2_1} -
({\xi^1_1}^2 + 1) {\xi^4_3}) {\xi^4_3});$\myquad
$ \mathbf{b} = \mathbf{J^6_5}= - (({\xi^4_3} {\xi^1_1}^2 + {\xi^4_3} - 2 {\xi^3_3} {\xi^2_1} {\xi^1_1} - 2 {\xi^2_1}) {\xi^4_3} + ({\xi^3_3}^2 + 1) {\xi^2_1}^2)/((({\xi^4_3} - {\xi^2_1}) {\xi^2_1} - ({\xi^1_1}^2 + 1)) {\xi^4_3} + ({\xi^3_3}^2 + 1) {\xi^2_1});$
 }
and the parameters are subject to the condition
\begin{eqnarray}
\label{condM10general1}
{\xi^2_1} {\xi^4_3}
({\xi^2_1} -{\xi^4_3})
 &\neq& 0 \\
(({\xi^4_3} - {\xi^2_1}) {\xi^2_1} - ({\xi^1_1}^2 + 1)) {\xi^4_3} + ({\xi^3_3}^2 +
1) {\xi^2_1}
 &\neq& 0.
 \end{eqnarray}
 Note that $b \neq 0$ since
$ \text{Num}(b) = -\| Y -Z\|^2$
where $ Y = {\xi^4_3} \left( \begin{smallmatrix} {\xi^1_1}\\1\end{smallmatrix} \right)
 \; , \;  Z = {\xi^2_1} \left( \begin{smallmatrix} {\xi^3_3}\\1\end{smallmatrix} \right).$
\par Now the automorphism group of
$ {M10}$ is comprised of the matrices
\begin{equation}
{\fontsize{8}{10}\selectfont
\label{automM10}
\Phi = \begin{pmatrix}
{b^1_1}&
{b^2_1} u&
0&
0&
0&
0\\
{b^2_1}&
 - {b^1_1} u&
0&
0&
0&
0\\
{b^3_1}&
{b^3_2}&
 - ({b^2_1}^2 + {b^1_1}^2) u&
0&
0&
0\\
{b^4_1}&
{b^4_2}&
0&
{b^2_1}^2 + {b^1_1}^2&
0&
0\\
{b^5_1}&
{b^5_2}&
{b^3_2} {b^1_1} - {b^3_1} {b^2_1} u + {b^4_1} {b^1_1} u + {b^4_2} {b^2_1}&
{b^5_4}&
 - ({b^2_1}^2 + {b^1_1}^2) {b^1_1} u&
({b^2_1}^2 + {b^1_1}^2) {b^2_1}\\
{b^6_1}&
{b^6_2}&
 - ({b^3_2} {b^2_1} + {b^3_1} {b^1_1} u + {b^4_1} {b^2_1} u - {b^4_2} {b^1_1})&
{b^6_4}&
({b^2_1}^2 + {b^1_1}^2) {b^2_1} u&
({b^2_1}^2 + {b^1_1}^2) {b^1_1}\end{pmatrix}
}
\end{equation}
%}
where $u^2=1$ and
$ {b^2_1}^2 + {b^1_1}^2 \neq 0.$
Equivalence by a suitable block-diagonal automorphism (\ref{automM10})
leads to the case $\xi^1_1=0$  and  we then can suppose $0< \xi^2_1 \leqslant 1.$
Equivalence by
\begin{equation}
{\fontsize{8}{10}\selectfont
\label{phiM1011}
\Phi = \begin{pmatrix}
1&
0&
0&
0&
0&
0\\
0&
1&
0&
0&
0&
0\\
0&
{\xi^3_1}/{\xi^2_1}&
1&
0&
0&
0\\
0&
{\xi^4_1}/{\xi^2_1}&
0&
1&
0&
0\\
0&
0&
{\xi^3_1}/{\xi^2_1}&
0&
1&
0\\
0&
0&
{\xi^4_1}/{\xi^2_1}&
0&
0&
1\end{pmatrix}
}
\end{equation}
leads to the case where  moreover
 $ {\xi^3_1}=0, {\xi^4_1}=0.$
Then, equivalence by
$$
{\fontsize{8}{10}\selectfont
\Phi = \begin{pmatrix}
1&
0&
0&
0&
0&
0\\
0&
1&
0&
0&
0&
0\\
0&
0&
1&
0&
0&
0\\
0&
0&
0&
1&
0&
0\\
{b^5_1}&
{b^5_2}&
0&
{\xi^5_3}/{\xi^4_3}&
1&
0\\
0&
0&
0&
{\xi^6_3}/{\xi^4_3}&
0&
1\end{pmatrix}
}
$$
with suitable $b^5_1,b^5_2$
leads to the case where moreover
 $ {\xi^3_1}=0, {\xi^4_1}=0,  {\xi^5_3}=0,  {\xi^6_1}=0 , {\xi^6_2}=0,  {\xi^6_3}=0 :$
\begin{equation}
\label{M10case1final}
 J({\xi^2_1},{\xi^3_3},{\xi^4_3}) =
 \text{diag} \left(
 \begin{pmatrix} 0& -\frac{1}{\xi^2_1}\\ {\xi^2_1}& 0  \end{pmatrix},
 \begin{pmatrix} \xi^3_3& -\frac{{\xi^3_3}^2+1}{\xi^4_3}\\ {\xi^4_3}& -\xi^3_3  \end{pmatrix},
 \begin{pmatrix} r& -\frac{1+r^2}{b}\\ b& -r  \end{pmatrix} \right)
\end{equation}
where
\begin{equation}
\label{condM10final1}
 0 < {\xi^2_1} \leqslant 1 ; \quad
   {\xi^2_1} {\xi^4_3} (  {\xi^2_1}-{\xi^4_3}) \neq0 ; \quad
({\xi^4_3} {\xi^2_1} - {\xi^2_1}^2 - 1) {\xi^4_3} + ({\xi^3_3}^2 + 1) {\xi^2_1}
\neq 0
\end{equation}
 and
$ r=(({\xi^2_1} + 1) ({\xi^2_1} - 1) {\xi^4_3} {\xi^3_3})/(({\xi^4_3} {\xi^2_1}
- {\xi^2_1}^2 - 1) {\xi^4_3} + ({\xi^3_3}^2 + 1) {\xi^2_1}),\\$
$ b= - (({\xi^4_3} - 2 {\xi^2_1}) {\xi^4_3} + ({\xi^3_3}^2 + 1) {\xi^2_1}^2)
/(({\xi^4_3} {\xi^2_1} - {\xi^2_1}^2 - 1) {\xi^4_3} + ({\xi^3_3}^2 + 1) {\xi^2_1}).\\$
Now suppose
 $J({\eta^2_1},{\eta^3_3},{\eta^4_3}) =\Phi^{-1} J({\xi^2_1},{\xi^3_3},{\xi^4_3}) \Phi$
where $\Phi$ is given in (\ref{automM10}) and the $\eta$'s and $\xi$'s satisfy
(\ref{condM10final1}).
Computing the matrix
 $J2=\Phi^{-1} J({\xi^2_1},{\xi^3_3},{\xi^4_3}) \Phi,$ one gets
$ {J2^1_1}=({b^2_1} {b^1_1} ({\xi^2_1}^2 - 1))/({\xi^2_1} ({b^2_1}^2 + {b^1_1}^2)
),$
$ {J2^2_1}= - ({b^2_1}^2 + {b^1_1}^2 {\xi^2_1}^2)/({\xi^2_1} u ({b^2_1}^2 + {b^1_1}^2)),$
$ {J2^3_3}={\xi^3_3},$
$ {J2^4_3}=-{\xi^4_3}u.$
From these formulae, we see that a necessary condition for equivalence
is that $ \eta^4_3=-u {\xi^4_3} \; (u = \pm 1) \;$ and $\eta^3_3 ={\xi^3_3}.$
As  $ u=1$ would change the sign of $ {\xi^2_1}$, we conclude that
$u=-1.$ Now to keep $J2^1_1=0$, one must have either $\xi^2_1=1$ or
$b^1_1b^2_1=0.$
If $\xi^2_1=1,$
or if $\xi^2_1<1$ and $ b^2_1=0,$
then $\eta^2_1=\xi^2_1$ and $\eta^3_3=\xi^3_3,\eta^4_3=\xi^4_3.$
If $\xi^2_1<1$ and $ b^1_1=0,$
then $\eta^2_1=1/\xi^2_1>1$ which is contradictory.
Hence
 $J({\eta^2_1},{\eta^3_3},{\eta^4_3})$ and $J({\xi^2_1},{\xi^3_3},{\xi^4_3})$ are not equivalent unless
$ {\eta^2_1}=  {\xi^2_1},   {\eta^3_3}
   ={\xi^3_3},
  {\eta^4_3}=   {\xi^4_3}.$
\par Commutation relations of $ \mathfrak{m} : $
$  [\tilde{x}_1,\tilde{x}_3]=
( - {\xi^4_3} {\xi^2_1} + 1) \tilde{x}_5  + {\xi^3_3} {\xi^2_1}\tilde{x}_6
 ;\; $
 $  [\tilde{x}_1,\tilde{x}_4]=
 {\xi^3_3} {\xi^2_1} \tilde{x}_5
  +
  \frac{{\xi^4_3}- {\xi^3_3}^2 {\xi^2_1} - {\xi^2_1}}{\xi^4_3}
  \, \tilde{x}_6
  ;$
 $  [\tilde{x}_2,\tilde{x}_3]=
 \frac{\xi^3_3}{\xi^2_1} \,  \tilde{x}_5
 +
 \frac{{\xi^4_3} - {\xi^2_1}}{\xi^2_1}\,
 \tilde{x}_6
 ;\;$
 $  [\tilde{x}_2,\tilde{x}_4]=
 \frac{{\xi^4_3} {\xi^2_1} - {\xi^3_3}^2 - 1}{{\xi^4_3} {\xi^2_1}}
 \, \tilde{x}_5
 -
 \frac{\xi^3_3}{\xi^2_1}\,
 \tilde{x}_6
 .$
\subsection{Case $ {\xi^4_3} = \xi^2_1 , \xi^3_3 = \xi^1_1.$}
\label{M10case21}
In that case one has necessarily $\xi^1_1=0.$
\begin{equation}
\label{M10case21general}
J = \begin{pmatrix}
0&
 - 1/{\xi^2_1}&
0&
0&
0&
0\\
\boxed{\xi^2_1}&
0&
0&
0&
0&
0\\
\boxed{\xi^3_1}&
{\xi^4_1}&
0&
 - 1/{\xi^2_1}&
0&
0\\
\boxed{\xi^4_1}&
 - {\xi^3_1}&
{\xi^2_1}&
0&
0&
0\\
*&
%( - ({\xi^6_2} - {\xi^6_1} {\xi^5_5} {\xi^2_1} + ({\xi^5_5} {\xi^4_1} + {\xi^3_1} {\xi^2_1}) {\xi^6_3} - {\xi^6_5} {\xi^5_3} {\xi^4_1}))/({\xi^6_5} {\xi^2_1})&
*&
%({\xi^6_2} {\xi^5_5} + {\xi^6_1} {\xi^2_1} + ({\xi^5_5} {\xi^3_1} {\xi^2_1} - {\xi^4_1}) {\xi^6_3}
% - {\xi^6_5} {\xi^5_3} {\xi^3_1} {\xi^2_1})/{\xi^6_5}&
\boxed{\xi^5_3}&
(({\xi^5_5}^2 + 1) {\xi^6_3} - {\xi^6_5} {\xi^5_5} {\xi^5_3})/({\xi^6_5} {\xi^2_1})&
\boxed{\xi^5_5}&
 - ({\xi^5_5}^2 + 1)/{\xi^6_5}\\
\boxed{\xi^6_1}&
\boxed{\xi^6_2}&
\boxed{\xi^6_3}&
 - ({\xi^6_5} {\xi^5_3} - {\xi^6_3} {\xi^5_5})/{\xi^2_1}&
\boxed{\xi^6_5}&
 - {\xi^5_5}\end{pmatrix}
 \end{equation}
where
{\fontsize{6}{10}\selectfont
$\mathbf{J^5_1}=
( - ({\xi^6_2} - {\xi^6_1} {\xi^5_5} {\xi^2_1} + ({\xi^5_5} {\xi^4_1} + {\xi^3_1} {\xi^2_1}) {\xi^6_3} - {\xi^6_5} {\xi^5_3} {\xi^4_1}))/({\xi^6_5} {\xi^2_1});\myquad
\mathbf{J^5_2}=
({\xi^6_2} {\xi^5_5} + {\xi^6_1} {\xi^2_1} + ({\xi^5_5} {\xi^3_1} {\xi^2_1} - {\xi^4_1}) {\xi^6_3}
 - {\xi^6_5} {\xi^5_3} {\xi^3_1} {\xi^2_1})/{\xi^6_5};
$}
and the parameters are subject to the condition
\begin{equation}
\label{condM10general21}
\xi^2_1=\pm 1 , \quad
 {\xi^6_5} \neq 0.
\end{equation}
As in the preceding case, by
equivalence by a suitable block-diagonal automorphism,
we  can suppose $\xi^5_5 =0, \; 0< \xi^6_5 \leqslant 1$ and then
equivalence by (\ref{phiM1011})
leads to the case where  moreover
 $ {\xi^3_1}=0, {\xi^4_1}=0.$
 Then equivalence by
$${\fontsize{8}{10}\selectfont
 \Phi = \begin{pmatrix}
1&
0&
0&
0&
0&
0\\
0&
1&
0&
0&
0&
0\\
0&
0&
1&
0&
0&
0\\
0&
0&
0&
1&
0&
0\\
({\xi^2_1} - {\xi^6_1})/{\xi^6_5}&
  {\xi^6_2} {\xi^2_1}/({\xi^6_5} {\xi^2_1})&
0&
{\xi^5_3}/{\xi^2_1}&
1&
0\\
0&
0&
0&
{\xi^6_3}/{\xi^2_1}&
0&
1\end{pmatrix}
}
$$
leads to the case where moreover
$ {\xi^5_3}=0,  {\xi^6_1}=0, {\xi^6_2}=0,  {\xi^6_3}=0 : $
\begin{equation}
\label{M10case21final}
 J({\xi^2_1},{\xi^6_5}) =
 \text{diag} \left(
 \begin{pmatrix} 0& -\frac{1}{\xi^2_1}\\ {\xi^2_1}& 0  \end{pmatrix},
 \begin{pmatrix} 0& -\frac{1}{\xi^2_1}\\ {\xi^2_1}& 0  \end{pmatrix},
 \begin{pmatrix} 0& -\frac{1}{\xi^6_5}\\ {\xi^6_5}& 0  \end{pmatrix} \right)
\end{equation}
and the parameters are subject to the condition
\begin{equation}
\label{condM10final21}
 {\xi^2_1}= \pm 1, \;  0< {\xi^6_5} \leqslant 1.
\end{equation}
As in the preceding case, we can see that
 $J({\eta^2_1},{\eta^6_5})$ and $J({\xi^2_1},{\xi^4_5})$ are not equivalent unless
$ {\eta^2_1}=  {\xi^2_1},     {\eta^6_5}=   {\xi^6_5}.$
Commutation relations of $ \mathfrak{m} : $
$ \mathfrak{m}$  is abelian.
 Since
$ \mathfrak{m}$  is abelian,
no $ J({\xi^2_1},{\xi^6_5})$ is equivalent to any
 $J({\xi^2_1},{\xi^3_3},{\xi^4_3})$ in  (\ref{M10case1final}).
\subsection{Case $ {\xi^4_3} = \xi^2_1 , \xi^3_3 \neq \xi^1_1.$}
\label{M10case22}
In that case one has necessarily $\xi^3_3\neq -\xi^1_1 $ as well.
%{\fontsize{8}{10} \selectfont
\begin{equation}
\label{M10case22general}
J = \begin{pmatrix}
\boxed{\xi^1_1}&
%(({\xi^3_3} - {\xi^1_1}) ({\xi^1_1}^2 + 1))/(({\xi^3_3} + {\xi^1_1}) {\xi^6_5})&
-\frac{{\xi^1_1}^2+1}{k}&
0&
0&
0&
0\\
k&
% - (({\xi^3_3} + {\xi^1_1}) {\xi^6_5})/({\xi^3_3} - {\xi^1_1})&
 - {\xi^1_1}&
0&
0&
0&
0\\
\boxed{\xi^3_1}&
% b&
\frac{(({\xi^3_3}^2 + 1) ({\xi^3_3} - {\xi^1_1}) {\xi^4_1} + ({\xi^3_3} + {\xi^1_1})^2 {\xi^6_5}
{\xi^3_1}) ({\xi^3_3} - {\xi^1_1})}{({\xi^3_3} + {\xi^1_1})^2 {\xi^6_5}^2}&
\boxed{\xi^3_3}&
-\frac{{\xi^3_3}^2+1}{k}&
%(({\xi^3_3}^2 + 1) ({\xi^3_3} - {\xi^1_1}))/(({\xi^3_3} + {\xi^1_1}) {\xi^6_5})&
0&
0\\
\boxed{\xi^4_1}&
 %c&
 - \frac{({\xi^3_3} + {\xi^1_1}) {\xi^6_5} {\xi^3_1} + ({\xi^3_3} - {\xi^1_1})^2
{\xi^4_1}}{({\xi^3_3} + {\xi^1_1}) {\xi^6_5}}&
%( - (({\xi^3_3} + {\xi^1_1}) {\xi^6_5} {\xi^3_1} + ({\xi^3_3} - {\xi^1_1})^2 {\xi^4_1}))/(({\xi^3_3} + {\xi^1_1}) {\xi^6_5})&
k&
 %- (({\xi^3_3} + {\xi^1_1}) {\xi^6_5})/({\xi^3_3} - {\xi^1_1})&
 - {\xi^3_3}&
0&
0\\
%(({\xi^6_5}^2 {\xi^6_1} {\xi^3_3} + {\xi^6_5}^2 {\xi^6_1} {\xi^1_1} + {\xi^6_5} {\xi^6_3} {\xi^4_1} {\xi^3_3} - {\xi^6_5} {\xi^6_3} {\xi^4_1} {\xi^1_1} + {\xi^6_5} {\xi^6_2} {\xi^3_3}^2 - {\xi^6_5} {\xi^6_2} {\xi^1_1}^2 - {\xi^6_3} {\xi^3_3}^2 {\xi^3_1} + 2 {\xi^6_3} {\xi^3_3} {\xi^3_1}
% {\xi^1_1} - {\xi^6_3} {\xi^3_1} {\xi^1_1}^2) ({\xi^3_3} + {\xi^1_1})^2 {\xi^6_5} - ({\xi^3_3}
%^2 + 1) ({\xi^3_3} - {\xi^1_1})^3 {\xi^6_3} {\xi^4_1} - (({\xi^3_3} - {\xi^1_1}) {\xi^5_3}
%{\xi^4_1} + ({\xi^1_1}^2 + 1) {\xi^6_1}) ({\xi^3_3} + {\xi^1_1}) ({\xi^3_3} - {\xi^1_1})^2
%{\xi^6_5})/(({\xi^3_3} + {\xi^1_1})^2 ({\xi^3_3} - {\xi^1_1})^2 {\xi^6_5}^2)&
*&
%(({\xi^6_5}^3 {\xi^6_2} {\xi^3_3}^3 + 3 {\xi^6_5}^3 {\xi^6_2} {\xi^3_3}^2 {\xi^1_1} + 3 {\xi^6_5}^3 {\xi^6_2} {\xi^3_3} {\xi^1_1}^2 + {\xi^6_5}^3 {\xi^6_2} {\xi^1_1}^3 - {\xi^6_5}^2
 %{\xi^6_3} {\xi^3_3}^3 {\xi^3_1} - {\xi^6_5}^2 {\xi^6_3} {\xi^3_3}^2 {\xi^3_1} {\xi^1_1} + {\xi^6_5}^2 {\xi^6_3} {\xi^3_3} {\xi^3_1} {\xi^1_1}^2 + {\xi^6_5}^2 {\xi^6_3} {\xi^3_1} {\xi^1_1}
%^3 - {\xi^6_5} {\xi^6_3} {\xi^4_1} {\xi^3_3}^4 + 2 {\xi^6_5} {\xi^6_3} {\xi^4_1} {\xi^3_3}^3
%{\xi^1_1} - 2 {\xi^6_5} {\xi^6_3} {\xi^4_1} {\xi^3_3} {\xi^1_1}^3 + {\xi^6_5} {\xi^6_3} {\xi^4_1}
 %{\xi^1_1}^4 - 2 {\xi^6_3} {\xi^3_3}^4 {\xi^3_1} {\xi^1_1} + 5 {\xi^6_3} {\xi^3_3}^3 {\xi^3_1} {\xi^1_1}^2 + {\xi^6_3} {\xi^3_3}^3 {\xi^3_1} - 3 {\xi^6_3} {\xi^3_3}^2 {\xi^3_1} {\xi^1_1}
%^3 - 3 {\xi^6_3} {\xi^3_3}^2 {\xi^3_1} {\xi^1_1} - {\xi^6_3} {\xi^3_3} {\xi^3_1} {\xi^1_1}^4
%+ 3 {\xi^6_3} {\xi^3_3} {\xi^3_1} {\xi^1_1}^2 + {\xi^6_3} {\xi^3_1} {\xi^1_1}^5 - {\xi^6_3} {\xi^3_1} {\xi^1_1}^3) ({\xi^3_3} + {\xi^1_1}) {\xi^6_5} - 2 ({\xi^3_3}^2 + 1) ({\xi^3_3} - {\xi^1_1})^4 {\xi^6_3} {\xi^4_1} {\xi^1_1} - (({\xi^3_3} + {\xi^1_1}) ({\xi^1_1}^2 + 1) {\xi^6_1} - ({\xi^3_3} - {\xi^1_1})^2 {\xi^5_3} {\xi^4_1}) ({\xi^3_3} + {\xi^1_1}) ({\xi^3_3} - {\xi^1_1})^3 {\xi^6_5} + ((2 {\xi^3_3} {\xi^1_1} + {\xi^1_1}^2 - 1) {\xi^6_2} + ({\xi^3_3} - {\xi^1_1}) {\xi^5_3} {\xi^3_1}) ({\xi^3_3} + {\xi^1_1})^2 ({\xi^3_3} - {\xi^1_1})^2 {\xi^6_5}^2)/
%(({\xi^3_3} + {\xi^1_1})^3 ({\xi^3_3} - {\xi^1_1})^2 {\xi^6_5}^3)&
*&
\boxed{\xi^5_3}&
 *&
%( - (({\xi^6_5}^3 {\xi^6_3} {\xi^3_3}^3 + 3 {\xi^6_5}^3 {\xi^6_3} {\xi^3_3}^2 {\xi^1_1} +
%3 {\xi^6_5}^3 {\xi^6_3} {\xi^3_3} {\xi^1_1}^2 + {\xi^6_5}^3 {\xi^6_3} {\xi^1_1}^3 - {\xi^6_5}^2 {\xi^5_3} {\xi^3_3}^4 + 2 {\xi^6_5}^2 {\xi^5_3} {\xi^3_3}^2 {\xi^1_1}^2 - {\xi^6_5}^
%2 {\xi^5_3} {\xi^1_1}^4 + 2 {\xi^6_5} {\xi^6_3} {\xi^3_3}^4 {\xi^1_1} - 2 {\xi^6_5} {\xi^6_3}
%{\xi^3_3}^3 {\xi^1_1}^2 - 2 {\xi^6_5} {\xi^6_3} {\xi^3_3}^3 - 2 {\xi^6_5} {\xi^6_3} {\xi^3_3}
%^2 {\xi^1_1}^3 + 2 {\xi^6_5} {\xi^6_3} {\xi^3_3}^2 {\xi^1_1} + 2 {\xi^6_5} {\xi^6_3} {\xi^3_3} {\xi^1_1}^4 + 2 {\xi^6_5} {\xi^6_3} {\xi^3_3} {\xi^1_1}^2 - 2 {\xi^6_5} {\xi^6_3} {\xi^1_1}
%^3 - {\xi^5_3} {\xi^3_3}^6 + 2 {\xi^5_3} {\xi^3_3}^5 {\xi^1_1} + 2 {\xi^5_3} {\xi^3_3}^4
%{\xi^1_1}^2 + {\xi^5_3} {\xi^3_3}^4 - 8 {\xi^5_3} {\xi^3_3}^3 {\xi^1_1}^3 - 4 {\xi^5_3} {\xi^3_3}^3 {\xi^1_1} + 7 {\xi^5_3} {\xi^3_3}^2 {\xi^1_1}^4 + 6 {\xi^5_3} {\xi^3_3}^2 {\xi^1_1}^2 - 2 {\xi^5_3} {\xi^3_3} {\xi^1_1}^5 - 4 {\xi^5_3} {\xi^3_3} {\xi^1_1}^3 + {\xi^5_3} {\xi^1_1}^4) ({\xi^3_3} + {\xi^1_1}) {\xi^6_5} + ({\xi^3_3}^2 + 1) ({\xi^3_3} - {\xi^1_1})^4 (
%{\xi^1_1}^2 + 1) {\xi^6_3}))/(({\xi^3_3} + {\xi^1_1})^3 ({\xi^3_3} - {\xi^1_1})^3 {\xi^6_5}
%^2)&
m&
%(({\xi^3_3} {\xi^1_1} - 1) ({\xi^3_3} - {\xi^1_1})^2 + ({\xi^3_3} + {\xi^1_1})^2 {\xi^6_5}%^
%2)/(({\xi^3_3} + {\xi^1_1}) ({\xi^3_3} - {\xi^1_1})^2)&
-\frac{m^2+1}{\xi^6_5}&
%( - (({\xi^6_5}^2 {\xi^3_3}^2 + 2 {\xi^6_5}^2 {\xi^3_3} {\xi^1_1} + {\xi^6_5}^2 {\xi^1_1}
%^2 + 2 {\xi^3_3}^3 {\xi^1_1} - 4 {\xi^3_3}^2 {\xi^1_1}^2 - 2 {\xi^3_3}^2 + 2 {\xi^3_3}
%{\xi^1_1}^3 + 4 {\xi^3_3} {\xi^1_1} - 2 {\xi^1_1}^2) ({\xi^3_3} + {\xi^1_1})^2 {\xi^6_5}^2
% + ({\xi^3_3}^2 + 1) ({\xi^3_3} - {\xi^1_1})^4 ({\xi^1_1}^2 + 1)))/(({\xi^3_3} + {\xi^1_1}%)^2 ({\xi^3_3} - {\xi^1_1})^4 {\xi^6_5})
\\
\boxed{\xi^6_1}&
\boxed{\xi^6_2}&
\boxed{\xi^6_3}&
 *&
%( - (({\xi^6_5} {\xi^6_3} {\xi^3_3} + {\xi^6_5} {\xi^6_3} {\xi^1_1} - {\xi^5_3} {\xi^3_3}^2 + 2
 %{\xi^5_3} {\xi^3_3} {\xi^1_1} - {\xi^5_3} {\xi^1_1}^2) ({\xi^3_3} + {\xi^1_1}) {\xi^6_5} - ({\xi^3_3}^2 + 1) ({\xi^3_3} - {\xi^1_1})^2 {\xi^6_3}))/(({\xi^3_3} + {\xi^1_1})^2 ({\xi^3_3}
%- {\xi^1_1}) {\xi^6_5})&
\boxed{\xi^6_5}&
 -m&
%( - (({\xi^3_3} {\xi^1_1} - 1) ({\xi^3_3} - {\xi^1_1})^2 + ({\xi^3_3} + {\xi^1_1})^2 {\xi^6_5}^2))/(({\xi^3_3} + {\xi^1_1}) ({\xi^3_3} - {\xi^1_1})^2)
\end{pmatrix}\end{equation}
%}
where
{\fontsize{8}{10} \selectfont
$\mathbf{k}=\mathbf{J^2_1}= - (({\xi^3_3} + {\xi^1_1}) {\xi^6_5})/({\xi^3_3} - {\xi^1_1});$\myquad
$ \mathbf{{J^5_1}}=(({\xi^6_5}^3 {\xi^6_1} {\xi^3_3}^2 + 2 {\xi^6_5}^3 {\xi^6_1} {\xi^3_3} {\xi^1_1} + {\xi^6_5}^3 {\xi^6_1} {\xi^1_1}^2 + {\xi^6_5}^2 {\xi^6_3} {\xi^4_1} {\xi^3_3}^2 -
{\xi^6_5}^2 {\xi^6_3} {\xi^4_1} {\xi^1_1}^2 + {\xi^6_5}^2 {\xi^6_2} {\xi^3_3}^3 + {\xi^6_5}
^2 {\xi^6_2} {\xi^3_3}^2 {\xi^1_1} - {\xi^6_5}^2 {\xi^6_2} {\xi^3_3} {\xi^1_1}^2 - {\xi^6_5}
^2 {\xi^6_2} {\xi^1_1}^3 - {\xi^6_5} {\xi^6_3} {\xi^3_3}^3 {\xi^3_1} + {\xi^6_5} {\xi^6_3} {\xi^3_3}^2 {\xi^3_1} {\xi^1_1} + {\xi^6_5} {\xi^6_3} {\xi^3_3} {\xi^3_1} {\xi^1_1}^2 - {\xi^6_5}
{\xi^6_3} {\xi^3_1} {\xi^1_1}^3 - {\xi^6_5} {\xi^6_1} {\xi^3_3}^3 {\xi^1_1} + {\xi^6_5} {\xi^6_1} {\xi^3_3}^2 {\xi^1_1}^2 + {\xi^6_5} {\xi^6_1} {\xi^3_3} {\xi^1_1}^3 - {\xi^6_5} {\xi^6_1}
{\xi^1_1}^4 - {\xi^6_5} {\xi^5_3} {\xi^4_1} {\xi^3_3}^3 + 3 {\xi^6_5} {\xi^5_3} {\xi^4_1} {\xi^3_3}^2 {\xi^1_1} - 3 {\xi^6_5} {\xi^5_3} {\xi^4_1} {\xi^3_3} {\xi^1_1}^2 + {\xi^6_5} {\xi^5_3}
{\xi^4_1} {\xi^1_1}^3 - {\xi^6_3} {\xi^4_1} {\xi^3_3}^4 + 3 {\xi^6_3} {\xi^4_1} {\xi^3_3}^3
{\xi^1_1} - 3 {\xi^6_3} {\xi^4_1} {\xi^3_3}^2 {\xi^1_1}^2 + {\xi^6_3} {\xi^4_1} {\xi^3_3} {\xi^1_1}^3) ({\xi^3_3} + {\xi^1_1}) + ({\xi^6_5} {\xi^6_1} {\xi^3_3} + {\xi^6_5} {\xi^6_1} {\xi^1_1}
 + {\xi^6_3} {\xi^4_1} {\xi^3_3} - {\xi^6_3} {\xi^4_1} {\xi^1_1}) ({\xi^3_3} {\xi^1_1} - 1) ({\xi^3_3} - {\xi^1_1})^2)/(({\xi^3_3} + {\xi^1_1})^2 ({\xi^3_3} - {\xi^1_1})^2 {\xi^6_5}^2)
;$
\myquad
$ \mathbf{{J^5_2}}=({\xi^6_5}^4 {\xi^6_2} {\xi^3_3}^4 + 4 {\xi^6_5}^4 {\xi^6_2} {\xi^3_3}^3
{\xi^1_1} + 6 {\xi^6_5}^4 {\xi^6_2} {\xi^3_3}^2 {\xi^1_1}^2 + 4 {\xi^6_5}^4 {\xi^6_2} {\xi^3_3} {\xi^1_1}^3 + {\xi^6_5}^4 {\xi^6_2} {\xi^1_1}^4 - {\xi^6_5}^3 {\xi^6_3} {\xi^3_3}^4
{\xi^3_1} - 2 {\xi^6_5}^3 {\xi^6_3} {\xi^3_3}^3 {\xi^3_1} {\xi^1_1} + 2 {\xi^6_5}^3 {\xi^6_3}
 {\xi^3_3} {\xi^3_1} {\xi^1_1}^3 + {\xi^6_5}^3 {\xi^6_3} {\xi^3_1} {\xi^1_1}^4 - {\xi^6_5}^2
 {\xi^6_3} {\xi^4_1} {\xi^3_3}^5 + {\xi^6_5}^2 {\xi^6_3} {\xi^4_1} {\xi^3_3}^4 {\xi^1_1} + 2
{\xi^6_5}^2 {\xi^6_3} {\xi^4_1} {\xi^3_3}^3 {\xi^1_1}^2 - 2 {\xi^6_5}^2 {\xi^6_3} {\xi^4_1}
{\xi^3_3}^2 {\xi^1_1}^3 - {\xi^6_5}^2 {\xi^6_3} {\xi^4_1} {\xi^3_3} {\xi^1_1}^4 + {\xi^6_5}
^2 {\xi^6_3} {\xi^4_1} {\xi^1_1}^5 + 2 {\xi^6_5}^2 {\xi^6_2} {\xi^3_3}^5 {\xi^1_1} + {\xi^6_5}^2 {\xi^6_2} {\xi^3_3}^4 {\xi^1_1}^2 - {\xi^6_5}^2 {\xi^6_2} {\xi^3_3}^4 - 4 {\xi^6_5}
^2 {\xi^6_2} {\xi^3_3}^3 {\xi^1_1}^3 - 2 {\xi^6_5}^2 {\xi^6_2} {\xi^3_3}^2 {\xi^1_1}^4 +
 2 {\xi^6_5}^2 {\xi^6_2} {\xi^3_3}^2 {\xi^1_1}^2 + 2 {\xi^6_5}^2 {\xi^6_2} {\xi^3_3} {\xi^1_1}^5 + {\xi^6_5}^2 {\xi^6_2} {\xi^1_1}^6 - {\xi^6_5}^2 {\xi^6_2} {\xi^1_1}^4 + {\xi^6_5}
^2 {\xi^5_3} {\xi^3_3}^5 {\xi^3_1} - {\xi^6_5}^2 {\xi^5_3} {\xi^3_3}^4 {\xi^3_1} {\xi^1_1} -
 2 {\xi^6_5}^2 {\xi^5_3} {\xi^3_3}^3 {\xi^3_1} {\xi^1_1}^2 + 2 {\xi^6_5}^2 {\xi^5_3} {\xi^3_3}^2 {\xi^3_1} {\xi^1_1}^3 + {\xi^6_5}^2 {\xi^5_3} {\xi^3_3} {\xi^3_1} {\xi^1_1}^4 - {\xi^6_5}^2 {\xi^5_3} {\xi^3_1} {\xi^1_1}^5 - 2 {\xi^6_5} {\xi^6_3} {\xi^3_3}^5 {\xi^3_1} {\xi^1_1}
+ 3 {\xi^6_5} {\xi^6_3} {\xi^3_3}^4 {\xi^3_1} {\xi^1_1}^2 + {\xi^6_5} {\xi^6_3} {\xi^3_3}^4
{\xi^3_1} + 2 {\xi^6_5} {\xi^6_3} {\xi^3_3}^3 {\xi^3_1} {\xi^1_1}^3 - 2 {\xi^6_5} {\xi^6_3} {\xi^3_3}^3 {\xi^3_1} {\xi^1_1} - 4 {\xi^6_5} {\xi^6_3} {\xi^3_3}^2 {\xi^3_1} {\xi^1_1}^4 + 2
{\xi^6_5} {\xi^6_3} {\xi^3_3} {\xi^3_1} {\xi^1_1}^3 + {\xi^6_5} {\xi^6_3} {\xi^3_1} {\xi^1_1}^6
- {\xi^6_5} {\xi^6_3} {\xi^3_1} {\xi^1_1}^4 - {\xi^6_5} {\xi^6_1} {\xi^3_3}^5 {\xi^1_1}^2 -
{\xi^6_5} {\xi^6_1} {\xi^3_3}^5 + {\xi^6_5} {\xi^6_1} {\xi^3_3}^4 {\xi^1_1}^3 + {\xi^6_5} {\xi^6_1} {\xi^3_3}^4 {\xi^1_1} + 2 {\xi^6_5} {\xi^6_1} {\xi^3_3}^3 {\xi^1_1}^4 + 2 {\xi^6_5} {\xi^6_1} {\xi^3_3}^3 {\xi^1_1}^2 - 2 {\xi^6_5} {\xi^6_1} {\xi^3_3}^2 {\xi^1_1}^5 - 2 {\xi^6_5} {\xi^6_1} {\xi^3_3}^2 {\xi^1_1}^3 - {\xi^6_5} {\xi^6_1} {\xi^3_3} {\xi^1_1}^6 - {\xi^6_5}
{\xi^6_1} {\xi^3_3} {\xi^1_1}^4 + {\xi^6_5} {\xi^6_1} {\xi^1_1}^7 + {\xi^6_5} {\xi^6_1} {\xi^1_1}^5 + {\xi^6_5} {\xi^5_3} {\xi^4_1} {\xi^3_3}^6 - 4 {\xi^6_5} {\xi^5_3} {\xi^4_1} {\xi^3_3}^5
 {\xi^1_1} + 5 {\xi^6_5} {\xi^5_3} {\xi^4_1} {\xi^3_3}^4 {\xi^1_1}^2 - 5 {\xi^6_5} {\xi^5_3}
{\xi^4_1} {\xi^3_3}^2 {\xi^1_1}^4 + 4 {\xi^6_5} {\xi^5_3} {\xi^4_1} {\xi^3_3} {\xi^1_1}^5 -
{\xi^6_5} {\xi^5_3} {\xi^4_1} {\xi^1_1}^6 - 2 {\xi^6_3} {\xi^4_1} {\xi^3_3}^6 {\xi^1_1} + 8 {\xi^6_3} {\xi^4_1} {\xi^3_3}^5 {\xi^1_1}^2 - 12 {\xi^6_3} {\xi^4_1} {\xi^3_3}^4 {\xi^1_1}^3 -
 2 {\xi^6_3} {\xi^4_1} {\xi^3_3}^4 {\xi^1_1} + 8 {\xi^6_3} {\xi^4_1} {\xi^3_3}^3 {\xi^1_1}^4
+ 8 {\xi^6_3} {\xi^4_1} {\xi^3_3}^3 {\xi^1_1}^2 - 2 {\xi^6_3} {\xi^4_1} {\xi^3_3}^2 {\xi^1_1}
^5 - 12 {\xi^6_3} {\xi^4_1} {\xi^3_3}^2 {\xi^1_1}^3 + 8 {\xi^6_3} {\xi^4_1} {\xi^3_3} {\xi^1_1}^4 - 2 {\xi^6_3} {\xi^4_1} {\xi^1_1}^5)/(({\xi^3_3} + {\xi^1_1})^3 ({\xi^3_3} - {\xi^1_1}
)^2 {\xi^6_5}^3);$
\myquad
$ \mathbf{{J^5_4}}=((({\xi^6_5}^3 {\xi^5_3} {\xi^3_3} + {\xi^6_5}^3 {\xi^5_3} {\xi^1_1} + {\xi^6_5} {\xi^5_3} {\xi^3_3}^3 - 2 {\xi^6_5} {\xi^5_3} {\xi^3_3}^2 {\xi^1_1} + {\xi^6_5} {\xi^5_3}
{\xi^3_3} {\xi^1_1}^2 - {\xi^6_3} {\xi^3_3}^2 + 2 {\xi^6_3} {\xi^3_3} {\xi^1_1} - {\xi^6_3} {\xi^1_1}^2) ({\xi^3_3} + {\xi^1_1}) + ({\xi^3_3} {\xi^1_1} - 1) ({\xi^3_3} - {\xi^1_1})^2 {\xi^6_5} {\xi^5_3}) ({\xi^3_3} + {\xi^1_1}) ({\xi^3_3} - {\xi^1_1})^2 - (({\xi^3_3} {\xi^1_1} - 1
) ({\xi^3_3} - {\xi^1_1})^2 + ({\xi^3_3} + {\xi^1_1})^2 {\xi^6_5}^2)^2 {\xi^6_3})/(({\xi^3_3} + {\xi^1_1})^3 ({\xi^3_3} - {\xi^1_1})^3 {\xi^6_5}^2);$
\myquad
$ \mathbf{m}=\mathbf{{J^5_5}}=(({\xi^3_3} {\xi^1_1} - 1) ({\xi^3_3} - {\xi^1_1})^2 + ({\xi^3_3} + {\xi^1_1}
)^2 {\xi^6_5}^2)/(({\xi^3_3} + {\xi^1_1}) ({\xi^3_3} - {\xi^1_1})^2);$
\myquad
$ \mathbf{{J^6_4}}=( - (({\xi^6_5}^2 {\xi^6_3} {\xi^3_3} + {\xi^6_5}^2 {\xi^6_3} {\xi^1_1} - {\xi^6_5} {\xi^5_3} {\xi^3_3}^2 + 2 {\xi^6_5} {\xi^5_3} {\xi^3_3} {\xi^1_1} - {\xi^6_5} {\xi^5_3}
{\xi^1_1}^2 - {\xi^6_3} {\xi^3_3}^3 + 2 {\xi^6_3} {\xi^3_3}^2 {\xi^1_1} - {\xi^6_3} {\xi^3_3}
 {\xi^1_1}^2) ({\xi^3_3} + {\xi^1_1}) + ({\xi^3_3} {\xi^1_1} - 1) ({\xi^3_3} - {\xi^1_1})^2
{\xi^6_3}))/(({\xi^3_3} + {\xi^1_1})^2 ({\xi^3_3} - {\xi^1_1}) {\xi^6_5});$
}
and the parameters are subject to the condition
\begin{equation}
\label{condM10general22}
{\xi^3_3} \neq \pm {\xi^1_1}, \; {\xi^2_1} {\xi^3_4} {\xi^6_5} \neq 0 .
\end{equation}
As in the preceding cases,
equivalence by a suitable block-diagonal automorphism
leads to the case $\xi^1_1=0.$
Equivalences by successively
$${\fontsize{8}{10}\selectfont
  \begin{pmatrix}
1&
0&
0&
0&
0&
0\\
0&
1&
0&
0&
0&
0\\
0&
 - {\xi^3_1}/{\xi^6_5}&
1&
0&
0&
0\\
0&
 - {\xi^4_1}/{\xi^6_5}&
0&
1&
0&
0\\
0&
0&
 - {\xi^3_1}/{\xi^6_5}&
0&
1&
0\\
0&
0&
 - {\xi^4_1}/{\xi^6_5}&
0&
0&
1\end{pmatrix}
%$$
\text{ and }
  \begin{pmatrix}
1&
0&
0&
0&
0&
0\\
0&
1&
0&
0&
0&
0\\
0&
0&
1&
0&
0&
0\\
0&
0&
0&
1&
0&
0\\
 - {\xi^6_1} {\xi^3_3}/({\xi^6_5} {\xi^3_3})&
- {\xi^6_5} {\xi^6_2} {\xi^3_3}/
({\xi^6_5}^2 {\xi^3_3})&
0&
 - {\xi^5_3}/{\xi^6_5}&
1&
0\\
0&
0&
0&
 - {\xi^6_3}/{\xi^6_5}&
0&
1\end{pmatrix}
}
$$
lead to the case where  moreover
 $ {\xi^3_1}=0, {\xi^4_1}=0$
and $ {\xi^5_3}=0, {\xi^6_1}=0, {\xi^6_2}=0,{\xi^6_3}=0 :$
\begin{equation}
{\fontsize{8}{10} \selectfont
J = \begin{pmatrix}
0&
1/{\xi^6_5}&
0&
0&
0&
0\\
 - {\xi^6_5}&
0&
0&
0&
0&
0\\
0&
0&
{\xi^3_3}&
({\xi^3_3}^2 + 1)/{\xi^6_5}&
0&
0\\
0&
0&
 - {\xi^6_5}&
 - {\xi^3_3}&
0&
0\\
0&
0&
0&
0&
({\xi^6_5} + 1) ({\xi^6_5} - 1)/{\xi^3_3}&
- (({\xi^6_5}^2 - 2) {\xi^6_5}^2 + {\xi^3_3}^2 + 1)/({\xi^6_5} {\xi^3_3}^2)\\
0&
0&
0&
0&
{\xi^6_5}&
- ({\xi^6_5} + 1) ({\xi^6_5} - 1)/{\xi^3_3}\end{pmatrix}
}
\end{equation}
Computing the matrix
 $J2=\Phi^{-1} J \Phi$
where $\Phi$ is given in (\ref{automM10}) and
$ {b^2_1}={b^1_1}=1,$ one gets
$ {J2^2_1}=\frac{1+{\xi^6_5}^2}{2 {\xi^6_5}} u,$
$ {J2^4_3}={\xi^6_5} u.$
Hence if $ {\xi^6_5}^2 \neq 1$ we are back to case 1 (\ref{M10case1}).
Finally, if $ {\xi^6_5}^2 =1,$ equivalence by
$ \Phi =  \text{diag}(1,-1,-1,1,-1,1)$
%$$ \Phi = \begin{pmatrix}1&0&0&0&0&0\\0&-1&0&0&0&0\\0&0&-1&0&0&0\\0& 0&0&1&0&0\\0&0&0&0&-1&
%0\\0&0&0&0&0&1\end{pmatrix}$$
leads to the case where  moreover
$ {\xi^6_5}=1 : $
\begin{equation}
\label{M10case22final}
 J({\xi^3_3}) =
 \text{diag} \left(
 \begin{pmatrix} 0& 1\\ -1& 0  \end{pmatrix},
 \begin{pmatrix} \xi^3_3& {\xi^3_3}^2+1\\ -1& -\xi^3_3  \end{pmatrix},
 \begin{pmatrix} 0& -1\\ 1& 0  \end{pmatrix} \right)
 \quad \quad ({\xi^3_3} \neq 0) .
\end{equation}
 $J({\xi^3_3})$  and $J({\eta^3_3})$ are non equivalent unless $\eta^3_3=\xi^3_3,$
 and
 $J({\xi^3_3})$ in (\ref{M10case22final}) is equivalent neither to
any  $J({\xi^2_1},{\xi^3_3},{\xi^4_3})$ in  (\ref{M10case1final}) nor
to any  $ J({\xi^2_1},{\xi^6_5})$ in (\ref{M10case21final}).
\par Commutation relations of $ \mathfrak{m} : $
 $  [\tilde{x}_1,\tilde{x}_3]= -  {\xi^3_3}  \tilde{x}_6  ;$
 $  [\tilde{x}_1,\tilde{x}_4]= -  {\xi^3_3}  \tilde{x}_5   -    {\xi^3_3}^2   \tilde{x}_6
;$
 $  [\tilde{x}_2,\tilde{x}_3]= - {\xi^3_3}  \tilde{x}_5  ;$
 $  [\tilde{x}_2,\tilde{x}_4]= -  {\xi^3_3}^2   \tilde{x}_5  + {\xi^3_3} \tilde{x}_6.$
\subsection{Conclusions.}
To solve the initial system comprised of all the torsion equations and the equation $J^2=-1$ in $\Rmath^{36}$ in the 3 cases 1
(\ref{M10case1}),
2.1 (\ref{M10case21}), 2.2 (\ref{M10case22}),
one has to complete first a set of common steps,  and then we are left with solving the
system $S$ of the remaining equations in the 12 variables
$\xi^1_1,\xi^2_1,\xi^3_1,\xi^3_3,\xi^4_1,\xi^4_3,\xi^5_3,\xi^5_5,
\xi^6_1,\xi^6_2,\xi^6_3,\xi^6_5$ in the open subset
$\xi^6_5 \xi^4_3 \xi^2_1 \neq 0$ of $\Rmath^{12}.$
Among these equations, we single out the 2 equations $13|6$ and $14|6$ which read:
 \begin{equation}
 \label{M10,f,g}
 \begin{cases}
 f=f_{13|6} =0\\ g=f_{14|6}=0
 \end{cases}
 \end{equation}
 where :
$ f_{13|6}=\xi^6_5 (\xi^3_3 +\xi^1_1) +\xi^5_5 (\xi^2_1 -\xi^4_3) -\xi^4_3 \xi^1_1 +\xi^3_3 \xi^2_1 $ and
$ f_{14|6}=(\xi^6_5 (\xi^4_3 \xi^2_1 -{\xi^3_3}^2 -1)+\xi^5_5 \xi^4_3(\xi^3_3 - \xi^1_1)
+\xi^4_3 \xi^3_3 \xi^1_1  + \xi^4_3 - ({\xi^3_3}^2+1) \xi^2_1 )/\xi^4_3.$
In each of the 3 cases, the remaining system is equivalent
to the system
 (\ref{M10,f,g}).
To conclude that
$\mathfrak{X}_{M10}$ is a 10-dimensional submanifold of $\Rmath^{36},$
it will be sufficient to prove that the
preceding system is of maximal rank 2  at any point of
$\mathfrak{X}_{M10},$
that is in each of the 3 cases some 2-jacobian doesn't vanish.
%\\ $\bullet$
In case 1, one has
$\frac{D(f,g)}{D(\xi^5_5,\xi^6_5)} = -\frac{1}{\xi^4_3} \, (
((\xi^4_3 -\xi^2_1)\xi^2_1 - ({\xi^1_1}^2+1))\xi^4_3 +({\xi^3_3}^2+1)\xi^2_1) \neq 0.$
%\\ $\bullet$
In case 2.1, one has
$\frac{D(f,g)}{D(\xi^1_1,\xi^2_1)} =
({\xi^6_5} - \xi^2_1)^2 +{\xi^5_5}^2 \neq 0$
 if $\xi^5_5 \neq 0.$
If $\xi^5_5 =0,$
$\frac{D(f,g)}{D(\xi^1_1,\xi^2_1)} =  ({\xi^6_5} -\xi^2_1)^2 $ and
$\frac{D(f,g)}{D(\xi^3_3,\xi^4_3)} =  ({\xi^6_5} +\xi^2_1)^2 .$
These 2-jacobians cannot simultaneously vanish.
%\\ $\bullet$
In case 2.2, one has
$\frac{D(f,g)}{D(\xi^2_1,\xi^5_5)} =
(\xi^3_3-\xi^1_1)(\xi^5_5+\xi^3_3) \neq 0$
 if $\xi^5_5 \neq -\xi^3_3.$
 If $\xi^5_5 = -\xi^3_3,$
$\frac{D(f,g)}{D(\xi^5_5,\xi^6_5)} =  {\xi^1_1}^2 -{\xi^3_3}^2 \neq 0.$
%\\
Hence the  system  (\ref{M10,f,g}) is of maximal rank 2  at any point of
$\mathfrak{X}_{M10},$
and $\mathfrak{X}_{M10}$ is a 10-dimensional submanifold of $\Rmath^{36}.$
\par
Any  CS is equivalent to one and only one of the following :
$J(\xi^2_1,\xi^3_3,\xi^4_3)$ in (\ref{M10case1final}) or
$J(\xi^2_1,\xi^6_5)$ in (\ref{M10case21final}) or
$J(\xi^3_3)$ in (\ref{M10case22final}).
%%%%%%%%%%%%%%%%%%%%%%%%%%%%%%%%%%%%%%%%%%%%%%%%%%%%%%%%%%%%%%%%%%%%%%%%%%%%%%%%%%%%%%%%%%
\subsection{}
%\subsubsection{}
\begin{equation*}
X_1 = \frac{\partial}{\partial x^1} - y^1 \, \frac{\partial}{\partial x^2}
 -x^2\, \frac{\partial}{\partial x^3}
-\left(y^2+ \frac{(y^1)^2}{2} \right)\, \frac{\partial}{\partial y^3}
\quad , \quad
X_2 = \frac{\partial}{\partial y^1} - y^2 \, \frac{\partial}{\partial x^3}
+x^2 \, \frac{\partial}{\partial y^3}.
\end{equation*}
%\subsubsection{Holomorphic functions.}
Let $G$ denote the group $G_0$ endowed with the left invariant  structure
of complex manifold  defined by
 $J({\xi^2_1},{\xi^3_3},{\xi^4_3})$ in
(\ref{M10case1final}) with conditions (\ref{condM10final1}).
%where
%$ 0 < {\xi^2_1} \leqslant 1,$
%  $ {\xi^2_1} {\xi^4_3} \neq 0,$
% $ {\xi^2_1}\neq {\xi^4_3} ,$
%and $ D=({\xi^4_3} {\xi^2_1} - {\xi^2_1}^2 - 1) {\xi^4_3} + ({\xi^3_3}^2 + 1) {\xi^2_1}
%\neq 0.$
Then $ H_{\Cmath}(G) =\{f \in C^{\infty}(G_0) \; ; \; \tilde{X}_j^{-} \, f = 0
\; \forall j=1,3,5\}.$
One has
\begin{eqnarray*}
\tilde{X}_1^{-}  &=&
2\, \frac{\partial}{\partial \overline{w^1}}
-y^1\, \frac{\partial}{\partial x^2}
-\left( x^2 +i\xi^2_1 y^2 \right) \, \frac{\partial}{\partial x^3}
-\left( \frac{(y^1)^2}{2} +y^2  -i\xi^2_1 x^2 \right) \,\frac{\partial}{\partial y^3}
\quad ,
\\
\tilde{X}_3^{-}  &=& 2  \frac{\partial}{\partial \overline{w^2}}
\quad , \quad
\tilde{X}_5^{-}  =2\, \frac{\partial}{\partial \overline{w^3}}
\quad ,
\end{eqnarray*}
where
\begin{equation*}
w^1 = x^1+ \frac{i}{\xi^2_1} y^1
\quad , \quad
w^2 = x^2 -\frac{\xi^3_3}{\xi^4_3} \, y^2  + \frac{i}{\xi^4_3}\, y^2
\quad , \quad
w^3 = x^3 - \frac{r}{b} \, y^3 +\frac{i}{b} \,   y^3 .
\end{equation*}
Then $f \in C^{\infty}(G_0)$ is in $ H_{\Cmath}(G)$ if and only if it is holomorphic
with respect to  $w^2$ and $w^3$ and satisfies the equation
\begin{multline*}
2\, \frac{\partial f}{\partial \overline{w^1}}
+i \xi^2_1 \frac{w^1-\overline{w^1}}{2i}\,  \frac{\partial f}{\partial {w^2}}
+ \frac{1}{8b} \, \left[ {\xi^2_1}^2 (r-i)
\left( 2\overline{w^1} w^1 -(\overline{w^1})^2 -{w^1}^2 \right)
\right.\\ \left.
-4  \overline{w^2} \left( i(c+\xi^2_1 r) +\xi^2_1 +b\right)
-4  {w^2} \left( i(-c+\xi^2_1 r) +\xi^2_1 +b\right)
\right] \,
\frac{\partial f}{\partial {w^3}}
=0.
\end{multline*}
where $c= i\xi^4_3\xi^2_1b -\xi^4_3 r +i\xi^4_3 +i\xi^3_3 \xi^2_1 r +\xi^3_3\xi^2_1 + b\xi^3_3 .$
The 3 functions
\begin{equation*}
\varphi^1 =w^1
\quad , \quad
\varphi^2 = w^2 -  \frac{i\xi^2_1}{4}
\left(w^1\overline{w^1} -\frac{(\overline{w^1})^2}{2}\right)
\quad , \quad
\end{equation*}
\begin{multline*}
\varphi^3 =w^3 +   \frac{1}{48b} \, \left(
(\overline{w^1})^3 \xi^2_1 (c +ib)
+3(\overline{w^1})^2 w^1 {\xi^2_1}^2 (i-r)
+12 \overline{w^1}\overline{w^2}(ic+i\xi^2_1 r+\xi^2_1+b)
\right.\\ \left.
+3\overline{w^1} {w^1}^2 \xi^2_1 (-c+2\xi^2_1 r-2i\xi^2_1-ib)
+12(\overline{w^1}+w^1)w^2(-ic+i\xi^2_1 r+\xi^2_1 + b)
\right)
 \end{multline*}
are holomorphic.
Let $F : G \rightarrow \Cmath^3$ defined
by
$F=(\varphi^1,\varphi^2,\varphi^3).$
$F$ is a global chart on $G.$
We determine now how the multiplication of $G$ looks like in that chart.
Let $a,x \in G $ with respective second kind canonical coordinates
$(x^1,y^1,x^2,y^2,x^3,y^3), (\alpha^1, \beta^1, \alpha^2, \beta^2, \alpha^3, \beta^3)$
as in  (\ref{x6general}).
With obvious notations,
%$a =[w^1_a, w^2_a,w^3_a],$
%$x =[w^1_x, w^2_x,w^3_x],$
%$a \, x =[w^1_{a x}, w^2_{a x},w^3_{a x}],$
%$a =[\varphi^1_a, \varphi^2_a,\varphi^3_a],$
%$x =[\varphi^1_x, \varphi^2_x,\varphi^3_x],$
%$a \, x =[\varphi^1_{a x}, \varphi^2_{a x},\varphi^3_{a x}].$
computations yield:
\begin{eqnarray*}
w^1_{a x} &=& w^1_a + w^1_x\\
w^2_{a x} &=& w^2_a + w^2_x- b^1x^1 \\
w^3_{a x} &=& w^3_a + w^3_x
-a^2 x^1 +\frac{1}{2} b^1(x^1)^2 -b^2 y^1 +\frac{i-r}{b}
\left( -b^2x^1 -\frac{1}{2} \, (b^1)^2 x^1  +y^1(a^2-b^1x^1)\right).
\end{eqnarray*}
We then get
\begin{equation*}
\varphi^1_{a x} = \varphi^1_a + \varphi^1_x
\quad , \quad
\varphi^2_{a x} = \varphi^2_a + \varphi^2_x -
    \frac{i \xi^2_1}{4} (2 \overline{\varphi^1_a} - \varphi^1_a) \, \varphi^1_x
\quad , \quad
\varphi^3_{a x} = \varphi^3_a +\varphi^3_x  +\chi^3(a,x)
\end{equation*}
where
\begin{multline*}
\chi^3(a,x)=
\frac{\xi^2_1}{16((\xi^4_3 -\xi^2_1)^2+ {\xi^3_3}^2{\xi^2_1}^2)} \,
D_1 \, (\varphi^1_x)^2
\\
+\frac{1}{16(\xi^4_3 -i{\xi^3_3}{\xi^2_1}-\xi^2_1)} \,
\left(  D_2 \xi^2_1 \,  \varphi^1_x
+8\xi^4_3(1-{\xi^2_1}^2) \,
\left(\overline{\varphi^1_a}  +{\varphi^1_a} \right) \, \varphi^2_x
\right)
\end{multline*}
with
\\
{\fontsize{8}{10} \selectfont
$D_1=
i\overline{\varphi^1_a} \, (7{\xi^4_3}^2{\xi^2_1}^2
-3{\xi^4_3}^2
+7i{\xi^4_3}{\xi^3_3}{\xi^2_1}^3
-7i{\xi^4_3}{\xi^3_3}{\xi^2_1}
-7{\xi^4_3}{\xi^2_1}^3
-{\xi^4_3}{\xi^2_1}
+4{\xi^3_3}^2{\xi^2_1}^2
+4{\xi^2_1}^2)
+
i{\varphi^1_a} \, (-4{\xi^4_3}^2{\xi^2_1}^2
+{\xi^4_3}^2
-4i{\xi^4_3}{\xi^3_3}{\xi^2_1}^3
+4i{\xi^4_3}{\xi^3_3}{\xi^2_1}
+4{\xi^4_3}{\xi^2_1}^3
+2{\xi^4_3}{\xi^2_1}
-3{\xi^3_3}^2{\xi^2_1}^2
-3{\xi^2_1}^2)
;\\
D_2=
-i(\overline{\varphi^1_a})^2
\, ({\xi^4_3}^2{\xi^2_1}
-2{\xi^4_3}{\xi^2_1}^2
+2{\xi^4_3}
+{\xi^3_3}^2{\xi^2_1}
-2i{\xi^3_3}{\xi^2_1}
-{\xi^2_1}  )
 +4i\overline{\varphi^1_a} \varphi^1_a
\, ({\xi^4_3}^2{\xi^2_1}
+{\xi^4_3}{\xi^2_1}^2
-{\xi^4_3}
+{\xi^3_3}^2{\xi^2_1}
-2i{\xi^3_3}{\xi^2_1}
-{\xi^2_1}  )
-8\overline{\varphi^2_a}
\, ({\xi^4_3}^2
+{\xi^3_3}^2
-2i{\xi^3_3}  -1  )
+8{\varphi^2_a}
\, ({\xi^4_3}^2
+{\xi^3_3}^2
-2{\xi^3_3} \xi^2_1 +1  )
-i{\varphi^1_a}^2
\, ({\xi^4_3}^2{\xi^2_1}
+3{\xi^4_3}{\xi^2_1}^2
-2{\xi^4_3}
+{\xi^3_3}^2{\xi^2_1}
-3i{\xi^3_3}{\xi^2_1}
-2{\xi^2_1}  ).
$
}
\par In the case of
 $J({\xi^2_1},{\xi^6_5})$ in (\ref{M10case21final})
 where $ {\xi^2_1}= \pm 1,$ and $ 0< {\xi^6_5} \leqslant 1,$  the preceding computations apply with
 $r=0, b=\xi^6_5, \xi^3_3=0, \xi^4_3=\xi^2_1.$
 The only difference is that we get now
\begin{equation*}
\chi^3(a,x)=
\frac{1}{16\xi^6_5} \,
\left(
D_1 \, (\varphi^1_x)^2
+
D_2 \,  \varphi^1_x  \right)
+\frac{\xi^6_5+\xi^2_1}{2\xi^6_5} \,
\left(\overline{\varphi^1_a}  +{\varphi^1_a} \right) \, \varphi^2_x
\end{equation*}
with
%{\fontsize{8}{10} \selectfont
$D_1=
-i\overline{\varphi^1_a} \,
(3{\xi^6_5}{\xi^2_1}+7)
+
i{\varphi^1_a} \,
({\xi^6_5}{\xi^2_1}+4)
;\,
D_2=
-i(\overline{\varphi^1_a})^2
\,
(3{\xi^6_5}{\xi^2_1}+1)
 -8i\overline{\varphi^1_a} \varphi^1_a
-8\overline{\varphi^2_a}
\,
({\xi^6_5} -  {\xi^2_1})
+8{\varphi^2_a}
\,
({\xi^6_5} +  {\xi^2_1})
+i{\varphi^1_a}^2
\,
({\xi^6_5}{\xi^2_1}+4)
\quad .$
%}
\par In the case of
 $J(\xi^3_3)$ in (\ref{M10case22final})
 where $ {\xi^3_3}\neq 0,$   the general computations apply with
 $r=0, b=1, \xi^4_3=\xi^2_1=-1.$
 The only difference is that we get now
%\begin{equation*}
$\chi^3(a,x)=
\frac{1}{16} \,
\left(
D_1 \, (\varphi^1_x)^2  +
D_2 \,  \varphi^1_x  \right)
$
%\end{equation*}
with
%{\fontsize{8}{10} \selectfont
$D_1=
-4i\overline{\varphi^1_a} + 3i{\varphi^1_a}
;\,
D_2=
(2i-\xi^3_3) \, (\overline{\varphi^1_a})^2
\,
 +4(\xi^3_3-2i) \, \overline{\varphi^1_a} \varphi^1_a
-8(i\xi^3_3 +2)\, \overline{\varphi^2_a}
+8i\xi^3_3 \,{\varphi^2_a}
+(3i-\xi^3_3) \, {\varphi^1_a}^2
.$
%}
%%%%%%%%%%%%%%%%%%%%%%%%%%%%%%%%%%%%%%%%%%%%%%%%%%%%%%%%%%%%%%%%%%%%%%%%%%%%%%%%%%%%%%%%%%
\section{Lie Algebra $ {M14_{\gamma}} (\gamma =\pm 1).$}
%%%%%%%%%%%%%%%%%%%%%%%%%%%%%%%%%%%%%%%%%%%%%%%%%%%%%%%%%%%%%%%%%%%%%%%%%%%%%%%%%%%%%%%%%%
Commutation relations for
$ {{M14_{\gamma}}}:$
$[x_1,x_3]=x_4$;
$[x_1,x_4]=x_6$;
$[x_2,x_3]=x_5$;
$[x_2,x_5]= \gamma x_6.$
$ {M14_{-1}}$ has no CS.
%\subsection{Case $ {M14_{1}}.$}
We consider the case $ {M14_{1}}.$
\begin{equation} \label{M14+1general}
 J = \begin{pmatrix}
0&
 - \frac{1}{\xi^2_1}&
0&
0&
0&
0\\
\boxed{{\xi^2_1}}&
0&
0&
0&
0&
0\\
\boxed{{\xi^3_1}}&
\boxed{{\xi^3_2}}&
* &
 - ({\xi^6_6} {\xi^3_5} - {\xi^6_5} {\xi^3_6}) {\xi^2_1}&
\boxed{{\xi^3_5}}&
\boxed{{\xi^3_6}}\\
\boxed{{\xi^4_1}}&
* &
 *  &
 - \frac{({\xi^6_6} {\xi^3_5} - {\xi^6_5} {\xi^3_6})^2}{\xi^3_6}&
 - \frac{(({\xi^6_5} {\xi^3_5} + 1) {\xi^3_6} - {\xi^6_6} {\xi^3_5}^2) {\xi^2_1}}{\xi^3_6}  &
({\xi^6_6} {\xi^3_5} - {\xi^6_5} {\xi^3_6}) {\xi^2_1}\\
\boxed{{\xi^5_1}}&
   *             &
    *              &
 - \frac{(({\xi^6_5} {\xi^3_5} - 1) {\xi^3_6} - {\xi^6_6} {\xi^3_5}^2) {\xi^2_1}}{\xi^3_6}
     &
 - \frac{{\xi^3_5}^2}{\xi^3_6}&
 - {\xi^3_5}\\
      *                &
 *&
  * &
- \frac{({\xi^6_6}^2 {\xi^3_5} - {\xi^6_6} {\xi^6_5} {\xi^3_6} + {\xi^3_5}) {\xi^2_1}}
{\xi^3_6}&
\boxed{{\xi^6_5}}&
\boxed{{\xi^6_6}}\end{pmatrix}
\end{equation}
where
{\fontsize{8}{10}\selectfont
$  \mathbf{J^3_3}=({\xi^6_5}^2 {\xi^3_6}^2 + {\xi^3_5}^2 + {\xi^6_6}^2 {\xi^3_5}^2 - (2 {\xi^6_5} {\xi^3_5} + 1) {\xi^6_6} {\xi^3_6})/{\xi^3_6};$
\myquad
$  \mathbf{J^4_2}=(({\xi^6_6} {\xi^3_5} - {\xi^6_5} {\xi^3_6}) {\xi^3_2} {\xi^2_1} + {\xi^5_1} {\xi^3_6} + {\xi^3_5} {\xi^3_1})/{\xi^3_6};$
\myquad
$  J^4_3=((({\xi^6_6}^2 {\xi^3_5}^3 - 3 {\xi^6_6} {\xi^6_5} {\xi^3_6} {\xi^3_5}^2 - {\xi^6_6} {\xi^3_6} {\xi^3_5} + 3 {\xi^6_5}^2 {\xi^3_6}^2 {\xi^3_5} + {\xi^6_5} {\xi^3_6}^2 + {\xi^3_5}^3) {\xi^6_6} - ({\xi^6_5}^3 {\xi^3_6}^2 + {\xi^6_5} {\xi^3_5}^2 + {\xi^3_5}) {\xi^3_6}) {\xi^2_1})/{\xi^3_6}^2;$
\myquad
$  \mathbf{J^5_2}=( - ({\xi^4_1} {\xi^3_6} + {\xi^3_5} {\xi^3_2} + {\xi^6_5} {\xi^3_6} {\xi^3_1} {\xi^2_1} - {\xi^6_6} {\xi^3_5} {\xi^3_1} {\xi^2_1}))/{\xi^3_6};$
\myquad
$  \mathbf{J^5_3}=( - (({\xi^6_6} {\xi^3_5} - 2 {\xi^6_5} {\xi^3_6}) {\xi^6_6} {\xi^3_5}^2 + {\xi^6_5}^2 {\xi^3_6}^2 {\xi^3_5} - {\xi^6_5} {\xi^3_6}^2 + {\xi^3_5}^3))/{\xi^3_6}^2;$
\myquad
$  \mathbf{J^6_1}=( - ({\xi^3_6} {\xi^3_2} {\xi^2_1} + {\xi^3_5}^2 {\xi^3_1} + {\xi^5_1} {\xi^3_6}
{\xi^3_5} + ({\xi^6_5} {\xi^3_1} + {\xi^4_1} {\xi^2_1}) {\xi^6_5} {\xi^3_6}^2 + (({\xi^6_6} {\xi^3_5} - 2 {\xi^6_5} {\xi^3_6}) {\xi^3_5} {\xi^3_1} - ({\xi^4_1} {\xi^3_5} {\xi^2_1} + {\xi^3_1})
{\xi^3_6}) {\xi^6_6}))/{\xi^3_6}^2;$
\myquad
$  \mathbf{J^6_2}=({\xi^6_6} {\xi^5_1} {\xi^3_5} + {\xi^6_6} {\xi^3_2} {\xi^2_1} - {\xi^6_5} {\xi^5_1}
{\xi^3_6} + {\xi^4_1} {\xi^3_5} {\xi^2_1} + {\xi^3_1})/({\xi^3_6} {\xi^2_1});$
\myquad
$  \mathbf{J^6_3}=({\xi^6_6}^3 {\xi^3_5}^2 - 2 {\xi^6_6}^2 {\xi^6_5} {\xi^3_6} {\xi^3_5} - {\xi^6_6}^2 {\xi^3_6} + {\xi^6_6} {\xi^6_5}^2 {\xi^3_6}^2 + {\xi^6_6} {\xi^3_5}^2 - {\xi^3_6})/{\xi^3_6}^2;$
}
and the parameters are subject to the condition
\begin{equation}
\label{condM14+1general}
\xi^2_1= \pm 1 \; ; \;
{\xi^3_6} \neq 0 .
 \end{equation}
\par Now the automorphism group of
$ {M14_1}$ is comprised of the matrices
$${\fontsize{8}{10}\selectfont
 \Phi = \begin{pmatrix}
{b^1_1}&
{b^2_1} u&
0&
0&
0&
0\\
{b^2_1}&
 - {b^1_1} u&
0&
0&
0&
0\\
0&
0&
{b^3_3}&
0&
0&
0\\
{b^4_1}&
{b^4_2}&
{b^4_3}&
{b^3_3} {b^1_1}&
{b^3_3} {b^2_1} u&
0\\
 - ({b^4_2} - {b^2_1} k) u&
{b^4_1} u - {b^1_1} k&
{b^5_3}&
{b^3_3} {b^2_1}&
 - {b^3_3} {b^1_1} u&
0\\
{b^6_1}&
{b^6_2}&
{b^6_3}&
{b^5_3} {b^2_1} + {b^4_3} {b^1_1}&
 - ({b^5_3} {b^1_1} - {b^4_3} {b^2_1}) u&
({b^2_1}^2 + {b^1_1}^2) {b^3_3}\end{pmatrix}
}
$$
where $({b^2_1}^2 + {b^1_1}^2) {b^3_3} \neq 0 $
and $ u= \pm 1, \; k \in \Rmath.$
Taking suitable values for the $b^i_j$'s,
equivalence by $\Phi$
leads to the case where
$ {\xi^3_1}= {\xi^3_2}={\xi^3_5}= {\xi^4_1}= {\xi^5_1}= {\xi^6_6}=0.$
Then equivalence by
$${\fontsize{8}{10}\selectfont
 \Phi = \begin{pmatrix}
1&
0&
0&
0&
0&
0\\
0&
 {\xi^2_1}&
0&
0&
0&
0\\
0&
0&
1&
0&
0&
0\\
 -  {\xi^2_1} /2&
0&
 - {\xi^6_5} {\xi^2_1}&
1&
0&
0\\
0&
 - 1/2&
0&
0&
 {\xi^2_1}&
0\\
 {\xi^6_5} /2&
0&
0&
 - {\xi^6_5} {\xi^2_1}&
0&
1/|{\xi^3_6}|\end{pmatrix}
}
$$
leads to the case where moreover
$ {\xi^6_5}=0,  {\xi^2_1}=1,  {\xi^3_6}^2=1 :$
 \selectfont
\begin{equation}
{\fontsize{8}{10}\selectfont
\label{M14+1final}
 J({\xi^3_6}) = \begin{pmatrix}
0&
-1&
0&
0&
0&
0\\
1&
0&
0&
0&
0&
0\\
0&
0&
0&
0&
0&
{\xi^3_6}\\
0&
0&
0&
0&
-1&
0\\
0&
0&
0&
1&
0&
0\\
0&
0&
 - 1/{\xi^3_6}&
0&
0&
0\end{pmatrix}
}
\quad , \quad
 {\xi^3_6}= \pm 1.
\end{equation}
The 2 matrices corresponding to $ {\xi^3_6}= \pm 1$ are not equivalent.
Commutation relations of $ \mathfrak{m} : $
$  [\tilde{x}_1,\tilde{x}_3]=\tilde{x}_4;$
$  [\tilde{x}_1,\tilde{x}_6]= - {\xi^3_6}  \tilde{x}_5  ;$
$  [\tilde{x}_2,\tilde{x}_3]=\tilde{x}_5;$
$  [\tilde{x}_2,\tilde{x}_6]= {\xi^3_6} \tilde{x}_4 .$
%\subsection{Conclusions.}
\par
From (\ref{M14+1general}), $\mathfrak{X}_{M14_{1}}$
is a submanifold of dimension 8 in $\Rmath^{36}.$
There are only 2 orbits, and  any CS
on $ {M14_{1}}$
is equivalent to one of the two non equivalent structures in
(\ref{M14+1final}).
%%%%%%%%%%%%%%%%%%%%%%%%%%%%%%%%%%%%%%%%%%%%%%%%%%%%%%%%%%%%%%%%%%%%%%%%%%%%%%%%%%%%%%%%%%
\subsection{}
%\subsubsection{}
\begin{equation*}
X_1 = \frac{\partial}{\partial x^1} - x^2 \, \frac{\partial}{\partial y^2}
-y^2\, \frac{\partial}{\partial y^3}
\quad , \quad
X_2 = \frac{\partial}{\partial y^1} - x^2 \, \frac{\partial}{\partial x^3}
-x^3 \, \frac{\partial}{\partial y^3}.
\end{equation*}
%\subsubsection{Holomorphic functions for $J.$}
Let $G$ denote the group $G_0$ endowed with the left invariant  structure
of complex manifold  defined by
 $J({\xi^3_6})$ in
(\ref{M14+1final}),
where
$ {\xi^3_6} = \pm 1.$
Then
$ H_{\Cmath}(G) =\{f \in C^{\infty}(G_0) \; ; \; \tilde{X}_j^{-} \, f = 0
\; \forall j=1,3,5\}.$
One has
\begin{equation*}
\tilde{X}_1^{-}  =
2\, \frac{\partial}{\partial \overline{w^1}}
-x^2\, \frac{\partial}{\partial y^2}
-ix^2\, \frac{\partial}{\partial x^3}
-iw^3\, \frac{\partial}{\partial y^3}
\quad , \quad
\tilde{X}_3^{-}  = 2  \frac{\partial}{\partial \overline{w^2}}
\quad , \quad
\tilde{X}_5^{-}  =2\, \frac{\partial}{\partial \overline{w^3}}
\end{equation*}
where
$ w^1 = x^1+ i y^1
\quad , \quad
w^2 = x^2 -i{\xi^3_6}\, y^3
\quad , \quad
w^3 = x^3 - i  y^2 . $
Then $f \in C^{\infty}(G_0)$ is in $ H_{\Cmath}(G)$ if and only if it is holomorphic
with respect to  $w^2$ and $w^3$ and satisfies the equation
\begin{equation*}
2\, \frac{\partial f}{\partial \overline{w^1}}
- \xi^3_6 w^3 \,  \frac{\partial f}{\partial {w^2}}
=0.
\end{equation*}
The 3 functions
$\varphi^1 =w^1
\quad , \quad
\varphi^2 = w^2 + \frac{\xi^3_6}{2}\, w^3\overline{w^1}
\quad , \quad
\varphi^3 =w^3 $
are holomorphic.
Let $F : G \rightarrow \Cmath^3$ defined
by
$F=(\varphi^1,\varphi^2,\varphi^3).$
$F$ is a global chart on $G.$
We determine now how the multiplication of $G$ looks like in that chart.
Let $a,x \in G $ with respective second kind canonical coordinates
$(x^1,y^1,x^2,y^2,x^3,y^3), (\alpha^1, \beta^1, \alpha^2, \beta^2, \alpha^3, \beta^3)$
as in  (\ref{x6general}).
With obvious notations,
%$a =[w^1_a, w^2_a,w^3_a],$
%$x =[w^1_x, w^2_x,w^3_x],$
%$a \, x =[w^1_{a x}, w^2_{a x},w^3_{a x}],$
%$a =[\varphi^1_a, \varphi^2_a,\varphi^3_a],$
%$x =[\varphi^1_x, \varphi^2_x,\varphi^3_x],$
%$a \, x =[\varphi^1_{a x}, \varphi^2_{a x},\varphi^3_{a x}].$
computations yield:
\begin{eqnarray*}
w^1_{a x} &=& w^1_a + w^1_x\\
w^2_{a x} &=& w^2_a + w^2_x-i \xi^3_6 ( -b^2x^1
+\frac{1}{2} \, a^2{x^1}^2 -a^3y^1 +\frac{1}{2} \, a^2{y^1}^2)\\
w^3_{a x} &=& w^3_a + w^3_x
-a^2 y^1 +ia^2x^1.
\end{eqnarray*}
We then get
\begin{equation*}
\varphi^1_{a x} = \varphi^1_a + \varphi^1_x
\quad , \quad
\varphi^2_{a x} = \varphi^2_a + \varphi^2_x +\chi^2(a,x)
\quad , \quad
\varphi^3_{a x} = \varphi^3_a +\varphi^3_x  +\chi^3(a,x)
\end{equation*}
where
\begin{eqnarray*}
\chi^2(a,x)&= &
\frac{1}{8} \, \varphi^1_x \left(2i\xi^3_6\overline{\varphi^1_a}
(\overline{\varphi^2_a}+ {\varphi^2_a})
+4\xi^3_6\overline{\varphi^3_a}
-i(\overline{\varphi^1_a})^2 {\varphi^3_a}
-i\,\overline{\varphi^1_a}\overline{\varphi^3_a} {\varphi^1_a}\right)
+\frac{\xi^3_6}{2}\,
\overline{\varphi^1_a}{\varphi^3_x} ;\\
\chi^3(a,x)&= &
\frac{i}{4} \, \varphi^1_x \left(
-3\xi^3_6 (\overline{\varphi^1_a} {\varphi^3_a} +{\varphi^1_a} \overline{\varphi^3_a})
+2({\varphi^2_a}+ \overline{\varphi^2_a})\right).
\end{eqnarray*}
%%%%%%%%%%%%%%%%%%%%%%%%%%%%%%%%%%%%%%%%%%%%%%%%%%%%%%%%%%%%%%%%%%%%%%%%%%%%%%%%%%%%%%%%%%
%%%%%%%%%%%%%%%%%%%%%%%%%%%%%%%%%%%%%%%%%%%%%%%%%%%%%%%%%%%%%%%%%%%%%%%%%%%%%%%%%%%%%%%%%%
\section{Lie Algebra $ {M18_{\gamma}} (\gamma =\pm 1).$}
%%%%%%%%%%%%%%%%%%%%%%%%%%%%%%%%%%%%%%%%%%%%%%%%%%%%%%%%%%%%%%%%%%%%%%%%%%%%%%%%%%%%%%%%%%
Commutation relations for
$ {{M18_{\gamma}}}:$
$[x_1,x_2]=x_3$;
$[x_1,x_3]=x_4$;
$[x_1,x_4]=x_6$;
$[x_2,x_3]=x_5$;
$[x_2,x_5]=\gamma x_6.$
$ {M18_{-1}}$ has no CS.
%\subsection{Case $ {M18_{1}}.$}
We consider the case $ M18_{1}.$
%{\fontsize{8}{10} \selectfont
%\begin{equation}
%\label{M18+1general}
%J = \text{same matrix as (\ref{M14+1general})}
Then $J$ is the same matrix as (\ref{M14+1general})
%\end{equation}
and the parameters are subject to the same condition
(\ref{condM14+1general}).
%\begin{equation}
%\label{condM18+1general}
%\text{same condition as (\ref{condM14+1general})}
% \end{equation}
This comes as no surprise, since the commutations relations of $M18_{\gamma}$ are simply
those  of $M14_{\gamma}$ plus $[x_1,x_2]=x_3,$ and any $J \in {\frak{X}}_{M18_{1}}$ has
$\xi^1_k=\xi^2_k= 0$ for $3 \leqslant k \leqslant 6$  and $\xi^1_1=\xi^2_2=0.$
\par Now the automorphism group of
$ {M18_{1}}$ is comprised of the matrices
$$
{\fontsize{8}{10}\selectfont
\Phi = \begin{pmatrix}
{b^1_1}&
{b^2_1} u&
0&
0&
0&
0\\
{b^2_1}&
 - {b^1_1} u&
0&
0&
0&
0\\
{b^3_1}&
{b^3_2}&
 - ({b^2_1}^2 + {b^1_1}^2) u&
0&
0&
0\\
{b^4_1}&
{b^4_2}&
{b^3_2} {b^1_1} - {b^3_1} {b^2_1} u&
 - ({b^2_1}^2 + {b^1_1}^2) {b^1_1} u&
 - ({b^2_1}^2 + {b^1_1}^2) {b^2_1}&
0\\
{b^5_1}&
{b^5_2}&
{b^3_2} {b^2_1} + {b^3_1} {b^1_1} u&
 - ({b^2_1}^2 + {b^1_1}^2) {b^2_1} u&
({b^2_1}^2 + {b^1_1}^2) {b^1_1}&
0\\
{b^6_1}&
{b^6_2}&
{b^4_2} {b^1_1} - {b^4_1} {b^2_1} u + {b^5_1} {b^1_1} u + {b^5_2} {b^2_1}&
({b^2_1}^2 + {b^1_1}^2) {b^3_2}&
 - ({b^2_1}^2 + {b^1_1}^2) {b^3_1}&
 - ({b^2_1}^2 + {b^1_1}^2)^2 u\end{pmatrix}
 } $$
where $ {b^2_1}^2 + {b^1_1}^2 \neq 0 $
and $ u= \pm 1.$
Taking suitable values for the $b^i_j$'s,
equivalence by $\Phi$ leads to the case where
 $ {\xi^3_1}=0, {\xi^3_2}=0,  {\xi^3_5}=0, {\xi^4_1}=0, {\xi^5_1}=0, {\xi^6_5}=0,
{\xi^6_6}=0.$
Then equivalence by
$ \Phi = \text{diag}(1, {\xi^2_1},{\xi^2_1}/|{\xi^3_6}|,
{\xi^2_1}/|{\xi^3_6}|,1/|{\xi^3_6}|,{\xi^2_1}/|{\xi^3_6}|^2)$
leads to the case where moreover
$  {\xi^2_1}=1,  {\xi^3_6}^2=1 ,$ that is the same matrix
 $J({\xi^3_6})$ as in (\ref{M14+1final}) with the same condition.
Again, the two matrices corresponding to $ {\xi^3_6}= \pm 1$ are not equivalent.
Commutation relations of $ \mathfrak{m}$ :
 $  [\tilde{x}_1,\tilde{x}_3]=\tilde{x}_4;$
 $  [\tilde{x}_1,\tilde{x}_6]= - {\xi^3_6}  \tilde{x}_5  ;$
 $  [\tilde{x}_2,\tilde{x}_3]=\tilde{x}_5;$
 $  [\tilde{x}_2,\tilde{x}_6]=   {\xi^3_6}  \tilde{x}_4   .$
%\subsection{Conclusions.}
From (\ref{M14+1general}), $\mathfrak{X}_{M18_{1}}$
is a submanifold of dimension 8 in $\Rmath^{36}.$
There are only 2 orbits, and  any
$ J \in \mathfrak{X}_{M18_{1}}$
is equivalent to one of the two non equivalent structures
in   (\ref{M14+1final}).
%%%%%%%%%%%%%%%%%%%%%%%%%%%%%%%%%%%%%%%%%%%%%%%%%%%%%%%%%%%%%%%%%%%%%%%%%%%%%%%%%%%%%%%%%%
\subsection{}
%\subsubsection{}
\begin{eqnarray*}
X_1 &=& \frac{\partial}{\partial x^1} - y^1 \, \frac{\partial}{\partial x^2}
-x^2\, \frac{\partial}{\partial y^2}
+\frac{1}{2}\, {y^1}^2 \, \frac{\partial}{\partial x^3}
-\left(y^2 +\frac{1}{6}\, {y^1}^3 \right)\, \frac{\partial}{\partial y^3}\\
X_2 &=& \frac{\partial}{\partial y^1} - x^2 \, \frac{\partial}{\partial x^3}
-x^3 \, \frac{\partial}{\partial y^3}.
\end{eqnarray*}
%\subsubsection{Holomorphic functions for $J.$}
Let $G$ denote the group $G_0$ endowed with the left invariant  structure
of complex manifold  defined by
 $J({\xi^3_6})$ in
(\ref{M14+1final}),
where
$ {\xi^3_6} = \pm 1.$
Then
$ H_{\Cmath}(G) =\{f \in C^{\infty}(G_0) \; ; \; \tilde{X}_j^{-} \, f = 0
\; \forall j=1,3,5\}.$
One has
\begin{eqnarray*}
\tilde{X}_1^{-}  &=&
2\, \frac{\partial}{\partial \overline{w^1}}
-y^1\, \frac{\partial}{\partial x^2}
-x^2\, \frac{\partial}{\partial y^2}
+\left(\frac{1}{2}\, {y^1}^2 -ix^2\right)\, \frac{\partial}{\partial x^3}
-\left(y^2 +\frac{1}{6}\, {y^1}^3 +ix^3\right)\, \frac{\partial}{\partial y^3}
\quad ,
\\
\tilde{X}_3^{-}  &=& 2  \frac{\partial}{\partial \overline{w^2}}
\quad ,
\tilde{X}_5^{-}  =2\, \frac{\partial}{\partial \overline{w^3}}
\quad ,
\end{eqnarray*}
where
$
w^1 = x^1+ i y^1
\quad , \quad
w^2 = x^2 -i{\xi^3_6}\, y^3
\quad , \quad
w^3 = x^3 - i  y^2 . $
Then $f \in C^{\infty}(G_0)$ is in $ H_{\Cmath}(G)$ if and only if it is holomorphic
with respect to  $w^2$ and $w^3$ and satisfies the equation
\begin{equation*}
2\, \frac{\partial f}{\partial \overline{w^1}}
+\left(-y^1 +\frac{i\xi^3_6}{6}\, {y^1}^3 -\xi^3_6 w^3\right)
\, \frac{\partial f}{\partial w^2}
+\frac{1}{2} \, {y^1}^2 \,  \frac{\partial f}{\partial {w^3}}
=0.
\end{equation*}
The 3 functions
\begin{eqnarray*}
\varphi^1 &=&w^1\\
\varphi^2 &=& w^2 + \frac{\xi^3_6}{2}\, w^3\overline{w^1}
+ \frac{1}{8}\, (w^1-\overline{w^1})^2 \left(i
- \frac{\xi^3_6}{48}\, (w^1-\overline{w^1})^2 \right)
+ \frac{\xi^3_6}{384}\, (\overline{w^1})^2 \left(
6{w^1}^2-8w^1\overline{w^1} +3(\overline{w^1})^2\right)\\
\varphi^3 &=&w^3 +\frac{1}{48}\,
\overline{w^1} \left(  3{w^1}^2 - 3w^1\overline{w^1} +(\overline{w^1})^2\right)
\end{eqnarray*}
are holomorphic.
Let $F : G \rightarrow \Cmath^3$ defined
by
$F=(\varphi^1,\varphi^2,\varphi^3).$
$F$ is a a global chart on $G.$
We determine now how the multiplication of $G$ looks like in that chart.
Let $a,x \in G $ with respective second kind canonical coordinates
$(x^1,y^1,x^2,y^2,x^3,y^3), (\alpha^1, \beta^1, \alpha^2, \beta^2, \alpha^3, \beta^3)$
as in  (\ref{x6general}).
With obvious notations,
computations yield:
\begin{eqnarray*}
w^1_{a x} &=& w^1_a + w^1_x\\
w^2_{a x} &=& w^2_a + w^2_x
\\ & &
-b^1x^1 +i \xi^3_6 \left( a^3y^1+b^2x^1
-\frac{1}{2} \, a^2{x^1}^2
+\frac{1}{6} \, {b^1}^3{x^1}
+\frac{1}{6} \, {b^1}{x^1}^3
+\frac{1}{2} \, {b^1}^2x^1{y^1}
-\frac{1}{2} \, {y^1}^2(a^2 -b^1{x^1})
\right)\\
w^3_{a x} &=& w^3_a + w^3_x
+\frac{1}{2} \, {b^1}^2x^1
-{y^1}(a^2 -b^1{x^1})
-i\left(
\frac{1}{2} \, b^1{x^1}^2 -a^2x^1\right).
\end{eqnarray*}
We then get
\begin{equation*}
\varphi^1_{a x} = \varphi^1_a + \varphi^1_x
\quad , \quad
\varphi^2_{a x} = \varphi^2_a + \varphi^2_x +\chi^2(a,x)
\quad , \quad
\varphi^3_{a x} = \varphi^3_a +\varphi^3_x  +\chi^3(a,x)
\end{equation*}
where
\begin{eqnarray*}
\chi^2(a,x)&= &
\frac{\xi^3_6}{32} \, (\varphi^1_x)^3 \left(\overline{\varphi^1_a} -{\varphi^1_a} \right)
+\frac{\xi^3_6}{64} \, (\varphi^1_x)^2 \left(4(\overline{\varphi^1_a})^2 -3({\varphi^1_a})^2 \right)
+\frac{1}{512}\,  D_2\, \varphi^1_x
+\frac{\xi^3_6}{2} \, \varphi^3_x \overline{\varphi^1_a} ;
\\
\chi^3(a,x)&= &
\frac{1}{16} \, (\varphi^1_x)^2 \left(  4\overline{\varphi^1_a} -3{\varphi^1_a} \right) +
\frac{1}{256} D_3\, \varphi^1_x ;
\end{eqnarray*}
where
\begin{multline*}
D_2= 16\xi^3_6 \left(
-(\overline{\varphi^1_a})^3
+8i\overline{\varphi^1_a} \overline{\varphi^2_a}
+2\overline{\varphi^1_a}(\varphi^1_a)^2
+8i\overline{\varphi^1_a}\varphi^2_a
+16\overline{\varphi^3_a}
-(\varphi^1_a)^3 \right)
\\
+i \left(
(\overline{\varphi^1_a})^5
-4(\overline{\varphi^1_a})^4 {\varphi^1_a}
+8(\overline{\varphi^1_a})^3 ({\varphi^1_a})^2
-4(\overline{\varphi^1_a})^2 ({\varphi^1_a})^3
\right. \\ \left.
-64(\overline{\varphi^1_a})^2 {\varphi^3_a}
-64\overline{\varphi^1_a} \overline{\varphi^3_a} {\varphi^1_a}
+\overline{\varphi^1_a} ({\varphi^1_a})^4
-256(\overline{\varphi^1_a} -{\varphi^1_a})\right)
\end{multline*}
\begin{multline*}
D_3= i\xi^3_6 \left(
(\overline{\varphi^1_a})^4
-4(\overline{\varphi^1_a})^3{\varphi^1_a}
+8(\overline{\varphi^1_a})^2(\varphi^1_a)^2
-4\overline{\varphi^1_a}(\varphi^1_a)^3
-64\overline{\varphi^1_a}\varphi^3_a
-64\overline{\varphi^3_a}\varphi^1_a
+({\varphi^1_a})^4 \right)
\\
-32(\overline{\varphi^1_a})^2
-16({\varphi^1_a})^2
+64\overline{\varphi^1_a}\varphi^1_a
+128i(\overline{\varphi^2_a}+\varphi^2_a)  .
\end{multline*}
%%%%%%%%%%%%%%%%%%%%%%%%%%%%%%%%%%%%%%%%%%%%%%%%%%%%%%%%%%%%%%%%%%%%%%%%%%%%%%%%%%%%%%%%%%
\section{Lie Algebra $ {M5}.$}
%%%%%%%%%%%%%%%%%%%%%%%%%%%%%%%%%%%%%%%%%%%%%%%%%%%%%%%%%%%%%%%%%%%%%%%%%%%%%%%%%%%%%%%%%%
Commutation relations for
$ {M5}:$
$[x_1,x_3]=x_5$;
$[x_1,x_4]=x_6$;
$[x_2,x_3]= - x_6$;
$[x_2,x_4]=x_5.$
This is the realification of the 3-dimensional complex Heisenberg Lie algebra
$\mathfrak{n} \;$
$[Z_1,Z_2]=Z_3$ we get by letting $x_1=Z_1,x_2=-iZ_1,x_3=Z_2,x_4=iZ_2,x_5=Z_3,x_6=iZ_3.$
\subsection{Case $ {\xi^2_4}^2 + {\xi^2_3}^2 \neq 0 .$}
\label{xi(2,4)}
%{\fontsize{8}{10}\selectfont
\begin{equation}
\label{globalchartxi(2,1)M5}
J = \begin{pmatrix}
\boxed{\xi^1_1}&
*&
\boxed{\xi^1_3}&
\boxed{\xi^1_4}&
0&
0\\
\boxed{\xi^2_1}&
*&
\boxed{\xi^2_3}&
\boxed{\xi^2_4}&
0&
0\\
*&
*&
*&
*&
0&
0\\
*&
*&
*&
*&
0&
0\\
*&
*&
*&
*&
\boxed{\xi^5_5}&
 - \frac{{\xi^5_5}^2 + 1}{\xi^6_5}\\
\boxed{\xi^6_1}&
\boxed{\xi^6_2}&
\boxed{\xi^6_3}&
\boxed{\xi^6_4}&
\boxed{\xi^6_5}&
 - {\xi^5_5}\end{pmatrix}
 \end{equation}
%}
where
{\fontsize{6}{10}\selectfont
$  \mathbf{J^1_2}=( - ((({\xi^6_5} {\xi^2_3} - {\xi^6_5} {\xi^1_4} - {\xi^5_5} {\xi^2_4}) {\xi^1_4} +
 ({\xi^2_3} - 2 {\xi^1_4}) {\xi^5_5} {\xi^1_3}) {\xi^6_5} - ({\xi^5_5}^2 + 1) ({\xi^2_4} + {\xi^1_3}) {\xi^1_3} - (({\xi^1_4}^2 + {\xi^1_3}^2) {\xi^2_1} - {\xi^2_3} {\xi^1_3} {\xi^1_1} -
{\xi^2_4} {\xi^1_4} {\xi^1_1}) {\xi^6_5}))/(({\xi^2_4} {\xi^1_3} - {\xi^2_3} {\xi^1_4}) {\xi^6_5})
;$
\myquad
$  \mathbf{J^2_2}=( - ((({\xi^6_5} {\xi^2_3} - {\xi^6_5} {\xi^1_4} - {\xi^5_5} {\xi^2_4} - {\xi^5_5}
{\xi^1_3}) {\xi^2_4} + ({\xi^2_3} - {\xi^1_4}) {\xi^5_5} {\xi^2_3}) {\xi^6_5} - ({\xi^5_5}^2 +
1) ({\xi^2_4} + {\xi^1_3}) {\xi^2_3} + (({\xi^2_4} {\xi^1_1} - {\xi^2_1} {\xi^1_4}) {\xi^2_4} +
({\xi^2_3} {\xi^1_1} - {\xi^2_1} {\xi^1_3}) {\xi^2_3}) {\xi^6_5}))/(({\xi^2_4} {\xi^1_3} - {\xi^2_3} {\xi^1_4}) {\xi^6_5});$
\myquad
$  \mathbf{J^3_1}=( - ((({\xi^2_1} {\xi^1_3} - {\xi^1_4} {\xi^1_1} - {\xi^2_3} {\xi^1_1}) {\xi^2_1} +
 ({\xi^1_1}^2 + 1) {\xi^2_4}) {\xi^6_5} + (({\xi^5_5}^2 + 1) ({\xi^2_4} + {\xi^1_3}) - ({\xi^2_3} - {\xi^1_4}) {\xi^6_5} {\xi^5_5}) {\xi^2_1}))/(({\xi^2_4} {\xi^1_3} - {\xi^2_3} {\xi^1_4})
 {\xi^6_5});$
\myquad
$  \mathbf{J^3_2}=( - (((({\xi^6_5}^2 {\xi^2_3}^2 - 2 {\xi^6_5}^2 {\xi^2_3} {\xi^1_4} + {\xi^6_5}^2 {\xi^1_4}^2 - 2 {\xi^6_5} {\xi^5_5} {\xi^2_4} {\xi^2_3} + 2 {\xi^6_5} {\xi^5_5} {\xi^2_4}
{\xi^1_4} - 2 {\xi^6_5} {\xi^5_5} {\xi^2_3} {\xi^1_3} + 2 {\xi^6_5} {\xi^5_5} {\xi^1_4} {\xi^1_3}
+ {\xi^6_5} {\xi^2_4} {\xi^2_1} {\xi^1_3} - {\xi^6_5} {\xi^2_3} {\xi^2_1} {\xi^1_4} - 2 {\xi^6_5}
{\xi^2_3} {\xi^1_3} {\xi^1_1} + {\xi^5_5}^2 {\xi^2_4}^2 + 2 {\xi^5_5}^2 {\xi^2_4} {\xi^1_3} -
 2 {\xi^5_5}^2 {\xi^2_3}^2 + 2 {\xi^5_5}^2 {\xi^2_3} {\xi^1_4} + {\xi^5_5}^2 {\xi^1_3}^2
- {\xi^5_5} {\xi^2_4}^2 {\xi^1_1} + {\xi^5_5} {\xi^2_4} {\xi^2_1} {\xi^1_4} - 2 {\xi^5_5} {\xi^2_3}^2 {\xi^1_1} + 2 {\xi^5_5} {\xi^2_3} {\xi^2_1} {\xi^1_3} + {\xi^5_5} {\xi^1_3}^2 {\xi^1_1} +
 {\xi^2_4}^2 + 2 {\xi^2_4} {\xi^1_3} - 2 {\xi^2_3}^2 + 2 {\xi^2_3} {\xi^1_4} + {\xi^1_3}^2
+ ({\xi^2_1} {\xi^1_3} + 2 {\xi^1_4} {\xi^1_1}) {\xi^6_5} {\xi^1_3}) {\xi^2_4} + (({\xi^2_3} -
{\xi^1_4}) {\xi^6_5} {\xi^5_5} - 2 {\xi^5_5}^2 {\xi^1_3} - 2 {\xi^1_3}) ({\xi^2_3} - {\xi^1_4})
 {\xi^2_3} - (2 {\xi^2_3}^2 {\xi^1_1} - 2 {\xi^2_3} {\xi^2_1} {\xi^1_3} + {\xi^2_1} {\xi^1_4}
{\xi^1_3}) {\xi^5_5} {\xi^1_3} + (2 {\xi^2_3}^3 {\xi^1_1} - 2 {\xi^2_3}^2 {\xi^2_1} {\xi^1_3}
- 2 {\xi^2_3}^2 {\xi^1_4} {\xi^1_1} + 3 {\xi^2_3} {\xi^2_1} {\xi^1_4} {\xi^1_3} - 2 {\xi^2_1}
{\xi^1_4}^2 {\xi^1_3}) {\xi^6_5}) {\xi^5_5} - ({\xi^2_4}^2 {\xi^1_1} - {\xi^2_4} {\xi^2_1} {\xi^1_4} - {\xi^2_4} {\xi^1_3} {\xi^1_1} + 2 {\xi^2_3}^2 {\xi^1_1} - 2 {\xi^2_3} {\xi^2_1} {\xi^1_3} + {\xi^2_1} {\xi^1_4} {\xi^1_3}) ({\xi^2_4} + {\xi^1_3})) {\xi^6_5} + ({\xi^5_5}^2 + 1)^2
({\xi^2_4} + {\xi^1_3})^2 {\xi^2_3} + (({\xi^2_3} {\xi^1_1} - {\xi^2_1} {\xi^1_3} - {\xi^1_4}
{\xi^1_1}) {\xi^2_4} + {\xi^2_1} {\xi^1_4}^2) ({\xi^2_3} - {\xi^1_4}) {\xi^6_5}^3 + ((({\xi^2_4} {\xi^1_1}^2 - {\xi^2_4} - {\xi^2_1} {\xi^1_4} {\xi^1_1} - ({\xi^1_1}^2 + 1) {\xi^1_3}) {\xi^2_3} - ({\xi^2_1} {\xi^1_3} {\xi^1_1} + {\xi^1_4} {\xi^1_1}^2 - {\xi^1_4}) {\xi^2_4} + ({\xi^2_1} {\xi^1_4} {\xi^1_3} + 2 {\xi^1_4}^2 {\xi^1_1} + {\xi^1_3}^2 {\xi^1_1}) {\xi^2_1}) {\xi^2_4} + ({\xi^2_3}^2 {\xi^1_1}^2 - 2 {\xi^2_3} {\xi^2_1} {\xi^1_3} {\xi^1_1} + {\xi^2_1}^2 {\xi^1_3}^2 + {\xi^2_1} {\xi^1_4} {\xi^1_3} {\xi^1_1} + {\xi^1_4}^2) {\xi^2_3} - ({\xi^1_4}^2 +
{\xi^1_3}^2) {\xi^2_1}^2 {\xi^1_4}) {\xi^6_5}^2))/(({\xi^2_4} {\xi^1_3} - {\xi^2_3} {\xi^1_4}
)^2 {\xi^6_5}^2);$
\myquad
$  \mathbf{J^3_3}=( - ((({\xi^5_5}^2 + 1) ({\xi^2_4} + {\xi^1_3}) - ({\xi^2_3} - {\xi^1_4}) {\xi^6_5} {\xi^5_5}) {\xi^2_3} - ({\xi^2_3}^2 {\xi^1_1} - {\xi^2_3} {\xi^2_1} {\xi^1_3} + {\xi^2_1}
{\xi^1_4} {\xi^1_3} - {\xi^2_4} {\xi^1_3} {\xi^1_1}) {\xi^6_5}))/(({\xi^2_4} {\xi^1_3} - {\xi^2_3}
 {\xi^1_4}) {\xi^6_5});$\\
$  J^3_4=( - ((({\xi^5_5}^2 + 1) ({\xi^2_4} + {\xi^1_3}) - ({\xi^2_3} - {\xi^1_4}) {\xi^6_5} {\xi^5_5}) {\xi^2_4} - ({\xi^2_4} {\xi^2_3} {\xi^1_1} - {\xi^2_4} {\xi^2_1} {\xi^1_3} - {\xi^2_4} {\xi^1_4} {\xi^1_1} + {\xi^2_1} {\xi^1_4}^2) {\xi^6_5}))/(({\xi^2_4} {\xi^1_3} - {\xi^2_3}
{\xi^1_4}) {\xi^6_5});$
$  \mathbf{J^4_1}= - (({\xi^2_1} {\xi^1_4} + {\xi^1_3} {\xi^1_1}) {\xi^2_1} - ({\xi^1_1}^2 + 1) {\xi^2_3} - {\xi^2_4} {\xi^2_1} {\xi^1_1} + ({\xi^2_4} + {\xi^1_3}) {\xi^5_5} {\xi^2_1} - ({\xi^2_3}
 - {\xi^1_4}) {\xi^6_5} {\xi^2_1})/({\xi^2_4} {\xi^1_3} - {\xi^2_3} {\xi^1_4});$
\myquad
$  \mathbf{J^4_2}=( - (((({\xi^2_3} - {\xi^1_4}) {\xi^6_5} - 2 {\xi^5_5} {\xi^1_3}) ({\xi^2_3} - {\xi^1_4}) {\xi^5_5} {\xi^2_3} + ({\xi^6_5}^2 {\xi^2_3}^2 - 2 {\xi^6_5}^2 {\xi^2_3} {\xi^1_4} +
 {\xi^6_5}^2 {\xi^1_4}^2 - 2 {\xi^6_5} {\xi^5_5} {\xi^2_4} {\xi^2_3} + 2 {\xi^6_5} {\xi^5_5}
{\xi^2_4} {\xi^1_4} - 2 {\xi^6_5} {\xi^5_5} {\xi^2_3} {\xi^1_3} + 2 {\xi^6_5} {\xi^5_5} {\xi^1_4}
{\xi^1_3} + {\xi^5_5}^2 {\xi^2_4}^2 + 2 {\xi^5_5}^2 {\xi^2_4} {\xi^1_3} - 2 {\xi^5_5}^2 {\xi^2_3}^2 + 2 {\xi^5_5}^2 {\xi^2_3} {\xi^1_4} + {\xi^5_5}^2 {\xi^1_3}^2 - 2 {\xi^5_5} {\xi^2_4}^2 {\xi^1_1} + 2 {\xi^5_5} {\xi^2_4} {\xi^2_1} {\xi^1_4} - 2 {\xi^5_5} {\xi^2_4} {\xi^1_3}
{\xi^1_1} + {\xi^5_5} {\xi^2_3} {\xi^2_1} {\xi^1_3} - 2 {\xi^5_5} {\xi^2_3} {\xi^1_4} {\xi^1_1} +
3 {\xi^5_5} {\xi^2_1} {\xi^1_4} {\xi^1_3}) {\xi^2_4} - ({\xi^2_3}^2 {\xi^2_1} - {\xi^2_3} {\xi^2_1} {\xi^1_4} + 2 {\xi^2_3} {\xi^1_3} {\xi^1_1} - 2 {\xi^2_1} {\xi^1_3}^2) {\xi^5_5} {\xi^1_4})
{\xi^6_5} + ({\xi^5_5}^2 + 1) ({\xi^5_5} {\xi^2_4} {\xi^2_3} + {\xi^5_5} {\xi^2_3} {\xi^1_3} -
{\xi^2_4} {\xi^2_3} {\xi^1_1} + {\xi^2_3} {\xi^2_1} {\xi^1_4} - {\xi^2_3} {\xi^1_3} {\xi^1_1} + {\xi^2_1} {\xi^1_3}^2) ({\xi^2_4} + {\xi^1_3}) + (2 ({\xi^2_4} {\xi^1_1} - {\xi^2_1} {\xi^1_4})
{\xi^2_4} + ({\xi^2_3} {\xi^1_1} - {\xi^2_1} {\xi^1_3}) ({\xi^2_3} + {\xi^1_4})) ({\xi^2_3} - {\xi^1_4}) {\xi^6_5}^2 - ((({\xi^2_3} {\xi^1_4} {\xi^1_1} - {\xi^2_1} {\xi^1_4} {\xi^1_3} + {\xi^1_4}^2 {\xi^1_1} + 2 {\xi^1_3}^2 {\xi^1_1}) {\xi^2_1} - ({\xi^1_1} + 1) ({\xi^1_1} - 1) {\xi^2_3} {\xi^1_3}) {\xi^2_3} - ({\xi^1_4}^2 + {\xi^1_3}^2) {\xi^2_1}^2 {\xi^1_3} - ({\xi^2_4}^
2 {\xi^1_1}^2 - 2 {\xi^2_4} {\xi^2_1} {\xi^1_4} {\xi^1_1} + {\xi^2_3}^2 {\xi^1_1}^2 - {\xi^2_3}^2 - {\xi^2_3} {\xi^2_1} {\xi^1_3} {\xi^1_1} + {\xi^2_3} {\xi^1_4} {\xi^1_1}^2 + {\xi^2_3}
{\xi^1_4} + {\xi^2_1}^2 {\xi^1_4}^2 - {\xi^2_1} {\xi^1_4} {\xi^1_3} {\xi^1_1} + {\xi^1_3}^2)
{\xi^2_4}) {\xi^6_5}))/(({\xi^2_4} {\xi^1_3} - {\xi^2_3} {\xi^1_4})^2 {\xi^6_5});$
\myquad
$  \mathbf{J^4_3}= - ({\xi^2_3} {\xi^2_1} {\xi^1_4} - {\xi^2_3} {\xi^1_3} {\xi^1_1} + {\xi^2_1} {\xi^1_3}^2 - {\xi^2_4} {\xi^2_3} {\xi^1_1} + ({\xi^2_4} + {\xi^1_3}) {\xi^5_5} {\xi^2_3} - ({\xi^2_3}
 - {\xi^1_4}) {\xi^6_5} {\xi^2_3})/({\xi^2_4} {\xi^1_3} - {\xi^2_3} {\xi^1_4});$
\myquad
$  \mathbf{J^4_4}=(({\xi^2_4} {\xi^1_1} - {\xi^2_1} {\xi^1_4}) {\xi^2_4} + ({\xi^2_3} {\xi^1_1} - {\xi^2_1} {\xi^1_3}) {\xi^1_4} - ({\xi^2_4} + {\xi^1_3}) {\xi^5_5} {\xi^2_4} + ({\xi^2_3} - {\xi^1_4}
) {\xi^6_5} {\xi^2_4})/({\xi^2_4} {\xi^1_3} - {\xi^2_3} {\xi^1_4});$
\myquad
$  \mathbf{J^5_1}=(((({\xi^5_5} - {\xi^1_1}) {\xi^6_1} - {\xi^6_2} {\xi^2_1}) ({\xi^2_4} {\xi^1_3} -
{\xi^2_3} {\xi^1_4}) - ({\xi^2_3} - {\xi^1_4}) {\xi^6_3} {\xi^5_5} {\xi^2_1} + (({\xi^2_1} {\xi^1_3} - {\xi^1_4} {\xi^1_1} - {\xi^2_3} {\xi^1_1}) {\xi^2_1} + ({\xi^1_1}^2 + 1) {\xi^2_4}) {\xi^6_3}) {\xi^6_5} + (({\xi^2_1} {\xi^1_4} + {\xi^1_3} {\xi^1_1}) {\xi^2_1} - ({\xi^1_1}^2 + 1)
{\xi^2_3} - {\xi^2_4} {\xi^2_1} {\xi^1_1}) {\xi^6_5} {\xi^6_4} - (({\xi^6_5} {\xi^2_3} - {\xi^6_5}
 {\xi^1_4} - {\xi^5_5} {\xi^2_4} - {\xi^5_5} {\xi^1_3}) {\xi^6_5} {\xi^6_4} - ({\xi^5_5}^2 + 1)
 ({\xi^2_4} + {\xi^1_3}) {\xi^6_3}) {\xi^2_1})/(({\xi^2_4} {\xi^1_3} - {\xi^2_3} {\xi^1_4}) {\xi^6_5}^2);$
\myquad
$  \mathbf{J^5_2}=(((({\xi^2_3} - {\xi^1_4}) {\xi^6_5} - 2 {\xi^5_5} {\xi^1_3}) ({\xi^2_3} - {\xi^1_4}) {\xi^6_5} {\xi^5_5} {\xi^2_3} + {\xi^6_5}^3 {\xi^2_4} {\xi^2_3}^2 - 2 {\xi^6_5}^3 {\xi^2_4} {\xi^2_3} {\xi^1_4} + {\xi^6_5}^3 {\xi^2_4} {\xi^1_4}^2 - 2 {\xi^6_5}^2 {\xi^5_5} {\xi^2_4}
^2 {\xi^2_3} + 2 {\xi^6_5}^2 {\xi^5_5} {\xi^2_4}^2 {\xi^1_4} - 2 {\xi^6_5}^2 {\xi^5_5} {\xi^2_4} {\xi^2_3} {\xi^1_3} + 2 {\xi^6_5}^2 {\xi^5_5} {\xi^2_4} {\xi^1_4} {\xi^1_3} + {\xi^6_5} {\xi^5_5}^2 {\xi^2_4}^3 + 2 {\xi^6_5} {\xi^5_5}^2 {\xi^2_4}^2 {\xi^1_3} - 2 {\xi^6_5} {\xi^5_5}^2 {\xi^2_4} {\xi^2_3}^2 + 2 {\xi^6_5} {\xi^5_5}^2 {\xi^2_4} {\xi^2_3} {\xi^1_4} + {\xi^6_5}
 {\xi^5_5}^2 {\xi^2_4} {\xi^1_3}^2 - 2 {\xi^6_5} {\xi^5_5} {\xi^2_4}^3 {\xi^1_1} + 2 {\xi^6_5} {\xi^5_5} {\xi^2_4}^2 {\xi^2_1} {\xi^1_4} - 2 {\xi^6_5} {\xi^5_5} {\xi^2_4}^2 {\xi^1_3} {\xi^1_1} + {\xi^6_5} {\xi^5_5} {\xi^2_4} {\xi^2_3} {\xi^2_1} {\xi^1_3} - 2 {\xi^6_5} {\xi^5_5} {\xi^2_4} {\xi^2_3} {\xi^1_4} {\xi^1_1} + 3 {\xi^6_5} {\xi^5_5} {\xi^2_4} {\xi^2_1} {\xi^1_4} {\xi^1_3} +
{\xi^5_5}^3 {\xi^2_4}^2 {\xi^2_3} + 2 {\xi^5_5}^3 {\xi^2_4} {\xi^2_3} {\xi^1_3} + {\xi^5_5}^
3 {\xi^2_3} {\xi^1_3}^2 - {\xi^5_5}^2 {\xi^2_4}^2 {\xi^2_3} {\xi^1_1} + {\xi^5_5}^2 {\xi^2_4} {\xi^2_3} {\xi^2_1} {\xi^1_4} - 2 {\xi^5_5}^2 {\xi^2_4} {\xi^2_3} {\xi^1_3} {\xi^1_1} + {\xi^5_5}^2 {\xi^2_4} {\xi^2_1} {\xi^1_3}^2 + {\xi^5_5} {\xi^2_4}^2 {\xi^2_3} + 2 {\xi^5_5} {\xi^2_4} {\xi^2_3} {\xi^1_3} + {\xi^5_5} {\xi^2_3} {\xi^1_3}^2 + ({\xi^2_3} {\xi^2_1} {\xi^1_4} - {\xi^2_3} {\xi^1_3} {\xi^1_1} + {\xi^2_1} {\xi^1_3}^2) {\xi^5_5}^2 {\xi^1_3} - ({\xi^2_3}^2 {\xi^2_1} - {\xi^2_3} {\xi^2_1} {\xi^1_4} + 2 {\xi^2_3} {\xi^1_3} {\xi^1_1} - 2 {\xi^2_1} {\xi^1_3}^2)
 {\xi^6_5} {\xi^5_5} {\xi^1_4} - ({\xi^2_4} {\xi^2_3} {\xi^1_1} - {\xi^2_3} {\xi^2_1} {\xi^1_4} +
{\xi^2_3} {\xi^1_3} {\xi^1_1} - {\xi^2_1} {\xi^1_3}^2) ({\xi^2_4} + {\xi^1_3})) {\xi^6_5} {\xi^6_4} + ({\xi^5_5}^2 + 1)^2 ({\xi^2_4} + {\xi^1_3})^2 {\xi^6_3} {\xi^2_3} + (2 ({\xi^2_4} {\xi^1_1} - {\xi^2_1} {\xi^1_4}) {\xi^2_4} + ({\xi^2_3} {\xi^1_1} - {\xi^2_1} {\xi^1_3}) ({\xi^2_3}
+ {\xi^1_4})) ({\xi^2_3} - {\xi^1_4}) {\xi^6_5}^3 {\xi^6_4} - (({\xi^6_2} {\xi^2_3} + {\xi^6_1}
 {\xi^1_3}) ({\xi^2_4} {\xi^1_3} - {\xi^2_3} {\xi^1_4}) - ({\xi^5_5} {\xi^2_4}^2 + {\xi^5_5} {\xi^2_4} {\xi^1_3} - 2 {\xi^5_5} {\xi^2_3}^2 + 2 {\xi^5_5} {\xi^2_3} {\xi^1_4} - {\xi^2_4}^2 {\xi^1_1} + {\xi^2_4} {\xi^2_1} {\xi^1_4} + {\xi^2_4} {\xi^1_3} {\xi^1_1} - 2 {\xi^2_3}^2 {\xi^1_1}
 + 2 {\xi^2_3} {\xi^2_1} {\xi^1_3} - {\xi^2_1} {\xi^1_4} {\xi^1_3}) {\xi^6_3}) ({\xi^5_5}^2 + 1
) ({\xi^2_4} + {\xi^1_3}) {\xi^6_5} + (({\xi^6_2} {\xi^2_4} + {\xi^6_1} {\xi^1_4}) ({\xi^2_4} {\xi^1_3} - {\xi^2_3} {\xi^1_4}) + ({\xi^2_3} - {\xi^1_4}) {\xi^6_3} {\xi^5_5} {\xi^2_4} + (({\xi^2_3} {\xi^1_1} - {\xi^2_1} {\xi^1_3} - {\xi^1_4} {\xi^1_1}) {\xi^2_4} + {\xi^2_1} {\xi^1_4}^2) {\xi^6_3}) ({\xi^2_3} - {\xi^1_4}) {\xi^6_5}^3 - ((({\xi^2_3} {\xi^1_4} {\xi^1_1} - {\xi^2_1} {\xi^1_4} {\xi^1_3} + {\xi^1_4}^2 {\xi^1_1} + 2 {\xi^1_3}^2 {\xi^1_1}) {\xi^2_1} - ({\xi^1_1} + 1
) ({\xi^1_1} - 1) {\xi^2_3} {\xi^1_3}) {\xi^2_3} - ({\xi^1_4}^2 + {\xi^1_3}^2) {\xi^2_1}^2
{\xi^1_3} - ({\xi^2_4}^2 {\xi^1_1}^2 - 2 {\xi^2_4} {\xi^2_1} {\xi^1_4} {\xi^1_1} + {\xi^2_3}^
2 {\xi^1_1}^2 - {\xi^2_3}^2 - {\xi^2_3} {\xi^2_1} {\xi^1_3} {\xi^1_1} + {\xi^2_3} {\xi^1_4} {\xi^1_1}^2 + {\xi^2_3} {\xi^1_4} + {\xi^2_1}^2 {\xi^1_4}^2 - {\xi^2_1} {\xi^1_4} {\xi^1_3} {\xi^1_1} + {\xi^1_3}^2) {\xi^2_4}) {\xi^6_5}^2 {\xi^6_4} - (((2 {\xi^5_5} {\xi^2_4} {\xi^2_3} -
2 {\xi^5_5} {\xi^2_4} {\xi^1_4} + 2 {\xi^5_5} {\xi^2_3} {\xi^1_3} - 2 {\xi^5_5} {\xi^1_4} {\xi^1_3} - {\xi^2_4} {\xi^2_1} {\xi^1_3} + {\xi^2_3} {\xi^2_1} {\xi^1_4} + 2 {\xi^2_3} {\xi^1_3} {\xi^1_1} - ({\xi^2_1} {\xi^1_3} + 2 {\xi^1_4} {\xi^1_1}) {\xi^1_3}) {\xi^2_4} - ({\xi^2_3} - {\xi^1_4})
^2 {\xi^5_5} {\xi^2_3} - (2 {\xi^2_3}^3 {\xi^1_1} - 2 {\xi^2_3}^2 {\xi^2_1} {\xi^1_3} - 2
{\xi^2_3}^2 {\xi^1_4} {\xi^1_1} + 3 {\xi^2_3} {\xi^2_1} {\xi^1_4} {\xi^1_3} - 2 {\xi^2_1} {\xi^1_4}^2 {\xi^1_3})) {\xi^6_3} {\xi^5_5} - (({\xi^2_4}^2 {\xi^1_1} - {\xi^2_4} {\xi^2_1} {\xi^1_4}
 + {\xi^2_3}^2 {\xi^1_1} - {\xi^2_3} {\xi^2_1} {\xi^1_3} - ({\xi^2_4}^2 - {\xi^2_3}^2 + 2
{\xi^2_3} {\xi^1_4}) {\xi^5_5}) {\xi^6_2} + ({\xi^2_4} {\xi^1_4} {\xi^1_1} + {\xi^2_3} {\xi^1_3}
{\xi^1_1} - {\xi^2_1} {\xi^1_4}^2 - {\xi^2_1} {\xi^1_3}^2 - ({\xi^2_4} {\xi^1_4} - {\xi^2_3}
{\xi^1_3} + 2 {\xi^1_4} {\xi^1_3}) {\xi^5_5}) {\xi^6_1}) ({\xi^2_4} {\xi^1_3} - {\xi^2_3} {\xi^1_4}) - ((({\xi^2_4} {\xi^1_1}^2 - {\xi^2_4} - {\xi^2_1} {\xi^1_4} {\xi^1_1} - ({\xi^1_1}^2 + 1
) {\xi^1_3}) {\xi^2_3} - ({\xi^2_1} {\xi^1_3} {\xi^1_1} + {\xi^1_4} {\xi^1_1}^2 - {\xi^1_4}) {\xi^2_4} + ({\xi^2_1} {\xi^1_4} {\xi^1_3} + 2 {\xi^1_4}^2 {\xi^1_1} + {\xi^1_3}^2 {\xi^1_1}) {\xi^2_1}) {\xi^2_4} + ({\xi^2_3}^2 {\xi^1_1}^2 - 2 {\xi^2_3} {\xi^2_1} {\xi^1_3} {\xi^1_1} + {\xi^2_1}^2 {\xi^1_3}^2 + {\xi^2_1} {\xi^1_4} {\xi^1_3} {\xi^1_1} + {\xi^1_4}^2) {\xi^2_3} - (
{\xi^1_4}^2 + {\xi^1_3}^2) {\xi^2_1}^2 {\xi^1_4}) {\xi^6_3}) {\xi^6_5}^2)/(({\xi^2_4} {\xi^1_3} - {\xi^2_3} {\xi^1_4})^2 {\xi^6_5}^3);$
\myquad
$  \mathbf{J^5_3}=( - ((({\xi^6_5} {\xi^2_3} - {\xi^6_5} {\xi^1_4} - {\xi^5_5} {\xi^2_4} - {\xi^5_5}
{\xi^1_3}) {\xi^6_5} {\xi^6_4} - ({\xi^5_5}^2 + 1) ({\xi^2_4} + {\xi^1_3}) {\xi^6_3}) {\xi^2_3}
 - ({\xi^2_3} {\xi^2_1} {\xi^1_4} - {\xi^2_3} {\xi^1_3} {\xi^1_1} + {\xi^2_1} {\xi^1_3}^2 - {\xi^2_4} {\xi^2_3} {\xi^1_1}) {\xi^6_5} {\xi^6_4} + (({\xi^6_2} {\xi^2_3} + {\xi^6_1} {\xi^1_3}) ({\xi^2_4} {\xi^1_3} - {\xi^2_3} {\xi^1_4}) - ({\xi^2_4} {\xi^1_3} - {\xi^2_3}^2) {\xi^6_3} {\xi^5_5} + ({\xi^2_3}^2 {\xi^1_1} - {\xi^2_3} {\xi^2_1} {\xi^1_3} + {\xi^2_1} {\xi^1_4} {\xi^1_3} - {\xi^2_4} {\xi^1_3} {\xi^1_1}) {\xi^6_3}) {\xi^6_5}))/(({\xi^2_4} {\xi^1_3} - {\xi^2_3} {\xi^1_4})
{\xi^6_5}^2);$
\myquad
$  \mathbf{J^5_4}=( - (({\xi^6_5} {\xi^2_4} {\xi^2_3} - {\xi^6_5} {\xi^2_4} {\xi^1_4} - {\xi^5_5} {\xi^2_4}^2 - 2 {\xi^5_5} {\xi^2_4} {\xi^1_3} + {\xi^5_5} {\xi^2_3} {\xi^1_4}) {\xi^6_5} {\xi^6_4} -
 ({\xi^5_5}^2 + 1) ({\xi^2_4} + {\xi^1_3}) {\xi^6_3} {\xi^2_4} + (({\xi^2_4} {\xi^1_1} - {\xi^2_1} {\xi^1_4}) {\xi^2_4} + ({\xi^2_3} {\xi^1_1} - {\xi^2_1} {\xi^1_3}) {\xi^1_4}) {\xi^6_5} {\xi^6_4} + (({\xi^2_4} {\xi^2_3} {\xi^1_1} - {\xi^2_4} {\xi^2_1} {\xi^1_3} - {\xi^2_4} {\xi^1_4} {\xi^1_1} + {\xi^2_1} {\xi^1_4}^2 + ({\xi^2_3} - {\xi^1_4}) {\xi^5_5} {\xi^2_4}) {\xi^6_3} + ({\xi^6_2} {\xi^2_4} + {\xi^6_1} {\xi^1_4}) ({\xi^2_4} {\xi^1_3} - {\xi^2_3} {\xi^1_4})) {\xi^6_5}))/((
{\xi^2_4} {\xi^1_3} - {\xi^2_3} {\xi^1_4}) {\xi^6_5}^2);$
}
and the parameters are subject to the condition
%\begin{equation}
%\label{condM5general1}
 ${\xi^6_5} ({\xi^2_4}^2 + {\xi^2_3}^2)(\xi^2_4\xi^1_3-\xi^2_3 \xi^1_4) \neq 0 .$
 %\end{equation}
\subsection{Case $ \xi^2_1  \neq 0 .$}
\subsubsection{Case $ \xi^2_1 \xi^2_4  \neq 0 .$}
%{\fontsize{8}{10}\selectfont
$$ J = \begin{pmatrix}
*&
*&
*&
*&
0&
0\\
\boxed{\xi^2_1}&
\boxed{\xi^2_2}&
\boxed{\xi^2_3}&
\boxed{\xi^2_4}&
0&
0\\
*&
*&
\boxed{\xi^3_3}&
\boxed{\xi^3_4}&
0&
0\\
*&
*&
*&
*&
0&
0\\
*&
*&
*&
*&
\boxed{\xi^5_5}&
 - \frac{{\xi^5_5}^2 + 1}{{\xi^6_5}}\\
\boxed{\xi^6_1}&
\boxed{\xi^6_2}&
\boxed{\xi^6_3}&
\boxed{\xi^6_4}&
\boxed{\xi^6_5}&
 - {\xi^5_5}\end{pmatrix}$$
%}
where
the parameters are subject to the condition
\begin{equation}
\label{condM5general2}
\xi^2_1 \xi^6_5 \xi^2_4 C_2\neq 0
 \end{equation}
 with $C_2=(({\xi^6_5} {\xi^3_4} {\xi^2_4}
 + {\xi^6_5} {\xi^2_4} {\xi^2_1} + {\xi^5_5} {\xi^3_4} {\xi^2_3}) {\xi^6_5} + ({\xi^5_5}^2 + 1)
 ({\xi^3_3} + {\xi^2_2}) {\xi^2_3} + ({\xi^2_4} {\xi^2_2} + {\xi^2_3} {\xi^2_1} + {\xi^3_3} {\xi^2_4}) {\xi^6_5} {\xi^5_5}) {\xi^2_4} - (({\xi^3_3} {\xi^2_4} {\xi^2_2} - {\xi^3_3} {\xi^2_3} {\xi^2_1} - {\xi^2_4}) {\xi^2_4} - ({\xi^2_4} {\xi^2_2} - {\xi^2_3} {\xi^2_1}) {\xi^3_4} {\xi^2_3}) {\xi^6_5}.$
 The starred $J^i_j$'s are smooth rational functions under condition
(\ref{condM5general2}). However we do not give them here, as some are huge
and we do not have explicit use of them in the rest of the paper.
We refer instead to (\cite{sc1}, pp. 133-136).
\subsubsection{Case $ {\xi^2_1} \xi^2_3 \neq 0 ,\xi^2_4 = 0 .$}
\label{2123neq0240}
%{\fontsize{8}{10}\selectfont
$$ J = \begin{pmatrix}
*&
*&
- \frac{({\xi^3_3} + {\xi^2_2}) {\xi^2_3}}{{\xi^2_1}}&
- \frac{{\xi^3_4} {\xi^2_3}}{{\xi^2_1}}&
0&
0\\
\boxed{\xi^2_1}&
\boxed{\xi^2_2}&
\boxed{\xi^2_3}&
0&
0&
0\\
*&
*&
\boxed{\xi^3_3}&
\boxed{\xi^3_4}&
0&
0\\
*&
*&
*&
*&
0&
0\\
*&
*&
*&
*&
\boxed{\xi^5_5}&
 - \frac{{\xi^5_5}^2 + 1}{{\xi^6_5}}\\
\boxed{\xi^6_1}&
\boxed{\xi^6_2}&
\boxed{\xi^6_3}&
\boxed{\xi^6_4}&
\boxed{\xi^6_5}&
 - {\xi^5_5}\end{pmatrix}$$
%}
where
{\fontsize{6}{10}\selectfont
$  \mathbf{J^1_1}=( - (({\xi^5_5}^2 + 1) ({\xi^3_3} + {\xi^2_2}) + {\xi^6_5} {\xi^2_2} {\xi^2_1} +
 ({\xi^3_4} {\xi^2_2} + {\xi^3_3} {\xi^2_1} + ({\xi^3_4} + {\xi^2_1}) {\xi^5_5}) {\xi^6_5}))/({\xi^6_5} {\xi^2_1});$
 \myquad
$  \mathbf{J^1_2}=( - (({\xi^3_3} + {\xi^2_2}) {\xi^2_2} - {\xi^3_4} {\xi^2_1} - ({\xi^3_3} + {\xi^2_2}) {\xi^5_5} - ({\xi^3_4} + {\xi^2_1}) {\xi^6_5}))/{\xi^2_1};$
 \myquad
$  \mathbf{J^3_1}=(({\xi^3_4} {\xi^2_2} + {\xi^3_3} {\xi^2_1} + ({\xi^3_4} + {\xi^2_1}) {\xi^5_5}) {\xi^6_5} + ({\xi^5_5}^2 + 1) ({\xi^3_3} + {\xi^2_2}))/({\xi^6_5} {\xi^2_3});$
 \myquad
$  \mathbf{J^3_2}=( - ({\xi^6_5} {\xi^3_4} + {\xi^6_5} {\xi^2_1} + {\xi^5_5} {\xi^3_3} + {\xi^5_5} {\xi^2_2} + {\xi^3_4} {\xi^2_1} - {\xi^3_3} {\xi^2_2} + 1))/{\xi^2_3};$
 \myquad
$  \mathbf{J^4_1}=(((({\xi^3_3} + {\xi^2_2}) {\xi^5_5} + 2 {\xi^6_5} {\xi^2_1} + 2 {\xi^5_5}^2 + 2
) ({\xi^3_3} + {\xi^2_2}) {\xi^5_5} {\xi^2_1} + {\xi^6_5}^2 {\xi^3_4} {\xi^2_1}^2 + {\xi^6_5}
^2 {\xi^2_1}^3 + {\xi^6_5} {\xi^5_5}^2 {\xi^3_4}^2 + 2 {\xi^6_5} {\xi^5_5}^2 {\xi^3_4} {\xi^2_1} + {\xi^6_5} {\xi^5_5}^2 {\xi^2_1}^2 + 2 {\xi^6_5} {\xi^5_5} {\xi^3_4}^2 {\xi^2_2} +
{\xi^6_5} {\xi^5_5} {\xi^3_4} {\xi^3_3} {\xi^2_1} + 3 {\xi^6_5} {\xi^5_5} {\xi^3_4} {\xi^2_2} {\xi^2_1} + 2 {\xi^5_5}^3 {\xi^3_4} {\xi^3_3} + 2 {\xi^5_5}^3 {\xi^3_4} {\xi^2_2} + 2 {\xi^5_5}^2
 {\xi^3_4} {\xi^3_3} {\xi^2_2} + 2 {\xi^5_5}^2 {\xi^3_4} {\xi^2_2}^2 + 2 {\xi^5_5} {\xi^3_4}
{\xi^3_3} + 2 {\xi^5_5} {\xi^3_4} {\xi^2_2}) {\xi^6_5} + ({\xi^5_5}^2 + 1)^2 ({\xi^3_3} + {\xi^2_2})^2 + (({\xi^3_3} + {\xi^2_2}) {\xi^2_1} + 2 {\xi^3_4} {\xi^2_2}) ({\xi^3_3} + {\xi^2_2}
) {\xi^6_5} + ({\xi^3_4}^2 {\xi^2_2}^2 + {\xi^3_4} {\xi^3_3} {\xi^2_2} {\xi^2_1} + {\xi^2_1}^
2 + ({\xi^2_2}^2 + {\xi^2_1}^2) {\xi^3_4} {\xi^2_1}) {\xi^6_5}^2)/({\xi^6_5}^2 {\xi^3_4}
{\xi^2_3} {\xi^2_1});$\\
$  J^4_2=( - ({\xi^6_5}^2 {\xi^5_5} {\xi^3_4}^2 + 2 {\xi^6_5}^2 {\xi^5_5} {\xi^3_4} {\xi^2_1} + {\xi^6_5}^2 {\xi^5_5} {\xi^2_1}^2 + {\xi^6_5}^2 {\xi^3_4}^2 {\xi^2_2} - {\xi^6_5}^
2 {\xi^2_2} {\xi^2_1}^2 + 2 {\xi^6_5} {\xi^5_5}^2 {\xi^3_4} {\xi^3_3} + 2 {\xi^6_5} {\xi^5_5}
^2 {\xi^3_4} {\xi^2_2} + 2 {\xi^6_5} {\xi^5_5}^2 {\xi^3_3} {\xi^2_1} + 2 {\xi^6_5} {\xi^5_5}^
2 {\xi^2_2} {\xi^2_1} + {\xi^6_5} {\xi^5_5} {\xi^3_4}^2 {\xi^2_1} + {\xi^6_5} {\xi^5_5} {\xi^3_4}
 {\xi^2_1}^2 - 2 {\xi^6_5} {\xi^5_5} {\xi^3_3} {\xi^2_2} {\xi^2_1} - 2 {\xi^6_5} {\xi^5_5} {\xi^2_2}^2 {\xi^2_1} + {\xi^6_5} {\xi^3_4}^2 {\xi^2_2} {\xi^2_1} - {\xi^6_5} {\xi^3_4} {\xi^3_3} {\xi^2_2}^2 + {\xi^6_5} {\xi^3_4} {\xi^3_3} - {\xi^6_5} {\xi^3_4} {\xi^2_2}^3 - {\xi^6_5} {\xi^3_4} {\xi^2_2} {\xi^2_1}^2 + {\xi^6_5} {\xi^3_4} {\xi^2_2} + {\xi^5_5}^3 {\xi^3_3}^2 + 2 {\xi^5_5}^3 {\xi^3_3} {\xi^2_2} + {\xi^5_5}^3 {\xi^2_2}^2 + {\xi^5_5}^2 {\xi^3_4} {\xi^3_3} {\xi^2_1} + {\xi^5_5}^2 {\xi^3_4} {\xi^2_2} {\xi^2_1} - {\xi^5_5}^2 {\xi^3_3}^2 {\xi^2_2} - 2 {\xi^5_5}^2 {\xi^3_3} {\xi^2_2}^2 - {\xi^5_5}^2 {\xi^2_2}^3 + {\xi^5_5} {\xi^3_3}^2 + 2 {\xi^5_5} {\xi^3_3} {\xi^2_2} + {\xi^5_5} {\xi^2_2}^2 + {\xi^3_4} {\xi^3_3} {\xi^2_1} + {\xi^3_4} {\xi^2_2} {\xi^2_1} - {\xi^3_3}^2 {\xi^2_2} - 2 {\xi^3_3} {\xi^2_2}^2 - {\xi^2_2}^3))/({\xi^6_5}
{\xi^3_4} {\xi^2_3} {\xi^2_1});$
 \myquad
$  \mathbf{J^4_3}=(({\xi^6_5} {\xi^3_4} {\xi^2_1} + {\xi^6_5} {\xi^2_1}^2 + {\xi^5_5} {\xi^3_4} {\xi^3_3} + {\xi^5_5} {\xi^3_4} {\xi^2_2} + 2 ({\xi^3_3} + {\xi^2_2}) {\xi^5_5} {\xi^2_1}) {\xi^6_5} +
 ({\xi^5_5}^2 + 1) ({\xi^3_3} + {\xi^2_2})^2 + ({\xi^2_2}^2 + {\xi^2_1}^2 + {\xi^3_3} {\xi^2_2}) {\xi^6_5} {\xi^3_4})/({\xi^6_5} {\xi^3_4} {\xi^2_1});$
 \myquad
$  \mathbf{J^4_4}=(({\xi^5_5}^2 + 1) ({\xi^3_3} + {\xi^2_2}) + ({\xi^5_5} {\xi^3_4} + {\xi^5_5} {\xi^2_1} + {\xi^3_4} {\xi^2_2}) {\xi^6_5})/({\xi^6_5} {\xi^2_1});$
 \myquad
$  \mathbf{J^5_1}=( - ((((({\xi^6_2} {\xi^2_1} - {\xi^6_1} {\xi^5_5}) {\xi^6_5} {\xi^2_3} + ({\xi^5_5}^2 + 1) ({\xi^3_3} + {\xi^2_2}) {\xi^6_3}) {\xi^2_1} + ({\xi^3_4} {\xi^2_2} + {\xi^3_3} {\xi^2_1} + ({\xi^3_4} + {\xi^2_1}) {\xi^5_5}) ({\xi^6_3} {\xi^2_1} - {\xi^6_1} {\xi^2_3}) {\xi^6_5})
{\xi^3_4} + ({\xi^3_4}^2 {\xi^2_2}^2 + {\xi^3_4} {\xi^3_3} {\xi^2_2} {\xi^2_1} + {\xi^2_1}^2
+ ({\xi^2_2}^2 + {\xi^2_1}^2) {\xi^3_4} {\xi^2_1}) {\xi^6_5} {\xi^6_4} - (({\xi^5_5}^2 + 1)
 ({\xi^3_3} + {\xi^2_2}) + {\xi^6_5} {\xi^2_2} {\xi^2_1}) {\xi^6_1} {\xi^3_4} {\xi^2_3}) {\xi^6_5}
 + (((({\xi^3_3} + {\xi^2_2}) {\xi^5_5} + 2 {\xi^6_5} {\xi^2_1} + 2 {\xi^5_5}^2 + 2) ({\xi^3_3} + {\xi^2_2}) {\xi^5_5} {\xi^2_1} + {\xi^6_5}^2 {\xi^3_4} {\xi^2_1}^2 + {\xi^6_5}^2 {\xi^2_1}^3 + {\xi^6_5} {\xi^5_5}^2 {\xi^3_4}^2 + 2 {\xi^6_5} {\xi^5_5}^2 {\xi^3_4} {\xi^2_1} +
{\xi^6_5} {\xi^5_5}^2 {\xi^2_1}^2 + 2 {\xi^6_5} {\xi^5_5} {\xi^3_4}^2 {\xi^2_2} + {\xi^6_5}
{\xi^5_5} {\xi^3_4} {\xi^3_3} {\xi^2_1} + 3 {\xi^6_5} {\xi^5_5} {\xi^3_4} {\xi^2_2} {\xi^2_1} + 2
{\xi^5_5}^3 {\xi^3_4} {\xi^3_3} + 2 {\xi^5_5}^3 {\xi^3_4} {\xi^2_2} + 2 {\xi^5_5}^2 {\xi^3_4}
 {\xi^3_3} {\xi^2_2} + 2 {\xi^5_5}^2 {\xi^3_4} {\xi^2_2}^2 + 2 {\xi^5_5} {\xi^3_4} {\xi^3_3} +
 2 {\xi^5_5} {\xi^3_4} {\xi^2_2}) {\xi^6_5} + ({\xi^5_5}^2 + 1)^2 ({\xi^3_3} + {\xi^2_2})^2
 + (({\xi^3_3} + {\xi^2_2}) {\xi^2_1} + 2 {\xi^3_4} {\xi^2_2}) ({\xi^3_3} + {\xi^2_2}) {\xi^6_5}
) {\xi^6_4}))/({\xi^6_5}^3 {\xi^3_4} {\xi^2_3} {\xi^2_1});$\\
$  J^5_2=(((({\xi^6_5} {\xi^3_4}^2 + 2 {\xi^6_5} {\xi^3_4} {\xi^2_1} + {\xi^6_5} {\xi^2_1}
^2 + 2 {\xi^5_5} {\xi^3_4} {\xi^3_3} + 2 {\xi^5_5} {\xi^3_4} {\xi^2_2} + {\xi^3_4}^2 {\xi^2_1}
 + {\xi^3_4} {\xi^2_1}^2 + 2 ({\xi^5_5} - {\xi^2_2}) ({\xi^3_3} + {\xi^2_2}) {\xi^2_1}) {\xi^5_5} + ({\xi^3_4} + {\xi^2_1}) ({\xi^3_4} - {\xi^2_1}) {\xi^6_5} {\xi^2_2}) {\xi^6_5} + ({\xi^5_5}
^2 + 1) ({\xi^5_5} {\xi^3_3} + {\xi^5_5} {\xi^2_2} + {\xi^3_4} {\xi^2_1} - {\xi^3_3} {\xi^2_2}
- {\xi^2_2}^2) ({\xi^3_3} + {\xi^2_2})) {\xi^6_4} + ({\xi^6_3} {\xi^2_1} - {\xi^6_1} {\xi^2_3})
 ({\xi^3_4} + {\xi^2_1}) {\xi^6_5}^2 {\xi^3_4} + ({\xi^3_4} {\xi^2_2} {\xi^2_1} - {\xi^3_3} {\xi^2_2}^2 + {\xi^3_3} - ({\xi^2_2}^2 + {\xi^2_1}^2 - 1) {\xi^2_2}) {\xi^6_5} {\xi^6_4} {\xi^3_4} - ((({\xi^3_4} {\xi^2_1} - {\xi^3_3} {\xi^2_2} - {\xi^2_2}^2 + ({\xi^3_3} + {\xi^2_2}) {\xi^5_5}) {\xi^6_1} - ({\xi^5_5} - {\xi^2_2}) {\xi^6_2} {\xi^2_1}) {\xi^2_3} - ({\xi^3_4} {\xi^2_1}
- {\xi^3_3} {\xi^2_2} + 1 + ({\xi^3_3} + {\xi^2_2}) {\xi^5_5}) {\xi^6_3} {\xi^2_1}) {\xi^6_5} {\xi^3_4})/({\xi^6_5}^2 {\xi^3_4} {\xi^2_3} {\xi^2_1});$
 \myquad
$  \mathbf{J^5_3}=( - ({\xi^6_5}^2 {\xi^6_4} {\xi^3_4} {\xi^2_1} + {\xi^6_5}^2 {\xi^6_4} {\xi^2_1}
^2 + {\xi^6_5} {\xi^6_4} {\xi^5_5} {\xi^3_4} {\xi^3_3} + {\xi^6_5} {\xi^6_4} {\xi^5_5} {\xi^3_4}
{\xi^2_2} + 2 {\xi^6_5} {\xi^6_4} {\xi^5_5} {\xi^3_3} {\xi^2_1} + 2 {\xi^6_5} {\xi^6_4} {\xi^5_5}
{\xi^2_2} {\xi^2_1} + {\xi^6_5} {\xi^6_4} {\xi^3_4} {\xi^3_3} {\xi^2_2} + {\xi^6_5} {\xi^6_4} {\xi^3_4} {\xi^2_2}^2 + {\xi^6_5} {\xi^6_4} {\xi^3_4} {\xi^2_1}^2 - {\xi^6_5} {\xi^6_3} {\xi^5_5} {\xi^3_4} {\xi^2_1} + {\xi^6_5} {\xi^6_3} {\xi^3_4} {\xi^3_3} {\xi^2_1} + {\xi^6_5} {\xi^6_2} {\xi^3_4} {\xi^2_3} {\xi^2_1} - {\xi^6_5} {\xi^6_1} {\xi^3_4} {\xi^3_3} {\xi^2_3} - {\xi^6_5} {\xi^6_1} {\xi^3_4} {\xi^2_3} {\xi^2_2} + {\xi^6_4} {\xi^5_5}^2 {\xi^3_3}^2 + 2 {\xi^6_4} {\xi^5_5}^2 {\xi^3_3} {\xi^2_2} + {\xi^6_4} {\xi^5_5}^2 {\xi^2_2}^2 + {\xi^6_4} {\xi^3_3}^2 + 2 {\xi^6_4} {\xi^3_3} {\xi^2_2} + {\xi^6_4} {\xi^2_2}^2))/({\xi^6_5}^2 {\xi^3_4} {\xi^2_1});$
 \myquad
$  \mathbf{J^5_4}=( - ({\xi^6_5} {\xi^6_4} {\xi^5_5} {\xi^3_4} + {\xi^6_5} {\xi^6_4} {\xi^3_4} {\xi^2_2} + {\xi^6_5} {\xi^6_3} {\xi^3_4} {\xi^2_1} - {\xi^6_5} {\xi^6_1} {\xi^3_4} {\xi^2_3} + {\xi^6_4}
{\xi^5_5}^2 {\xi^3_3} + {\xi^6_4} {\xi^5_5}^2 {\xi^2_2} + {\xi^6_4} {\xi^3_3} + {\xi^6_4} {\xi^2_2}))/({\xi^6_5}^2 {\xi^2_1});$
}
and the parameters are subject to the conditions
%\begin{equation}
%\label{condM5general3}
$\xi^2_1 \xi^6_5 \xi^2_3 \xi^3_4 \neq 0.$
% \end{equation}
\subsubsection{Case $\xi^2_1 \neq 0, \xi^2_3 =\xi^2_4=0, \xi^2_2+\xi^3_3 \neq 0 .$}
%{\fontsize{8}{10}\selectfont
$$ J = \begin{pmatrix}
-\xi^2_2&
-\frac{{\xi^2_2}^2+1}{\xi^2_1}&
0&
0&
0&
0\\
\boxed{\xi^2_1}&
\boxed{\xi^2_2}&
0&
0&
0&
0\\
\boxed{\xi^3_1}&
- \frac{{\xi^4_1} {\xi^3_4} - {\xi^3_3} {\xi^3_1} + {\xi^3_1} {\xi^2_2}}{{\xi^2_1}}&
\boxed{\xi^3_3}&
\boxed{\xi^3_4}&
0&
0\\
\boxed{\xi^4_1}&
\frac{{\xi^4_1} {\xi^3_4} {\xi^3_3} + {\xi^4_1} {\xi^3_4} {\xi^2_2} + {\xi^3_3}^2 {\xi^3_1} + {\xi^3_1}}{{\xi^3_4} {\xi^2_1}}&
-\frac{{\xi^3_3}^2+1}{\xi^3_4}&
-\xi^3_3&
0&
0\\
*&
*&
*&
*&
*&
*\\
\boxed{\xi^6_1}&
\boxed{\xi^6_2}&
\boxed{\xi^6_3}&
\boxed{\xi^6_4}&
*&
*
 \end{pmatrix}$$
%}
where
{\fontsize{6}{10}\selectfont
$  \mathbf{J^5_1}=({\xi^6_4} {\xi^4_1} {\xi^3_4}^2 {\xi^2_1} + {\xi^6_4} {\xi^4_1} {\xi^3_4} {\xi^2_2}
^2 + {\xi^6_4} {\xi^4_1} {\xi^3_4} {\xi^2_1}^2 + {\xi^6_4} {\xi^4_1} {\xi^3_4} + {\xi^6_4} {\xi^4_1} {\xi^3_3}^2 {\xi^2_1} + {\xi^6_4} {\xi^4_1} {\xi^2_1} + {\xi^6_3} {\xi^3_4}^2 {\xi^3_1}
{\xi^2_1} + {\xi^6_3} {\xi^3_4} {\xi^3_1} {\xi^2_2}^2 + {\xi^6_3} {\xi^3_4} {\xi^3_1} {\xi^2_1}^
2 + {\xi^6_3} {\xi^3_4} {\xi^3_1} + {\xi^6_3} {\xi^3_3}^2 {\xi^3_1} {\xi^2_1} + {\xi^6_3} {\xi^3_1} {\xi^2_1} + {\xi^6_2} {\xi^3_4}^2 {\xi^2_1}^2 + {\xi^6_2} {\xi^3_4} {\xi^2_2}^2 {\xi^2_1}
+ {\xi^6_2} {\xi^3_4} {\xi^2_1}^3 + {\xi^6_2} {\xi^3_4} {\xi^2_1} + {\xi^6_2} {\xi^3_3}^2 {\xi^2_1}^2 + {\xi^6_2} {\xi^2_1}^2 - {\xi^6_1} {\xi^3_4} {\xi^3_3} {\xi^2_2}^2 + {\xi^6_1} {\xi^3_4} {\xi^3_3} {\xi^2_1}^2 - {\xi^6_1} {\xi^3_4} {\xi^3_3} - {\xi^6_1} {\xi^3_4} {\xi^2_2}^3 -
{\xi^6_1} {\xi^3_4} {\xi^2_2} {\xi^2_1}^2 - {\xi^6_1} {\xi^3_4} {\xi^2_2} - 2 {\xi^6_1} {\xi^3_3}
^2 {\xi^2_2} {\xi^2_1} - 2 {\xi^6_1} {\xi^2_2} {\xi^2_1})/
({\xi^3_4}^2 {\xi^2_1}^2 - 2 {\xi^3_4} \xi^2_1({\xi^3_3} {\xi^2_2} -1)+ ({\xi^2_2}^2 +1)( {\xi^3_3}^2 + 1))
;$
\myquad
$  \mathbf{J^5_2}=({\xi^6_4} {\xi^4_1} {\xi^3_4}^3 {\xi^3_3} {\xi^2_1} + {\xi^6_4} {\xi^4_1} {\xi^3_4}
^3 {\xi^2_2} {\xi^2_1} + {\xi^6_4} {\xi^4_1} {\xi^3_4}^2 {\xi^3_3} {\xi^2_2}^2 + {\xi^6_4} {\xi^4_1} {\xi^3_4}^2 {\xi^3_3} {\xi^2_1}^2 + {\xi^6_4} {\xi^4_1} {\xi^3_4}^2 {\xi^3_3} + {\xi^6_4} {\xi^4_1} {\xi^3_4}^2 {\xi^2_2}^3 + {\xi^6_4} {\xi^4_1} {\xi^3_4}^2 {\xi^2_2} {\xi^2_1}^2
 + {\xi^6_4} {\xi^4_1} {\xi^3_4}^2 {\xi^2_2} + {\xi^6_4} {\xi^4_1} {\xi^3_4} {\xi^3_3}^3 {\xi^2_1} + {\xi^6_4} {\xi^4_1} {\xi^3_4} {\xi^3_3}^2 {\xi^2_2} {\xi^2_1} + {\xi^6_4} {\xi^4_1} {\xi^3_4} {\xi^3_3} {\xi^2_1} + {\xi^6_4} {\xi^4_1} {\xi^3_4} {\xi^2_2} {\xi^2_1} + {\xi^6_4} {\xi^3_4}^2
 {\xi^3_3}^2 {\xi^3_1} {\xi^2_1} + {\xi^6_4} {\xi^3_4}^2 {\xi^3_1} {\xi^2_1} + {\xi^6_4} {\xi^3_4} {\xi^3_3}^2 {\xi^3_1} {\xi^2_2}^2 + {\xi^6_4} {\xi^3_4} {\xi^3_3}^2 {\xi^3_1} {\xi^2_1}^2
 + {\xi^6_4} {\xi^3_4} {\xi^3_3}^2 {\xi^3_1} + {\xi^6_4} {\xi^3_4} {\xi^3_1} {\xi^2_2}^2 + {\xi^6_4} {\xi^3_4} {\xi^3_1} {\xi^2_1}^2 + {\xi^6_4} {\xi^3_4} {\xi^3_1} + {\xi^6_4} {\xi^3_3}^4
{\xi^3_1} {\xi^2_1} + 2 {\xi^6_4} {\xi^3_3}^2 {\xi^3_1} {\xi^2_1} + {\xi^6_4} {\xi^3_1} {\xi^2_1}
 - {\xi^6_3} {\xi^4_1} {\xi^3_4}^4 {\xi^2_1} - {\xi^6_3} {\xi^4_1} {\xi^3_4}^3 {\xi^2_2}^2 -
{\xi^6_3} {\xi^4_1} {\xi^3_4}^3 {\xi^2_1}^2 - {\xi^6_3} {\xi^4_1} {\xi^3_4}^3 - {\xi^6_3} {\xi^4_1} {\xi^3_4}^2 {\xi^3_3}^2 {\xi^2_1} - {\xi^6_3} {\xi^4_1} {\xi^3_4}^2 {\xi^2_1} - {\xi^6_3} {\xi^3_4}^3 {\xi^3_3} {\xi^3_1} {\xi^2_1} + {\xi^6_3} {\xi^3_4}^3 {\xi^3_1} {\xi^2_2} {\xi^2_1} - {\xi^6_3} {\xi^3_4}^2 {\xi^3_3} {\xi^3_1} {\xi^2_2}^2 - {\xi^6_3} {\xi^3_4}^2 {\xi^3_3}
{\xi^3_1} {\xi^2_1}^2 - {\xi^6_3} {\xi^3_4}^2 {\xi^3_3} {\xi^3_1} + {\xi^6_3} {\xi^3_4}^2 {\xi^3_1} {\xi^2_2}^3 + {\xi^6_3} {\xi^3_4}^2 {\xi^3_1} {\xi^2_2} {\xi^2_1}^2 + {\xi^6_3} {\xi^3_4}^2 {\xi^3_1} {\xi^2_2} - {\xi^6_3} {\xi^3_4} {\xi^3_3}^3 {\xi^3_1} {\xi^2_1} + {\xi^6_3} {\xi^3_4} {\xi^3_3}^2 {\xi^3_1} {\xi^2_2} {\xi^2_1} - {\xi^6_3} {\xi^3_4} {\xi^3_3} {\xi^3_1} {\xi^2_1}
 + {\xi^6_3} {\xi^3_4} {\xi^3_1} {\xi^2_2} {\xi^2_1} + 2 {\xi^6_2} {\xi^3_4}^3 {\xi^2_2} {\xi^2_1}^2 - {\xi^6_2} {\xi^3_4}^2 {\xi^3_3} {\xi^2_2}^2 {\xi^2_1} + {\xi^6_2} {\xi^3_4}^2 {\xi^3_3} {\xi^2_1}^3 - {\xi^6_2} {\xi^3_4}^2 {\xi^3_3} {\xi^2_1} + {\xi^6_2} {\xi^3_4}^2 {\xi^2_2}^
3 {\xi^2_1} + {\xi^6_2} {\xi^3_4}^2 {\xi^2_2} {\xi^2_1}^3 + {\xi^6_2} {\xi^3_4}^2 {\xi^2_2}
{\xi^2_1} - {\xi^6_1} {\xi^3_4}^3 {\xi^2_2}^2 {\xi^2_1} - {\xi^6_1} {\xi^3_4}^3 {\xi^2_1} -
{\xi^6_1} {\xi^3_4}^2 {\xi^2_2}^4 - {\xi^6_1} {\xi^3_4}^2 {\xi^2_2}^2 {\xi^2_1}^2 - 2 {\xi^6_1} {\xi^3_4}^2 {\xi^2_2}^2 - {\xi^6_1} {\xi^3_4}^2 {\xi^2_1}^2 - {\xi^6_1} {\xi^3_4}^2
- {\xi^6_1} {\xi^3_4} {\xi^3_3}^2 {\xi^2_2}^2 {\xi^2_1} - {\xi^6_1} {\xi^3_4} {\xi^3_3}^2 {\xi^2_1} - {\xi^6_1} {\xi^3_4} {\xi^2_2}^2 {\xi^2_1} - {\xi^6_1} {\xi^3_4} {\xi^2_1})/
({\xi^3_4} {\xi^2_1}
({\xi^3_4}^2 {\xi^2_1}^2 - 2 {\xi^3_4} \xi^2_1({\xi^3_3} {\xi^2_2} -1)+ ({\xi^2_2}^2 +1)( {\xi^3_3}^2 + 1)))
;$
\myquad
$  \mathbf{J^5_3}=( - {\xi^6_4} {\xi^3_4}^2 {\xi^3_3}^2 {\xi^2_1} - {\xi^6_4} {\xi^3_4}^2 {\xi^2_1} - {\xi^6_4} {\xi^3_4} {\xi^3_3}^2 {\xi^2_2}^2 - {\xi^6_4} {\xi^3_4} {\xi^3_3}^2 {\xi^2_1}^
2 - {\xi^6_4} {\xi^3_4} {\xi^3_3}^2 - {\xi^6_4} {\xi^3_4} {\xi^2_2}^2 - {\xi^6_4} {\xi^3_4} {\xi^2_1}^2 - {\xi^6_4} {\xi^3_4} - {\xi^6_4} {\xi^3_3}^4 {\xi^2_1} - 2 {\xi^6_4} {\xi^3_3}^2
{\xi^2_1} - {\xi^6_4} {\xi^2_1} + {\xi^6_3} {\xi^3_4}^3 {\xi^3_3} {\xi^2_1} + {\xi^6_3} {\xi^3_4}
^3 {\xi^2_2} {\xi^2_1} + 2 {\xi^6_3} {\xi^3_4}^2 {\xi^3_3} {\xi^2_1}^2 + {\xi^6_3} {\xi^3_4}
{\xi^3_3}^3 {\xi^2_1} - {\xi^6_3} {\xi^3_4} {\xi^3_3}^2 {\xi^2_2} {\xi^2_1} + {\xi^6_3} {\xi^3_4} {\xi^3_3} {\xi^2_1} - {\xi^6_3} {\xi^3_4} {\xi^2_2} {\xi^2_1})/
 ({\xi^3_4}
({\xi^3_4}^2 {\xi^2_1}^2 - 2 {\xi^3_4} \xi^2_1({\xi^3_3} {\xi^2_2} -1)+ ({\xi^2_2}^2 +1)( {\xi^3_3}^2 + 1)))
 ;$
\myquad
$  \mathbf{J^5_4}=( - {\xi^6_4} {\xi^3_4}^2 {\xi^3_3} {\xi^2_1} + {\xi^6_4} {\xi^3_4}^2 {\xi^2_2}
{\xi^2_1} - 2 {\xi^6_4} {\xi^3_4} {\xi^3_3} {\xi^2_2}^2 - 2 {\xi^6_4} {\xi^3_4} {\xi^3_3} - {\xi^6_4} {\xi^3_3}^3 {\xi^2_1} - {\xi^6_4} {\xi^3_3}^2 {\xi^2_2} {\xi^2_1} - {\xi^6_4} {\xi^3_3}
{\xi^2_1} - {\xi^6_4} {\xi^2_2} {\xi^2_1} + {\xi^6_3} {\xi^3_4}^3 {\xi^2_1} + {\xi^6_3} {\xi^3_4}
^2 {\xi^2_2}^2 + {\xi^6_3} {\xi^3_4}^2 {\xi^2_1}^2 + {\xi^6_3} {\xi^3_4}^2 + {\xi^6_3} {\xi^3_4} {\xi^3_3}^2 {\xi^2_1} + {\xi^6_3} {\xi^3_4} {\xi^2_1})/
({\xi^3_4}^2 {\xi^2_1}^2 - 2 {\xi^3_4} \xi^2_1({\xi^3_3} {\xi^2_2} -1)+ ({\xi^2_2}^2 +1)( {\xi^3_3}^2 + 1))
 ;$
 \myquad
$  \mathbf{J^5_5}=( - {\xi^3_4}^2 {\xi^2_2} {\xi^2_1} + {\xi^3_4} {\xi^3_3} {\xi^2_2}^2 - {\xi^3_4}
 {\xi^3_3} {\xi^2_1}^2 + {\xi^3_4} {\xi^3_3} + {\xi^3_3}^2 {\xi^2_2} {\xi^2_1} + {\xi^2_2} {\xi^2_1})/
  ({\xi^3_4}^2 {\xi^2_1} + {\xi^3_4} ({\xi^2_2}^2 + {\xi^2_1}^2 + 1) +
({\xi^3_3}^2 +1)  {\xi^2_1});$
\myquad
$  \mathbf{J^5_6}=({\xi^3_4}^2 {\xi^2_2}^2 + {\xi^3_4}^2 + 2 {\xi^3_4} {\xi^3_3} {\xi^2_2} {\xi^2_1} + 2 {\xi^3_4} {\xi^2_1} + {\xi^3_3}^2 {\xi^2_1}^2 + {\xi^2_1}^2)/
  ({\xi^3_4}^2 {\xi^2_1} + {\xi^3_4} ({\xi^2_2}^2 + {\xi^2_1}^2 + 1) +
({\xi^3_3}^2 +1)  {\xi^2_1});$
$  \mathbf{J^6_5}=( - {\xi^3_4}^2 {\xi^2_1}^2 + 2 {\xi^3_4} {\xi^3_3} {\xi^2_2} {\xi^2_1} - 2 {\xi^3_4} {\xi^2_1} - {\xi^3_3}^2 {\xi^2_2}^2 - {\xi^3_3}^2 - {\xi^2_2}^2 - 1)/
  ({\xi^3_4}^2 {\xi^2_1} + {\xi^3_4} ({\xi^2_2}^2 + {\xi^2_1}^2 + 1) +
({\xi^3_3}^2 +1)  {\xi^2_1});$
\myquad
$  \mathbf{J^6_6}=({\xi^3_4}^2 {\xi^2_2} {\xi^2_1} - {\xi^3_4} {\xi^3_3} {\xi^2_2}^2 + {\xi^3_4} {\xi^3_3} {\xi^2_1}^2 - {\xi^3_4} {\xi^3_3} - {\xi^3_3}^2 {\xi^2_2} {\xi^2_1} - {\xi^2_2} {\xi^2_1})/
  ({\xi^3_4}^2 {\xi^2_1} + {\xi^3_4} ({\xi^2_2}^2 + {\xi^2_1}^2 + 1) +
({\xi^3_3}^2 +1)  {\xi^2_1});$
}
and the parameters are subject to the condition
\begin{equation}
\label{condM5general4}
\xi^2_1 \xi^3_4 (\xi^3_3 +\xi^2_2)
 ({\xi^3_4}^2 {\xi^2_1} + {\xi^3_4} ({\xi^2_2}^2 + {\xi^2_1}^2 +1) + ({\xi^3_3}^2 +1)
 {\xi^2_1})\neq 0 .
 \end{equation}
 (Note that
 ${\xi^3_4}^2 {\xi^2_1}^2 - 2{\xi^3_4}\xi^2_1(\xi^3_3\xi^2_2-1)+
 ({\xi^2_2}^2 +1) ({\xi^3_3}^2 +1) \neq 0$ is automatic from $\xi^3_3+\xi^2_2 \neq 0.$)
\subsubsection{Case $ {\xi^2_1}\neq 0, {\xi^2_3} = {\xi^2_4} =0, {\xi^3_3}= -{\xi^2_2}, {\xi^3_4}=- {\xi^2_1}. $}
%{\fontsize{4}{6}\selectfont
%{\fontsize{8}{10}\selectfont
$$ J = \begin{pmatrix}
0&
-\frac{1}{\xi^2_1}&
0&
0&
0&
0\\
{\xi^2_1}&
0&
0&
0&
0&
0\\
{\xi^3_1}&
{\xi^4_1}&
0&
 - {\xi^2_1}&
0&
0\\
{\xi^4_1}&
 - {\xi^3_1}&
\frac{1}{\xi^2_1}&
0&
0&
0\\
 \frac{- {\xi^6_4} {\xi^4_1} - {\xi^6_3} {\xi^3_1} - {\xi^6_2} {\xi^2_1} + {\xi^6_1} {\xi^5_5}
 }{\xi^6_5}&
\frac{{\xi^6_4} {\xi^3_1} {\xi^2_1} - {\xi^6_3} {\xi^4_1} {\xi^2_1} + {\xi^6_2} {\xi^5_5} {\xi^2_1} +{\xi^6_1}}{{\xi^6_5} {\xi^2_1}}&
\frac{ - {\xi^6_4} + {\xi^6_3} {\xi^5_5} {\xi^2_1}}{{\xi^6_5} {\xi^2_1}}&
\frac{{\xi^6_4} {\xi^5_5} + {\xi^6_3} {\xi^2_1}}{{\xi^6_5}}&
{\xi^5_5}&
 - \frac{{\xi^5_5}^2 + 1}{{\xi^6_5}}\\
{\xi^6_1}&
{\xi^6_2}&
{\xi^6_3}&
{\xi^6_4}&
{\xi^6_5}&
 - {\xi^5_5}\end{pmatrix}$$
%}
where
%{\fontsize{6}{10}\selectfont
%$  J^1_1=0;$\\
%$  J^1_2=( - 1)/{\xi^2_1};$\\
%$  J^1_3=0;$\\
%$  J^1_4=0;$\\
%$  J^1_5=0;$\\
%$  J^1_6=0;$\\
%$  J^2_1={\xi^2_1};$\\
%$  J^2_2=0;$\\
%$  J^2_3=0;$\\
%$  J^2_4=0;$\\
%$  J^2_5=0;$\\
%$  J^2_6=0;$\\
%$  J^3_1={\xi^3_1};$\\
%$  J^3_2={\xi^4_1};$\\
%$  J^3_3=0;$\\
%$  J^3_4= - {\xi^2_1};$\\
%$  J^3_5=0;$\\
%$  J^3_6=0;$\\
%$  J^4_1={\xi^4_1};$\\
%$  J^4_2= - {\xi^3_1};$\\
%$  J^4_3=1/{\xi^2_1};$\\
%$  J^4_4=0;$\\
%$  J^4_5=0;$\\
%$  J^4_6=0;$\\
%$  J^5_1=( - {\xi^6_4} {\xi^4_1} - {\xi^6_3} {\xi^3_1} - {\xi^6_2} {\xi^2_1} + {\xi^6_1} {%\xi^5_5})/{\xi^6_5};$\\
%$  J^5_2=({\xi^6_4} {\xi^3_1} {\xi^2_1} - {\xi^6_3} {\xi^4_1} {\xi^2_1} + {\xi^6_2} {\xi^5_5}
%{\xi^2_1} + {\xi^6_1})/({\xi^6_5} {\xi^2_1});$\\
%$  J^5_3=( - {\xi^6_4} + {\xi^6_3} {\xi^5_5} {\xi^2_1})/({\xi^6_5} {\xi^2_1});$\\
%$  J^5_4=({\xi^6_4} {\xi^5_5} + {\xi^6_3} {\xi^2_1})/{\xi^6_5};$\\
%$  J^5_5={\xi^5_5};$\\
%$  J^5_6= - ({\xi^5_5}^2 + 1)/{\xi^6_5};$\\
%$  J^6_1={\xi^6_1};$\\
%$  J^6_2={\xi^6_2};$\\
%$  J^6_3={\xi^6_3};$\\
%$  J^6_4={\xi^6_4};$\\
%$  J^6_5={\xi^6_5};$\\
%$  J^6_6= - {\xi^5_5};$\\
%}
%\begin{equation}
%\label{condM5general5}
${\xi^2_1}=\pm 1.$
 %\end{equation}
\subsubsection{Case $ {\xi^2_1}\neq 0, {\xi^2_3} = {\xi^2_4} =0, {\xi^3_3}= -{\xi^2_2}, {\xi^3_4}\neq - {\xi^2_1}. $}
%{\fontsize{8}{10}\selectfont
$$ J = \begin{pmatrix}
-\xi^2_2&
-\frac{{\xi^2_2}^2+1}{\xi^2_1}&
0&
0&
0&
0\\
\boxed{\xi^2_1}&
\boxed{\xi^2_2}&
0&
0&
0&
0\\
\boxed{\xi^3_1}&
- \frac{{\xi^4_1} {\xi^3_4} - 2{\xi^3_1} {\xi^2_2}}{{\xi^2_1}}&
-{\xi^2_2}&
\boxed{\xi^3_4}&
0&
0\\
\boxed{\xi^4_1}&
\frac{{\xi^3_1} ({\xi^2_2}^2 + 1)}{{\xi^3_4} {\xi^2_1})}&
- \frac{{\xi^2_2}^2 + 1}{{\xi^3_4}}&
{\xi^2_2}&
0&
0\\
*&
*&
*&
*&
\frac{{\xi^2_2} ( - {\xi^3_4} + {\xi^2_1})}{{\xi^3_4} + {\xi^2_1}}&
\frac{{\xi^3_4}^2 {\xi^2_2}^2 + {\xi^3_4}^2 - 2 {\xi^3_4} {\xi^2_2}^2 {\xi^2_1} +
2 {\xi^3_4} {\xi^2_1} + {\xi^2_2}^2 {\xi^2_1}^2 + {\xi^2_1}^2}{{\xi^3_4}^2 {\xi^2_1} + {\xi^3_4} {\xi^2_2}^2 + {\xi^3_4} {\xi^2_1}^2 + {\xi^3_4} + {\xi^2_2}^2 {\xi^2_1} + {\xi^2_1}}\\
\boxed{\xi^6_1}&
\boxed{\xi^6_2}&
\boxed{\xi^6_3}&
\boxed{\xi^6_4}&
- \frac{{\xi^3_4} {\xi^2_1} + {\xi^2_2}^2 + 1}{{\xi^3_4} + {\xi^2_1}}&
\frac{{\xi^2_2} ({\xi^3_4} - {\xi^2_1})}{{\xi^3_4} + {\xi^2_1}}
 \end{pmatrix}$$
%}
where
{\fontsize{6}{10}\selectfont
$  \mathbf{J^5_1}=({\xi^6_4} {\xi^4_1} {\xi^3_4} + {\xi^6_4} {\xi^4_1} {\xi^2_1} + {\xi^6_3} {\xi^3_4}
{\xi^3_1} + {\xi^6_3} {\xi^3_1} {\xi^2_1} + {\xi^6_2} {\xi^3_4} {\xi^2_1} + {\xi^6_2} {\xi^2_1}^2
 - 2 {\xi^6_1} {\xi^2_2} {\xi^2_1})/({\xi^3_4} {\xi^2_1} + {\xi^2_2}^2 + 1);$
 \myquad
$  \mathbf{J^5_2}=({\xi^6_4} {\xi^3_4} {\xi^3_1} {\xi^2_2}^2 + {\xi^6_4} {\xi^3_4} {\xi^3_1} + {\xi^6_4} {\xi^3_1} {\xi^2_2}^2 {\xi^2_1} + {\xi^6_4} {\xi^3_1} {\xi^2_1} - {\xi^6_3} {\xi^4_1} {\xi^3_4}^3 - {\xi^6_3} {\xi^4_1} {\xi^3_4}^2 {\xi^2_1} + 2 {\xi^6_3} {\xi^3_4}^2 {\xi^3_1} {\xi^2_2}
 + 2 {\xi^6_3} {\xi^3_4} {\xi^3_1} {\xi^2_2} {\xi^2_1} + 2 {\xi^6_2} {\xi^3_4}^2 {\xi^2_2} {\xi^2_1} - {\xi^6_1} {\xi^3_4}^2 {\xi^2_2}^2 - {\xi^6_1} {\xi^3_4}^2 - {\xi^6_1} {\xi^3_4} {\xi^2_2}^2 {\xi^2_1} - {\xi^6_1} {\xi^3_4} {\xi^2_1})/({\xi^3_4} {\xi^2_1} ({\xi^3_4} {\xi^2_1} + {\xi^2_2}^2 + 1));$
 \myquad
$  \mathbf{J^5_3}=( - {\xi^6_4} {\xi^3_4} {\xi^2_2}^2 - {\xi^6_4} {\xi^3_4} - {\xi^6_4} {\xi^2_2}^2
 {\xi^2_1} - {\xi^6_4} {\xi^2_1} - 2 {\xi^6_3} {\xi^3_4} {\xi^2_2} {\xi^2_1})/({\xi^3_4} ({\xi^3_4} {\xi^2_1} + {\xi^2_2}^2 + 1));$\myquad
$  \mathbf{J^5_4}=({\xi^3_4} (2 {\xi^6_4} {\xi^2_2} + {\xi^6_3} {\xi^3_4} + {\xi^6_3} {\xi^2_1}))/({\xi^3_4} {\xi^2_1} + {\xi^2_2}^2 + 1);$
}
and the parameters are subject to the condition
%\begin{equation}
%\label{condM5general6}
${\xi^2_1} \xi^3_4 (\xi^3_4+\xi^2_1)({\xi^2_2}^2 +\xi^3_4\xi^2_1+1) \neq 0.
$
 %\end{equation}
\subsubsection{Conclusions for the case $\xi^2_1\neq 0. $}
\label{xi(2,1)}
In each of the 5 subcases in the case $\xi^2_1 \neq 0,$ after completing a set of
common steps, one is left with solving the two equations $14|5$ and $14|6$ in the
14 variables
$\xi^2_1,\xi^2_2,\xi^2_3,\xi^2_4,\xi^3_1,\xi^3_3,\xi^3_4,\xi^4_1,\xi^5_5,
\xi^6_1,\xi^6_2,\xi^6_3,\xi^6_4,\xi^6_5$ in the open subset
$\xi^6_5 (\xi^3_4\xi^2_1 -\xi^3_1\xi^2_4) \neq 0$ of $\Rmath^{14}.$
That is, the initial system comprised of all the torsion equations
and the equation $J^2=-1$
 in $\Rmath^{36}$ is reduced after the common steps to  a system equivalent to the
 2 mentioned equations, which reads
 \begin{equation}
 \label{f,g}
 \begin{cases}
 f =0\\ g=0
 \end{cases}
 \end{equation}
 where :
$ f=\xi^6_5\xi^5_5\xi^3_4{\xi^2_1}^2
 +\xi^6_5\xi^5_5{\xi^2_1}^3
 +\xi^6_5\xi^4_1\xi^3_4\xi^2_4{\xi^2_1}
 -\xi^6_5\xi^4_1\xi^3_1{\xi^2_4}^2
 +\xi^6_5\xi^3_4\xi^2_2{\xi^2_1}^2
 +\xi^6_5\xi^3_3\xi^3_1\xi^2_4{\xi^2_1}
 +\xi^6_5\xi^3_3{\xi^2_1}^3
 -\xi^6_5{\xi^3_1}^2\xi^2_4{\xi^2_3}
 -\xi^6_5{\xi^3_1}\xi^2_4{\xi^2_2}\xi^2_1
 -\xi^6_5{\xi^3_1}\xi^2_3{\xi^2_1}^2
 +{\xi^5_5}^2\xi^3_3{\xi^2_1}^2
 +{\xi^5_5}^2\xi^2_2{\xi^2_1}^2
 +\xi^3_3{\xi^2_1}^2
 +\xi^2_2{\xi^2_1}^2$
 and
 $g=
 \xi^6_5\xi^3_4{\xi^2_1}^2
 +\xi^6_5{\xi^2_1}^3
 +{\xi^5_5}\xi^3_3{\xi^2_1}^2
 +{\xi^5_5}\xi^2_2{\xi^2_1}^2
 -\xi^4_1\xi^3_4\xi^2_3{\xi^2_1}
 +\xi^4_1\xi^3_1\xi^2_4{\xi^2_3}
 +\xi^3_4{\xi^2_1}^3
 -\xi^3_3\xi^3_1\xi^2_3{\xi^2_1}
 -\xi^3_3\xi^2_2{\xi^2_1}^2
 +{\xi^3_1}^2{\xi^2_3}^2
 -\xi^3_1\xi^2_4{\xi^2_1}^2
 +\xi^3_1\xi^2_3\xi^2_2{\xi^2_1}
 +{\xi^2_1}^2.$
Hence, if $\mathfrak{X}_{\xi^2_1 \neq 0}$ denotes the subset of
$\mathfrak{X}_{M5}$ such that $\xi^2_1 \neq 0,$ to conclude that
$\mathfrak{X}_{\xi^2_1 \neq 0}$ is a 12-dimensional submanifold of $\Rmath^{36},$
it will be sufficient to prove that the
preceding system is of maximal rank 2, that is in each of the subcases some 2-jacobian doesn't vanish.
\\ $\bullet$
First, one has
$\frac{D(f,g)}{D(\xi^4_1,\xi^3_1)} = -\frac{1}{{\xi^2_1}^2} \, (\xi^3_4 \xi^2_1 -\xi^3_1\xi^2_4)({\xi^2_4}^2+{\xi^2_3}^2)$
hence this 2-jacobian doesn't vanish if $\xi^2_4 \neq 0$ or $\xi^2_3 \neq 0.$
\\ $\bullet$
Suppose  $\xi^2_4=\xi^2_3 =0.$ Then
$\frac{D(f,g)}{D(\xi^5_5,\xi^6_5)} =
-\frac{1}{{\xi^6_5}^2} \, \left( (\xi^6_5 (\xi^3_4 +\xi^2_1)+ \xi^5_5(\xi^3_3 +\xi^2_2))^2
+({\xi^3_3}+{\xi^2_2})^2\right) $
hence this 2-jacobian doesn't vanish if $\xi^3_3 + \xi^2_2 \neq 0.$
\\ $\bullet$
Suppose  $\xi^2_4=\xi^2_3 =0$ and $\xi^3_3 =-\xi^2_2.$ Then
$\frac{D(f,g)}{D(\xi^5_5,\xi^6_5)} =
 (\xi^3_4 +\xi^2_1)^2 $
hence this 2-jacobian doesn't vanish if $\xi^3_4 + \xi^2_1 \neq 0.$
\\ $\bullet$
Suppose  $\xi^2_4=\xi^2_3 =0$ and $\xi^3_3 =-\xi^2_2$ and $\xi^3_4 =-\xi^2_1.$ Then
the equation $f=0$ reads $-2\xi^2_2\xi^2_1=0,$ hence $\xi^2_2 =0.$
Then the equation $g=0$ reads ${\xi^2_2}^2-{\xi^2_1}^2+1=0,$ hence $\xi^2_1 =\pm 1.$
Then
$\frac{D(f,g)}{D(\xi^2_2,\xi^3_3)} = -2\xi^5_5\xi^2_1$ and $
\frac{D(f,g)}{D(\xi^2_3,\xi^2_4)} =  {\xi^4_1}^2 +{\xi^3_1}^2$
hence if $\xi^5_5$ or $\xi^4_1$ or $\xi^3_1$ is $\neq 0,$
one of these  2-jacobians doesn't vanish.
On the other hand, if $\xi^5_5 = \xi^4_1 = \xi^3_1 =0,$
$\frac{D(f,g)}{D(\xi^2_1,\xi^3_3)} = \xi^2_1(1-{\xi^6_5}^2)$
hence if $\xi^6_5 \neq \pm 1,$
this 2-jacobian doesn't vanish. Suppose now that $\xi^6_5 = \pm 1.$
Then
$\frac{D(f,g)}{D(\xi^2_1,\xi^2_2)} =2{\xi^6_5} (\xi^2_1-{\xi^6_5})$
hence if $\xi^6_5 \neq \xi^2_1,$
this 2-jacobian doesn't vanish. Suppose finally that $\xi^6_5 = \xi^2_1.$
Then $\frac{D(f,g)}{D(\xi^3_3,\xi^3_4)} =4.$
This ends the proof that the system
(\ref{f,g}) is of maximal rank.
\subsection{}
From  \ref{xi(2,4)}, one has that
the subset  $\mathfrak{X}_{{\xi^2_4}^2 +{\xi^2_3}^2 \neq 0}$
of
$\mathfrak{X}_{M5}$ such that ${\xi^2_4}^2 +{\xi^2_3}^2 \neq 0,$
is a 12-dimensional submanifold of $\Rmath^{36}$ with the global chart
(\ref{globalchartxi(2,1)M5}).
Now if $\xi^2_4=\xi^2_3=0,$ then necessarily $\xi^2_1 \neq 0.$ Hence
$$\mathfrak{X}_{M5} =
\mathfrak{X}_{{\xi^2_4}^2 +{\xi^2_3}^2 \neq 0}
\cup \mathfrak{X}_{\xi^2_1 \neq 0}$$
and we conclude from   \ref{xi(2,1)} that $\mathfrak{X}_{M5}$
is a 12-dimensional submanifold of $\Rmath^{36}.$
A local chart  in a neighborhood of the canonical CS $J_0$
appears in  (\cite{sc1}, 6.5, p.155).
%%%%%%%%%%%%%%%%%%%%%%%%%%%%%%%%%%%%%%%%%%%%%%%%%%%%%%%%%%%%%%%%%%%%%%%%%%%%%%%%%%%%%%%%%%
\subsection{Equivalence.}
Due to the number of cases that had to be considered in the preceding computations,
we'll tackle the equivalence problem in a slightly different way, mixing equations solving
 and reduction by equivalence. First we give explicit computation of
all $\Phi \in \text{Aut}({M5}).$
\subsubsection{}
\begin{lemma}
\label{lemmaautM5}
$\text{Aut}({M5})$
is comprised of all those (real) matrices
of the following form :
\begin{equation}
{\fontsize{8}{10}\selectfont
\label{AutM5}
 \Phi = \begin{pmatrix}
\Phi_{\dagger} &
\begin{matrix} 0&0\\0&0\\0&0&\\0&0\end{matrix}
\\
\begin{matrix} b^5_1&b^5_2&b^5_3&b^5_4\\
b^6_1&b^6_2&b^6_3&b^6_4\end{matrix} &
\begin{matrix}
H&Ku\\
K&-Hu
\end{matrix}
 \end{pmatrix}
}
\end{equation}
 with
 \begin{equation}
{\fontsize{8}{10}\selectfont
 \label{phi0}
\Phi_{\dagger} = \begin{pmatrix}
b^1_1 & b^2_1 u & b^1_3 &-b^2_3u \\
b^2_1 & -b^1_1 u & b^2_3 & b^1_3u \\
b^3_1 & -b^4_1 u & b^3_3 & b^4_3u \\
b^4_1 & b^3_1 u & b^4_3 & -b^3_3u
 \end{pmatrix}
 }
 \end{equation}
% \Phi = \begin{pmatrix}
%b^1_1 & b^2_1 u & b^1_3 &-b^2_3u &0&0\\
%b^2_1 & -b^1_1 u & b^2_3 & b^1_3u &0&0\\
%b^3_1 & -b^4_1 u & b^3_3 & b^4_3u &0&0\\
%b^4_1 & b^3_1 u & b^4_3 & -b^3_3u &0&0\\
%b^5_1  & b^5_2 & b^5_3 & b^5_4&H&Ku\\
%b^6_1  & b^6_2 & b^6_3 & b^6_4&K&-Hu
% \end{pmatrix}
% }
% \end{equation}
where $u=\pm 1$ and
$H +i K = \det(\bv{w}_1,\bv{w}_3) \neq 0$ with
 $\bf{w}_1 =\left( \begin{smallmatrix} b^1_1 -i b^2_1\\ b^3_1 +i b^4_1  \end{smallmatrix}\right)$
and $\bf{w}_3 =\left(\begin{smallmatrix} b^1_3 -i b^2_3\\ b^3_3 +i b^4_3
\end{smallmatrix}\right).$
\end{lemma}
\begin{proof}
Let $\Phi=(b^i_j) \in \text{Aut}({M5}).$
Since $\Phi$ leaves   the $2^{\text{d}}$ central derivative $\mathcal{C}^2 ({M5})$
invariant, $b^i_5 =b^i_6=0$ for $ 1 \leqslant i \leqslant 4.$
Denote by $ij|k$ the equation obtained by projecting on $x_k$ the equation
$\left[\Phi(x_i),\Phi(x_j)\right] - \Phi ([x_i,x_j])=0.$
Equations $13|5,13|6,14|5,14|6$ yield
$b^5_5 =b^4_3b^2_1 -b^4_1 b^2_3 +b^3_3b^1_1 -b^3_1b^1_3; \;
b^6_5 =b^4_3b^1_1 -b^4_1 b^1_3 -b^3_3b^2_1 +b^3_1b^2_3; \;
b^5_6 =b^4_4b^2_1 -b^4_1 b^2_4 +b^3_4b^1_1 -b^3_1b^1_4; \;
b^6_6 =b^4_4b^1_1 -b^4_1 b^1_4 -b^3_4b^2_1 +b^3_1b^2_4.$
Now equations $12|5$ and $12|6$ read respectively
$\Delta^{1,3} + \Delta^{2,4} =0$ and
$\Delta^{1,4} - \Delta^{2,3} =0,$ with $\Delta^{i,j}$
the minors formed with the $1^{\text{st}}$ and $2^{\text{d}}$ columns and the rows indicated by the indices
 in  the matrix
${\Phi}_{\dagger} = {\left( b^i_j\right)}_{1\leqslant i,j \leqslant 4}.$
% \begin{equation*}
%{\Phi}_{\dagger} = \begin{pmatrix}
%b^1_1 & b^1_2  & b^1_3 &b^1_4 \\
%b^2_1 & b^2_2  & b^2_3 & b^2_4 \\
%b^3_1 & b^3_2  & b^3_3 & b^3_4 \\
%b^4_1 & b^4_2  & b^4_3 & b^4_4
% \end{pmatrix}
% \, .
% \end{equation*}
 If we introduce  for $1 \leqslant j \leqslant 4\quad $
 ${\bf{w}}_j =\begin{pmatrix} b^1_j -i b^2_j\\ b^3_j +i b^4_j\\  \end{pmatrix},$
 then
$\det(\bv{w}_1,\bv{w}_2) =\Delta^{1,3} + \Delta^{2,4} +i(\Delta^{1,4} - \Delta^{2,3}),$
hence equations $12|5$ and $12|6$ are equivalent to the single complex equation
$ \det(\bv{w}_1,\bv{w}_2) = 0,$ i.e. to the existence of $z=\alpha+i\beta \in \Cmath$ such that
 ${\bf{w}}_2 = z\, {\bf{w}}_1.$
In the same way,
equations $34|5$ and $34|6$ are equivalent
to the existence of $w=\gamma+i\delta \in \Cmath$ such that
 ${\bf{w}}_4 = w\, {\bf{w}}_3.$
 No, if we introduce $h=\alpha+\gamma , k = \beta+\delta,$ the system
 $23|5, 23|6$ reads
 \begin{equation}
 \label{autM5systhk}
\begin{cases}
h\, H
- k \, K  =0\\
k\, H
+h \,K
=0
 \end{cases}
 \end{equation}
where $H ={\Delta^{\prime}}^{2,4} +{\Delta^{\prime}}^{1,3}, \;
K={\Delta^{\prime}}^{1,4} -{\Delta^{\prime}}^{2,3},$
 ${\Delta^{\prime}}^{i,j}$
the minors formed with the $1^{\text{st}}$ and $3^{\text{d}}$ columns and the rows indicated by the indices
 in  the matrix  $\Phi_{\dagger}.$ Since
$H+iK =\det(\bv{w}_1,\bv{w}_3),$
the case $H=K=0$
would imply
$\det(\bv{w}_1,\bv{w}_3)=0$ which in turn leads to $\det \Phi =0.$ Hence
$(H,K)$  is a non trivial solution to the system
 (\ref{autM5systhk}). As its determinant is $h^2 + k^2$ we conclude that $h=k=0,$
 i.e. $\gamma=- \alpha, \delta = -\beta.$
 Now we are left only with  the system of equations $24|5,24|6.$ It reads
 \begin{equation}
 \label{autM5systHK}
\begin{cases}
(\alpha^2-\beta^2+1) H -2\alpha\beta K=0\\
2\alpha \beta H +(\alpha^2-\beta^2+1) K =0.
 \end{cases}
 \end{equation}
Again, as $(H,K)$  is a non trivial solution one has $\alpha\beta =0, \alpha^2-\beta^2+1=0$
i.e. $\alpha =0, \beta=\pm 1.$
Then we get  (\ref{AutM5}) with $u=\beta.$
\end{proof}
The subgroup $\text{Aut}({\mathfrak{n}}) \subset \text{Aut}({M5})$
of complex automorphisms of the complex Heisenberg Lie algebra $\frak{n}$
is the subgroup
comprised of all those matrices
in (\ref{AutM5})
for which $u=-1$ and $b^6_2=-b^5_1,b^5_2=b^6_1, b^6_4=b^5_3, b^5_4=-b^6_3.$
\subsubsection{}
Now, looking for CSs $J,$ after completing a set of general steps, one is left
to find solutions of the torsion equations and $J^2=-1$ of the following form:
$$ J = \begin{pmatrix}
\xi^1_1&
\xi^1_2&
\xi^1_3&
\xi^1_4&
0&
0\\
\xi^2_1&
\xi^2_2&
\xi^2_3&
\xi^2_4&
0&
0\\
\xi^3_1&
\xi^3_2&
\xi^3_3&
\xi^3_4&
0&
0\\
\xi^4_1&
\xi^4_2&
\xi^4_3&
\xi^4_4&
0&
0\\
\xi^5_1&
\xi^5_2&
\xi^5_3&
\xi^5_4&
\xi^5_5&
-\frac{{\xi^5_5}^2+1}{\xi^6_5}\\
{\xi^6_1}&
{\xi^6_2}&
{\xi^6_3}&
{\xi^6_4}&
{\xi^6_5}&
-\xi^5_5
 \end{pmatrix}$$
 where the $\xi^5_j$'s are certain linear expressions in the $\xi^6_k$'s ( $1 \leqslant j,k
 \leqslant 4$) and $\xi^4_4 = -\sum_{1}^3 \xi^j_j.$
Take the following $\Phi \in  \text{Aut}({M5}) :$
$${\fontsize{8}{10}\selectfont
 \Phi = \begin{pmatrix}
1&
0&
0&
0&
0&
0\\
0&
1&
0&
0&
0&
0\\
0&
0&
1&
0&
0&
0\\
0&
0&
0&
1&
0&
0\\
b^5_1&
b^5_2&
b^5_3&
b^5_4&
1&
0\\
{b^6_1}&
{b^6_2}&
{b^6_3}&
{b^6_4}&
0&
1\end{pmatrix}
}
$$
with
$ {b^5_j}=
\frac{1}{{\xi^6_5}} \,
({-{\xi^6_j}+ \sum_{k=1}^{4} {b^6_k} {\xi^k_j} +{\xi^5_5} {b^6_j}}).$
Then equivalence by $\Phi$ leads to the case where
$ \xi^5_j=\xi^6_j = 0\;  \forall j \; 1 \leqslant j \leqslant 4. $
Now consider the submatrix
\begin{equation}
\label{Jdagger}
J_{\dagger}=
\begin{pmatrix}
\xi^1_1&
\xi^1_2&
\xi^1_3&
\xi^1_4
\\
\xi^2_1&
\xi^2_2&
\xi^2_3&
\xi^2_4
\\
\xi^3_1&
\xi^3_2&
\xi^3_3&
\xi^3_4
\\
\xi^4_1&
\xi^4_2&
\xi^4_3&
\xi^4_4
 \end{pmatrix}=
\begin{pmatrix}
A&B\\C&D
 \end{pmatrix}
 .
 \end{equation}
We get 2 cases 1 and 2 below. Before we proceed further, we record a lemma.
\begin{lemma}
\label{reduclemma}
Let $J$ be a CS
on ${M5}$ of the following form :
\begin{equation*}
 J= \text{diag} \left(  J_{\dagger},
\begin{pmatrix}
\xi^5_5&
-\frac{{\xi^5_5}^2+1}{\xi^6_5}\\
{\xi^6_5}& -\xi^5_5
 \end{pmatrix} \right)
% J= \begin{pmatrix}  J_{\dagger}&
% \begin{matrix}0&0\\0&0\\0&0\\0&0\end{matrix} \\
% \begin{matrix}0&0&0&0\\0&0&0&0\end{matrix} &
% \begin{pmatrix}
%\xi^5_5&
%-\frac{{\xi^5_5}^2+1}{\xi^6_5}\\
%{\xi^6_5}& -\xi^5_5
% \end{pmatrix}
% \end{pmatrix}
 \end{equation*}
 where $J_{\dagger} $ is given in (\ref{Jdagger}).
\\ (i)
 Suppose that $B$ and $C$ are not simultaneously zero.
Then $J$ is equivalent to a CS of the same form for which $B \neq 0.$
\\ (ii)
 Suppose
$ B = \begin{pmatrix}\xi^1_3& \xi^1_4\\ \xi^2_3&\xi^2_4 \end{pmatrix} \neq 0.$
Then $J$ is equivalent to a CS of the same form for which $\xi^2_3 =1,
\xi^2_4=0.$
\\ (iii)
 Suppose
$ B = \begin{pmatrix}\xi^1_3& \xi^1_4\\ 1& 0\end{pmatrix}.$
Then $J$ is equivalent to a CS of the same form having the same $B$ and for which $\xi^2_1 =\xi^2_2=0.$
\\ (iv) Suppose $B=C=0$ and $\xi^2_2=-\xi^1_1, \xi^4_4 =-\xi^3_3.$
Then $J$ is equivalent to a CS of the same form for which $B=C=0$ and  $\xi^1_1 =\xi^3_3=0.$
\end{lemma}
\begin{proof}
(i)
The matrix $J$ has the form
$\begin{pmatrix}A&B&0\\C&D&0\\0&0&E \end{pmatrix}\, .$
Take
$\Phi_{\dagger}= \begin{pmatrix}0&I\\I&0 \end{pmatrix}$
and
$\Phi= \text{diag} \left( \Phi_{\dagger},\begin{pmatrix} -1&0\\0&1\end{pmatrix}
\right) \, .$
From Lemma \ref{lemmaautM5},
$\Phi \in \text{Aut}({M5}).$ Now
$$(\Phi_{\dagger})^{-1} \begin{pmatrix}A&B\\C&D\end{pmatrix} \Phi_{\dagger} =
 \begin{pmatrix}0&I\\I&0\end{pmatrix}
 \begin{pmatrix}A&B\\C&D\end{pmatrix}
 \begin{pmatrix}0&I\\I&0\end{pmatrix}
 =\begin{pmatrix}D&C\\B&A\end{pmatrix}.$$
 Hence if $B$ and $C$ are not simultaneously zero, one may suppose $B \neq 0.$
\\(ii)
Suppose first  $\xi^2_3=\xi^2_4=0.$ Then, since  $B\neq 0,$ the first row of $B$ is not zero.
Consider
$\Phi_{\dagger}= \begin{pmatrix}U&0\\0&U \end{pmatrix}$
and
$\Phi= \text{diag}\left( \Phi_{\dagger},\begin{pmatrix} 1&0\\0&-1\end{pmatrix} \right)$
with
$ U = \begin{pmatrix}0&1\\ 1&0\end{pmatrix}.$
Then   $\Phi \in \text{Aut}({M5}).$ Now
$$(\Phi_{\dagger})^{-1} \begin{pmatrix}A&B\\C&D\end{pmatrix} \Phi_{\dagger} =
 \begin{pmatrix}U&0\\0&U\end{pmatrix}
 \begin{pmatrix}A&B\\C&D\end{pmatrix}
 \begin{pmatrix}U&0\\0&U\end{pmatrix}
 =\begin{pmatrix}UAU&UBU\\UCU&UDU\end{pmatrix}.$$
 where
 $UBU =
\begin{pmatrix}0&1\\ 1&0\end{pmatrix}
\begin{pmatrix}\xi^1_3& \xi^1_4\\ \xi^2_3&\xi^2_4 \end{pmatrix}
\begin{pmatrix}0&1\\ 1&0\end{pmatrix}
 =
\begin{pmatrix}\xi^2_4& \xi^2_3\\ \xi^1_4&\xi^1_3 \end{pmatrix} .$
Hence we are led to the case where the second row of $B$ doesn't vanish.
Consider that case now.
Introduce
$\Phi_{\dagger}= \begin{pmatrix}I&0\\0&V\end{pmatrix}$
and
$\Phi= \text{diag}\left( \Phi_{\dagger},\begin{pmatrix} \alpha&-\beta\\ \beta&\alpha\end{pmatrix}\right) $
with
$ V = \begin{pmatrix}\alpha&-\beta\\ \beta &\alpha\end{pmatrix} \neq 0.$
Then $\Phi \in \text{Aut}({M5}).$ Now
$$(\Phi_{\dagger})^{-1} \begin{pmatrix}A&B\\C&D\end{pmatrix} \Phi_{\dagger} =
 \begin{pmatrix}I&0\\0&V^{-1}\end{pmatrix}
 \begin{pmatrix}A&B\\C&D\end{pmatrix}
 \begin{pmatrix}I&0\\0&V\end{pmatrix}
 =\begin{pmatrix}A&BV\\V^{-1}C&V^{-1}DV\end{pmatrix}.$$
 where
 $$BV =
\begin{pmatrix}\xi^1_3& \xi^1_4\\ \xi^2_3&\xi^2_4 \end{pmatrix}
\begin{pmatrix}\alpha&-\beta\\ \beta &\alpha\end{pmatrix}
 =
\begin{pmatrix}*& *\\ {\xi^\prime}^2_3 &{\xi^\prime}^2_4 \end{pmatrix}
$$
with ${\xi^\prime}^2_3 = \alpha \xi^2_3 + \beta \xi^2_4\; , \;
{\xi^\prime}^2_4\ = \alpha \xi^2_4  - \beta \xi^2_3 .$
We want $ {\xi^\prime}^2_3 =1 $ and ${\xi^\prime}^2_4 =0.$
This is a Cramer system in $\alpha, \beta$ since ${(\xi^2_3)}^2 + (\xi^2_4)^2 \neq 0$, hence it has a nontrivial solution.
Hence we are reduced to the case where $\xi^2_3=1$ and $\xi^2_4=0.$
\\ (iii) Suppose $ B = \begin{pmatrix}\xi^1_3& \xi^1_4\\ 1&0 \end{pmatrix}.$
Consider
$\Phi_{\dagger}= \begin{pmatrix}I&0\\T&I \end{pmatrix}$
and
$\Phi= \text{diag}\left( \Phi_{\dagger},\begin{pmatrix} 1&0\\0&1\end{pmatrix}\right) $
with
$ T = \begin{pmatrix}\alpha& \beta\\ \beta&-\alpha  \end{pmatrix}\neq 0.$
Then $\Phi \in \text{Aut}({M5}).$ Now
\begin{eqnarray*}
(\Phi_{\dagger})^{-1} \begin{pmatrix}A&B\\C&D\end{pmatrix} \Phi_{\dagger} &=&
 \begin{pmatrix}I&0\\-T&I\end{pmatrix}
 \begin{pmatrix}A&B\\C&D\end{pmatrix}
 \begin{pmatrix}I&0\\T&I\end{pmatrix} \\
 &=& \begin{pmatrix}A&B\\-TA+C&-TB+D\end{pmatrix}
 \begin{pmatrix}I&0\\T&I\end{pmatrix}
 \\ &=&\begin{pmatrix}A+BT&B\\ * & * \end{pmatrix}.
\end{eqnarray*}
 We want $A+BT= \begin{pmatrix}*& *\\ 0&0 \end{pmatrix} ,$ that is
$ BT = \begin{pmatrix}*& *\\ -\xi^2_1&-\xi^2_2 \end{pmatrix}.$
As $BT =
\begin{pmatrix}\xi^1_3& \xi^1_4\\ 1&0 \end{pmatrix}
\begin{pmatrix}\alpha& \beta\\ \beta&-\alpha  \end{pmatrix}
=
\begin{pmatrix}*& *\\ \alpha&\beta\end{pmatrix},$
one has to let $\alpha = -\xi^2_1$ and $\beta = - \xi^2_2,$
which is possible if one doesn't already has $\xi^2_1=0, \; \xi^2_2=0.$
Hence we are reduced to the case where $\xi^2_1=0$ and $\xi^2_2=0.$
\\ (iv)
The matrix $J$ has the form
$\begin{pmatrix}A&0&0\\0&D&0\\0&0&E \end{pmatrix}$
with moreover $\xi^2_2=-\xi^1_1$ and $\xi^4_4=-\xi^3_3.$
Consider
$\Phi_{\dagger}= \begin{pmatrix}U&0\\0&V \end{pmatrix}$
and ($\alpha^2+ \beta^2 \neq 0, \; \gamma^2 +\delta^2 \neq 0$)
$\Phi= \text{diag} \left( \Phi_{\dagger},\begin{pmatrix} H&-K\\K&H\end{pmatrix} \right)$
where $ U = \begin{pmatrix}\alpha& -\beta\\ \beta&\alpha  \end{pmatrix}
\; , \;  V = \begin{pmatrix}\gamma& -\delta\\ \delta&\gamma  \end{pmatrix}\, ,$
and $H +iK= (\alpha-i\beta)(\gamma+i\delta).$
Then   $\Phi \in \text{Aut}({M5}).$ Now
$$(\Phi_{\dagger})^{-1} \begin{pmatrix}A&0\\0&D\end{pmatrix} \Phi_{\dagger} =
 \begin{pmatrix}U^{-1}&0\\0&V^{-1}\end{pmatrix}
 \begin{pmatrix}A&0\\0&D\end{pmatrix}
 \begin{pmatrix}U&0\\0&V\end{pmatrix}
 =\begin{pmatrix}U^{-1}AU&0\\0&V^{-1}D\end{pmatrix}.$$
One has
 $$ U^{-1}AU = \frac{1}{\alpha^2+\beta^2}
\begin{pmatrix}{\xi^\prime}^1_1& {\xi^\prime}^1_2\\ {\xi^\prime}^2_1&{\xi^\prime}^2_2  \end{pmatrix},$$
with
\begin{eqnarray*}
{\xi^\prime}^1_1&=&
\alpha^2 \xi^1_1 + \alpha\beta (\xi^1_2+\xi^2_1) +\beta^2\xi^2_2\\
{\xi^\prime}^2_2&=&
\alpha^2 \xi^2_2 - \alpha\beta (\xi^1_2+\xi^2_1) +\beta^2\xi^1_1.
\end{eqnarray*}
Since  $\xi^1_1 =-\xi^2_2,$ one has ${\xi^\prime}^1_1= -{\xi^\prime}^2_2.$
Then the discriminant
of the equation  ${\xi^\prime}^2_2=0$ is the sum of two squares,
hence there exist $\alpha, \beta \in \Rmath$ with $\beta = 1 $ such that
${\xi^\prime}^1_1=0$  and ${\xi^\prime}^2_2=0.$
Then
 $$ U^{-1}AU = \frac{1}{\alpha^2+\beta^2}
\begin{pmatrix}0& {\xi^\prime}^1_2\\ {\xi^\prime}^2_1&0 \end{pmatrix}.$$
Similarly,
there exist $\gamma, \delta \in \Rmath$ with $\delta = 1 $ such that
 $$ V^{-1}DV = \frac{1}{\gamma^2+\delta^2}
\begin{pmatrix}0& {\xi^\prime}^3_4\\ {\xi^\prime}^4_3&0 \end{pmatrix}.$$
Hence we are reduced to the case where $\xi^1_1 = \xi^3_3 =0.$
\end{proof}
\subsubsection{Case 1.}
In the present case 1, we suppose that $B$ and $C$
in (\ref{Jdagger})
are not both $0.$
Then by equivalence
 (Lemma \ref{reduclemma} (i)),  we may suppose $B \neq 0.$
 Since $ B \neq 0, $ by equivalence (Lemma \ref{reduclemma} (ii)), we may assume
 $ {\xi^2_3}=1,$  $ {\xi^2_4}=0.$
 Again  by equivalence (Lemma \ref{reduclemma} (iii)), we may assume
 without altering $ B $  that
$ {\xi^2_1}=0,$
 $ {\xi^2_2}=0.$
 Then solving all equations, we   finally  get   the matrix
\begin{equation}
{\fontsize{8}{10}\selectfont
\label{JM5case1}
J(\xi^1_3,\xi^1_4,\xi^5_5,\xi^6_5)
  = \begin{pmatrix}
a&
 - {\xi^6_5} {\xi^1_4} + {\xi^6_5} - {\xi^5_5} {\xi^1_3}&
{\xi^1_3}&
{\xi^1_4}&
0&
0\\
0&
0&
1&
0&
0&
0\\
0&
-1&
0&
0&
0&
0\\
%( - (({\xi^6_5} {\xi^5_5}^2 {\xi^1_4}^2 - 2 {\xi^6_5} {\xi^5_5}^2 {\xi^1_4} + {\xi^6_5} {\xi^5_5}^2 + {\xi^6_5} + 2 {\xi^5_5}^3 {\xi^1_4} {\xi^1_3} - 2 {\xi^5_5}^3 {\xi^1_3} + 2 {\xi^5_5} {\xi^1_4} {\xi^1_3} - 2 {\xi^5_5} {\xi^1_3}) {\xi^6_5} + ({\xi^5_5}^2 + 1)^2 {\xi^1_3}^
%2))/({\xi^6_5}^2 {\xi^1_4})&
b&
%(({\xi^6_5} {\xi^5_5} {\xi^1_4}^2 - 2 {\xi^6_5} {\xi^5_5} {\xi^1_4} + {\xi^6_5} {\xi^5_5} + 2
%{\xi^5_5}^2 {\xi^1_4} {\xi^1_3} - 2 {\xi^5_5}^2 {\xi^1_3} + {\xi^1_4} {\xi^1_3}) {\xi^6_5} + (
%{\xi^5_5}^2 + 1) {\xi^5_5} {\xi^1_3}^2)/({\xi^6_5} {\xi^1_4})&
c&
%(({\xi^6_5} {\xi^1_4} - {\xi^6_5} - {\xi^5_5} {\xi^1_4} {\xi^1_3} + 2 {\xi^5_5} {\xi^1_3}) {\xi^6_5} - ({\xi^5_5}^2 + 1) {\xi^1_3}^2)/({\xi^6_5} {\xi^1_4})&
d&
%( - (({\xi^5_5}^2 + 1) {\xi^1_3} + ({\xi^1_4} - 1) {\xi^6_5} {\xi^5_5}))/{\xi^6_5}&
-a&
0&
0\\
0&
0&
0&
0&
{\xi^5_5}&
 - ({\xi^5_5}^2 + 1)/{\xi^6_5}\\
0&
0&
0&
0&
{\xi^6_5}&
 - {\xi^5_5}\end{pmatrix}
 }
 \end{equation}
 where : \\
${a} ={J^1_1}= (({\xi^5_5}^2 + 1) {\xi^1_3} + ({\xi^1_4} - 1) {\xi^6_5} {\xi^5_5})/{\xi^6_5}; \myquad
{b}={J^4_1}=  -\frac{1+a^2}{\xi^1_4};\myquad
{c}={J^4_2}=
(({\xi^6_5} {\xi^5_5} ({\xi^1_4} - 1)^2
+ 2  {\xi^5_5}^2 {\xi^1_3} ({\xi^1_4}-1) + {\xi^1_4} {\xi^1_3}) {\xi^6_5} + (
{\xi^5_5}^2 + 1) {\xi^5_5} {\xi^1_3}^2)/({\xi^6_5} {\xi^1_4});\myquad
{d}={J^4_3}=
(({\xi^6_5} ({\xi^1_4} - 1) + {\xi^5_5} {\xi^1_3} ( 2- {\xi^1_4})) {\xi^6_5} - ({\xi^5_5}^2 + 1) {\xi^1_3}^2)/({\xi^6_5} {\xi^1_4});\myquad
$
and the parameters are subject to the condition
\begin{equation}
\label{condM5equivcase1}
{\xi^1_4} {\xi^6_5}  \neq 0 .
 \end{equation}
Note also the following formulae from $J^2 =-1 :$
\begin{equation}
\label{eqd}
d=\frac{1}{\xi^1_4} \, (\xi^6_5 (\xi^1_4 -1) +\xi^1_3(\xi^5_5 -a))\quad ,\quad
c=ad + (1+a^2) \, \frac{\xi^1_3}{\xi^1_4}.
\end{equation}
\par Commutation relations of $ \mathfrak{m} : $
 $  [\tilde{x}_1,\tilde{x}_2]= a\tilde{x}_5 +d\tilde{x}_6; \;$
$  [\tilde{x}_1,\tilde{x}_3]= (b+1) \tilde{x}_5 -c \tilde{x}_6 ; \; $
$  [\tilde{x}_2,\tilde{x}_3]=
\frac{1}{{\xi^6_5} {\xi^1_4}}\,
(({\xi^6_5}( {\xi^1_4}-1)  +  {\xi^5_5}{\xi^1_3})^2 +{\xi^1_3}^2) \,
({\xi^5_5} \tilde{x}_5 + {\xi^6_5} \tilde{x}_6)  ;$
\; $  [\tilde{x}_2,\tilde{x}_4]=   (1- {\xi^1_4}) \tilde{x}_5 + {\xi^1_3} \tilde{x}_6 ; \;$
$  [\tilde{x}_3,\tilde{x}_4]=   a \tilde{x}_5 +
({\xi^6_5} ({\xi^1_4}-1) +\xi^5_5\xi^1_3)\, \tilde{x}_6 .$
%%%%%%%%%%%%%%%%%%%%%%%%%%%%%%%%%%%%%%%%%%%%%%%%%%%%%%%%%%%%%%%%%%%%%%%%%%%%%%%%%%%%
\par
$ \mathfrak{m}$ is abelian if and only if
\begin{equation}
\label{M5equivaabelian}
    {\xi^1_3}=0\quad ,\quad {\xi^1_4}=1.
     \end{equation}
Suppose now that (\ref{M5equivaabelian}) holds, and denote $ \alpha = {\xi^5_5}, \beta = {\xi^6_5} \neq 0 ,$
$ J(\alpha,\beta) = J(0,1,\xi^5_5,\xi^6_5).$ Then :
\begin{equation}
\label{J(0,beta)}
{\fontsize{8}{10}\selectfont
 J(\alpha,\beta)
= \begin{pmatrix}
0&
0&
0&
1&
0&
0\\
0&
0&
1&
0&
0&
0\\
0&
-1&
0&
0&
0&
0\\
-1&
0&
0&
0&
0&
0\\
0&
0&
0&
0&
\alpha&
 - \frac{\alpha^2 + 1}{\beta}\\
0&
0&
0&
0&
\beta&
 - \alpha\end{pmatrix}
 }
 \end{equation}
where $ \beta  \neq 0 .$
A direct study shows  that any $ J(\alpha,\beta) $
is equivalent to a unique $ J(0,\gamma) ,
\; 0 < \gamma \leqslant 1, $
and any two such $ J(0,\gamma) , \; 0 < \gamma \leqslant 1 $
are not equivalent unless $ \gamma = \gamma ' .$
%\begin{equation}
%{\fontsize{8}{10}\selectfont
%\label{J(0,beta)}
%J(0,\beta) = \begin{pmatrix}
%0&
%0&
%0&
%1&
%0&
%0\\
%0&
%0&
%1&
%0&
%0&
%0\\
%0&
%-1&
%0&
%0&
%0&
%0\\
%-1&
%0&
%0&
%0&
%0&
%0\\
%0&
%0&
%0&
%0&
%0&
%- \frac{1}{\beta}\\
%0&
%0&
%0&
%0&
%\beta&
%0\end{pmatrix}.
%}
%\end{equation}
\subsubsection{}
Let $J_1 =J(\xi^1_3,\xi^1_4,\xi^5_5,\xi^6_5)$ and
 $J_2 =J(\eta^1_3,\eta^1_4,\eta^5_5,\eta^6_5)$ as in (\ref{JM5case1}).
 It is clear that $J_1 \cong J_2$ if and only if there exists   a matrix $\Phi$
 such that
 \begin{equation}
 \label{existphi}
\Phi= \begin{pmatrix}\Phi_{\dagger}&0\\0& \begin{pmatrix} H&Ku\\K&-Hu\end{pmatrix}\end{pmatrix}
\end{equation}
where $u= \pm 1$, $\Phi_{\dagger}, H, K$ are defined in Lemma \ref{lemmaautM5}, $H^2 + K^2 \neq 0$ and
$J_2 = \Phi^{-1} J_1 \Phi .$
This equation implies that
\begin{equation}
\label{CN}
\begin{pmatrix}\xi^5_5& -\frac{{\xi^5_5}^2+1}{ \xi^6_5}\\ \xi^6_5&-\xi^5_5 \end{pmatrix}
\begin{pmatrix} H&Ku\\K&-Hu\end{pmatrix}
=\begin{pmatrix} H&Ku\\K&-Hu\end{pmatrix}
\begin{pmatrix}\eta^5_5& -\frac{{\eta^5_5}^2+1}{ \eta^6_5}\\ \eta^6_5&-\eta^5_5 \end{pmatrix}.
\end{equation}
Equation (\ref{CN}) has non trivial solutions in $H,K$ if and only if the following condition holds :
\begin{equation}
\label{CNSxi55xi65}
\frac{{\eta^5_5}^2 +{\eta^6_5}^2 +1}{\eta^6_5} =
-u\, \frac{{\xi^5_5}^2 +{\xi^6_5}^2 +1}{\xi^6_5} .
\end{equation}
However, though that condition is also sufficient for the existence of $\Phi$ in
 (\ref{existphi})
in the abelian case where
$\xi^1_3= \eta^1_3=0, \xi^1_4 = \eta^1_4 = 1,$ it is no longer sufficient in the nonabelian case.
For example, take $J(1,1,1,1)$ and $J(1,1,1,\eta^6_5),$ with $\eta^6_5 \neq 1 :$ then
(\ref{CNSxi55xi65}) holds if only if $u=-1, \eta^6_5 =2$ or $u=1, \eta^6_5 =-2,-1,$
and in neither  case does equivalence occur;
hence $J(1,1,1,1) \not \cong J(1,1,1,\eta^6_5)$ if $\eta^6_5 \neq 1.$
However $J(1,1,0,1)  \cong J(1,1,0,-1).$
One also has $J(\xi^1_3,1,0,1)) \cong J(\eta^1_3,\eta^1_4,0,1) \Leftrightarrow \eta^1_3 = \pm \xi^1_3 \; ; \; \eta^1_4 =\xi^1_4.$
In general, if $\eta^1_4 \neq \xi^1_4,$ for
$J(\xi^1_3,\xi^1_4,\xi^5_5,\xi^6_5)$ and
 $J(\eta^1_3,\eta^1_4,\eta^5_5,\eta^6_5)$ to be equivalent, it would be necessary that
 $P\xi^5_5 +Q =0$ where $P,Q$ are certain huge polynomials in the other variables.
 We simply
 conjecture here that equivalence implies
 $\xi^1_4 = \eta ^1_4$ and  $\eta^1_3 = \pm \xi^1_3,$
 and leave open the equivalence problem in the nonabelian case.
\subsubsection{Case 2.}
We now suppose $ B=C=0.$ Then necessarily $\xi^2_2=-\xi^1_1, \xi^4_4 = -\xi^3_3.$
By equivalence
 (Lemma \ref{reduclemma} (iv)),  we may suppose $ {\xi^1_1}=0 , {\xi^3_3}=0 .$
 \\
\textbf{Case 2.1.}
Suppose   ${\xi^5_5} = 0$.
If $ {\xi^4_3} {\xi^2_1} \neq 1 $ one gets then the matrix
\begin{equation}
\label{cas21xi(2,1)xi(4,1)neq0}
J(\xi^2_1,\xi^4_3)
 = \text{diag}\left(
 \begin{pmatrix}0& - 1/{\xi^2_1} \\ {\xi^2_1}&0\end{pmatrix},
 \begin{pmatrix}0& - {\xi^4_3} \\ 1/{\xi^4_3}&0\end{pmatrix},
 \begin{pmatrix}0& (\xi^4_3\xi^2_1- 1)/(\xi^4_3-{\xi^2_1}) \\
 (-\xi^4_3+\xi^2_1)/(\xi^4_3{\xi^2_1}-1) &0\end{pmatrix} \right)
\end{equation}
with the conditions
\begin{equation}
\label{cas21xi(2,1)xi(4,1)neq0cond}
\xi^2_1,\xi^4_3 \neq 0 \; ; \; \xi^4_3 \neq \xi^2_1, \frac{1}{\xi^2_1}.
\end{equation}
If $ {\xi^4_3} {\xi^2_1} = 1 ,$ we get
\begin{equation}
\label{cas21xi(2,1)xi(4,1)0}
 J  = \text{diag}\left(
 \begin{pmatrix}0& - 1/{\xi^2_1} \\ {\xi^2_1}&0\end{pmatrix},
 \begin{pmatrix}0& - {\xi^2_1} \\ 1/{\xi^2_1}&0\end{pmatrix},
 \begin{pmatrix}0& - 1/{\xi^6_5} \\ {\xi^6_5}&0\end{pmatrix} \right)
\end{equation}
with the conditions
\begin{equation}
\label{cas21xi(2,1)xi(4,1)=1cond}
\xi^2_1= \pm 1 \; ; \; \xi^6_5 \neq 0.
\end{equation}
 \\
\textbf{Case 2.2.}
Suppose   ${\xi^5_5} \neq  0$.
%We finally suppose $ B=C=0 $ and  ${\xi^5_5} \neq 0$.
%By equivalence,
% (Lemma \ref{reduclemma} (iv)),  we may suppose $ {\xi^1_1}=0 , {\xi^3_3}=0 .$
Then we get:
\begin{equation}
\label{case22}
 J  = \text{diag}\left(
 \begin{pmatrix}0& - 1/{\xi^2_1} \\ {\xi^2_1}&0\end{pmatrix},
 \begin{pmatrix}0& - 1/{\xi^2_1} \\ {\xi^2_1}&0\end{pmatrix},
 \begin{pmatrix}\xi^5_5& -({\xi^5_5}^2+ 1)/{\xi^6_5} \\ {\xi^6_5}&-\xi^5_5\end{pmatrix} \right)
 \end{equation}
with the conditions
\begin{equation}
\label{cas22cond}
\xi^2_1= \pm 1 \; ; \; \xi^5_5 \xi^6_5 \neq 0.
\end{equation}
\subsubsection{}
Computing intertwining automorphisms, one can prove that for all
$\xi^2_1,\xi^4_3,\eta^2_1,\eta^4_3 \in \Rmath $
satisfying conditions (\ref{cas21xi(2,1)xi(4,1)neq0cond}),
$J(\eta^2_1,\eta^4_3) \cong J(\xi^2_1,\xi^4_3) $
if and only if there exists $u = \pm 1$ such that one of the following is satisfied :
$$\left(\eta^2_1 = u\xi^2_1 \text{ or } \frac{u}{\xi^2_1}\right) \text{ and } \left(
\eta^4_3 = u\xi^4_3 \text{ or } \frac{u}{\xi^4_3}\right)$$
or
$$\left(\eta^2_1 = u\xi^4_3 \text{ or } \frac{u}{\xi^4_3} \right)\text{ and }
\left(
\eta^4_3 = u\xi^2_1 \text{ or } \frac{u}{\xi^2_1}\right).$$
Conditions (\ref{cas21xi(2,1)xi(4,1)neq0cond}) are preserved by the transformations.
For example, the canonical CS
$J_0=J(-1,1)$ and its opposite $-J_0$ are equivalent.
In fact, on has :
\begin{lemma}
Let $\Omega$ denote the $\text{Aut}({M5})$-orbit of the
canonical CS $J_0.$ Then $\Omega$ is the 4-dimensional space
comprised of the matrices
\begin{equation*}
{\fontsize{8}{10}\selectfont
%\label{M5canonicalorbit}
J = \begin{pmatrix}
0&
 - 1/{\xi^2_1}&
0&
0&
0&
0\\
{\xi^2_1}&
0&
0&
0&
0&
0\\
0&
0&
0&
1/{\xi^2_1}&
0&
0\\
0&
0&
 - {\xi^2_1}&
0&
0&
0\\
{\xi^6_2}&
 - {\xi^6_1}&
 - {\xi^6_4}&
{\xi^6_3}&
0&
1/{\xi^2_1}\\
{\xi^6_1}&
{\xi^6_2}&
{\xi^6_3}&
{\xi^6_4}&
 - {\xi^2_1}&
0\end{pmatrix}
}
\end{equation*}
where $\xi^2_1 = \pm1 , \xi^6_1,\xi^6_2,\xi^6_3,\xi^6_4 \in  \Rmath.$
\end{lemma}
\par
On the other hand, the $J$ in (\ref{cas21xi(2,1)xi(4,1)0}) appears simply as a limiting
case when $\xi^5_5 \rightarrow 0$ of
the structure $J(\xi^2_1,\xi^5_5,\xi^6_5)$ defined in (\ref{case22}) with $\xi^2_1 = \pm 1
\; , \; \xi^5_5 \xi^6_5 \neq 0.$
\subsubsection{}
${M5}$ is a complex algebra for the CS
$J(\xi^2_1,\xi^4_3)$
in (\ref{cas21xi(2,1)xi(4,1)neq0}) if and only if $\xi^2_1 =-1, \xi^4_3 =1$ or
$\xi^2_1 =1, \xi^4_3 =-1,$ \textit{i.e.}
$J(\xi^2_1,\xi^4_3)$ is the canonical CS
$J_0$ or its opposite  respectively.
Since ${M5}$ is not a complex algebra for the CS
$J(\xi^1_3,\xi^1_4,\xi^5_5,\xi^6_5)$ in (\ref{JM5case1}), the latter
is not equivalent to $J_0.$
\begin{lemma}
Suppose
$J(\xi^2_1,\xi^4_3) \not \in
\{J_0,-J_0\}.$
Then $J(\xi^2_1,\xi^4_3) $ is equivalent to some CS in case 1, \textit{i.e.}
there exist $\xi^1_3,\xi^1_4,\xi^5_5,\xi^6_5$ such that
$J(\xi^2_1,\xi^4_3)\cong J(\xi^1_3,\xi^1_4,\xi^5_5,\xi^6_5).$
\end{lemma}
\begin{proof}
Take
$\Phi_{\dagger}= \begin{pmatrix}I&S\\0&I \end{pmatrix}$
and
\begin{equation}
\label{philemma2}
\Phi= \begin{pmatrix}\Phi_{\dagger}&0\\0& \begin{pmatrix} 1&0\\0&1\end{pmatrix}\end{pmatrix}
\quad \text{ where } \quad
%\label{Slemma2}
 S = \begin{pmatrix}\alpha& \beta\\ \beta&-\alpha  \end{pmatrix} \quad (\alpha, \beta \in \Rmath).
\end{equation}
Then $\Phi \in \text{Aut}({M5}).$
Denote
${J(\xi^2_1,\xi^4_3)}_{\dagger} =\begin{pmatrix}A&0\\0&D\end{pmatrix}$
where
$A =\begin{pmatrix}0&-\frac{1}{\xi^2_1}\\ \xi^2_1&0\end{pmatrix}, \;
D =\begin{pmatrix}0&-\frac{1}{\xi^4_3}\\ \xi^4_3&0\end{pmatrix}.$
Then
$$(\Phi_{\dagger})^{-1} {J(\xi^2_1,\xi^4_3)}_{\dagger} \,  \Phi_{\dagger} =
 \begin{pmatrix}I&-S\\0&I\end{pmatrix}
\begin{pmatrix}A&0\\0&D\end{pmatrix}
 \begin{pmatrix}I&S\\0&I\end{pmatrix}
 =
\begin{pmatrix}A&AS-SD\\0&D\end{pmatrix} .$$
Now
$$AS-SD = \begin{pmatrix}-\beta( \frac{1}{\xi^2_1} + \xi^4_3)
&\alpha ( \frac{1}{\xi^2_1}  +\frac{1}{\xi^4_3}) \\
\alpha ({\xi^2_1}  +{\xi^4_3})
&\beta ({\xi^2_1}  +\frac{1}{\xi^4_3})\end{pmatrix}
.$$
The CS $\Phi^{-1} J(\xi^2_1,\xi^4_3) \,  \Phi$
is of type 1 if $AS- SD \neq 0.$
For this to hold, just choose
$\alpha =1, \beta =0$ if $\xi^2_1 +\xi^4_3 \neq 0$ and
$\alpha =0, \beta =1$ if $\xi^2_1 +\xi^4_3 = 0,$  noting in this latter case that
$\xi^2_1 +\frac{1}{\xi^4_3} = \frac{{(\xi^2_1)}^2 -1}{\xi^2_1} \neq 0 $ since $\xi^2_1 \neq \pm 1$ as
$J(\xi^2_1,\xi^4_3) \not \in \{J_0,-J_0\}.$
\end{proof}
\begin{lemma}
Let $J$ be either the CS defined in
 (\ref{case22}) with $\xi^2_1 = \pm 1 \; , \; \xi^5_5, \xi^6_5 \neq 0$ or
the one in (\ref{cas21xi(2,1)xi(4,1)0}).
Then
$J\cong J(0,\beta)$ for some $\beta \; (0< \beta \leqslant 1),$
where $J(0,\beta)$ is defined in (\ref{J(0,beta)}).
\end{lemma}
\begin{proof}
In both cases, $\mathfrak{m}$ is abelian, and
$J= \begin{pmatrix}A&0&0\\0&A&0\\0&0&E \end{pmatrix}$
where $A = \begin{pmatrix} 0&-\frac{1}{\xi^2_1}\\ \xi^2_1&0\end{pmatrix} \; , \;
E = \begin{pmatrix} \xi^5_5&-\frac{{\xi^5_5}^2+1}{\xi^6_5}\\ \xi^6_5&\xi^5_5\end{pmatrix}$
with $\xi^6_5 \neq 0,$  $\xi^2_1 = \pm 1,$
 $\xi^5_5 = 0$ in case
(\ref{cas21xi(2,1)xi(4,1)0})
and  $\xi^5_5 \neq 0$ in case
(\ref{case22}).
Take $\Phi, S$ as in (\ref{philemma2}).
Denote
${J}_{\dagger} =\begin{pmatrix}A&0\\0&A\end{pmatrix}.$
Then
$$(\Phi_{\dagger})^{-1} {J(\xi^2_1,\xi^4_3)}_{\dagger} \,  \Phi_{\dagger} =
 \begin{pmatrix}I&-S\\0&I\end{pmatrix}
\begin{pmatrix}A&0\\0&A\end{pmatrix}
 \begin{pmatrix}I&S\\0&I\end{pmatrix}
 =
\begin{pmatrix}A&AS-SA\\0&A\end{pmatrix} .$$
Now
$AS-SA =
\begin{pmatrix} -2\beta \xi^2_1 &2\alpha \xi^2_1 \\ 2\alpha \xi^2_1 &2\beta \xi^2_1\end{pmatrix}
=\begin{pmatrix} 0 &1\\ 1 &0\end{pmatrix} $
if we take $\alpha = \frac{\xi^2_1}{2}, \beta = 0.$
Hence $\Phi^{-1} J\,  \Phi$ is  a CS of type 1, and we readily see
that
$J\cong J(0,1,\xi^5_5,\xi^6_5),$ hence
$J\cong J(0,\beta)$ for some  $\beta \; , \; 0 < \beta \leqslant 1.$
\end{proof}
To summarize, we have shown the following :
\begin{theorem}
Any  CS on the Lie algebra
${M5}$ is equivalent to either the canonical CS $J_0$ or
some CS
$J(\xi^1_3,\xi^1_4,\xi^5_5,\xi^6_5)$  defined in (\ref{JM5case1}).
\end{theorem}
\subsubsection{}
$G_0$ is here the complex 3-dimensional Heisenberg group, considered as a real Lie group,
\textit{i.e.}
the real Lie group comprised of the matrices
\begin{equation}
\label{x}
 x=\begin{pmatrix}
1&x^1+iy^1&x^3+iy^3\\0&1&x^2+iy^2\\0&0&1\end{pmatrix}
\quad (x^k,y^k \in \Rmath \quad \forall k=1,2,3).\end{equation}
We here depart from the second kind coordinates to use the natural coordinates
defined by (\ref{x}). Then
the matrix $x$ in (\ref{x}) is
$x=\exp{(x^2 x_3+y^2 x_4+x^3 x_5+y^3 x_6)}\exp{(x^1 x_1-y^1 x_2)},$
and
\begin{eqnarray*}
%X_1 &=& \frac{\partial}{\partial x^1}\\
X_1 = \frac{\partial}{\partial x^1}\; &,& \;
X_2 = -\frac{\partial}{\partial y^1}\\
%X_2 &=& -\frac{\partial}{\partial y^1}\\
X_3 = \frac{\partial}{\partial x^2} +x^1 \frac{\partial}{\partial x^3} +y^1 \frac{\partial}{\partial y^3}  \; &,& \;
%X_3 &=& \frac{\partial}{\partial x^2} +x^1 \frac{\partial}{\partial x^3} +y^1 \frac{\partial}{\partial y^3} \\
%X_4 &=& \frac{\partial}{\partial y^2} -y^1 \frac{\partial}{\partial x^3}+x^1 \frac{\partial}{\partial y^3} \\
X_4 = \frac{\partial}{\partial y^2} -y^1 \frac{\partial}{\partial x^3}+x^1 \frac{\partial}{\partial y^3} \\
%X_5 &=& \frac{\partial}{\partial x^3}\\
%X_6 &=& \frac{\partial}{\partial y^3}.
X_5 = \frac{\partial}{\partial x^3}\; &,&\;
X_6 = \frac{\partial}{\partial y^3}.
\end{eqnarray*}
\subsubsection{Holomorphic functions for $J(\xi^1_3,\xi^1_4,\xi^5_5,\xi^6_5).$}
Let $G$ denote the group $G_0$ endowed with the left invariant  structure
of complex manifold  defined by
$J = J(\xi^1_3,\xi^1_4,\xi^5_5,\xi^6_5)$  in (\ref{JM5case1}).
One easily checks (with formulae (\ref{eqd}))  that
%\begin{eqnarray*}
\begin{equation*}
%\tilde{X}_6^{-} &=& -i \, \frac{1-i\xi^5_5}{\xi^6_5} \, \tilde{X}_5^{-} \\
\tilde{X}_6^{-} = -i \, \frac{1-i\xi^5_5}{\xi^6_5} \, \tilde{X}_5^{-} \; ; \;
%\tilde{X}_1^{-} &=& -i \, \frac{1+ia}{\xi^1_4} \, \tilde{X}_4^{-} \\
\tilde{X}_1^{-} = -i \, \frac{1+ia}{\xi^1_4} \, \tilde{X}_4^{-} \; ; \;
%i \, \tilde{X}_2^{-} - \tilde{X}_3^{-} &=&  -i \,\left( i d + \frac{(1+ia)\xi^1_3}{\xi^1_4}\right) \, \tilde{X}_4^{-}
i \, \tilde{X}_2^{-} - \tilde{X}_3^{-} =  -i \,\left( i d + \frac{(1+ia)\xi^1_3}{\xi^1_4}\right) \, \tilde{X}_4^{-} .
\end{equation*}
Hence
$ H_{\Cmath}(G) =\{f \in C^{\infty}(G_0) \; ; \; \tilde{X}_j^{-} \, f = 0
\; \forall j=1,3,5\}.$
Consider first the equation
\begin{equation}
\label{X5f}
\tilde{X}_5^{-} \, f = 0.
\end{equation}
One has
\begin{equation}
\label{X5f2}
\tilde{X}_5^{-} = X_5 +i(\xi^5_5 X_5 +\xi^6_5 X_6)
=\frac{\partial}{\partial x^3} +i \, \left( \xi^5_5 \, \frac{\partial}{\partial x^3} + \xi^6_5 \, \frac{\partial}{\partial y^3}\right)
=\frac{\partial}{\partial u^3} +i \, \frac{\partial}{\partial v^3}
=2\, \frac{\partial}{\partial \overline{w^3}}
\end{equation}
where
$
u^3 = x^3 -\frac{\xi^5_5}{\xi^6_5} \, y^3
\; ; \; v^3 = \frac{y^3}{\xi^6_5} \; ; \;
w^3 = u^3 +i v^3.
$
Equation (\ref{X5f}) simply means that $f$ is holomorphic with respect to $w^3.$
Consider now the equation
\begin{equation}
\label{X1f}
\tilde{X}_1^{-} \, f = 0.
\end{equation}
One has
\begin{equation*}
\tilde{X}_1^{-} = X_1 +i(aX_1 +b X_4)
= \frac{\partial }{\partial x^1} +i\left(
a\, \frac{\partial }{\partial x^1} +b \, \frac{\partial }{\partial y^2}\right)
+ib\left(
-y^1 \,\frac{\partial }{\partial x^3}
+x^1 \,\frac{\partial }{\partial y^3}\right).
\end{equation*}
We suppose that $f$ satisfies
equation (\ref{X5f}),  \textit{i.e.} $f$ is holomorphic with respect to $w^3.$
Hence
$$    \frac{\partial f}{\partial {x^3}} =   \frac{\partial f}{\partial {u^3}}
=   \frac{\partial f}{\partial {w^3}}$$
$$    \frac{\partial f}{\partial {y^3}} =  \frac{1}{\xi^6_5} \frac{\partial f}{\partial {v^3}}
-\frac{\xi^5_5}{\xi^6_5} \frac{\partial f}{\partial {u^3}}
=\frac{1}{\xi^6_5} (i-\xi^5_5) \frac{\partial f}{\partial {w^3}}.$$
Then equation (\ref{X1f}) reads
$$\frac{\partial f}{\partial x^1} +i\left(
a\, \frac{\partial f}{\partial x^1} +b \, \frac{\partial f}{\partial y^2}\right)
-b\left(
iy^1
+\frac{x^1}{\xi^6_5} (1+i\xi^5_5)\right) \frac{\partial f}{\partial {w^3}}=0$$
that is
\begin{equation}
\label{X1f=0}
2\, \frac{\partial f}{\partial \overline{w^1}} -\frac{b}{2}\, \left(
2iy^1
+\frac{1+i\xi^5_5}{\xi^6_5} ((1-ia)w^1 +(1+ia)\overline{w^1})\right)
\frac{\partial f}{\partial {w^3}}=0
\end{equation}
where
$u^1 = x^1 -\frac{a}{b} \, y^2
\; ; \; v^1 = \frac{y^2}{b}\ \; ; \;
w^1 = u^1 +i v^1.$
Finally, we turn to
the last equation
\begin{equation}
\label{X3f}
\tilde{X}_3^{-} \, f = 0.
\end{equation}
One has
\begin{eqnarray*}
\tilde{X}_3^{-} &=& i\xi^1_3  X_1 +iX_2 +X_3 +id X_4 \\
&=& i \xi^1_3\frac{\partial }{\partial x^1} -i\, \frac{\partial }{\partial y^1}
+\left( \frac{\partial }{\partial x^2}+ x^1 \,\frac{\partial }{\partial x^3}
+y^1 \,\frac{\partial }{\partial y^3}\right)
+id\, \left( \frac{\partial }{\partial y^2} -y^1 \,\frac{\partial }{\partial x^3}
+x^1 \,\frac{\partial }{\partial y^3}\right).
\end{eqnarray*}
Since we suppose $f$  holomorphic with respect to $w^3,$
equation (\ref{X3f}) then reads
\begin{equation*}
\frac{\partial f}{\partial x^2}
-i\, \frac{\partial f}{\partial y^1}
+ i\,  \xi^1_3\frac{\partial f}{\partial x^1}
+id\,  \frac{\partial f}{\partial y^2}
+\left(  x^1 -idy^1+(y^1+idx^1)\, \frac{i-\xi^5_5}{\xi^6_5} \right)
\frac{\partial f}{\partial {w^3}}=0
\end{equation*}
that is
%\begin{equation}
%\begin{split}
\begin{multline}
\label{X3f=0}
%2\, & \frac{\partial f}{\partial \overline{w^2}}
2\,  \frac{\partial f}{\partial \overline{w^2}}
+ \left(i \left(\xi^1_3-\frac{ad}{b}\right)-\frac{d}{b} \right)  \frac{\partial f}{\partial {w^1}}
+ \left(i \left(\xi^1_3-\frac{ad}{b}\right)+\frac{d}{b} \right) \frac{\partial f}{\partial \overline{w^1}}
\\
+ \frac{1}{2}\left[ \left((1-ia)w^1 +(1+ia)\overline{w^1})\right)
\left(1-\frac{d(1+i\xi^5_5)}{\xi^6_5}\right)
\right.
\\
\left.
-(w^2 -\overline{w^2})\left(-d+\frac{1+i\xi^5_5}{\xi^6_5}\right)
\right]
\frac{\partial f}{\partial {w^3}}=0
%\end{split}
%\end{equation}
\end{multline}
where
$w^2 = x^2 -i y^1.$
Now equation (\ref{X1f=0}) reads
\begin{equation}
\label{X1f=02}
2\, \frac{\partial f}{\partial \overline{w^1}} -\frac{b}{2}\, \left(
-w^2 +\overline{w^2}
+\frac{1+i\xi^5_5}{\xi^6_5} ((1-ia)w^1 +(1+ia)\overline{w^1})\right)
\frac{\partial f}{\partial {w^3}}=0.
\end{equation}
From equations
(\ref{X5f2}),
(\ref{X1f=02}),
(\ref{X3f=0}),
one readily sees that the functions $\varphi^1$ and $\varphi^2$ defined by
 \begin{eqnarray}
 \label{phi1}
 \varphi^1 &=& 2w^1 -
\left(i \left(\xi^1_3-\frac{ad}{b}\right)-\frac{d}{b} \right) \overline{w^2}\\
 \label{phi2}
 \varphi^2 &=& w^2
 \end{eqnarray}
 are holomorphic on $G.$
We look for a holomorphic function which depends on $w^3.$
For any $C^{\infty}$-function $\psi(w^1,w^2,\overline{w^2}),$ \textit{i.e.}
$\psi$ doesn't depend on $w^3,\overline{w^3},\overline{w^1},$  the following function
$f_1$ is a solution  of equations (\ref{X5f}) and (\ref{X1f=02}) :
%\begin{equation}
%\begin{split}
\begin{multline}
f_1 = w^3 + \frac{b}{4} \left[ -w^2\,  \overline{w^1} + \overline{w^2} \, \overline{w^1}
+\frac{1+i\xi^5_5}{\xi^6_5} \left(
(1-ia){w^1}\overline{w^1} + (1+
 ia)
%(1-ia){w^1}\overline{w^1} + (1+ia)
\frac{(\overline{w^1})^2}{2}
\right) \right]
 \\
 + \psi(w^1,w^2,\overline{w^2}) .
 %+ \psi(w^1,w^2,\overline{w^2}) .
\end{multline}
%\end{split}
%\end{equation}
We want to  choose $\psi$ such that $f_1$ is a solution of
(\ref{X3f=0}) as well.
First, we have :
\begin{eqnarray*}
\frac{\partial f_1}{\partial {w^1}} &=& \frac{b}{4\xi^6_5} \, (1-ia)(1+i\xi^5_5) \, \overline{w^1}
+ \frac{\partial \psi}{\partial {w^1}} \\
 \frac{\partial f_1}{\partial \overline{w^1}} &=&\frac{b}{4}\, \left(
-w^2 +\overline{w^2}
+\frac{1+i\xi^5_5}{\xi^6_5} ((1-ia)w^1 +(1+ia)\overline{w^1})\right)
 \\ \frac{\partial f_1}{\partial \overline{w^2}} &=&\frac{b}{4}\, \overline{w^1}
+ \frac{\partial \psi}{\partial \overline{w^2}} \\
 \\ \frac{\partial f_1}{\partial {w^3}} &=&1 \; .
\end{eqnarray*}
Introducing these values in (\ref{X3f=0})
we find that $f_1$ is a solution to (\ref{X3f=0}) if and only if
$$ N \, \overline{w^1}
 - M(w^2-\overline{w^2}) + \Lambda w^1
+ \left(i \left(\xi^1_3-\frac{ad}{b}\right)-\frac{d}{b} \right) \frac{\partial \psi}{\partial {w^1}}+
2\,  \frac{\partial \psi}{\partial \overline{w^2}}
 = 0$$
 where
 \begin{eqnarray*}
N&=&\frac{1}{\xi^6_5} \left( ib\xi^1_3 -b\xi^1_3\xi^5_5
+b\xi^6_5 +(1+ia)\xi^6_5 -d(1+ia)(1+i\xi^5_5)\right)\\
 M& =&
 \frac{1}{4} (ib\xi^1_3 +(1-ia)d)
 +\frac{1}{2} \left(\frac{1}{\xi^6_5} -d  +i\frac{\xi^5_5}{\xi^6_5}\right) \\
 \Lambda &=&
 \frac{(1+i\xi^5_5)(1-ia)}{4\xi^6_5}\left(ib\xi^1_3 +(1-ia)d \right)
 + \frac{1-ia}{2}\left(1-\frac{d(1+i\xi^5_5)}{\xi^6_5}  \right).
 \end{eqnarray*}
 A computation shows that $N$ is actually equal to $0.$
 Hence
$f_1$ is a solution to (\ref{X3f=0}) if and only if
\begin{equation}
\label{eqpsi}
- M(w^2-\overline{w^2}) + \Lambda w^1
+ \left(i \left(\xi^1_3-\frac{ad}{b}\right)-\frac{d}{b} \right) \frac{\partial \psi}{\partial {w^1}}+
2\,  \frac{\partial \psi}{\partial \overline{w^2}}
 = 0.
 \end{equation}
 Note that
$$i \left(\xi^1_3-\frac{ad}{b}\right)-\frac{d}{b} = 0 \Leftrightarrow \xi^1_3=d=0
\Leftrightarrow \xi^1_3=0, \xi^1_4 =1.$$
Hence, in the nonabelian case where one doesn't have simulteanously
$\xi^1_3=0, \xi^1_4 =1,$
 for  (\ref{eqpsi}) to hold, it is sufficient to have
 \begin{eqnarray*}
\frac{\partial \psi}{\partial {w^1}}&=& -\frac{\Lambda}{i \left(\xi^1_3-\frac{ad}{b}\right)-\frac{d}{b}} \, w^1\\
\frac{\partial \psi}{\partial \overline{w^2}}&=&
  \frac{M}{2}(w^2-\overline{w^2})
 \end{eqnarray*}
 which gives a solution
 \begin{equation*}
 \psi =
-\frac{\Lambda}{i \left(\xi^1_3-\frac{ad}{b}\right)-\frac{d}{b}} \, \frac{(w^1)^2}{2}+
  \frac{M}{2}\left(w^2\overline{w^2}-\frac{(\overline{w^2})^2}{2}\right).
 \end{equation*}
 In the abelian case, one can take
 \begin{equation*}
 \psi =
-\frac{\Lambda}{2} \,
w^1 \overline{w^2}+
  \frac{M}{2}\left(w^2\overline{w^2}-\frac{(\overline{w^2})^2}{2}\right).
 \end{equation*}
We finally get the holomorphic function $f_1 :$ in the nonabelian case,
\begin{multline}
\label{phi3}
f_1 = w^3 + \frac{b}{4} \left[ -w^2\,  \overline{w^1} + \overline{w^2} \, \overline{w^1}
+\frac{1+i\xi^5_5}{\xi^6_5} \left(
(1-ia){w^1}\overline{w^1} + (1+
 ia)
%(1-ia){w^1}\overline{w^1} + (1+ia)
\frac{(\overline{w^1})^2}{2}
\right) \right]
 \\
-\frac{\Lambda}{i \left(\xi^1_3-\frac{ad}{b}\right)-\frac{d}{b}} \, \frac{(w^1)^2}{2}+
  \frac{M}{2}\left(w^2\overline{w^2}-\frac{(\overline{w^2})^2}{2}\right).
\end{multline}
In the abelian case,
\begin{multline}
f_1 = w^3 + \frac{b}{4} \left[ -w^2\,  \overline{w^1} + \overline{w^2} \, \overline{w^1}
+\frac{1+i\xi^5_5}{\xi^6_5} \left(
(1-ia){w^1}\overline{w^1} + (1+ ia)
\frac{(\overline{w^1})^2}{2}
\right) \right]
 \\
-\frac{\Lambda}{2} \,
w^1\overline{w^2}+
  \frac{M}{2}\left(w^2\overline{w^2}-\frac{(\overline{w^2})^2}{2}\right).
\end{multline}
Note that in the abelian case, one can take $\xi^5_5 = 0, \xi^6_5 = \beta \; , \;
0 < \beta \leqslant 1,$ and then $a=c=d=0,\linebreak b=-1, M=\frac{1}{2\beta}, \Lambda=\frac{1}{2},$ hence
\begin{multline}
\label{phi3abelienbeta}
f_1 = w^3 - \frac{1}{4} \left[ -w^2\,  \overline{w^1} + \overline{w^2} \, \overline{w^1}
+\frac{1}{\beta} \left(
{w^1}\overline{w^1} +
\frac{(\overline{w^1})^2}{2}
\right) \right]
 \\
-\frac{1}{4}w^1 \overline{w^2}+
  \frac{1}{4\beta}\left(w^2\overline{w^2}-\frac{(\overline{w^2})^2}{2}\right).
\end{multline}
In both abelian and nonabelian cases, let $F : G \rightarrow \Cmath^3$ defined
by
$F=(\varphi^1,\varphi^2,\varphi^3)$
where $\varphi^1, \varphi^2$ are defined in
(\ref{phi1}), (\ref{phi2}) and $\varphi^3 = f_1.$
$F$ is a global chart on $G.$
We determine now how the multiplication of $G$ looks like in the chart $F.$
Recall first the formulae :
\begin{equation}
\label{w1-2-3}
w^1 =\left(x^1 -\frac{a}{b}y^2\right) +i \frac{y^2}{b}\\
\quad , \quad
w^2 =x^2 -iy^1
\quad , \quad
w^3= \left(x^3 - \frac{\xi^5_5}{\xi^6_5}\right) +i \, \frac{y^3}{\xi^6_5}.
\end{equation}
Let $a,x \in G $ with respective canonical coordinates
$(x^1,y^1,x^2,y^2,x^3,y^3), (\alpha^1, \beta^1, \alpha^2, \beta^2, \alpha^3, \beta^3)$
as in  (\ref{x}).
With obvious notations,
%$\alpha =[w^1_\alpha, w^2_\alpha,w^3_\alpha],$
%$x =[w^1_x, w^2_x,w^3_x],$
%$\alpha \, x =[w^1_{\alpha x}, w^2_{\alpha x},w^3_{\alpha x}],$
%$\alpha =[\varphi^1_\alpha, \varphi^2_\alpha,\varphi^3_\alpha],$
%$x =[\varphi^1_x, \varphi^2_x,\varphi^3_x],$
%$\alpha \, x =[\varphi^1_{\alpha x}, \varphi^2_{\alpha x},\varphi^3_{\alpha x}].$
%We want to compute
%$\varphi^1_{\alpha x}, \varphi^2_{\alpha x},\varphi^3_{\alpha x}.$
%First, compute
%$w^1_{\alpha x}, w^2_{\alpha x},w^3_{\alpha x}.$
by matrix multiplication and  (\ref{w1-2-3}) one gets
\begin{equation*}
w^1_{\alpha x} = w^1_\alpha + w^1_x
\quad , \quad
w^2_{\alpha x} = w^2_\alpha + w^2_x
\quad , \quad
w^3_{\alpha x} = w^3_\alpha + w^3_x + \chi(\alpha,x)
\end{equation*}
where $\chi(\alpha,x) =\alpha^1 x^2 -\beta^1 y^2 +\frac{i-\xi^5_5}{\xi^6_5} \, (\alpha^1 y^2 + \beta^1 x^2).$
Then from (\ref{phi1}), (\ref{phi2}) :
\begin{equation*}
\varphi^1_{\alpha x} = \varphi^1_\alpha + \varphi^1_x
\quad , \quad
\varphi^2_{\alpha x} = \varphi^2_\alpha + \varphi^2_x.
\end{equation*}
To get
$\varphi^3_{\alpha x},$ we just make the substitutions
$w^1\rightarrow w^1_\alpha + w^1_x, \;
w^2\rightarrow w^2_\alpha + w^2_x, \;
w^3\rightarrow w^3_\alpha + w^3_x + \chi(\alpha,x)$
in (\ref{phi3}) (we consider here the nonabelian case).
Now, let
$\Delta = \varphi^3_{\alpha x} - \varphi^3_\alpha - \varphi^3_x .$
Computations give:
\begin{equation*}
\Delta =
\frac{1}{8\xi^6_5} \left(
C_1  \varphi^1_x
+C_2 \xi^6_5 \varphi^2_x \right)
\end{equation*}
where \\
%\begin{multline*}
$$C_1=
b(1+i\xi^5_5) \overline{\varphi^1_\alpha}
+(b-1+ia)\xi^6_5 \overline{\varphi^2_\alpha}
%\\ +\frac{b(a+i)\xi^6_5}{b\xi^1_3 - d(a-i)} {\varphi^1_\alpha}
 +\frac{b(a+i)\xi^6_5}{b\xi^1_3 - d(a-i)} {\varphi^1_\alpha}
+\left((\xi^5_5 -i)(b\xi^1_3 - da-id) -b\xi^6_5\right) {\varphi^2_\alpha}
\quad ,
$$
%\end{multline*}
\begin{multline*}
C_2 =
(1-b+ia)\overline{\varphi^1_\alpha}
\\
+\left(\frac{4(1+i\xi^5_5)}{\xi^6_5} +(a+i(b+1))\xi^1_3 -\frac{d}{b}(1+a^2) -d(1+ia)
\right) \overline{\varphi^2_\alpha}
\\
+(1-ai) {\varphi^1_\alpha}
\\
+\left(-\frac{2(1+i\xi^5_5)}{\xi^6_5} +(a+i(b-1))\xi^1_3
-\frac{d}{b}(1+a^2) +d(1-ia)
\right)
{\varphi^2_\alpha} .
\end{multline*}
%The fact that left translations are holomorphic, though the multiplication isn't, appears
%clearly on these formulae.
For example, in the case of $J(1,1,1,1),$
\begin{multline*}
\Delta =
\left (
-\frac{3-i}{4}\,  \overline{\varphi^2_\alpha}
-\frac{5+5i}{8} \, \overline{\varphi^1_\alpha}
+\frac{3+4i}{8}\,  {\varphi^2_\alpha}
+\frac{7+i}{16} \, {\varphi^1_\alpha}\right)
\varphi^1_x
\\
+
\left (
\frac{3+i}{4} \, \overline{\varphi^2_\alpha}
+\frac{3+i}{4} \, \overline{\varphi^1_\alpha}
-\frac{1+3i}{4} \, {\varphi^2_\alpha}
+\frac{1-2i}{8}\,  {\varphi^1_\alpha}\right)
\varphi^2_x \; .
\end{multline*}
Finally, in the abelian case  $J(0,1,0,\beta),$ one has to use
(\ref{phi3abelienbeta}) and one gets :
\begin{equation*}
%\begin{multline}
\Delta = -
\left (
\frac{1}{4}\,  \overline{\varphi^2_\alpha}
+\frac{1}{8\beta} \, \overline{\varphi^1_\alpha}
-\frac{1}{8}\,  {\varphi^2_\alpha}
+\frac{1}{16\beta} \, {\varphi^1_\alpha}\right)
\varphi^1_x
+
\left (
\frac{1}{2\beta} \, \overline{\varphi^2_\alpha}
+\frac{1}{4} \, \overline{\varphi^1_\alpha}
-\frac{1}{4\beta} \, {\varphi^2_\alpha}
+\frac{1}{8}\,  {\varphi^1_\alpha}\right)
\varphi^2_x \; .
%\end{multline}
\end{equation*}
%%%%%%%%%%%%%%%%%%%%%%%%%%%%%%%%%%%%%%%%%%%%%%%%%%%%%%%%%%%%%%%%%%%%%%%%%%%%%%%%%%%%
%%%%%%%%%%%%%%%%%%%%%%%%%%%%%%%%%%%%%%%%%%%%%%%%%%%%%%%%%%%%%%%%%%%%%%%%%%%%%%%%%%%%%%%%%%

\end{document}